\documentclass{amsart}

\usepackage{graphicx}
\usepackage{longtable}

\begin{document}

\title[Covering Finite Square Lattices]{Finite Square Lattice Vertex Cover by a Baseline Set Defined With a Minimum Sublattice}

\author{Richard J. Mathar}
\urladdr{http://www.strw.leidenuniv.nl/~mathar}
\address{Leiden Observatory, Leiden University, P.O. Box 9513, 2300 RA Leiden, The Netherlands}


\date{\today}

\begin{abstract}
Each straight infinite line defined by two vertices of a finite
square point lattice contains (covers) these two points and a---possibly empty---subset
of points that happen to be collinear to these.
This work documents vertex subsets of minimum order such that
the sum of the infinite straight lines associated with the edges of their
complete subgraph covers the entire set of vertices (nodes). This is an abstraction
to the problem of sending a light signal to all stations (receivers) in
a square array with a minimum number of stations also equipped
with transmitters to redirect the light to other transmitters.
\end{abstract}

\maketitle
\section{Geometry of the Optimization Problem}
Given a $(n+1)\times (n+1)$ regular square lattice,
we aim at a quantitative assessment of the following optimization
problem considered by Alon
\cite{AlonGFA1}:
a subset of vertices of minimum size (order) $t$ is to be found such that 
each of the $(n+1)^2$ vertices of the square lattice lies on at least one of the
infinite lines defined by connecting all pairs of the $t$ points.
The number of these lines is generally ${t \choose 2}=t(t-1)/2$, once
a subset of vertices has been selected, but may
actually be smaller if some of the vertices of that subset are collinear.

The requirement that solutions cover points by straight lines puts this
in the category of geometric optimizations rather than 
reduction problems in graph theory which deal with the topology and
connectivity of the full square lattice. It also defines this as a model for
distribution of a signal via classical optical (straight) light from
one of the vertices to all the others, supposed that (i) most of them
are of some kind of semi-transparent receivers, and (ii) those that we
selected as a sublattice (transmitters) can in addition dispatch
(redirect) the signal
to others of their kind. Our optimization calls for the cheapest solution,
which minimizes the number of---supposedly more versatile and more expensive---transmitters.
It also minimizes broadcast times if the signal is distributed along
some Eulerian trail between the transmitters.

Solutions to this task are not necessarily unique. First, there is a generic
multiplicity of solutions associated with the symmetry group 4mm
of the square \cite{TerzibaschPSS133}: there are up to 8 equivalent
(congruential) geometries,
when any setup with a sublattice of $t(n)$ vertices is rotated by
multiples of 90 degrees (operations $\delta_4^1$, $\delta_4^2$, $\delta_4^3$)
or flipped into its mirror image (operations $m_x$ or $m_y$ flip
along the horizontal or vertical mid axis, $m_d$ or $m_{d'}$
along the main or secondary diagonal). The actual
multiplicity depends on how far the sublattice itself shows the full
symmetry of the square or `breaks' this symmetry. Second, as we shall
see, sets of incongruential solutions may display the same
minimum count $t(n)$ of vertices to span the lines.

The results are put into two categories: Section \ref{sec.congr}
summarizes minimum solutions for the cases of small lattices $n=3$
to $n=6$ for which the entire configuration space has been scanned:
the corresponding $t(n)$ are proven to be minimum at those
individual $n$. One representative of each of the incongruential
solutions will be drawn for $n=3$ to $5$, and some incongruential solutions for $n=6$.
Section \ref{sec.boun} addresses larger $n$ in the range $n=7$ to $110$,
where solutions have been searched by a mixture of heuristics;
those small $t(n)$ actually found establish only upper bounds
since even smaller $t(n)$ may have evaded detection.

\section{Optimum Solutions}\label{sec.congr}
The graphs show the general square lattice vertices as grey squares,
the vertices of the spanning sublattice in Red---darker in black\&white printouts.

All base lines between the vertices of the red sublattice are drawn,
which allows a visual verification of the matching requirement
that each vertex is hit by at least one line.
(Unlike some associated `attacking Queens' problems on chess boards,
the slope of the lines is not constrained to multiples of 45 degrees.)

\begin{figure}[htb]
\includegraphics[scale=0.25]{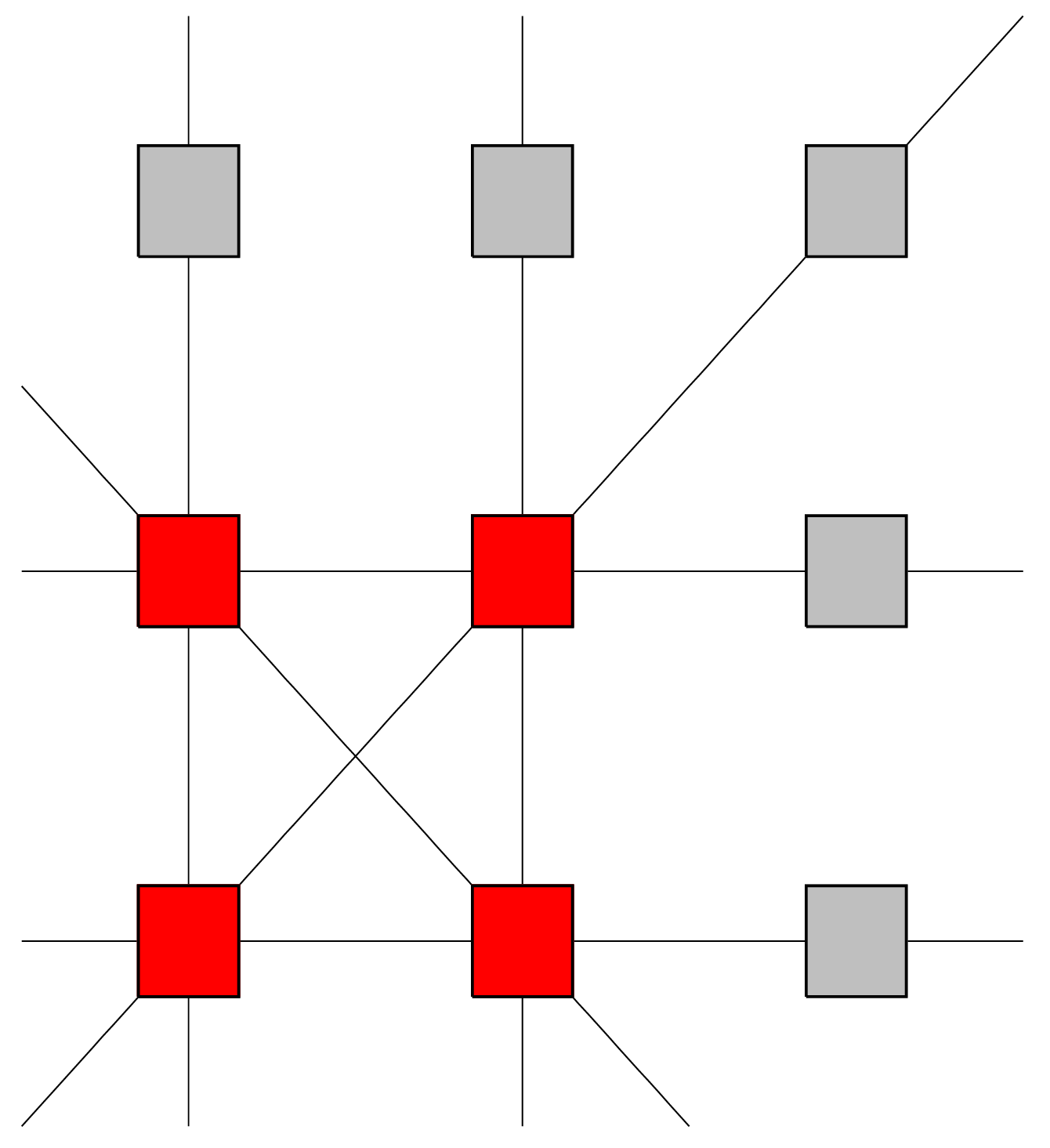}
\includegraphics[scale=0.25]{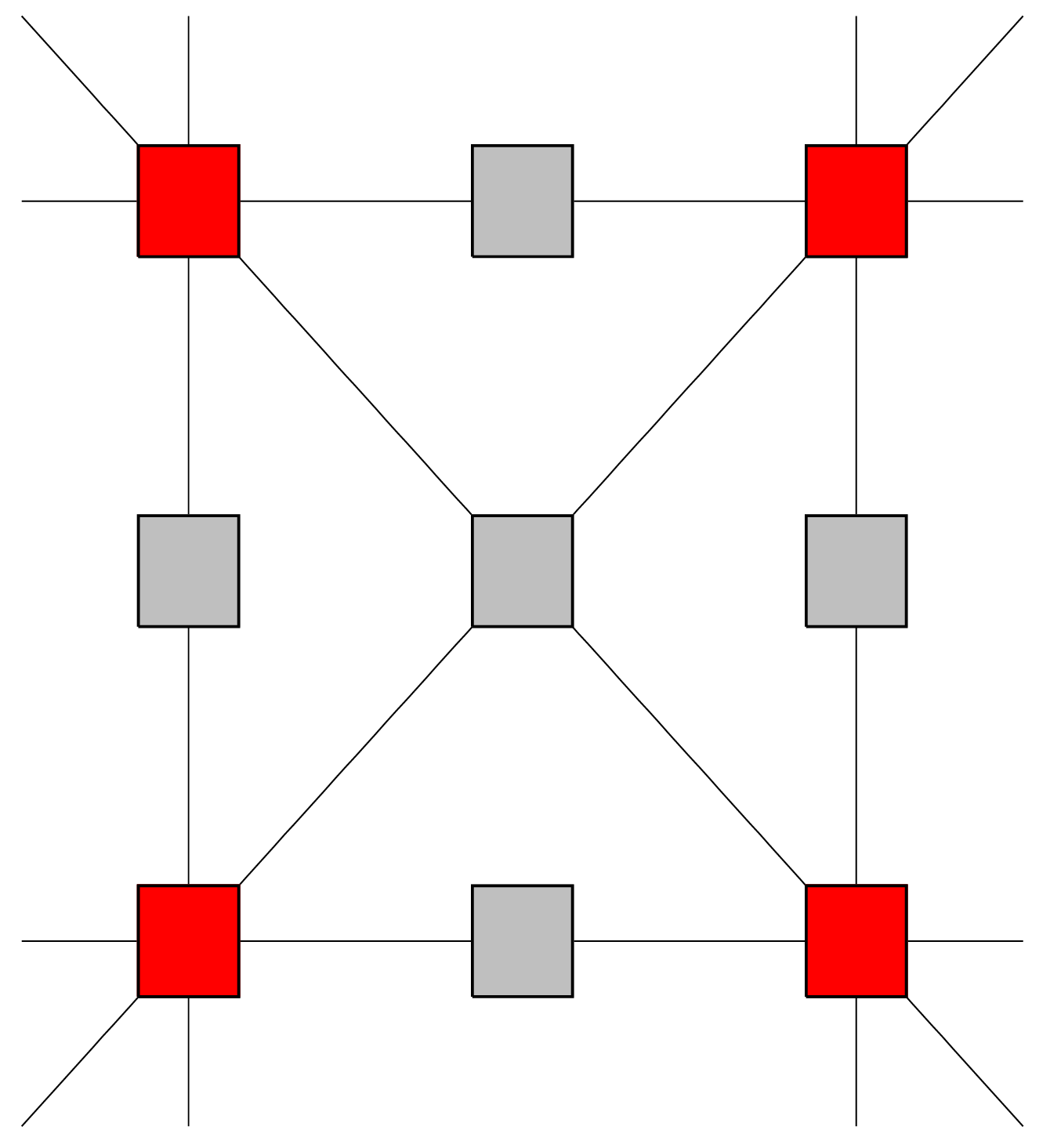}
\caption{$t(2)=4$.}
\label{fig.n2}
\end{figure}

\begin{figure}[htb]
\includegraphics[scale=0.3]{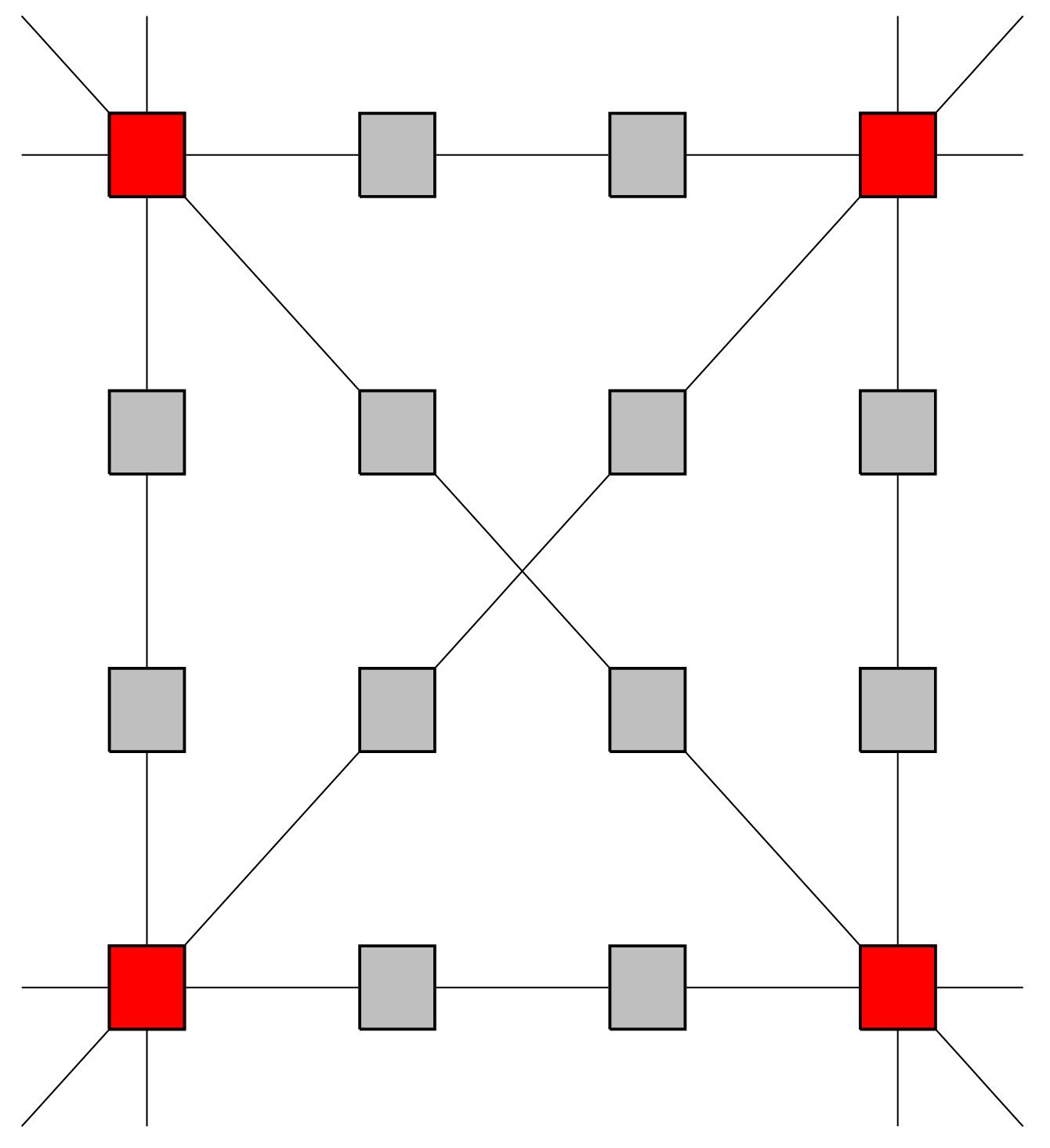}
\includegraphics[scale=0.3]{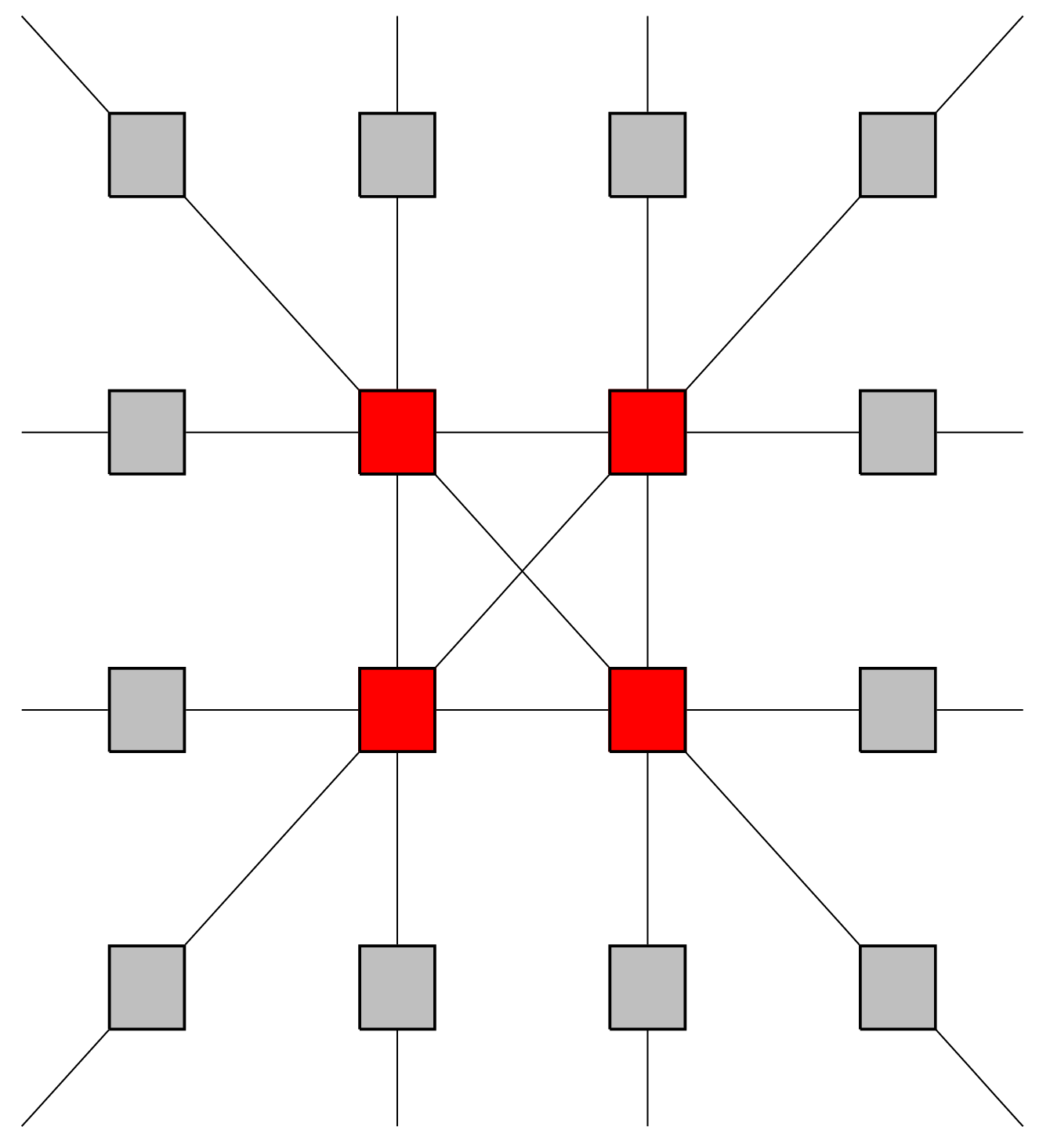}
\caption{$t(3)=4$.}
\label{fig.n3}
\end{figure}

The obvious two incongruential solutions for $n=2$ with sublattices
of $t(2)=4$ are shown in Figure \ref{fig.n2}, and the two incongruential
solutions for $n=3$ with $t(3)=4$ in Figure \ref{fig.n3}.

\begin{figure}[h]
\includegraphics[scale=0.3]{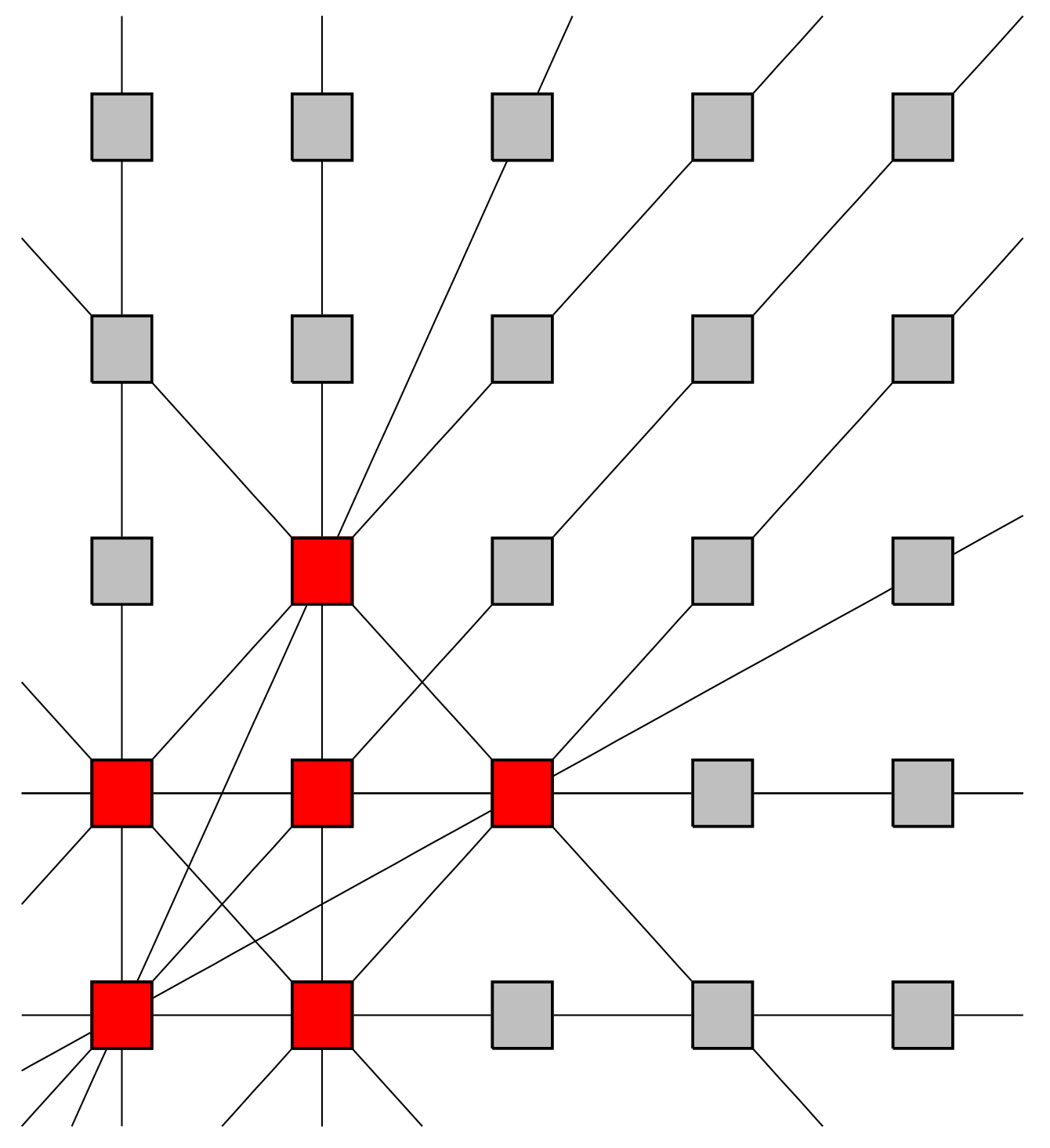}
\includegraphics[scale=0.3]{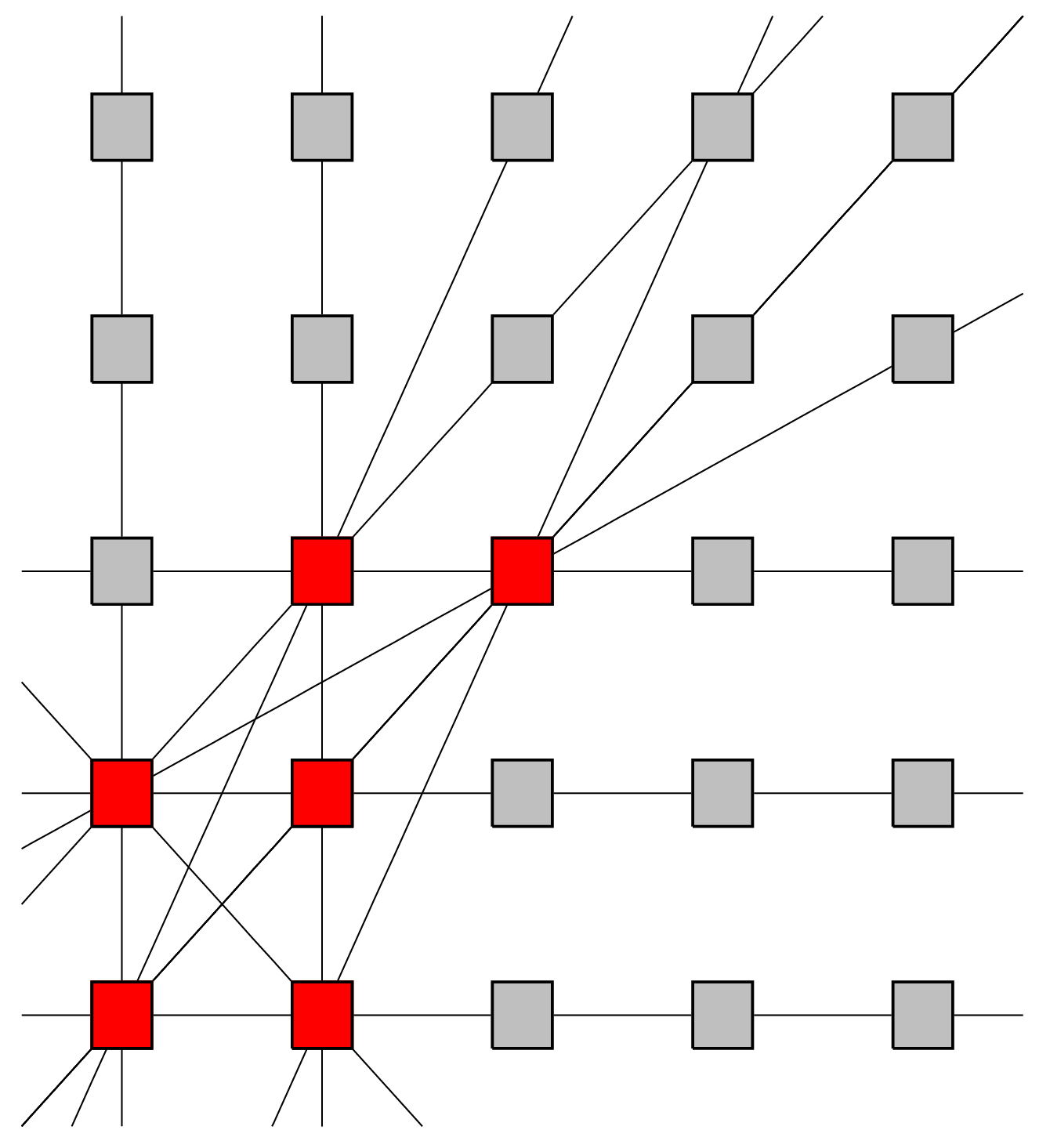}
\includegraphics[scale=0.3]{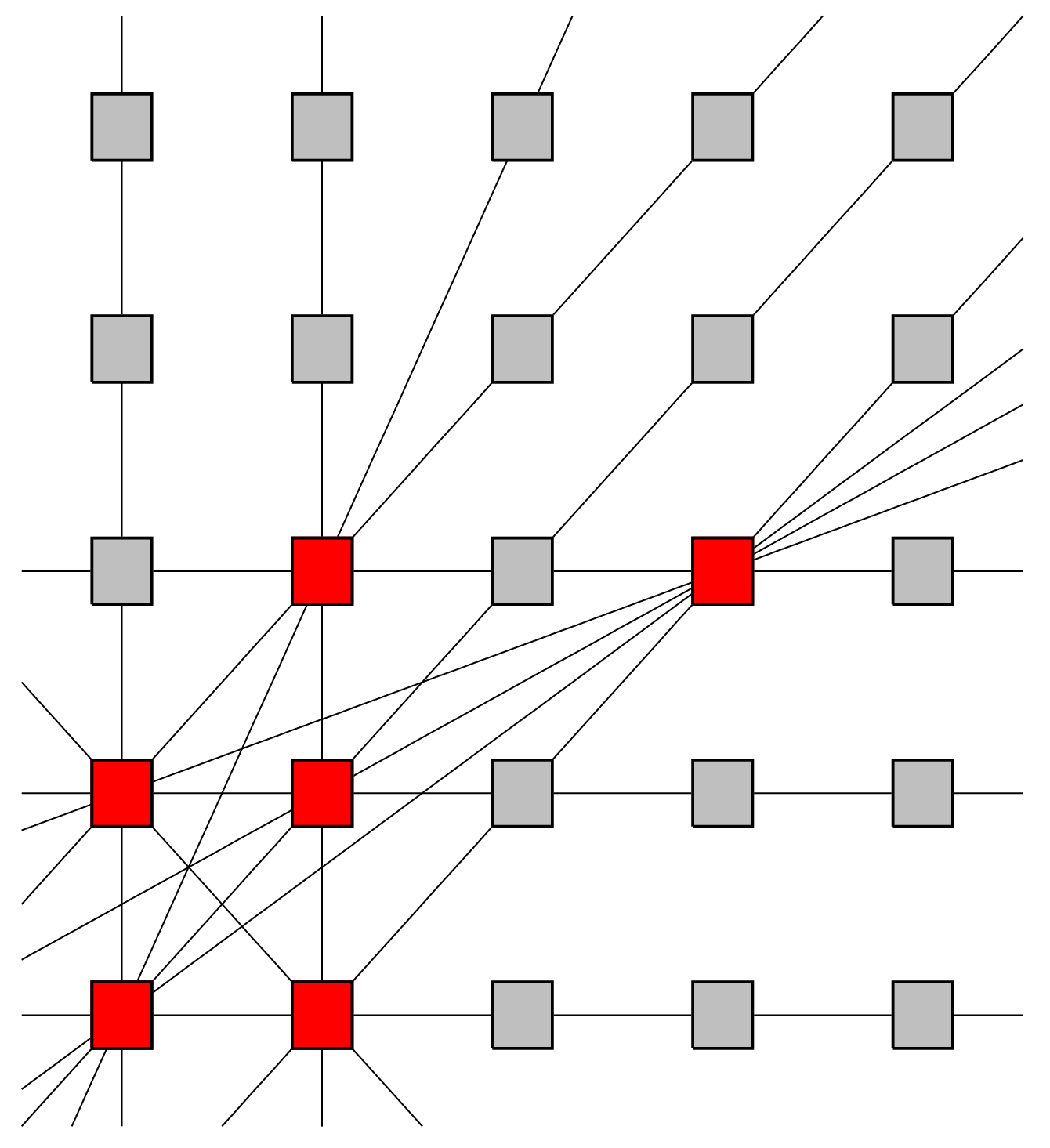}
\includegraphics[scale=0.3]{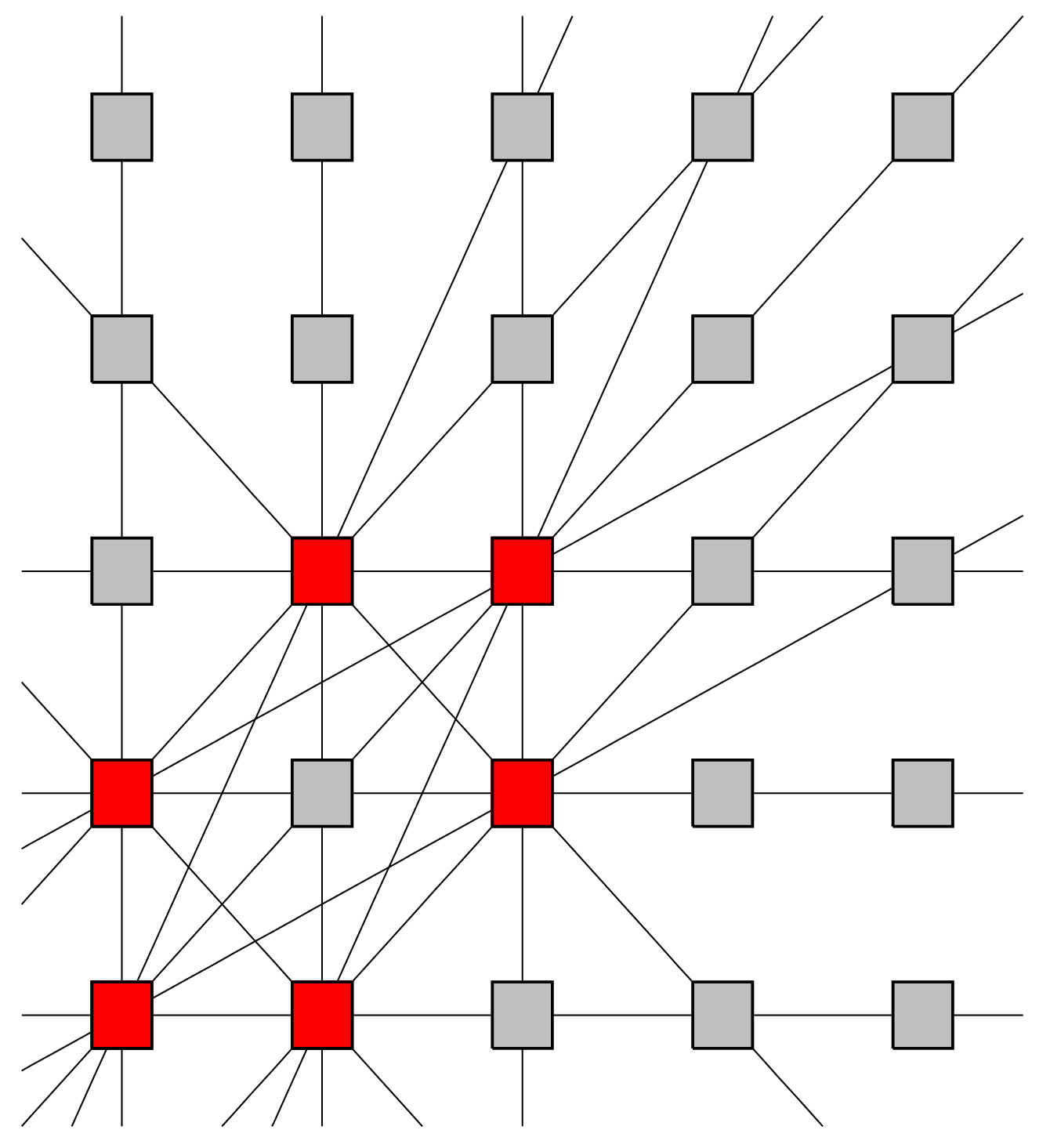}
\includegraphics[scale=0.3]{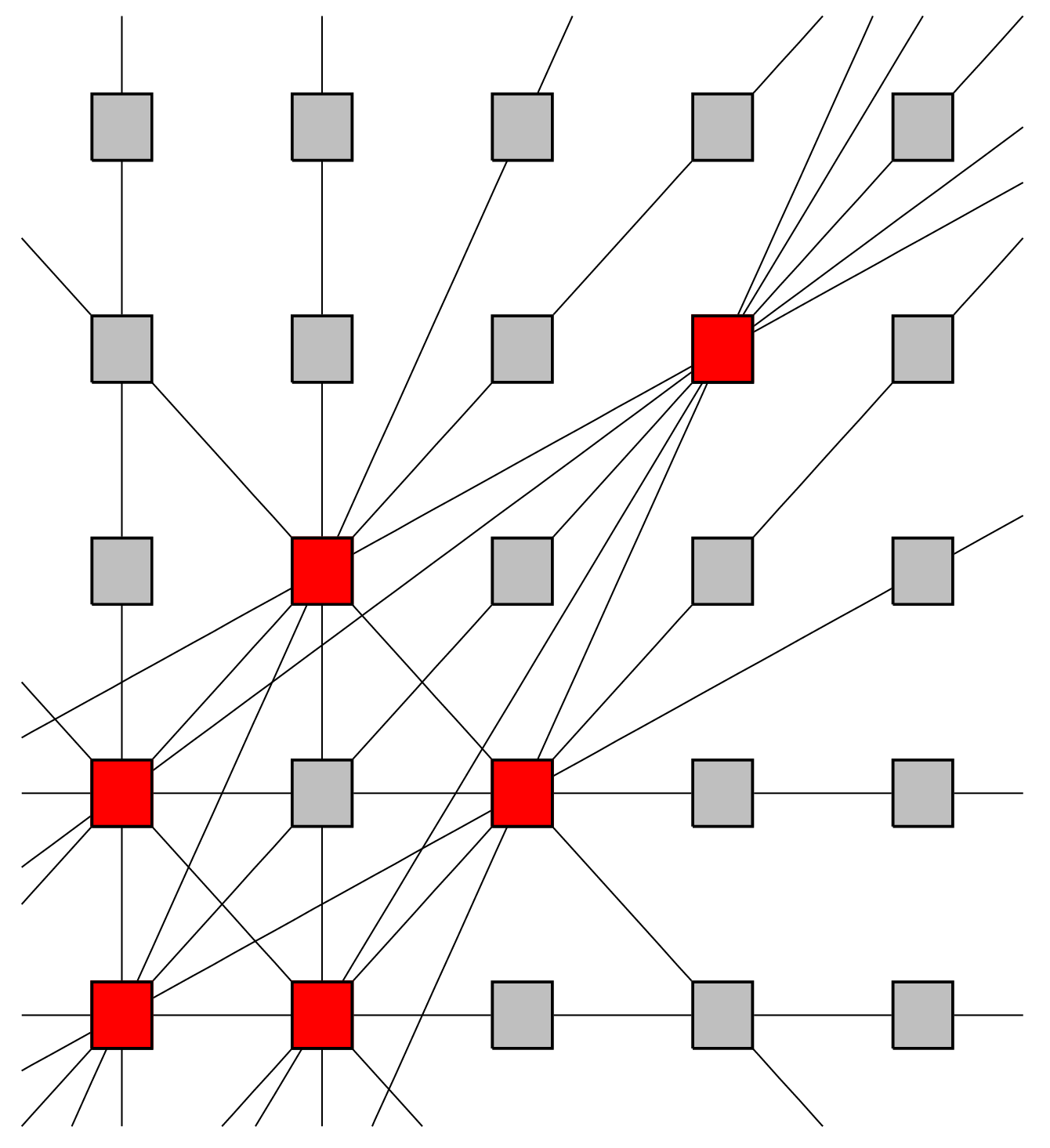}
\includegraphics[scale=0.3]{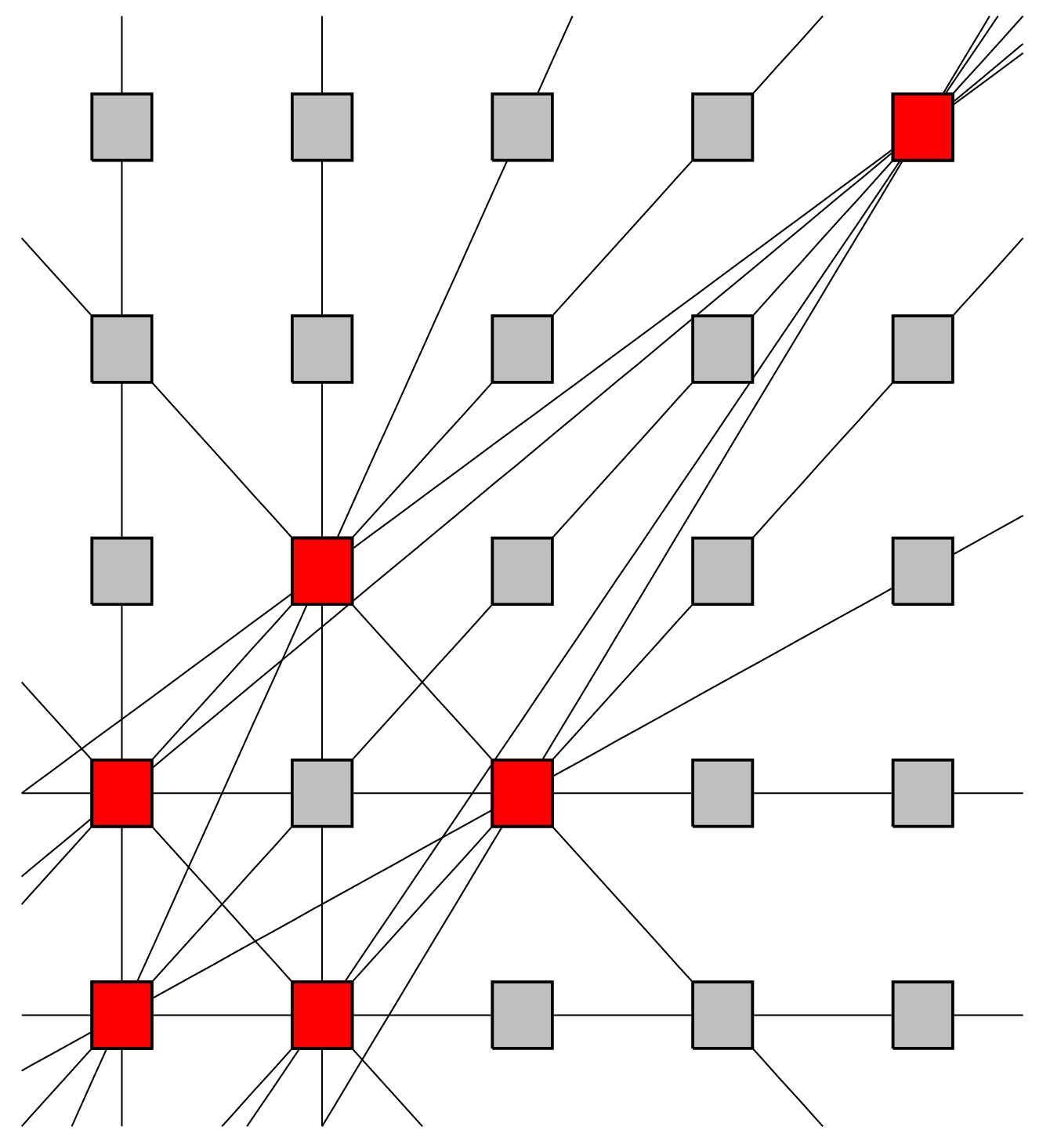}
\includegraphics[scale=0.3]{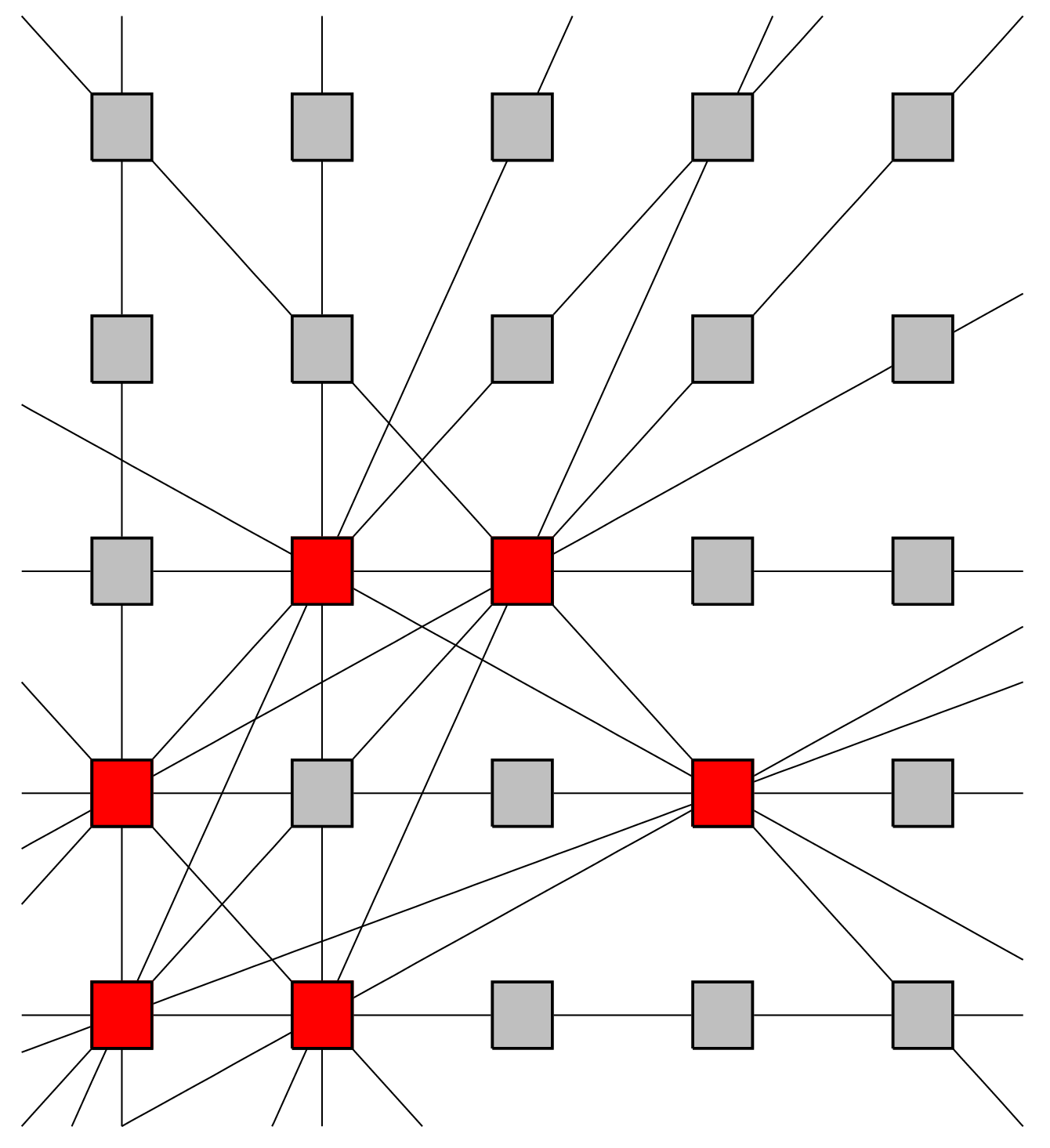}
\includegraphics[scale=0.3]{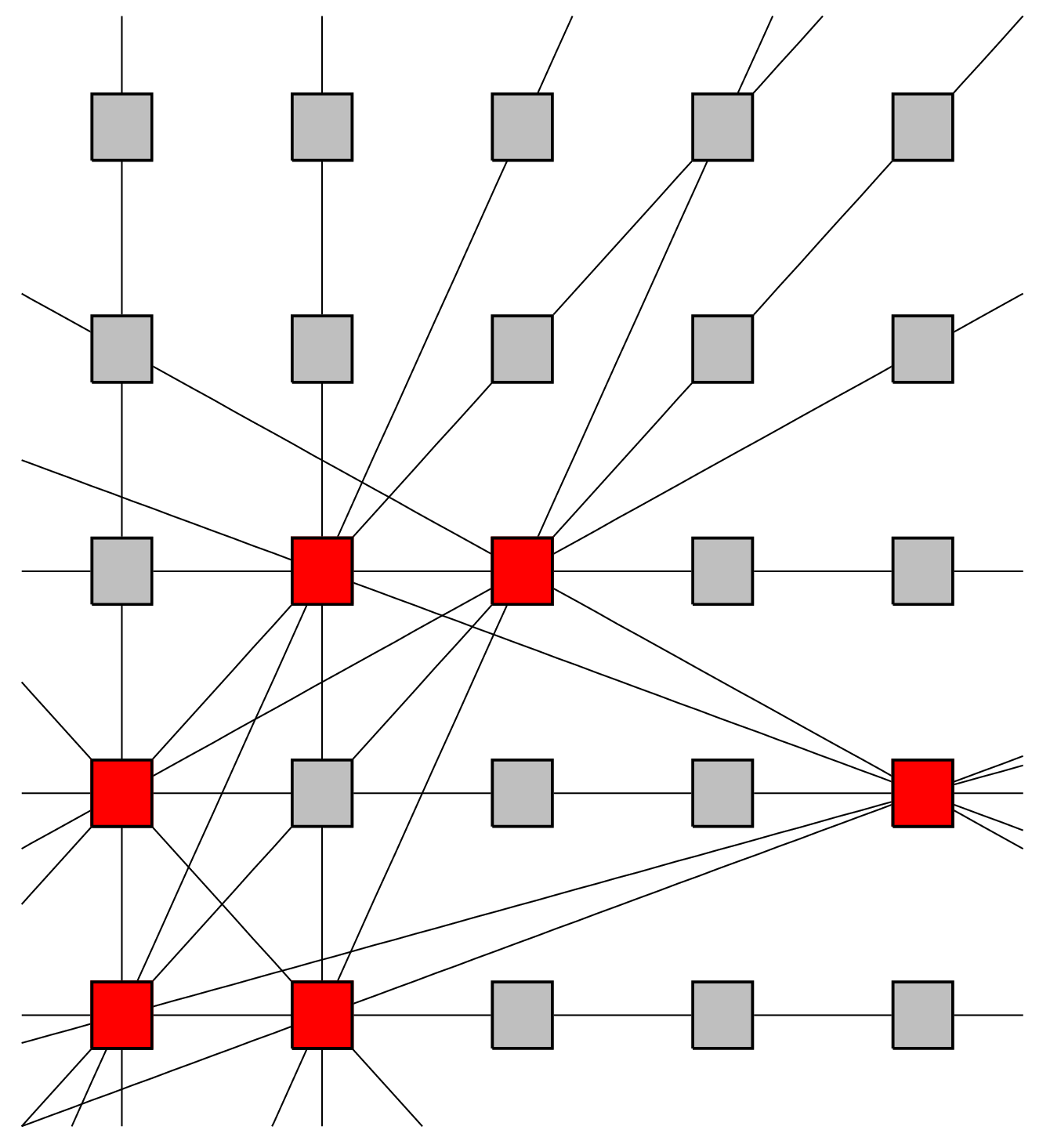}
\includegraphics[scale=0.3]{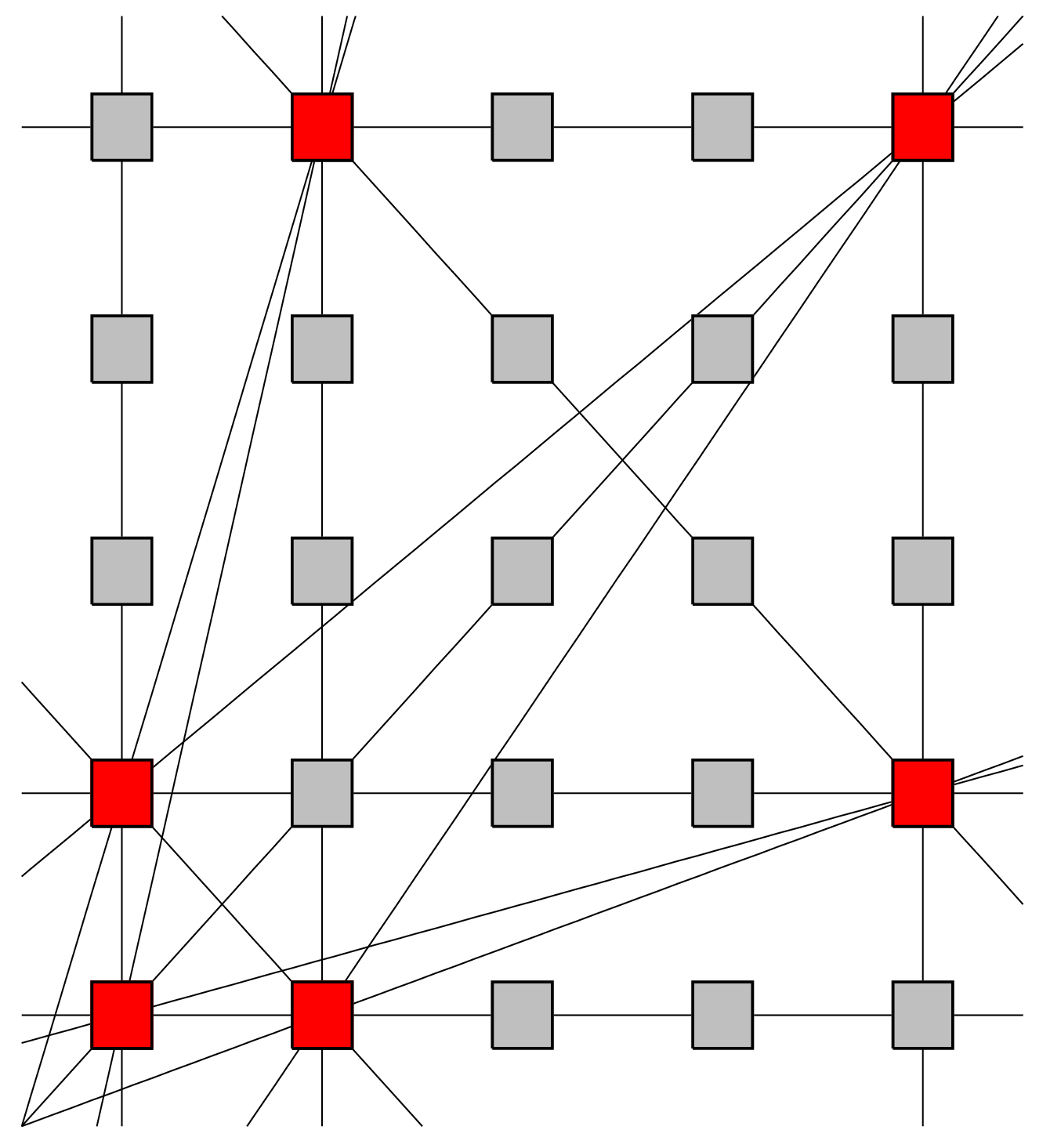}
\includegraphics[scale=0.3]{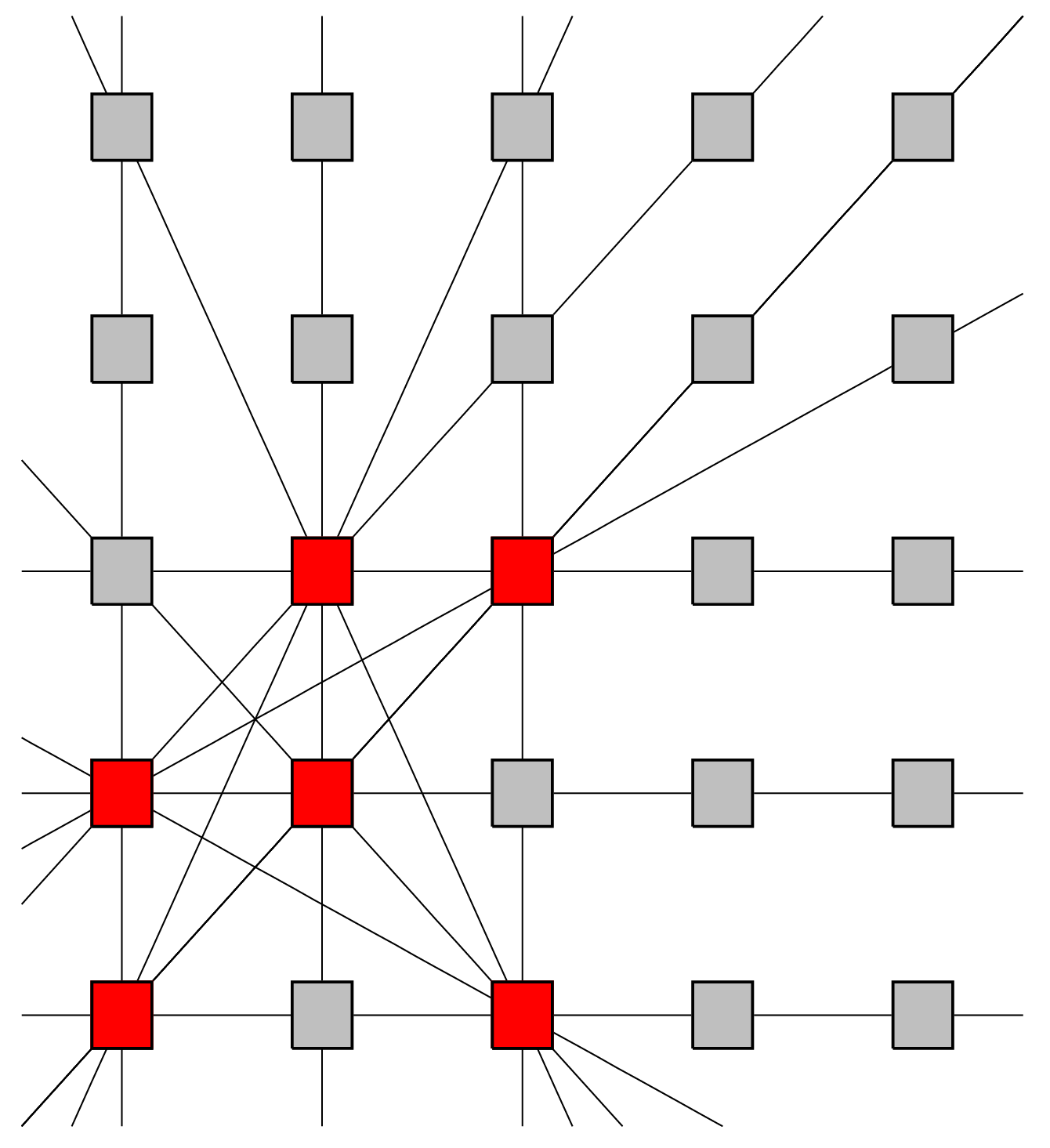}
\includegraphics[scale=0.3]{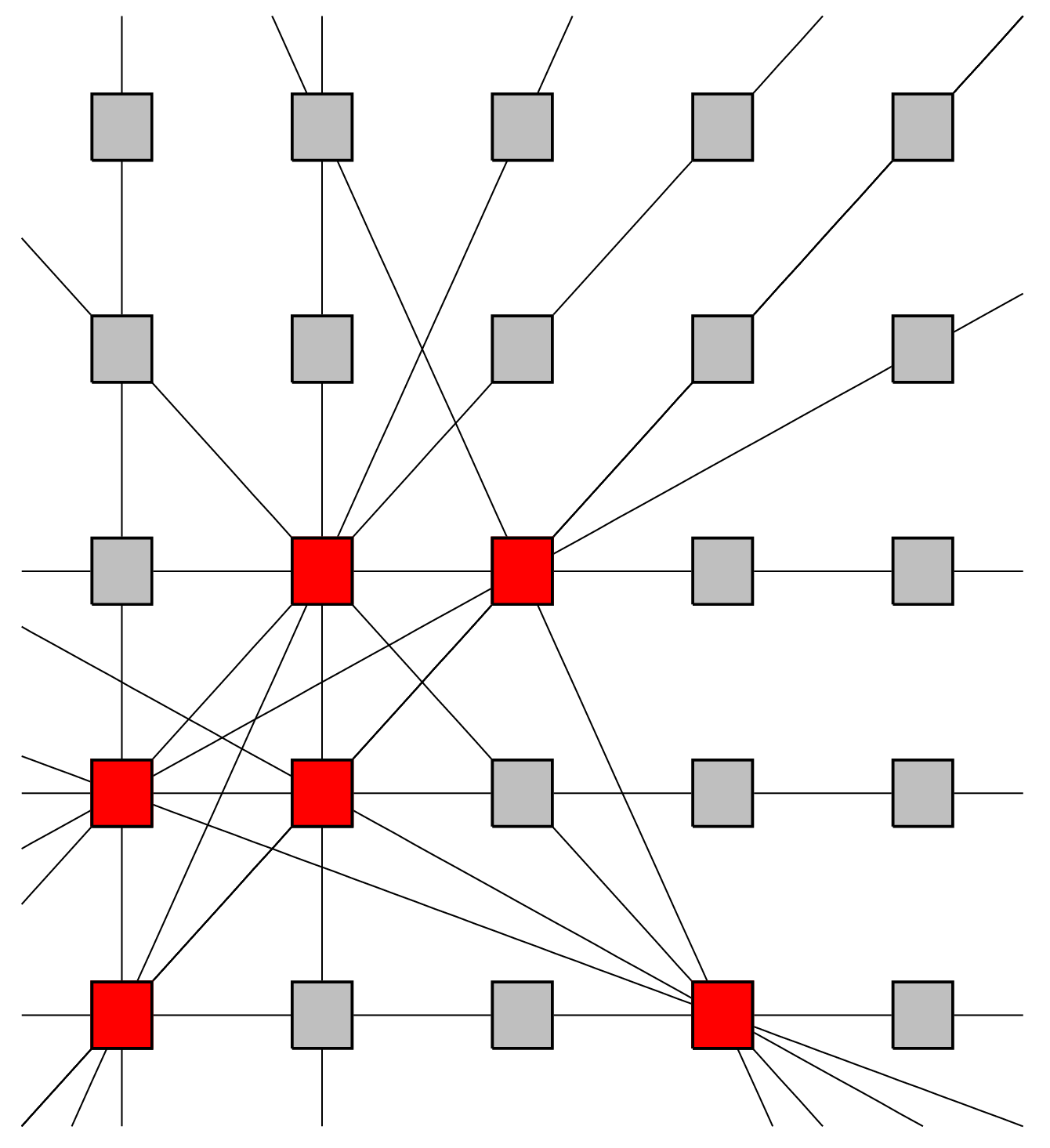}
\includegraphics[scale=0.3]{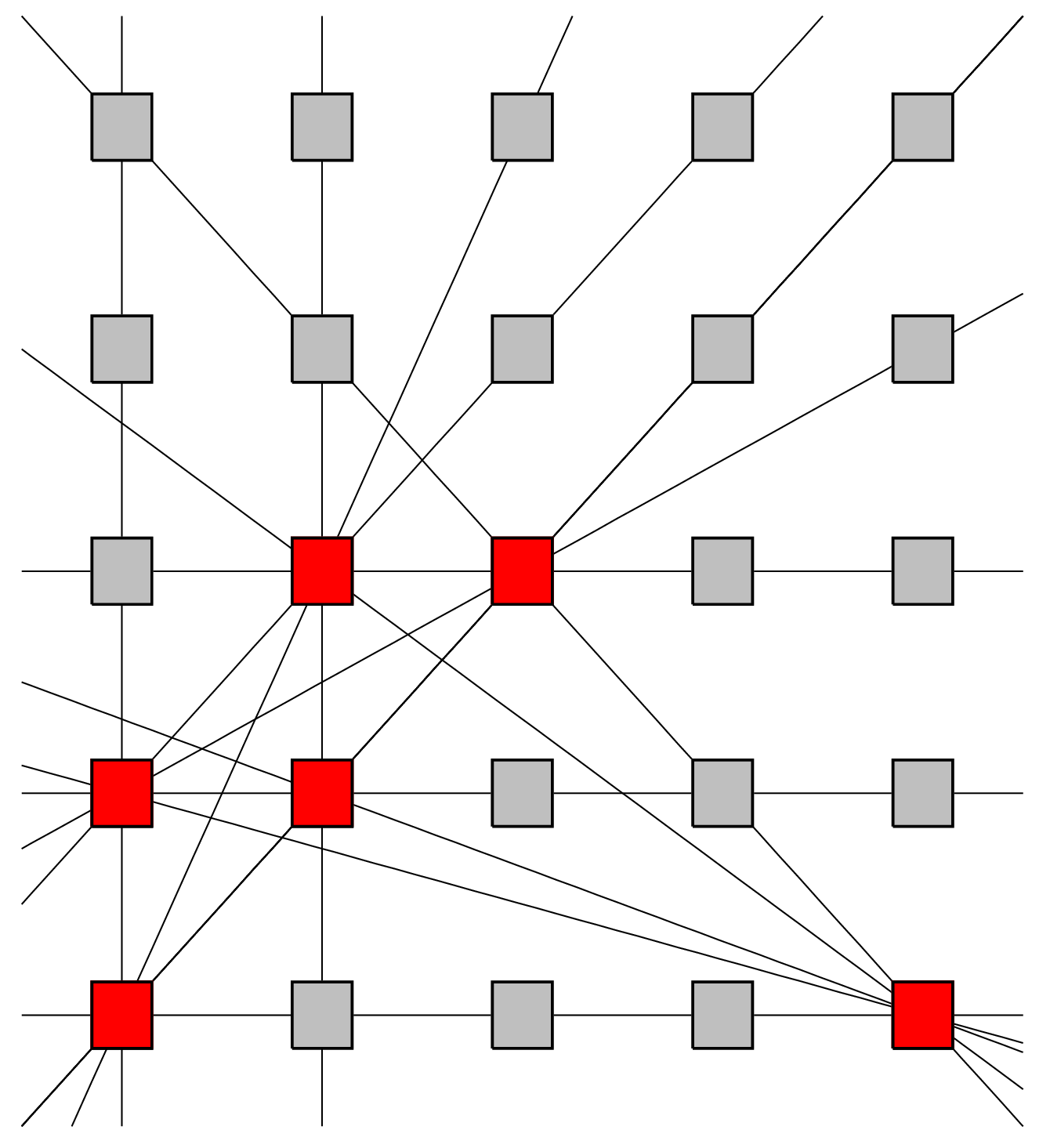}
\caption{$t(4)=6$.}
\label{fig.n4_1}
\end{figure}

\begin{figure}
\includegraphics[scale=0.3]{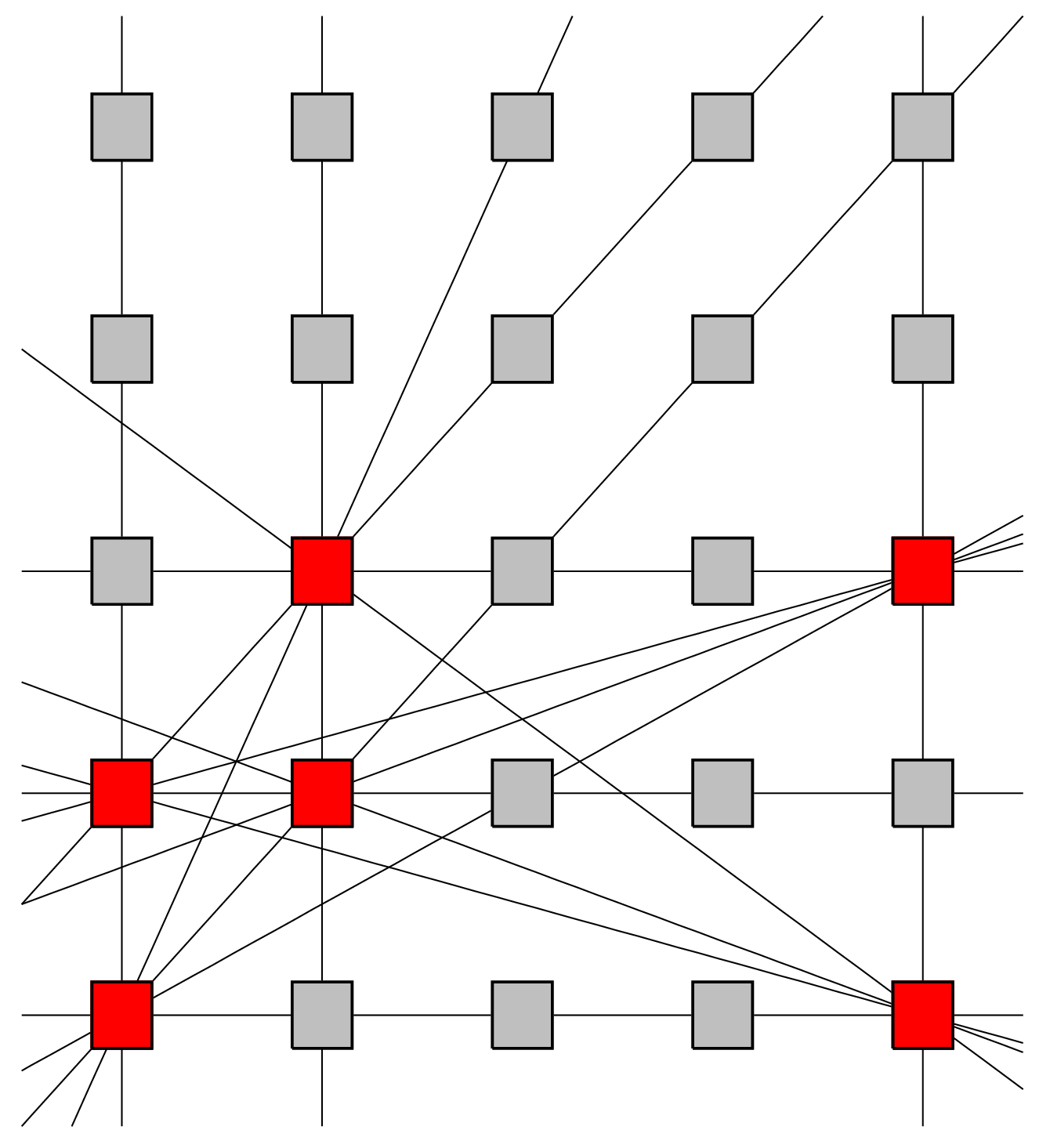}
\includegraphics[scale=0.3]{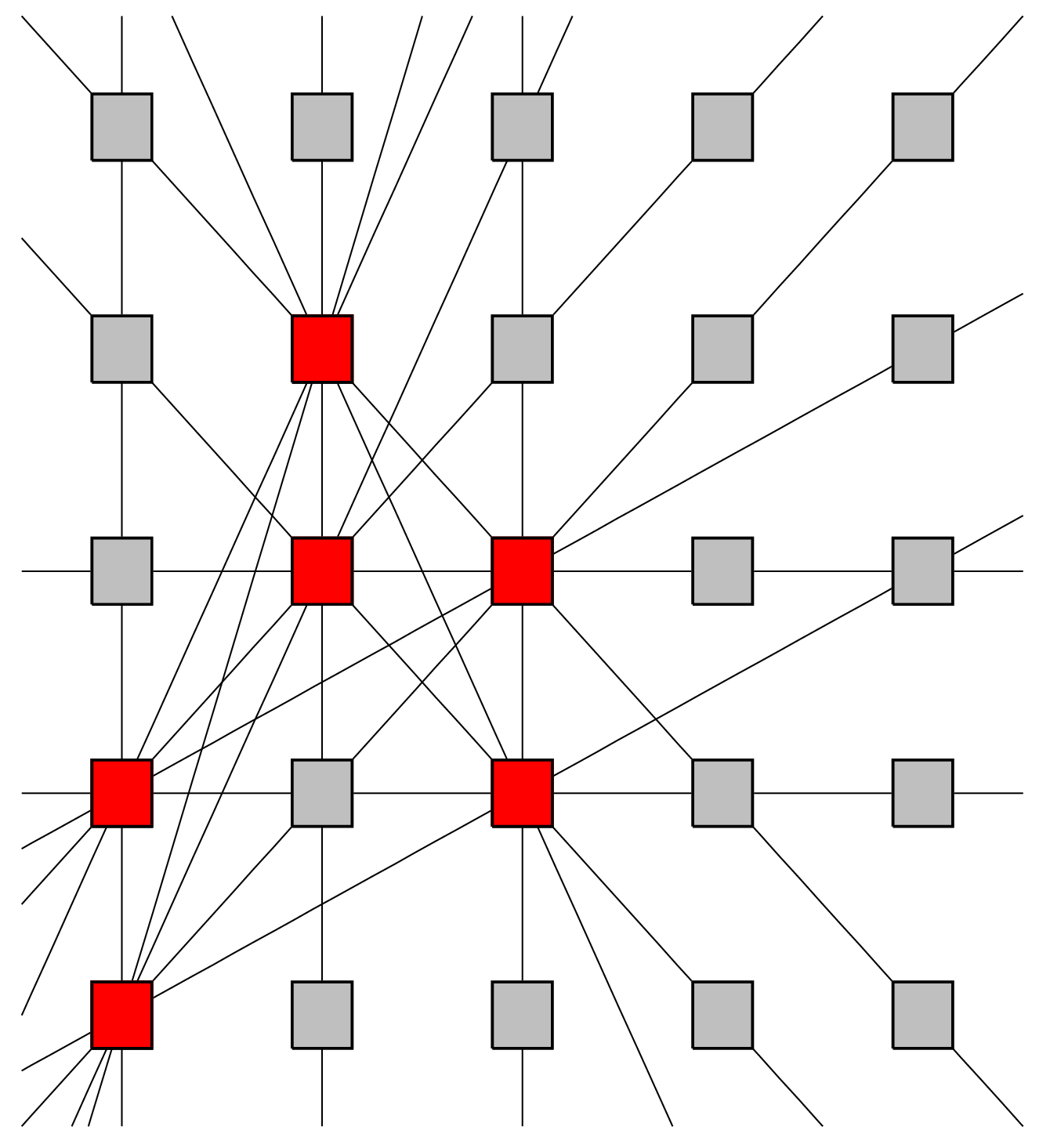}
\includegraphics[scale=0.3]{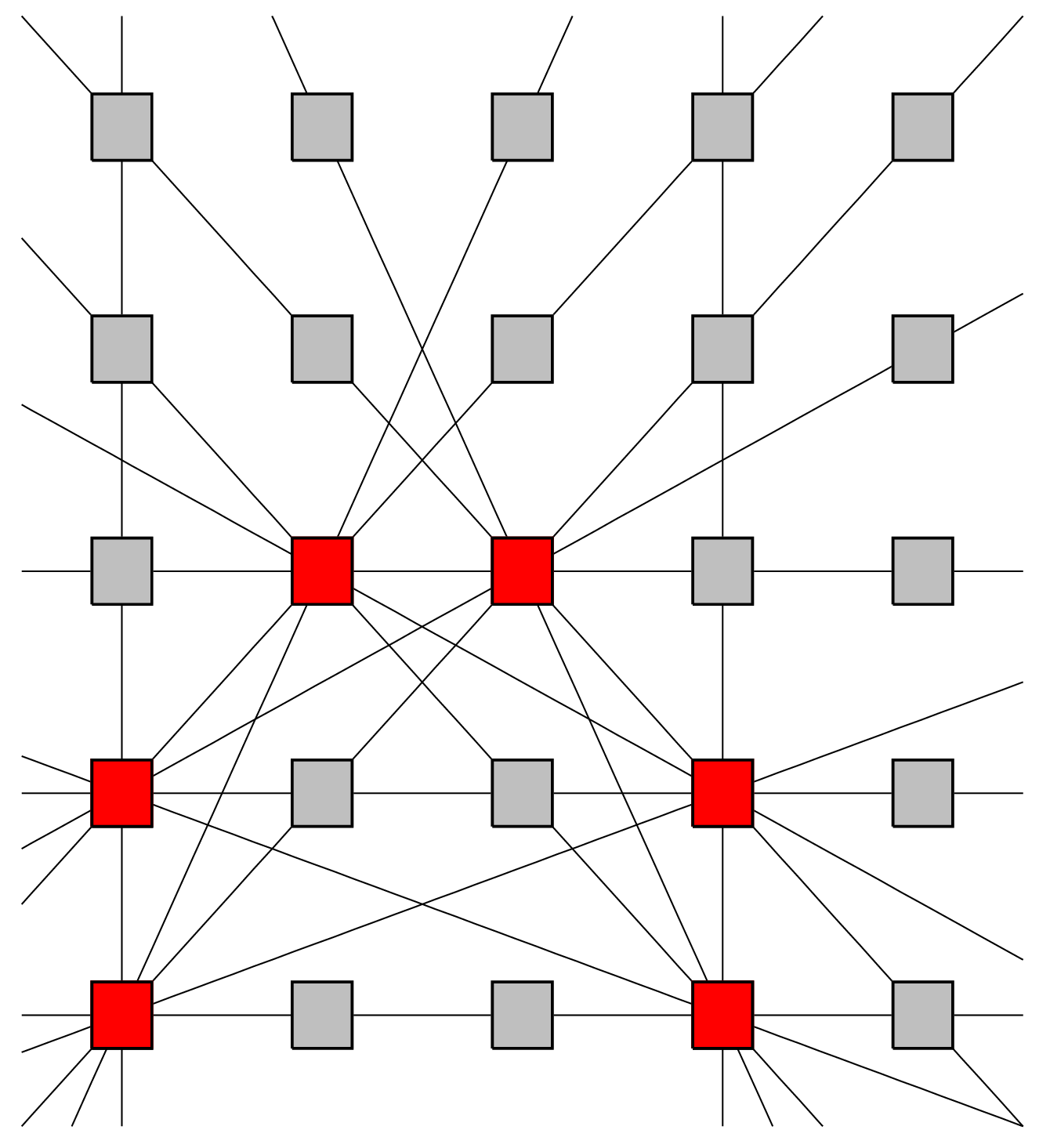}
\includegraphics[scale=0.3]{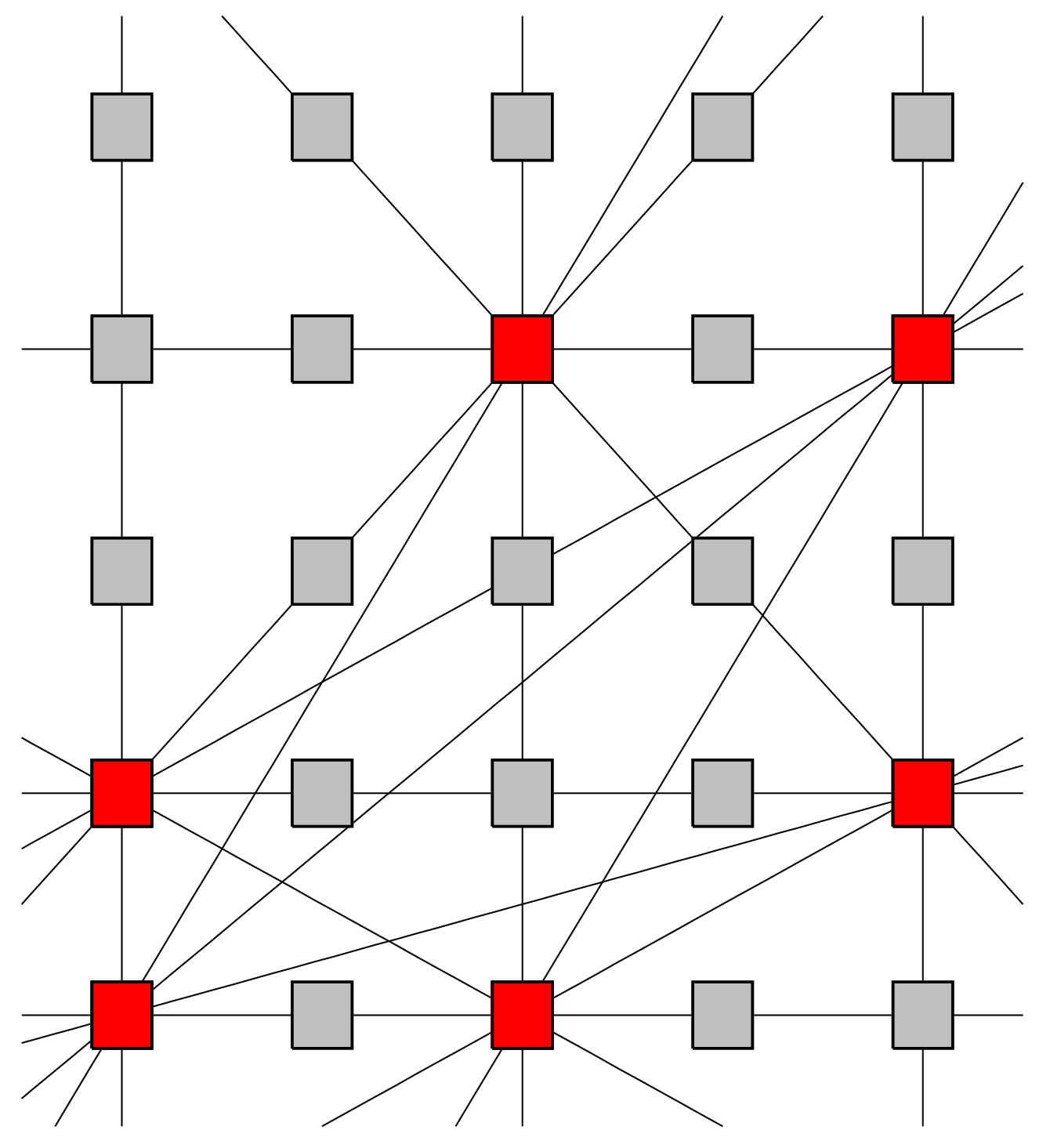}
\includegraphics[scale=0.3]{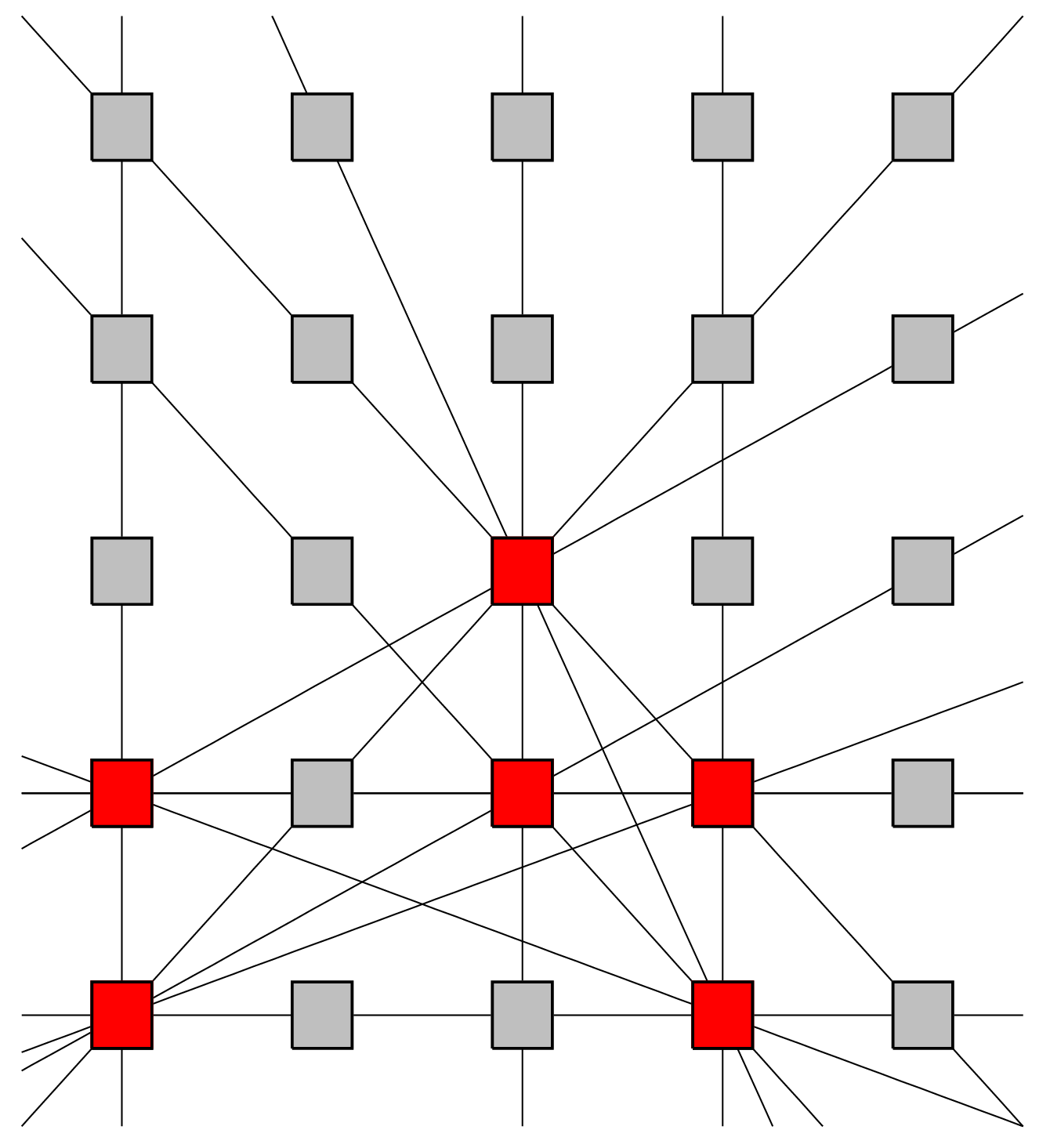}
\includegraphics[scale=0.3]{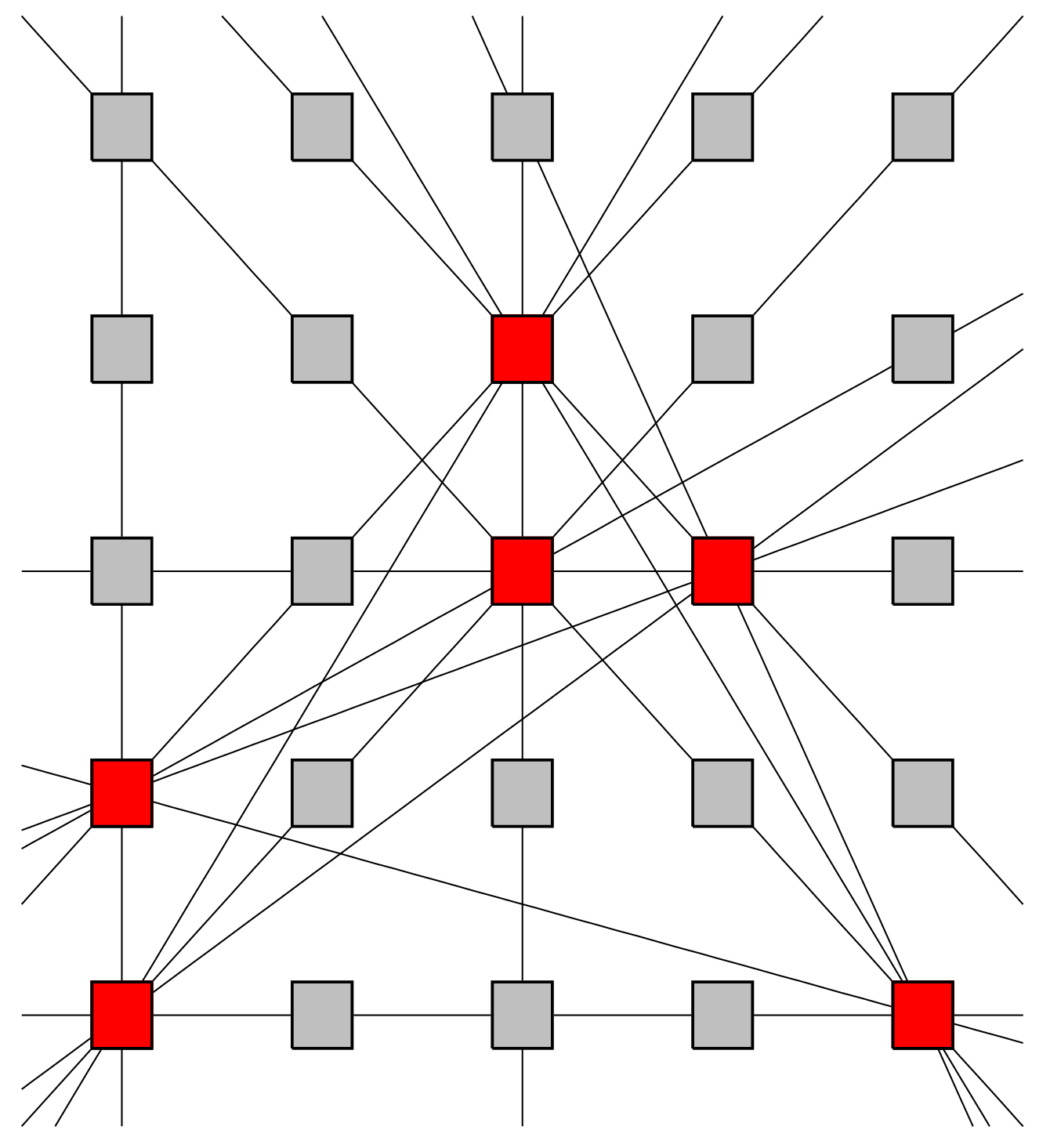}
\includegraphics[scale=0.3]{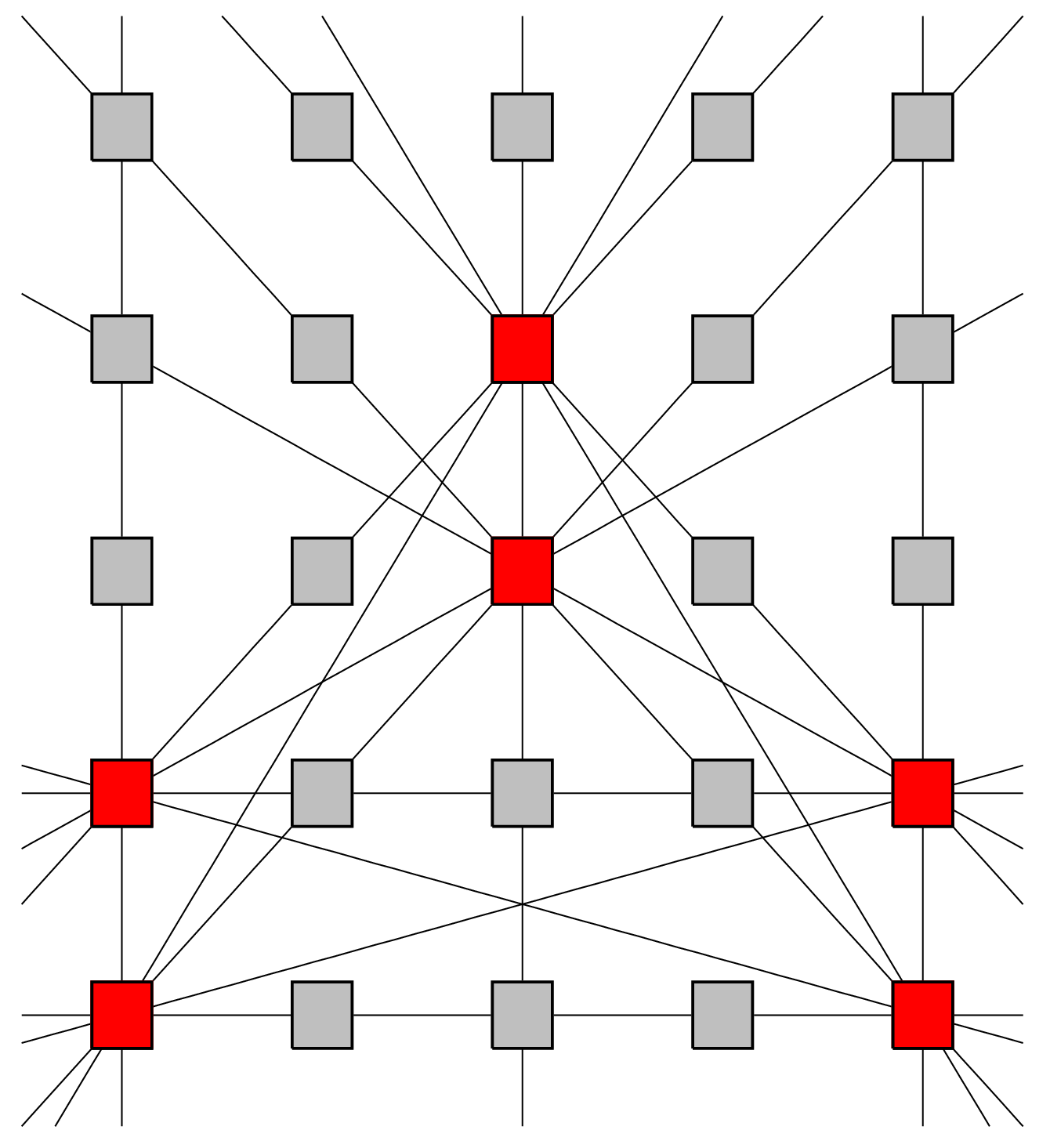}
\includegraphics[scale=0.3]{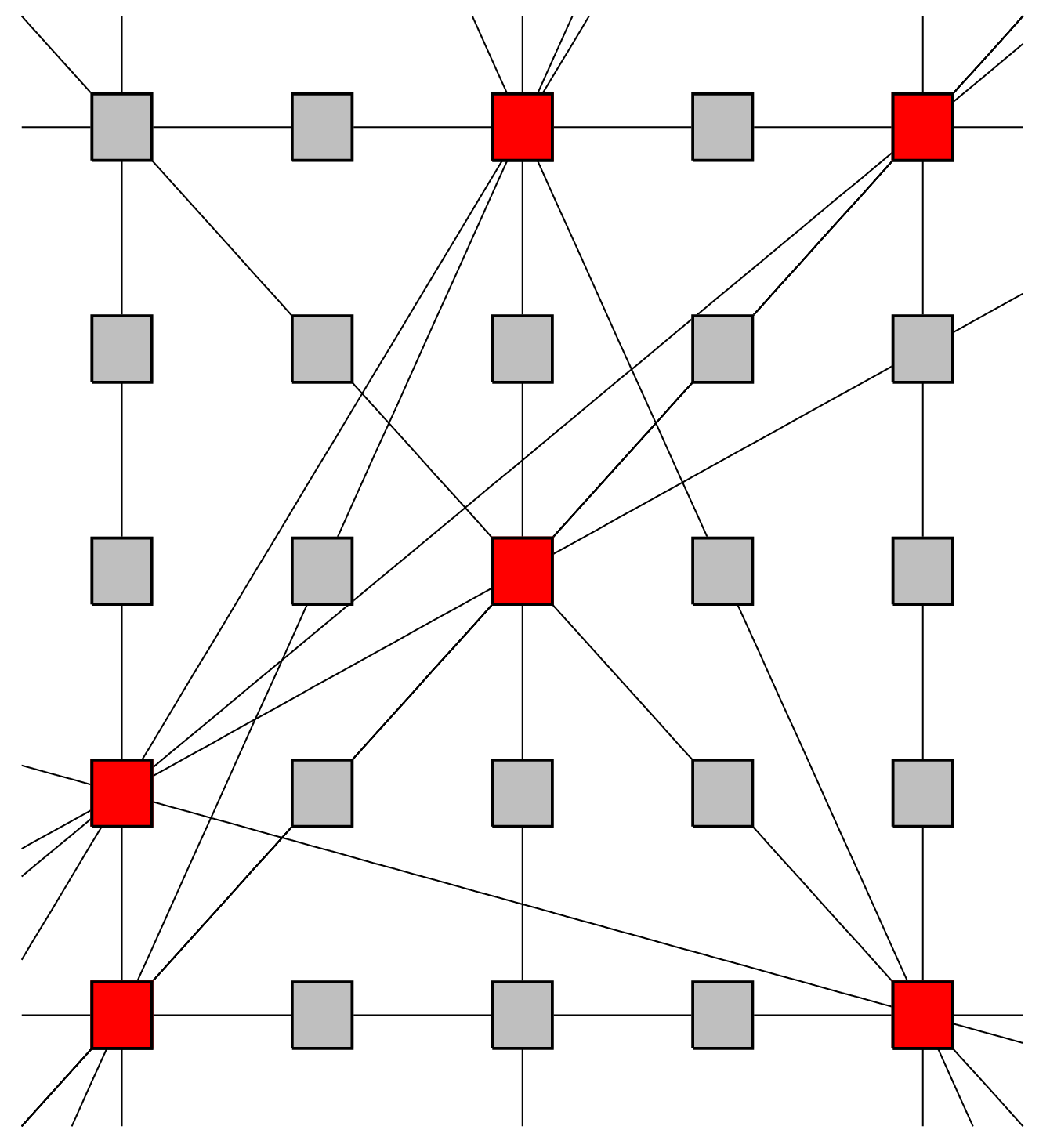}
\includegraphics[scale=0.3]{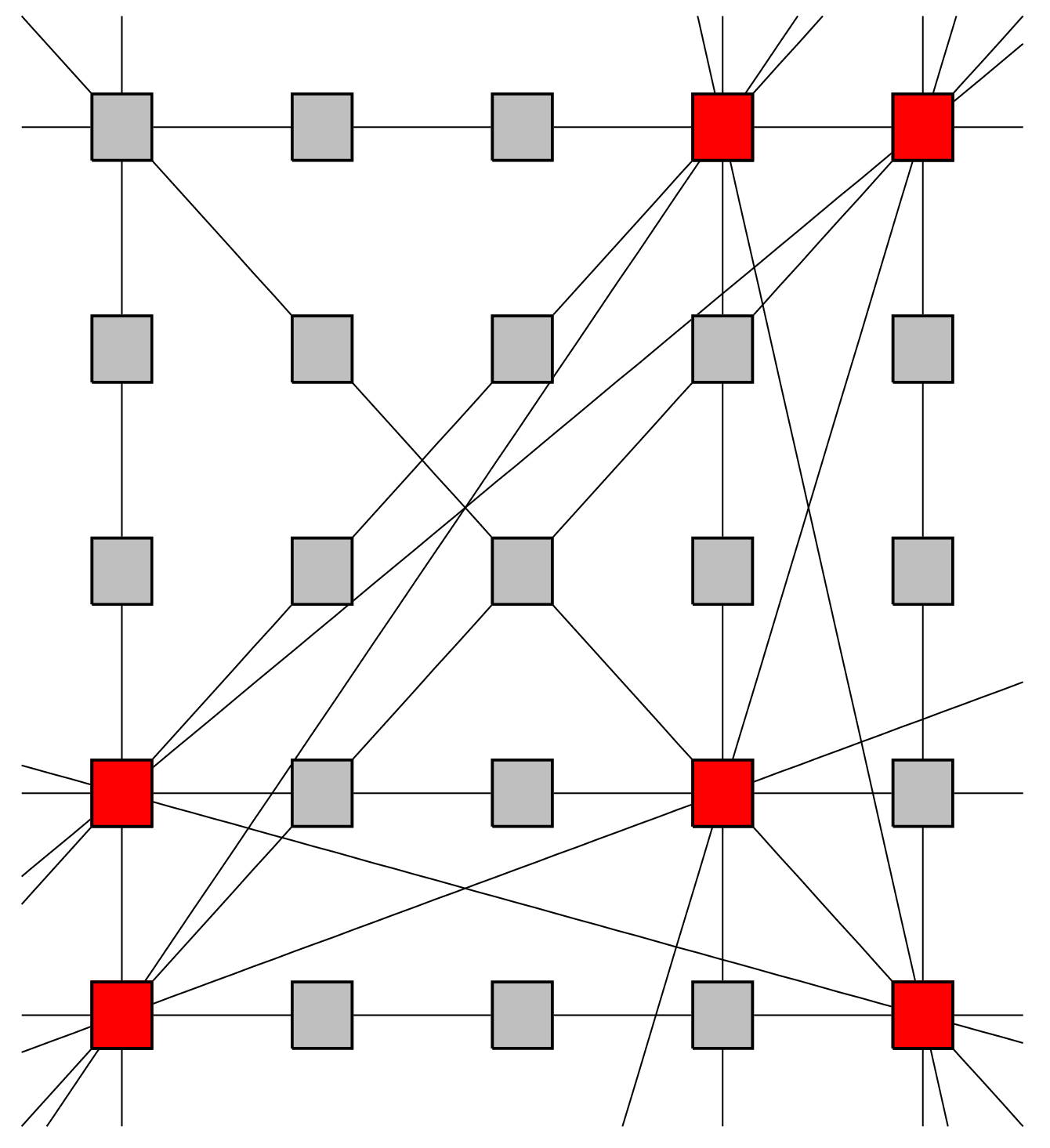}
\includegraphics[scale=0.3]{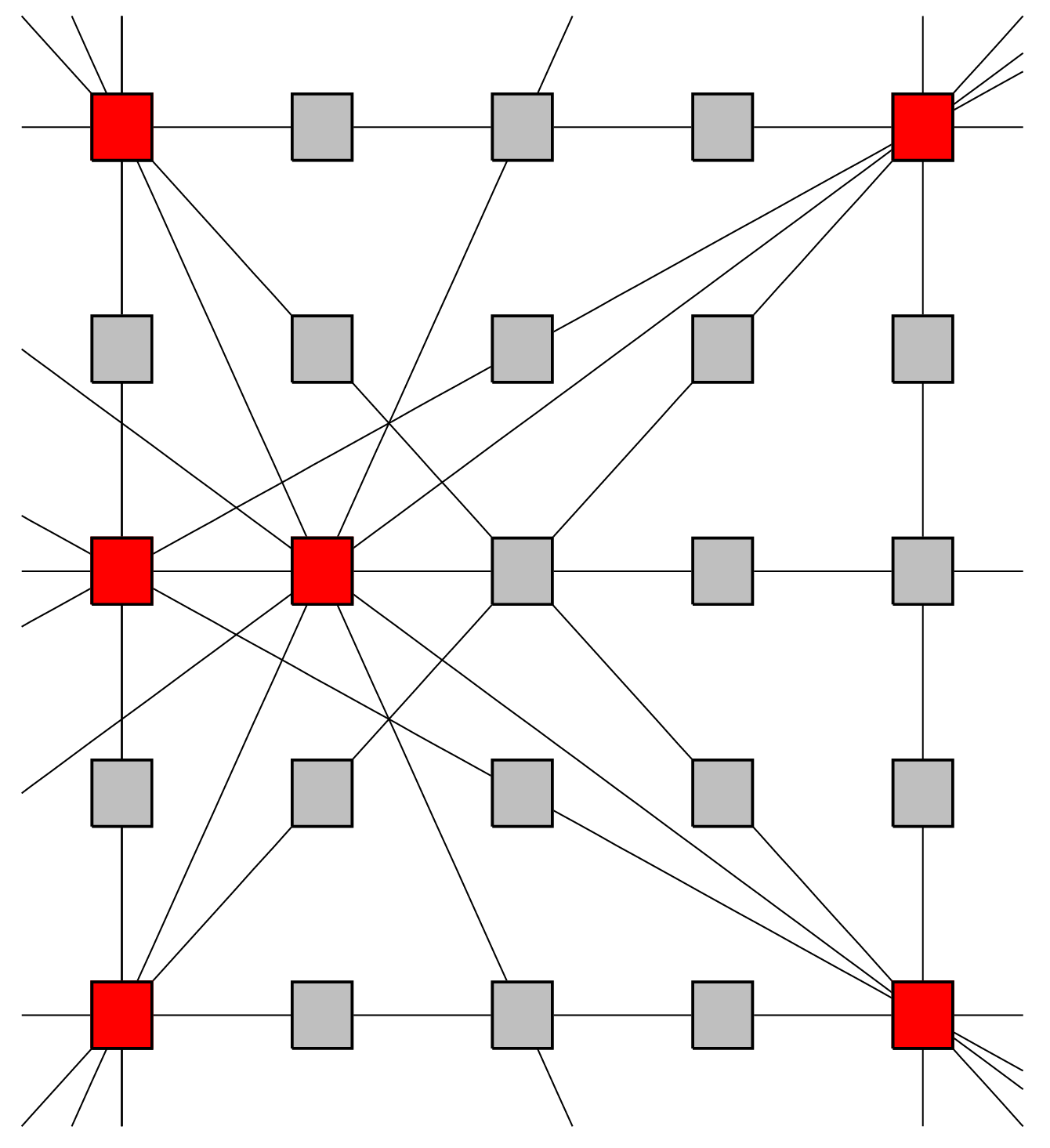}
\includegraphics[scale=0.3]{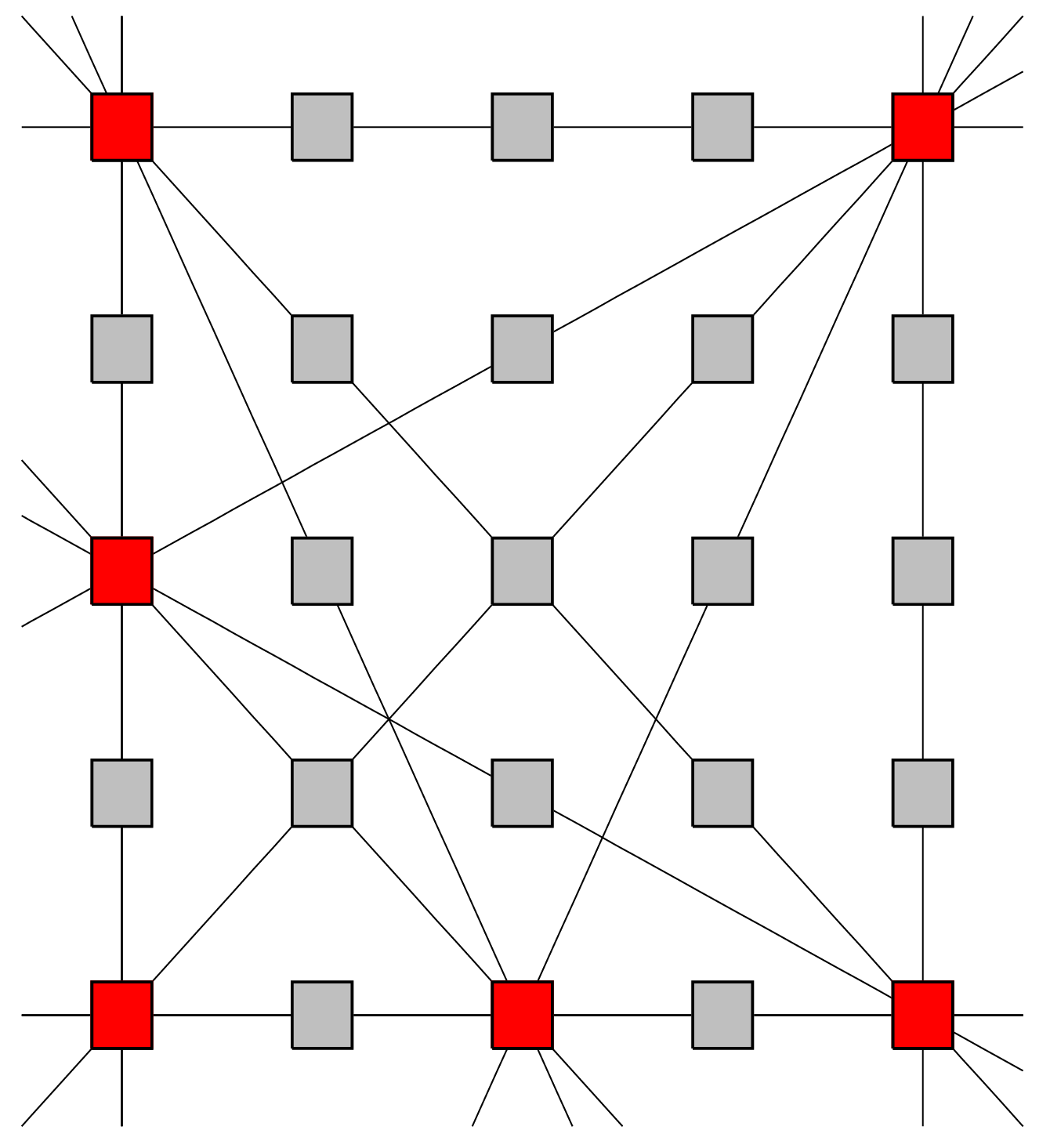}
\includegraphics[scale=0.3]{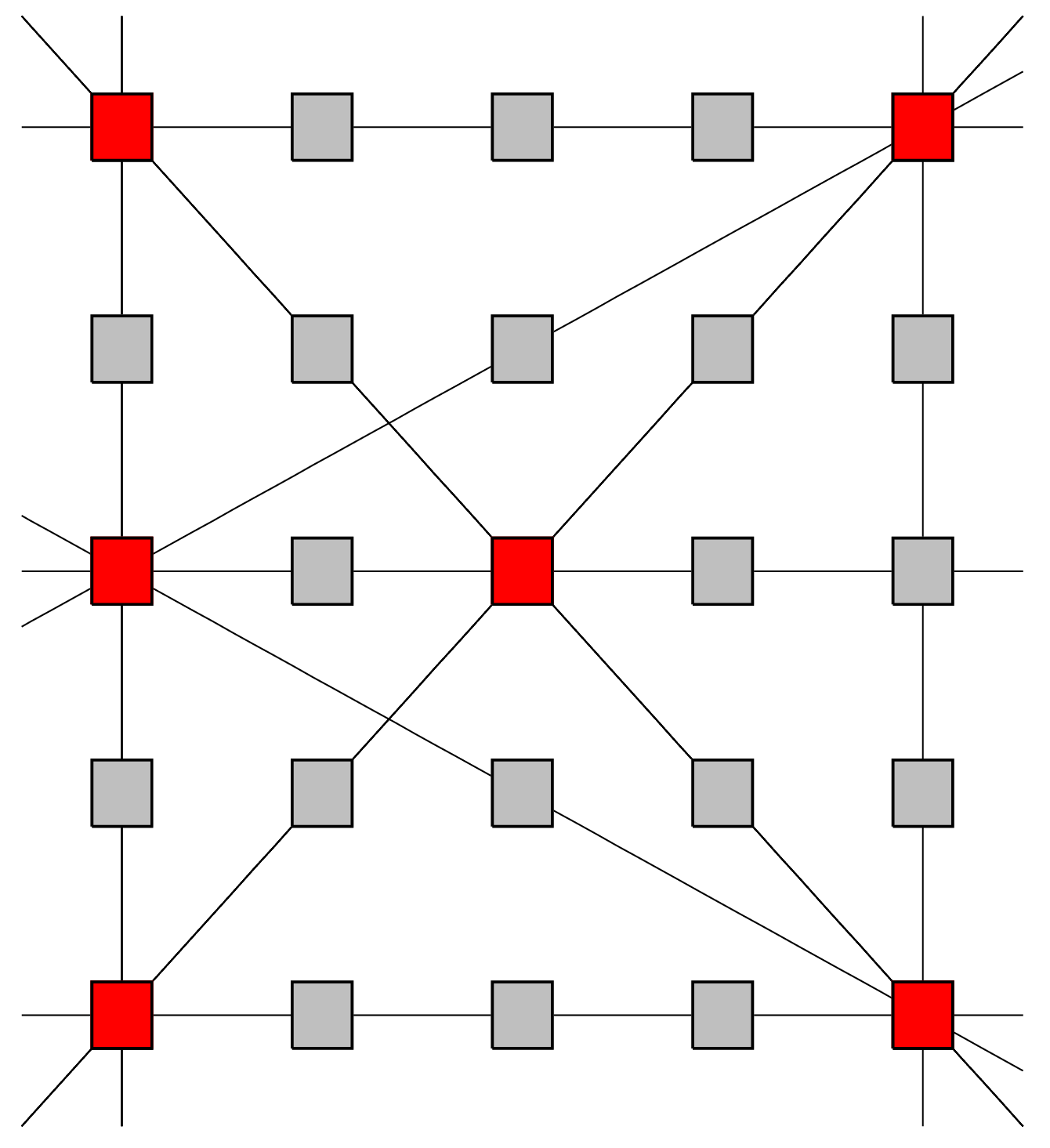}
\caption{$t(4)=6$ (continued).}
\label{fig.n4_2}
\end{figure}

\begin{figure}
\includegraphics[scale=0.3]{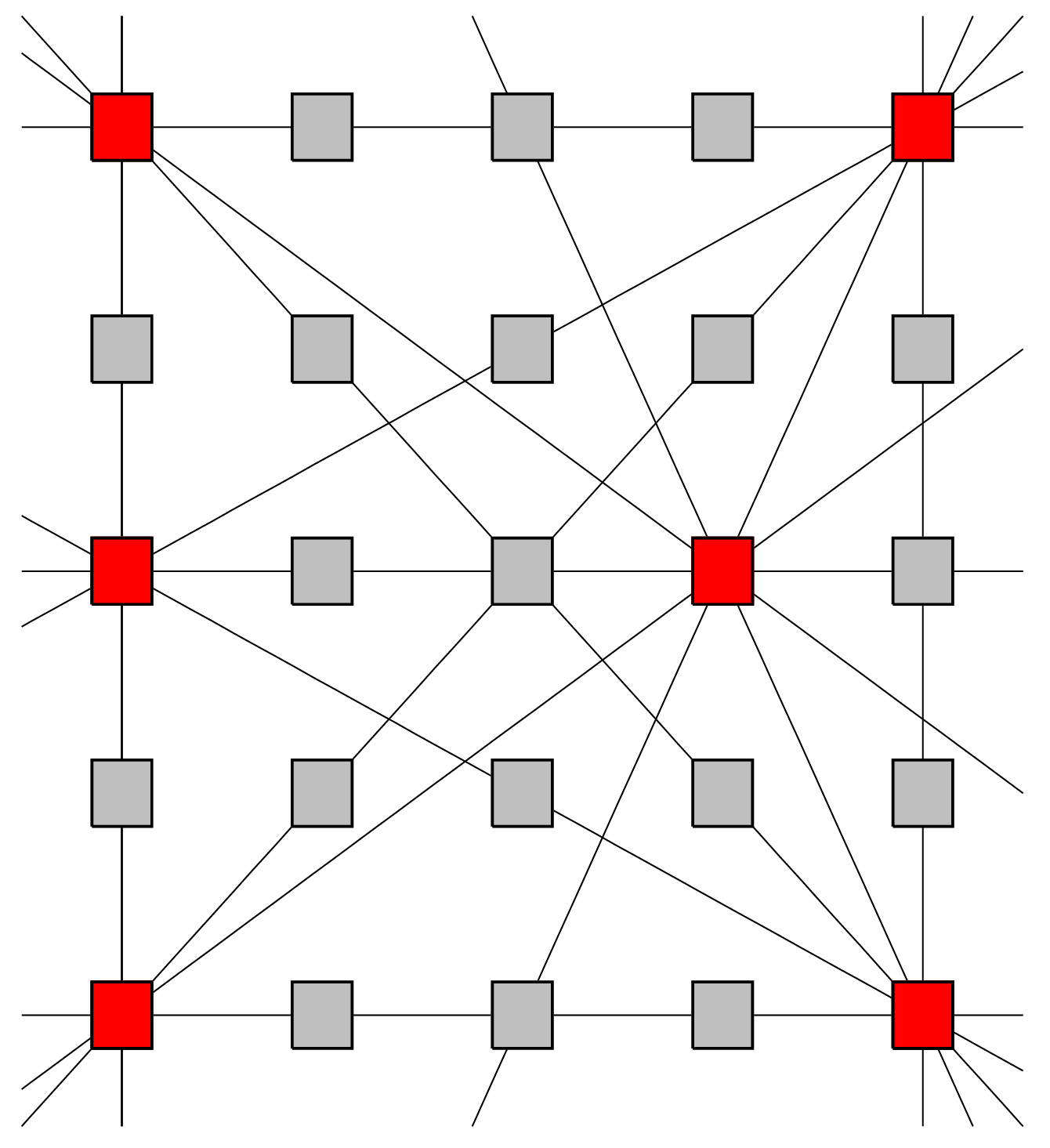}
\includegraphics[scale=0.3]{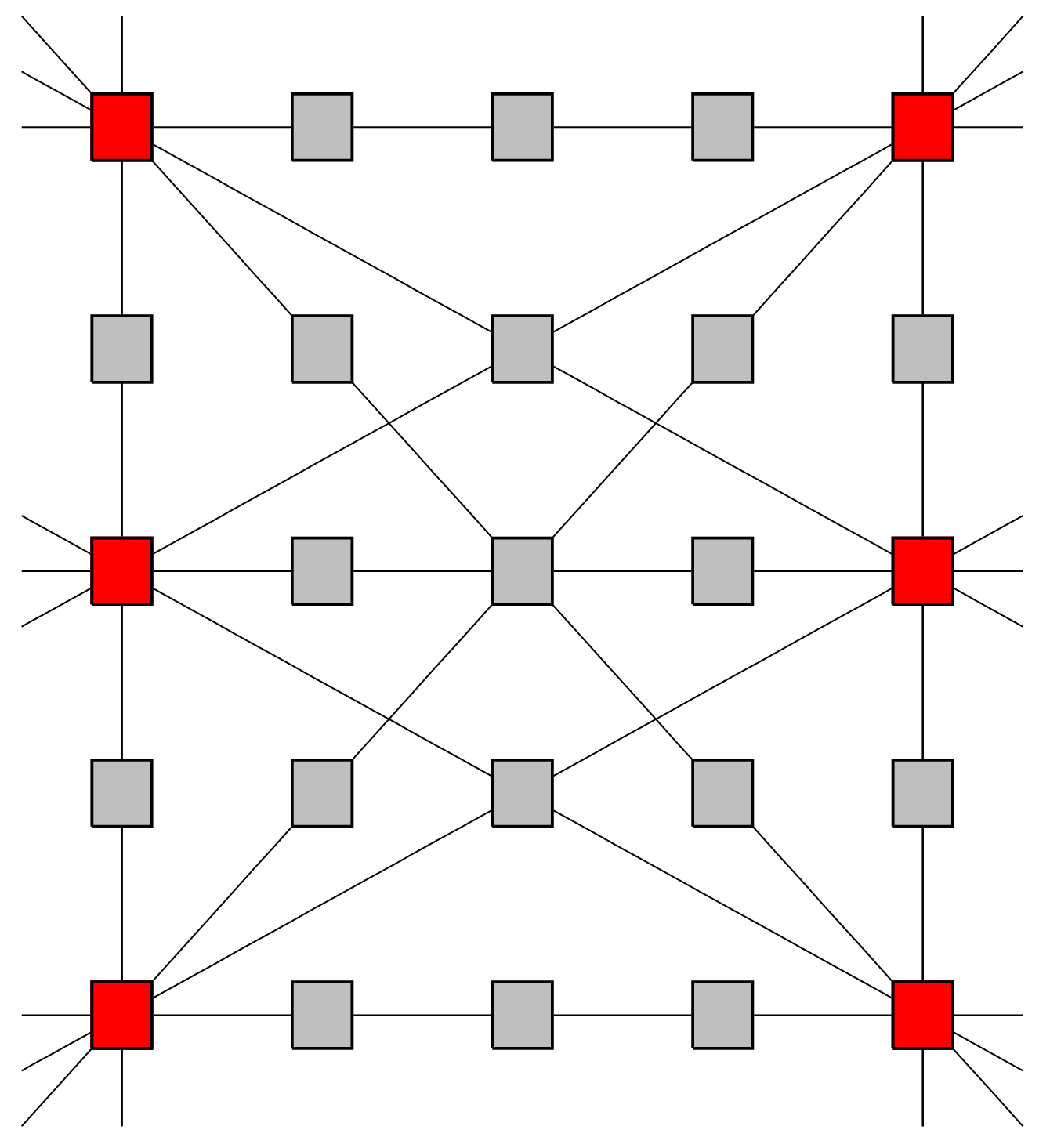}
\includegraphics[scale=0.3]{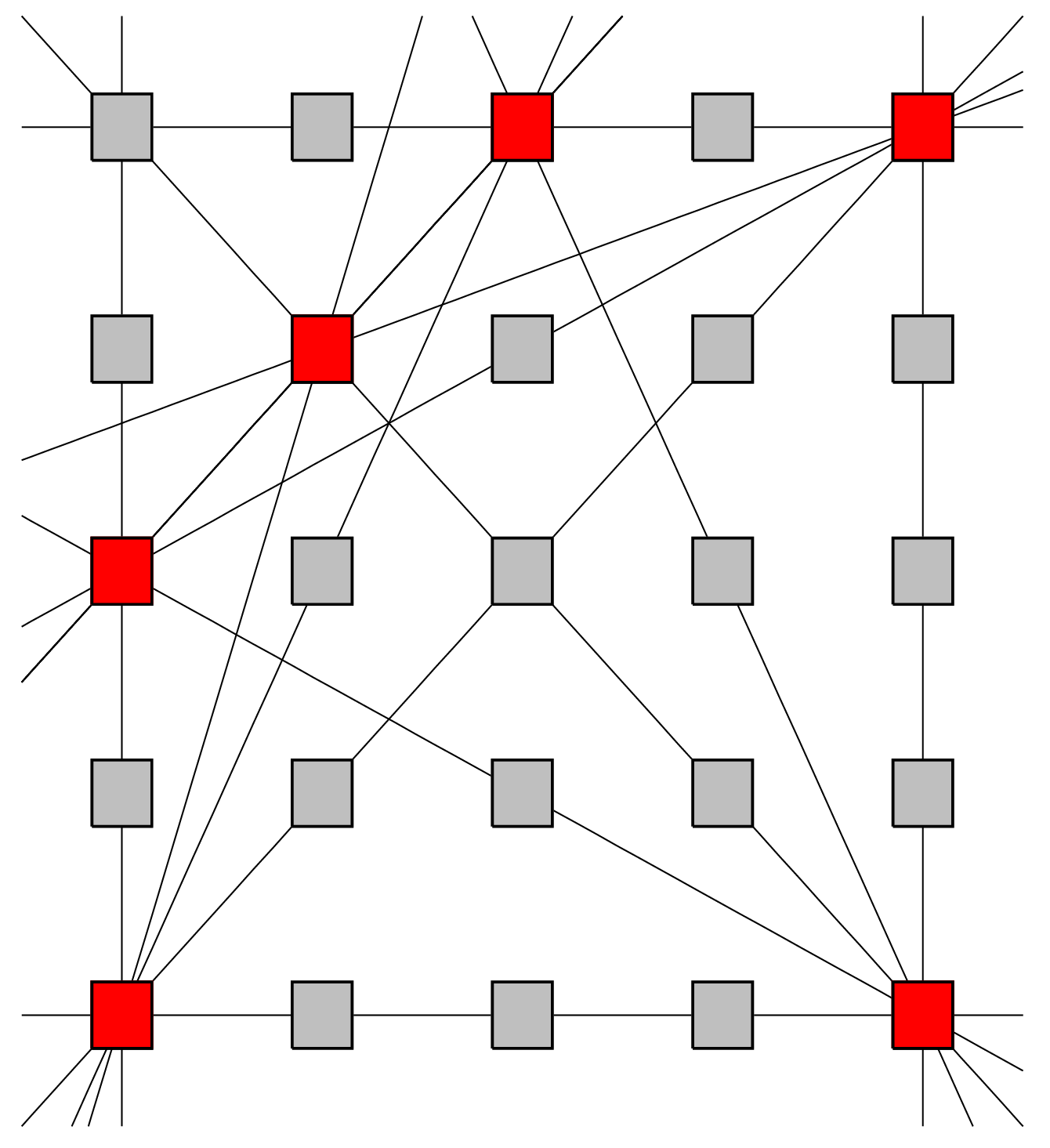}
\includegraphics[scale=0.3]{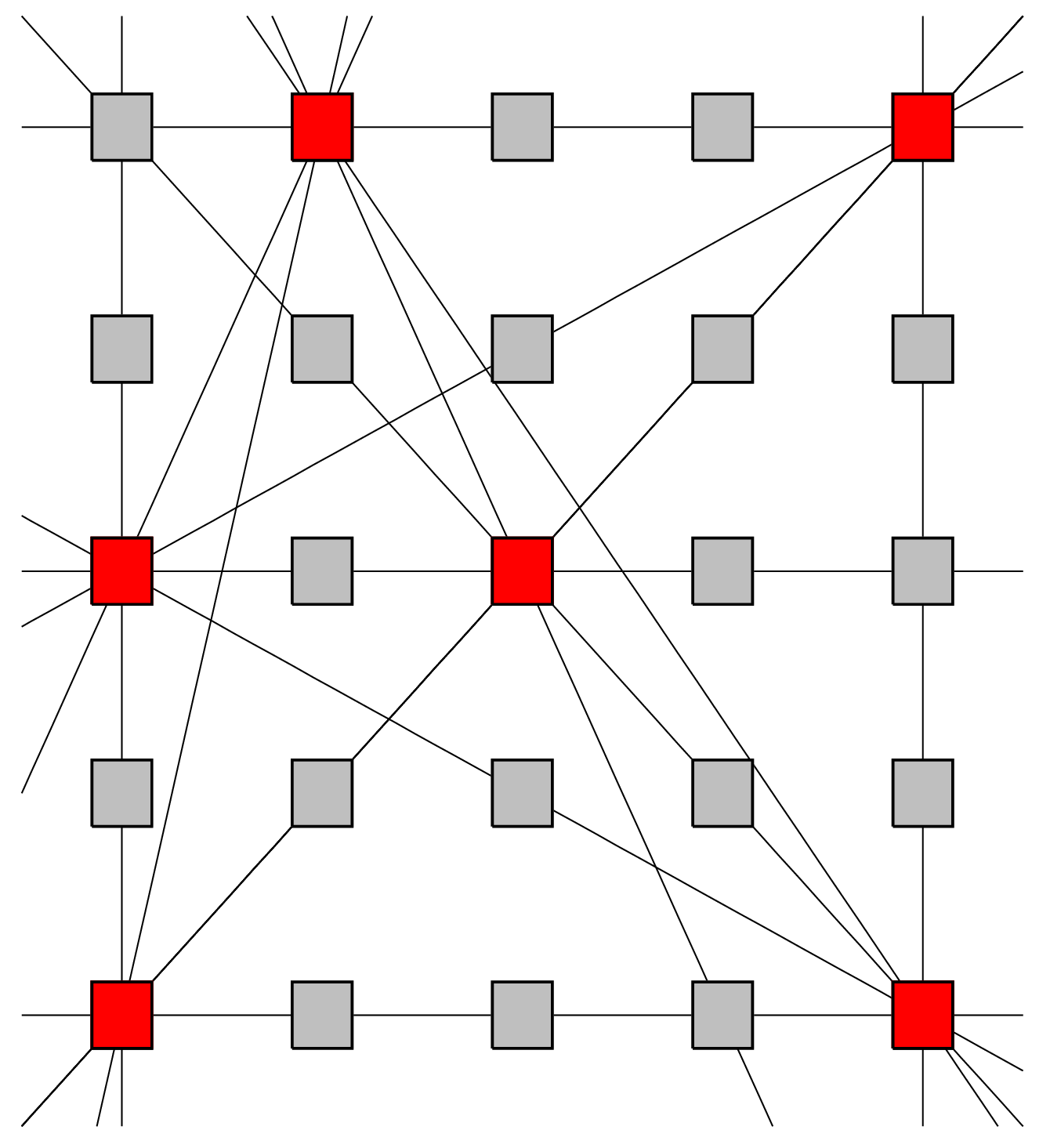}
\includegraphics[scale=0.3]{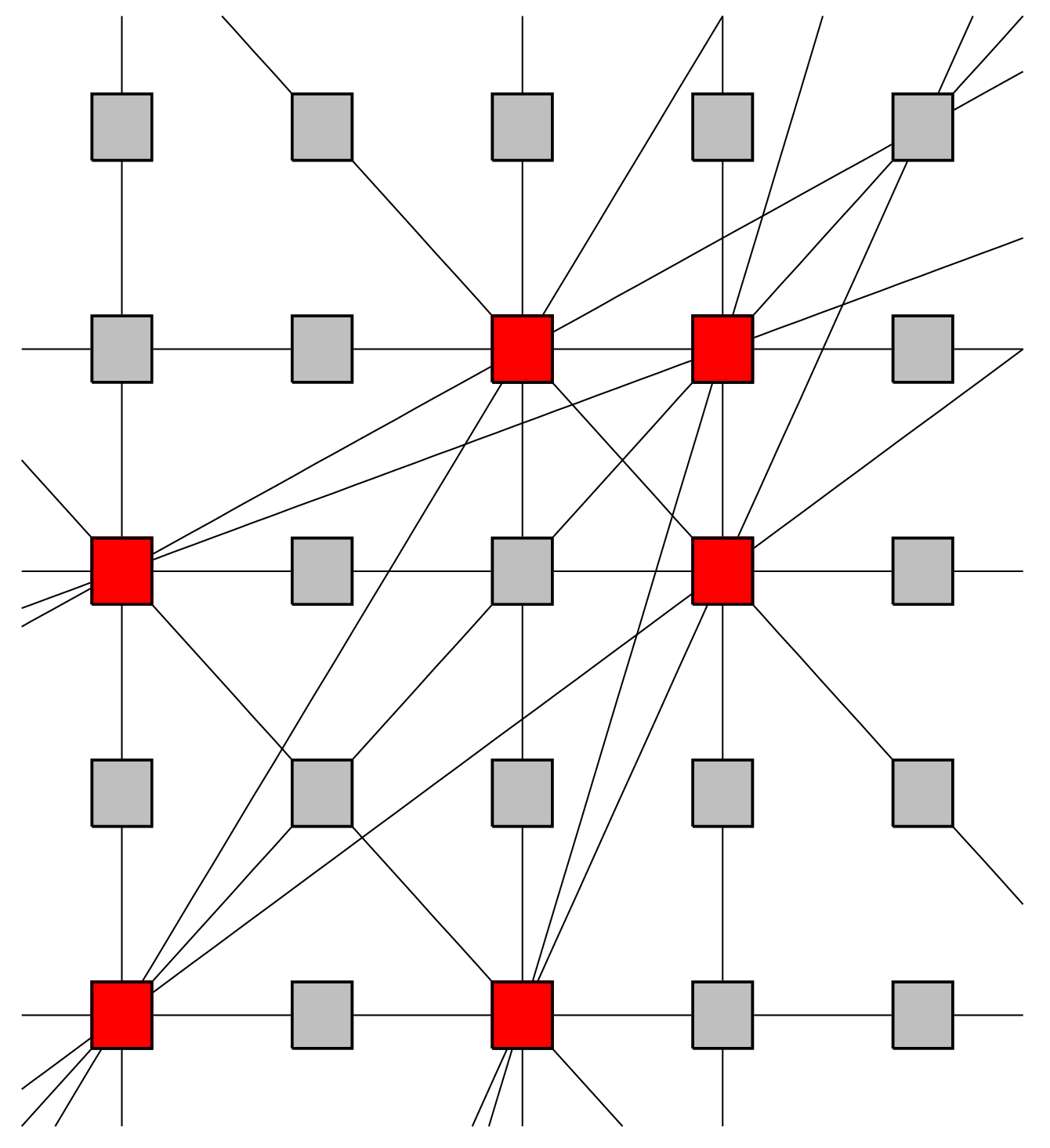}
\includegraphics[scale=0.3]{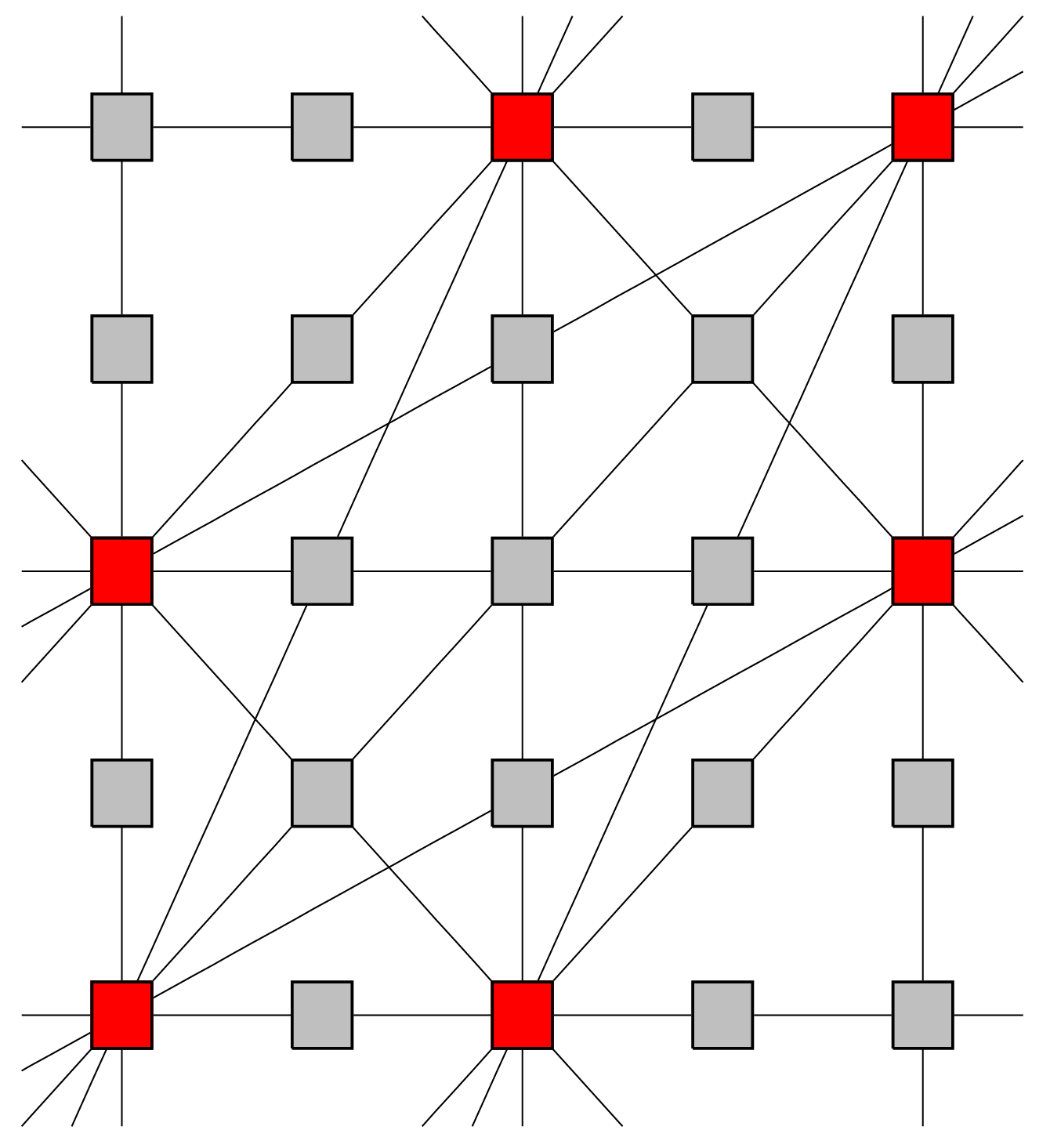}
\includegraphics[scale=0.3]{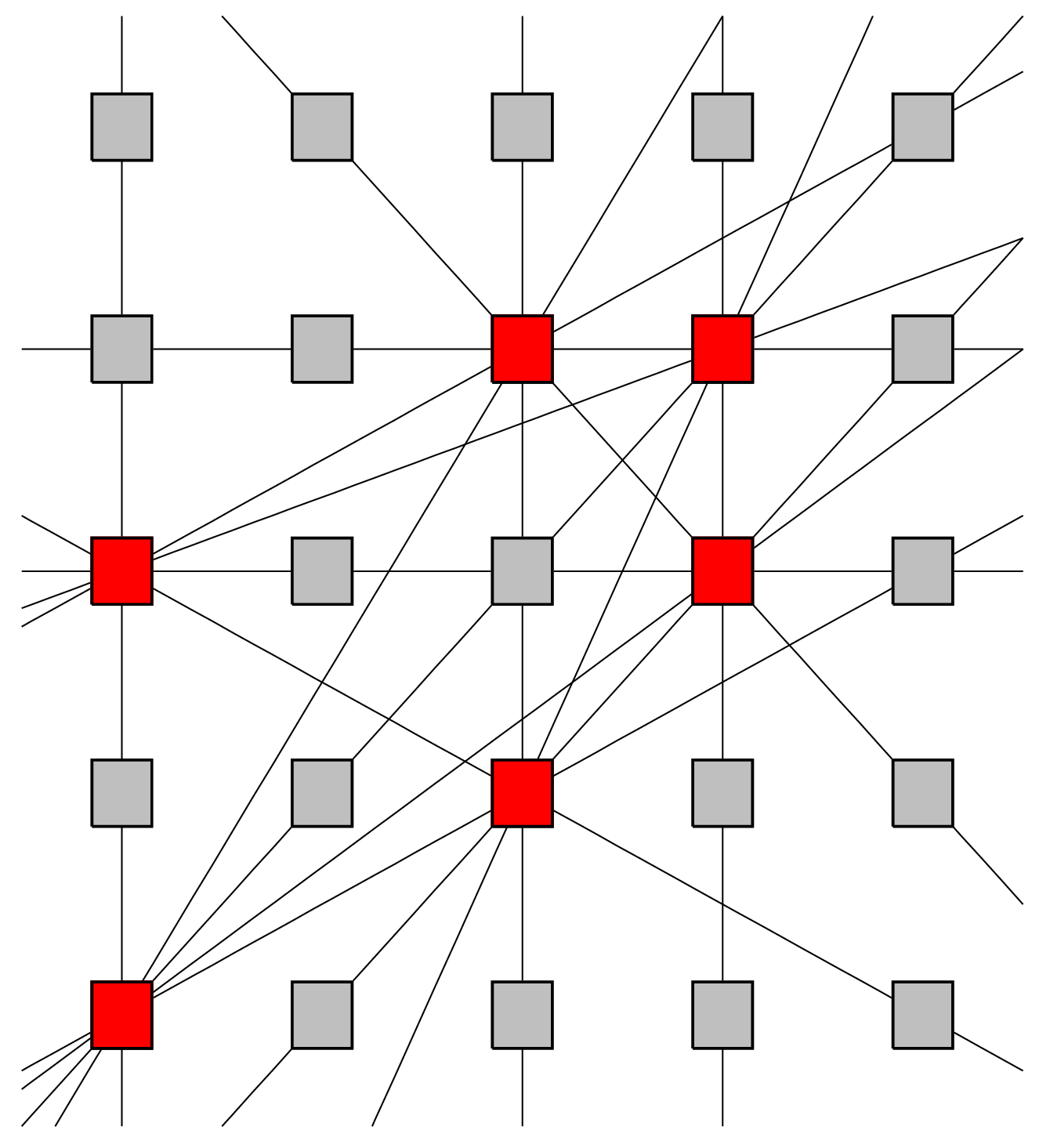}
\includegraphics[scale=0.3]{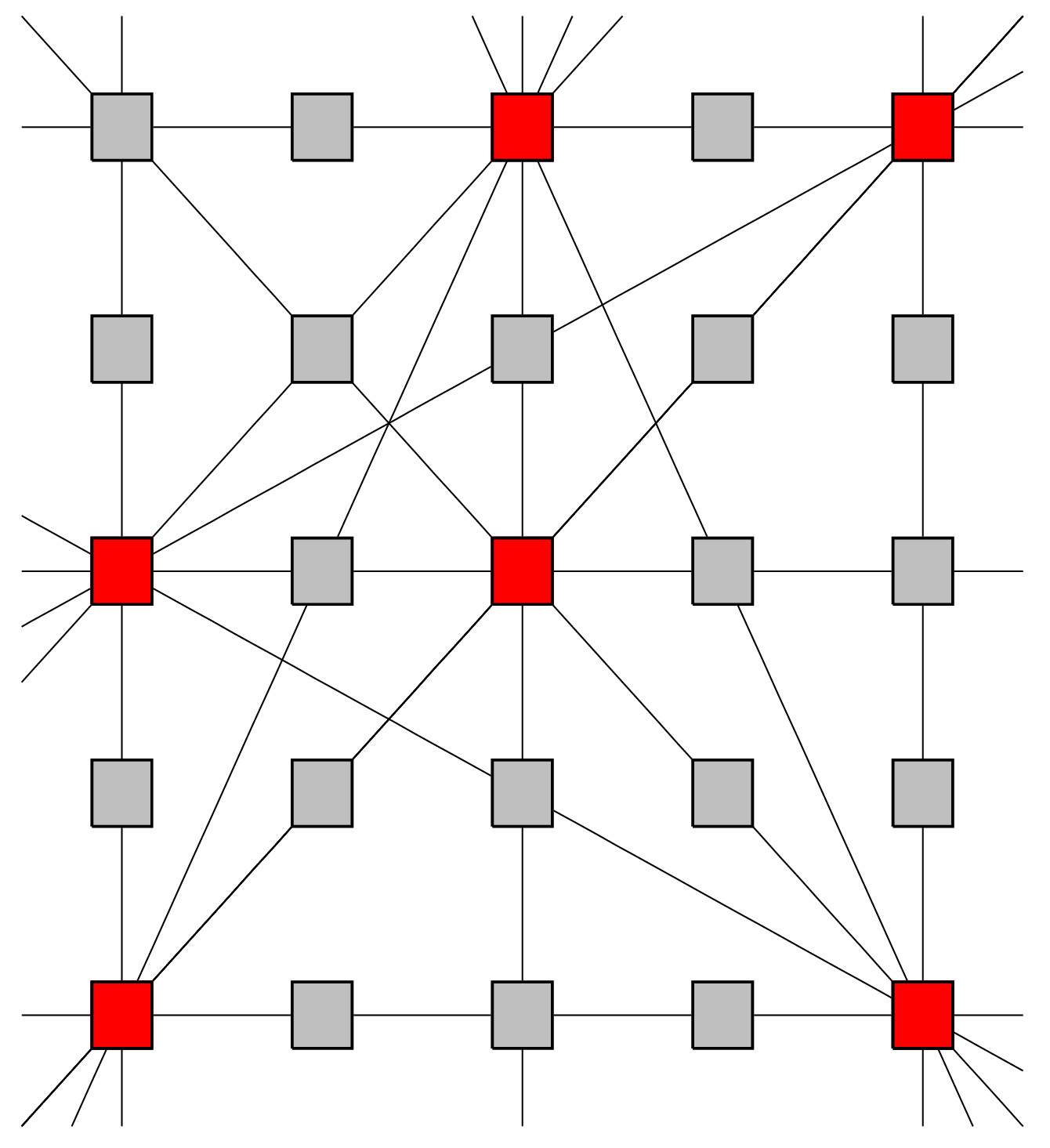}
\includegraphics[scale=0.3]{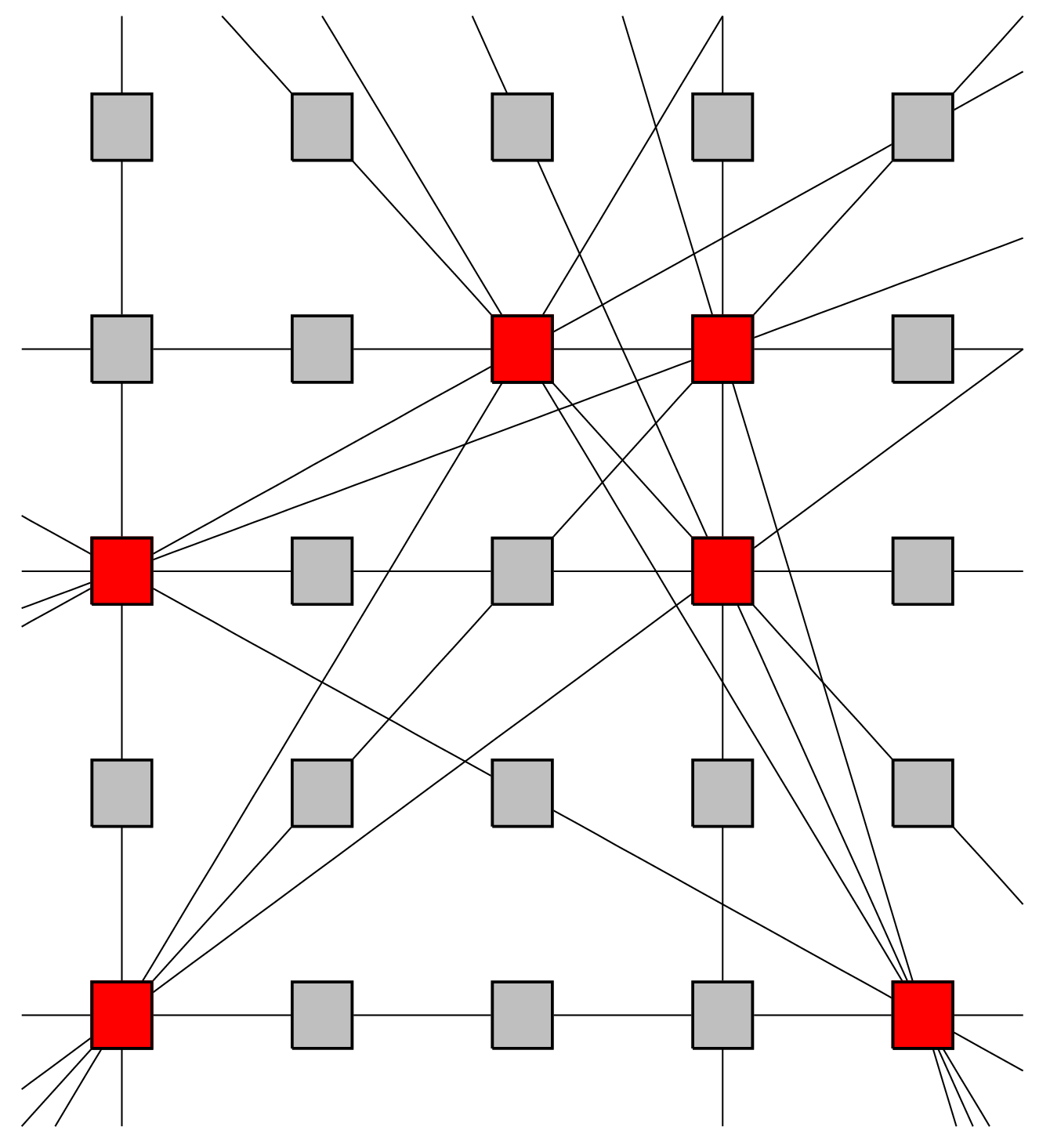}
\includegraphics[scale=0.3]{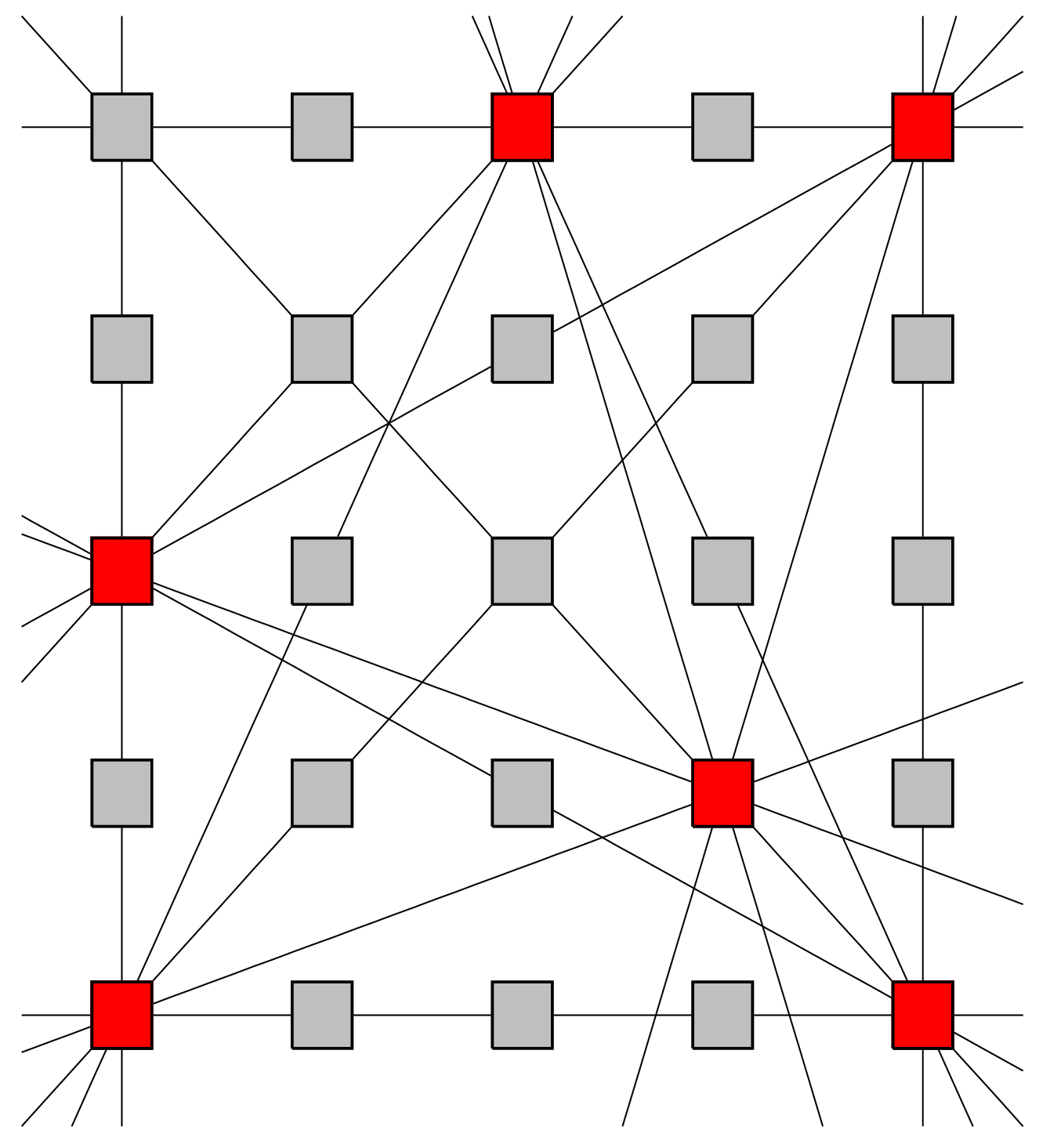}
\includegraphics[scale=0.3]{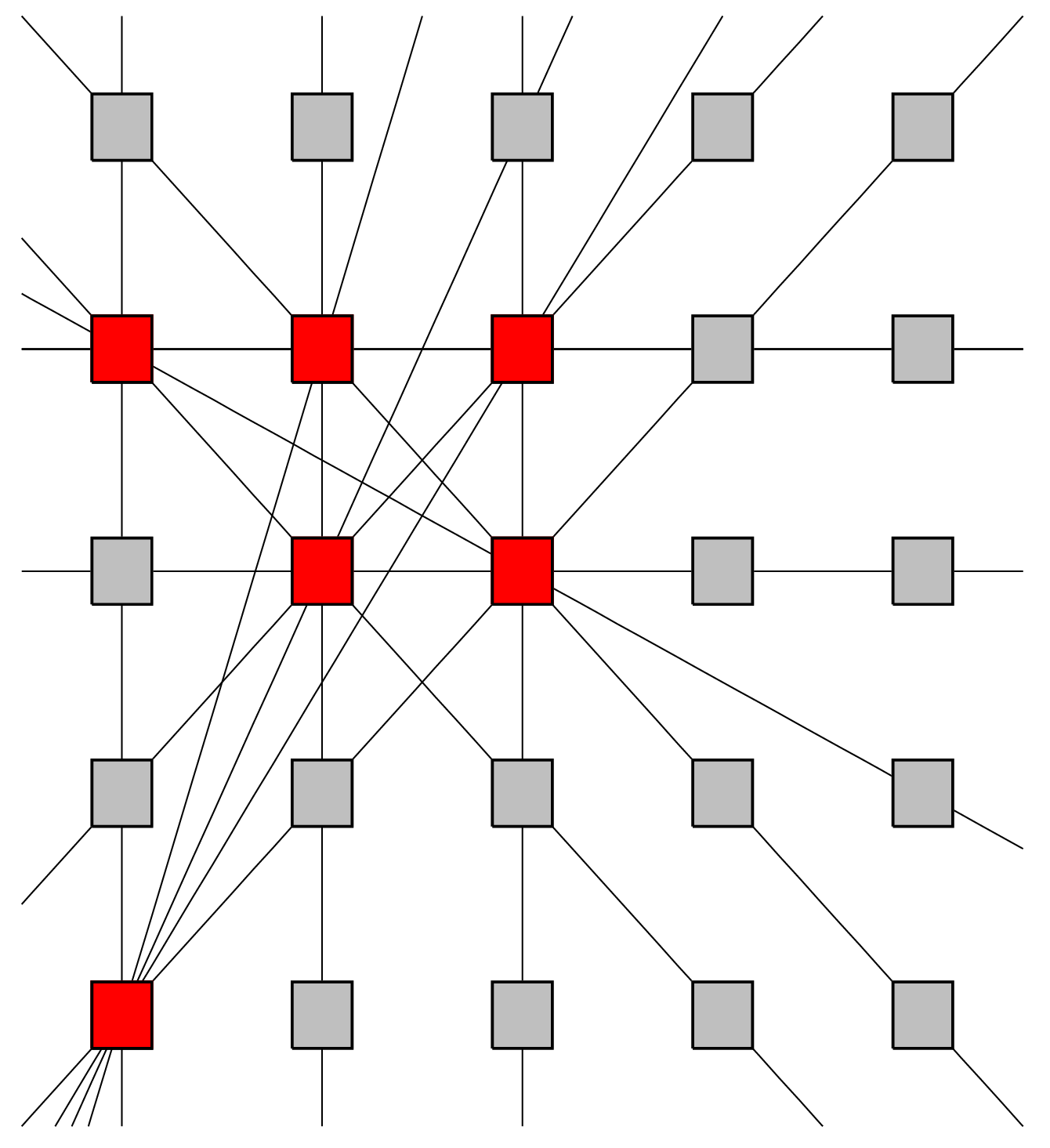}
\includegraphics[scale=0.3]{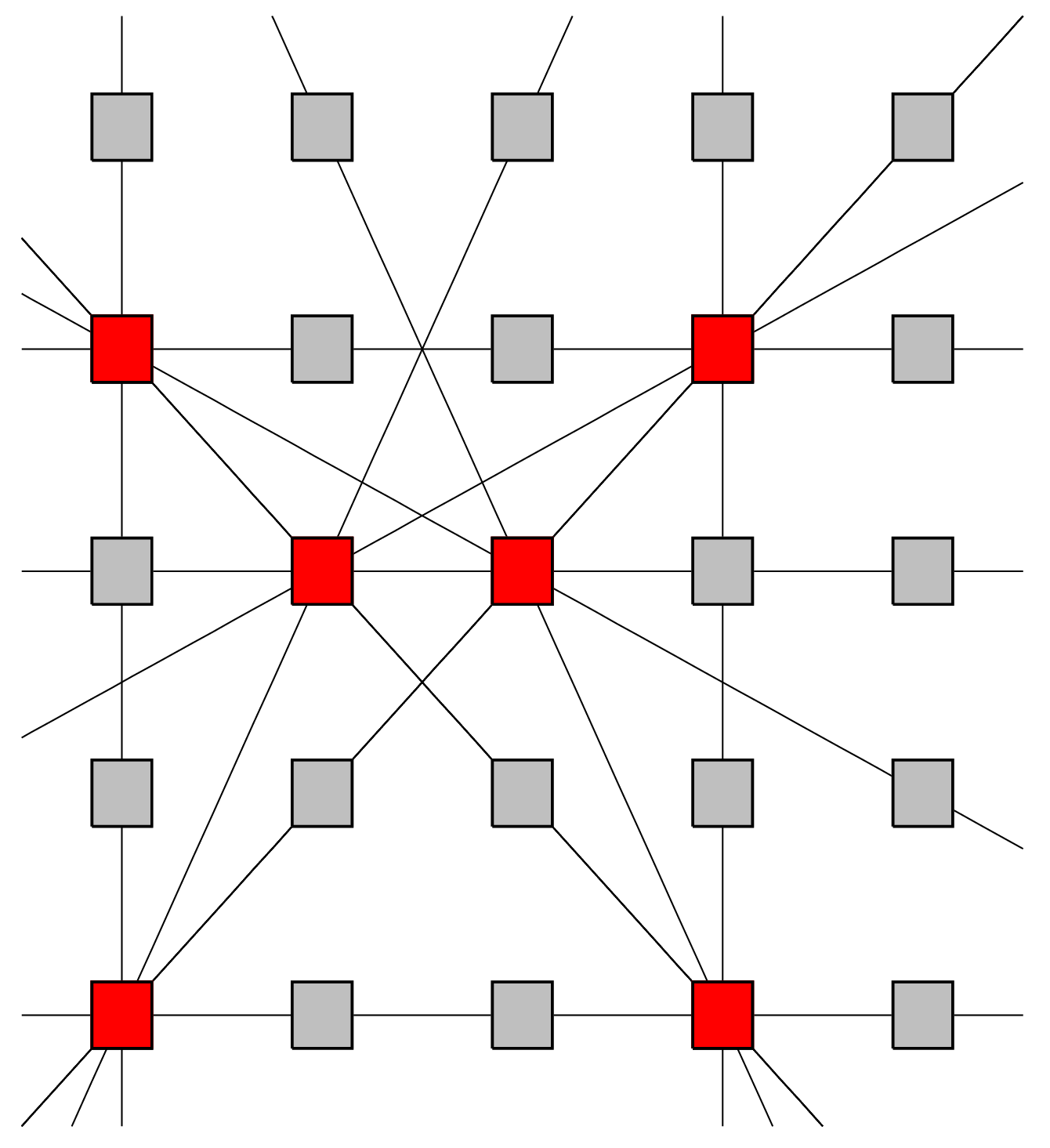}
\caption{$t(4)=6$ (continued).}
\label{fig.n4_3}
\end{figure}

\begin{figure}
\includegraphics[scale=0.3]{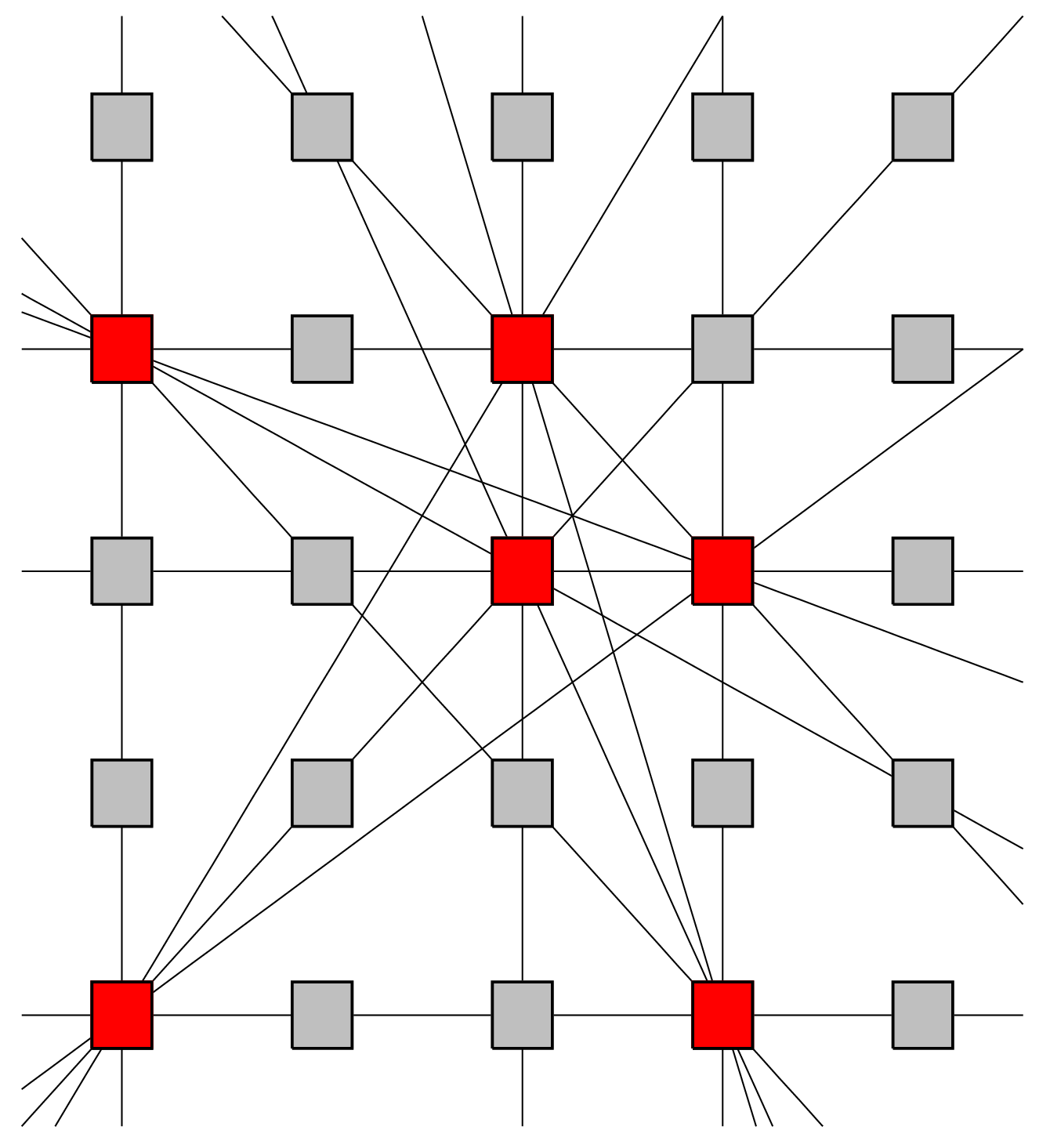}
\includegraphics[scale=0.3]{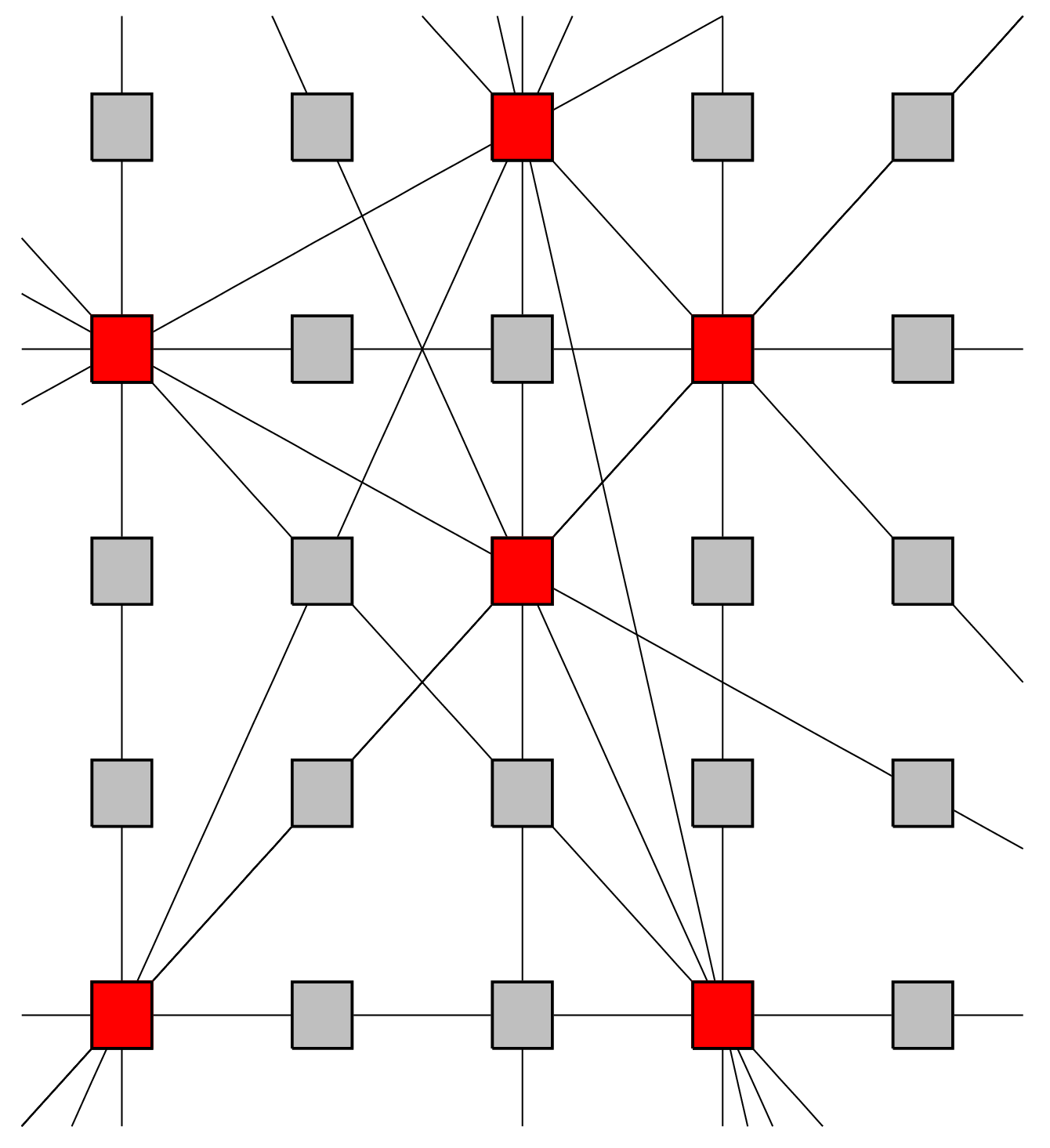}
\includegraphics[scale=0.3]{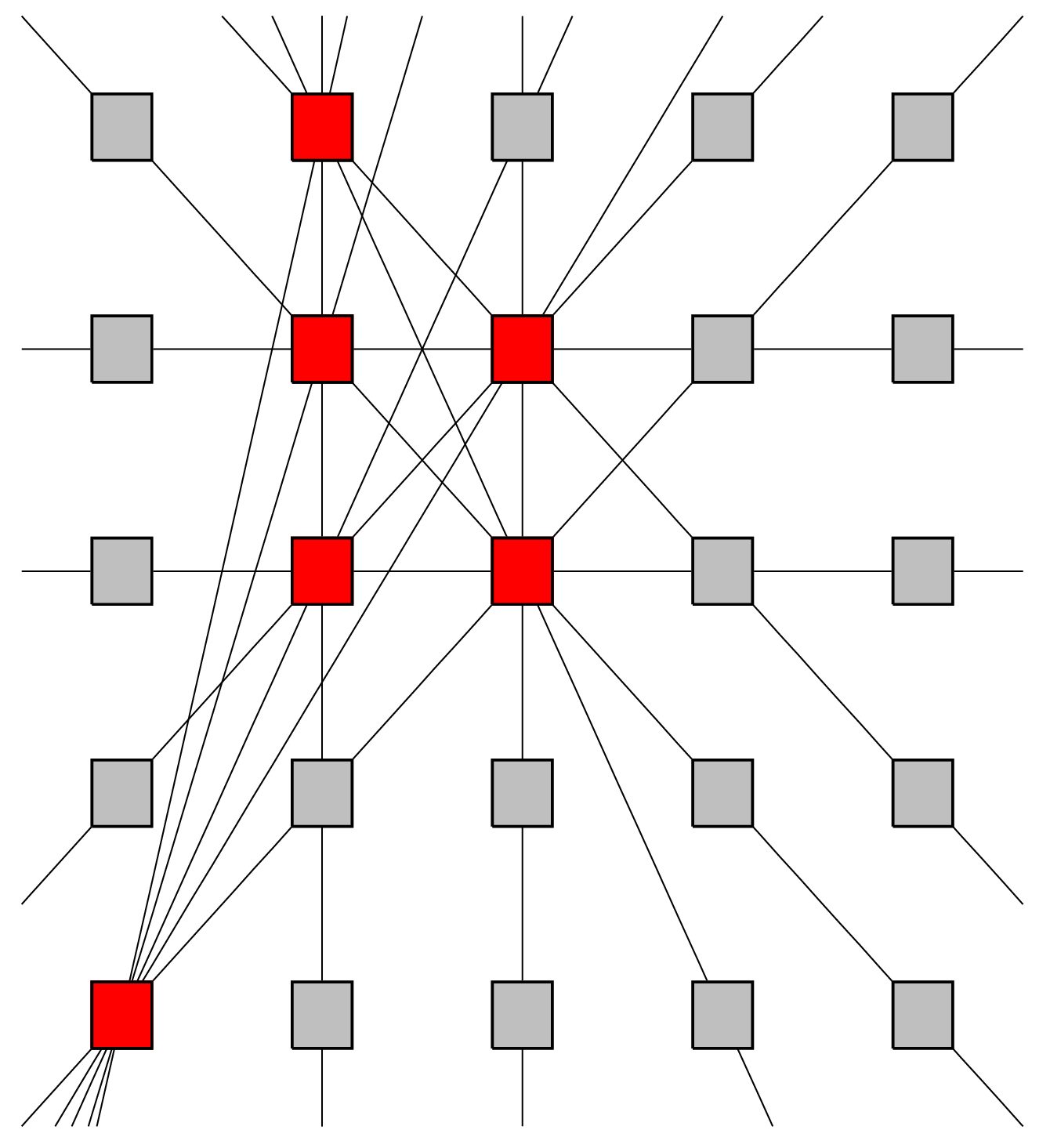}
\includegraphics[scale=0.3]{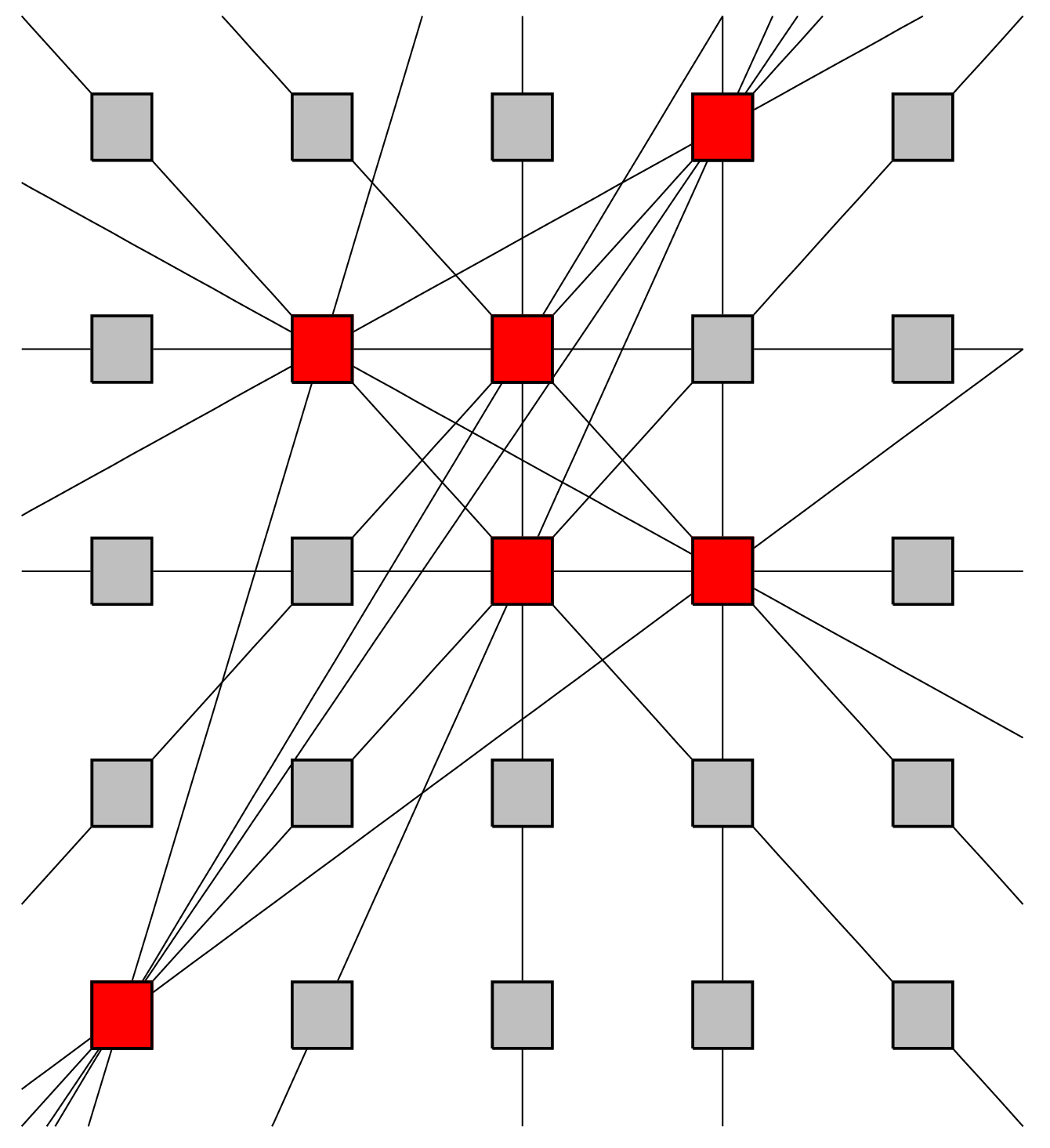}
\includegraphics[scale=0.3]{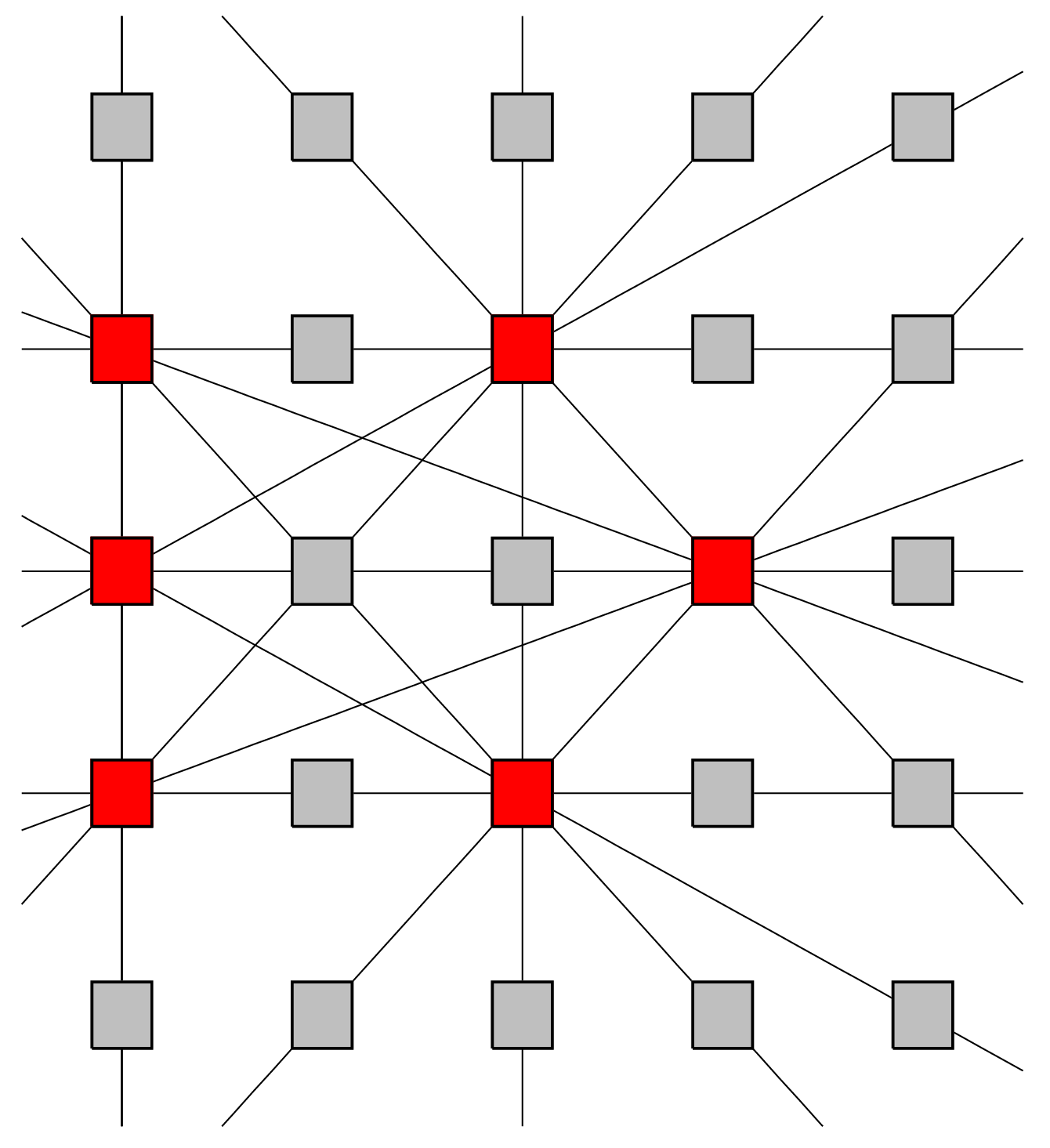}
\includegraphics[scale=0.3]{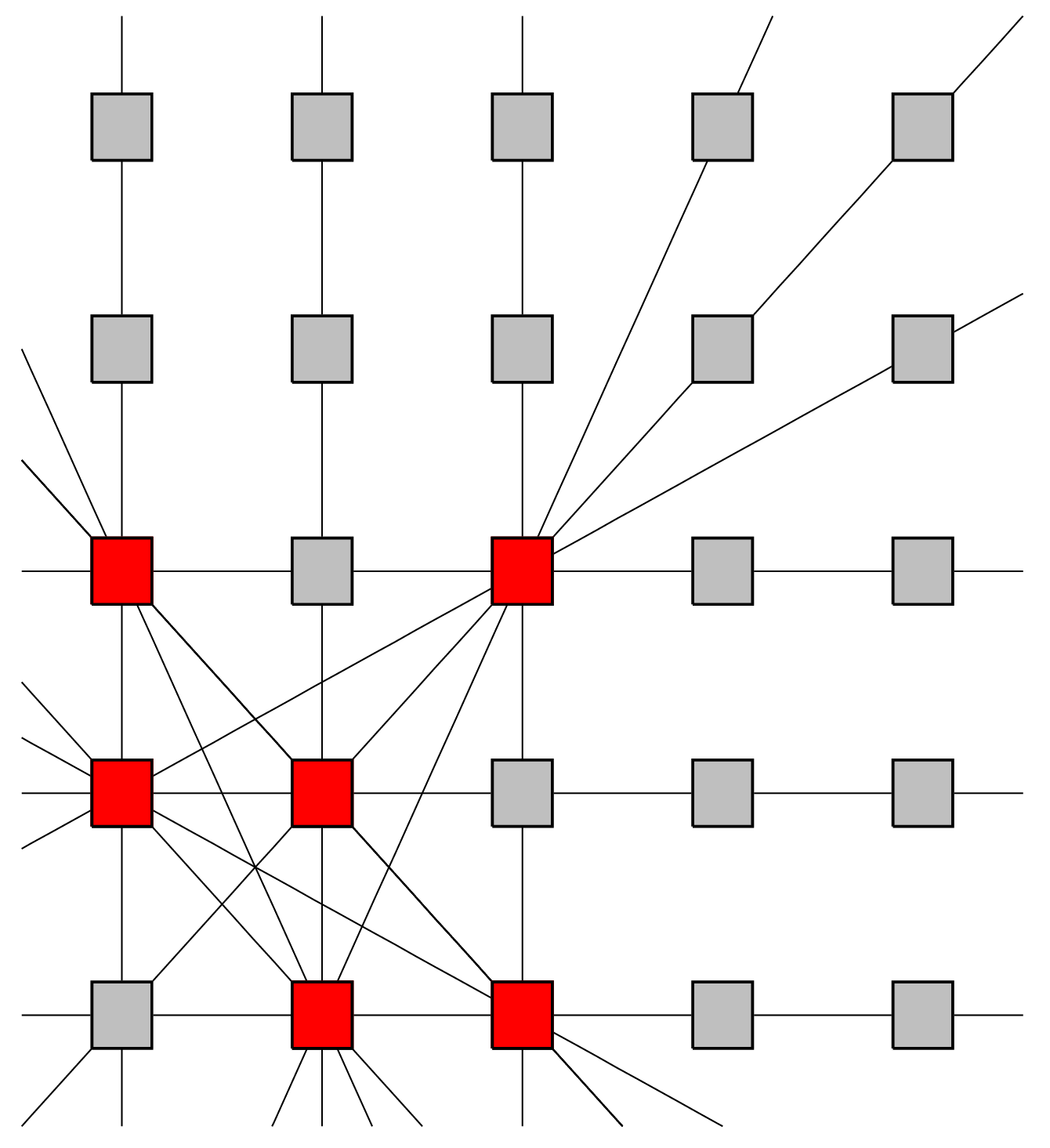}
\includegraphics[scale=0.3]{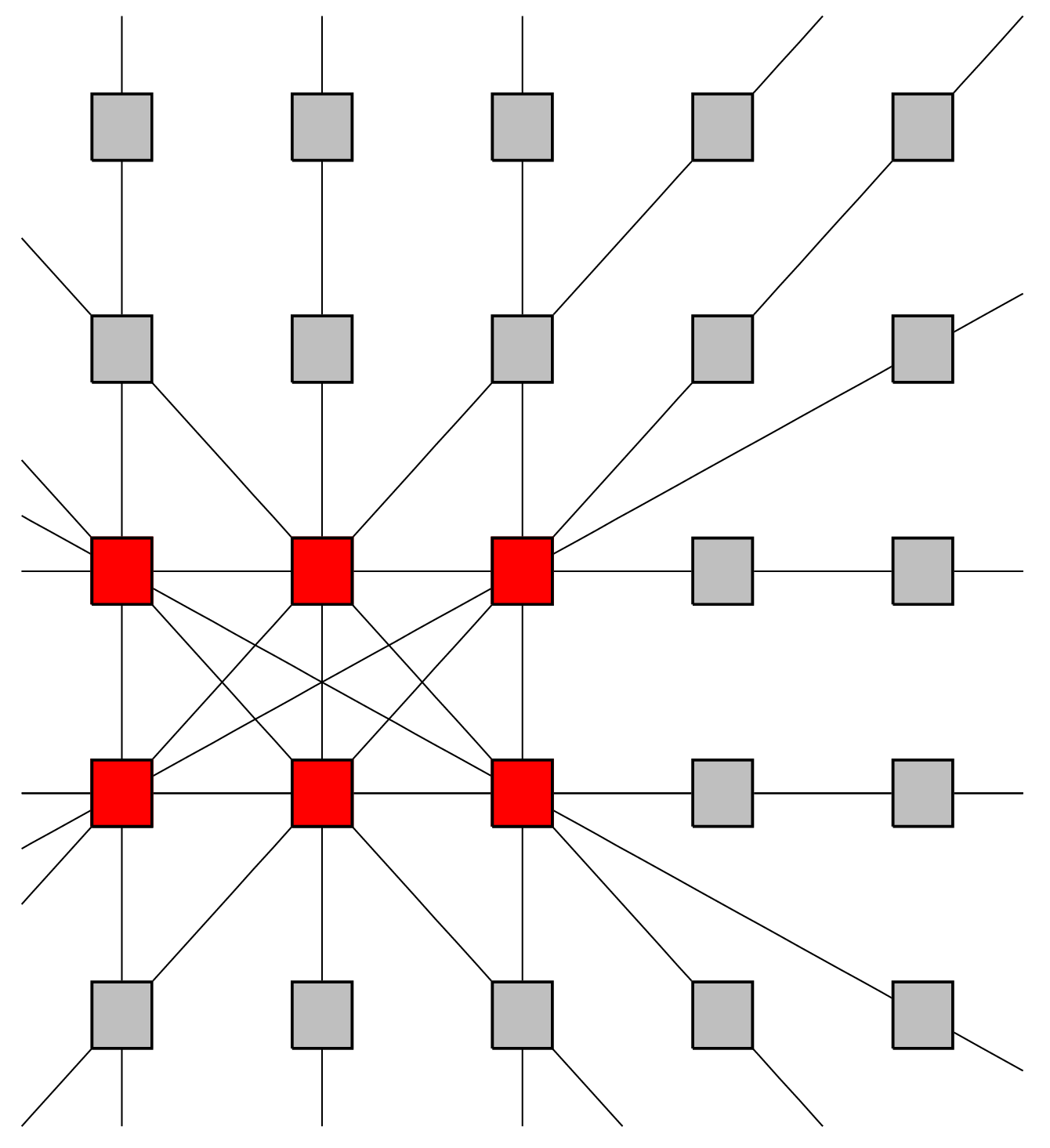}
\includegraphics[scale=0.3]{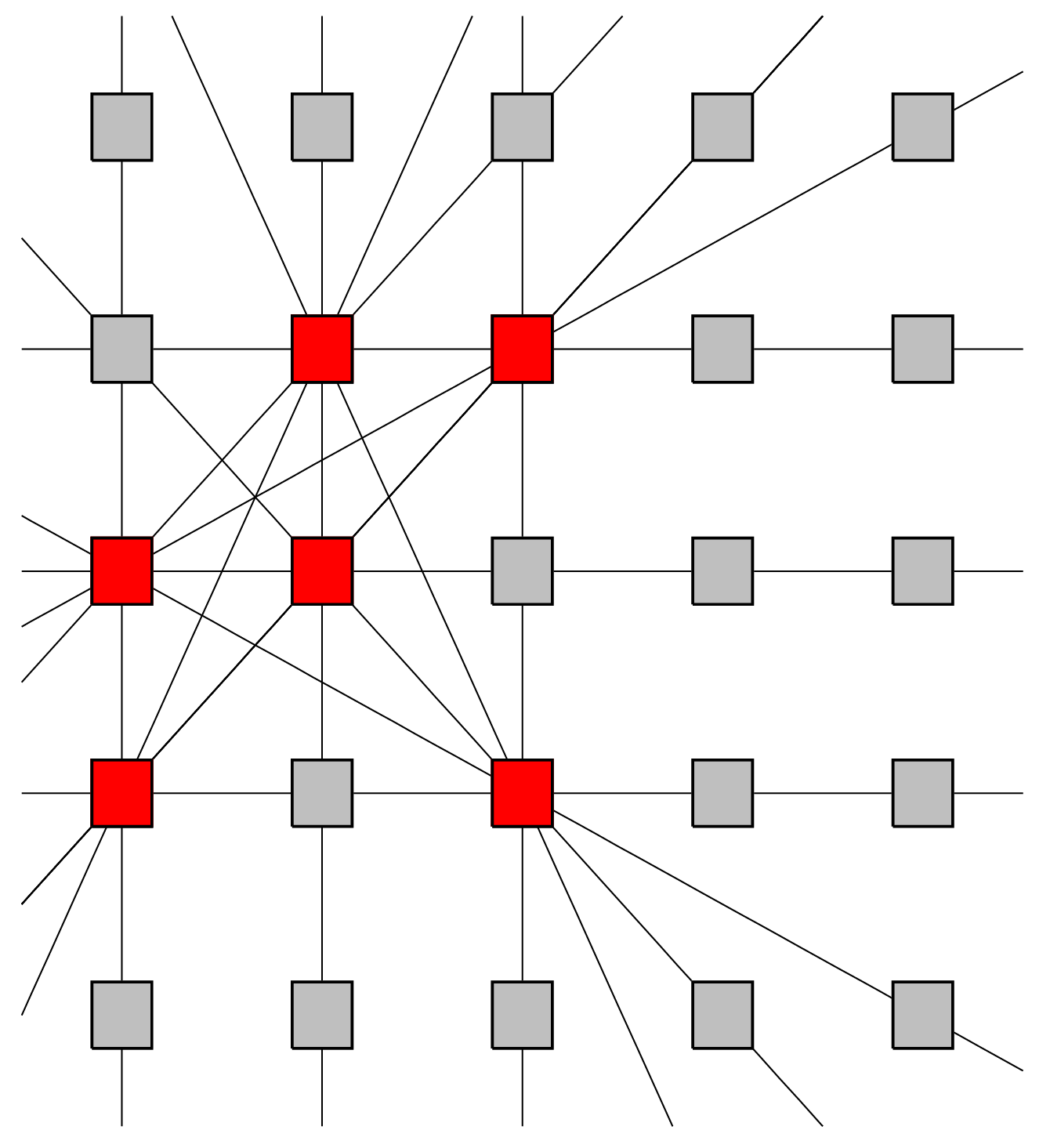}
\includegraphics[scale=0.3]{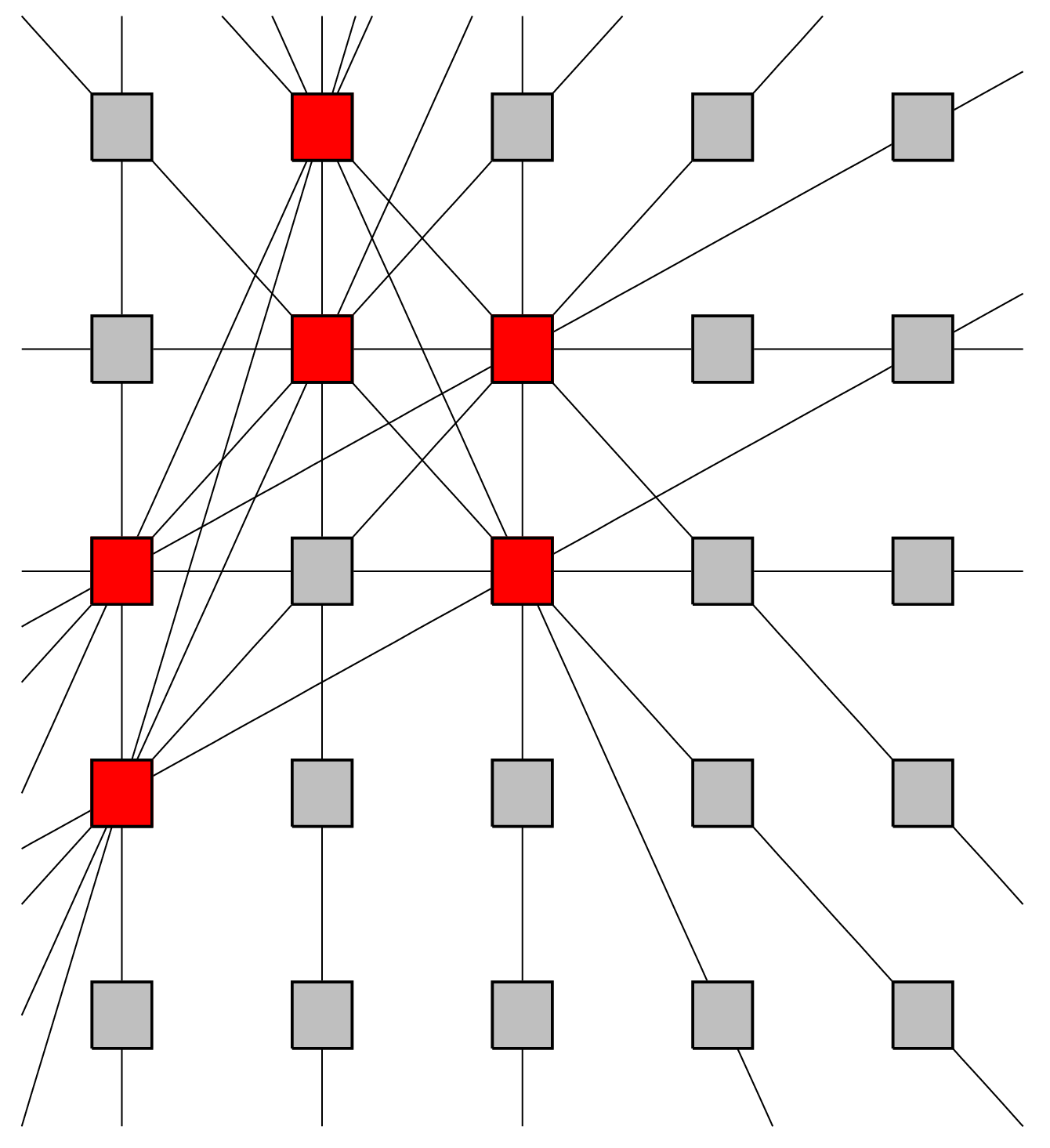}
\includegraphics[scale=0.3]{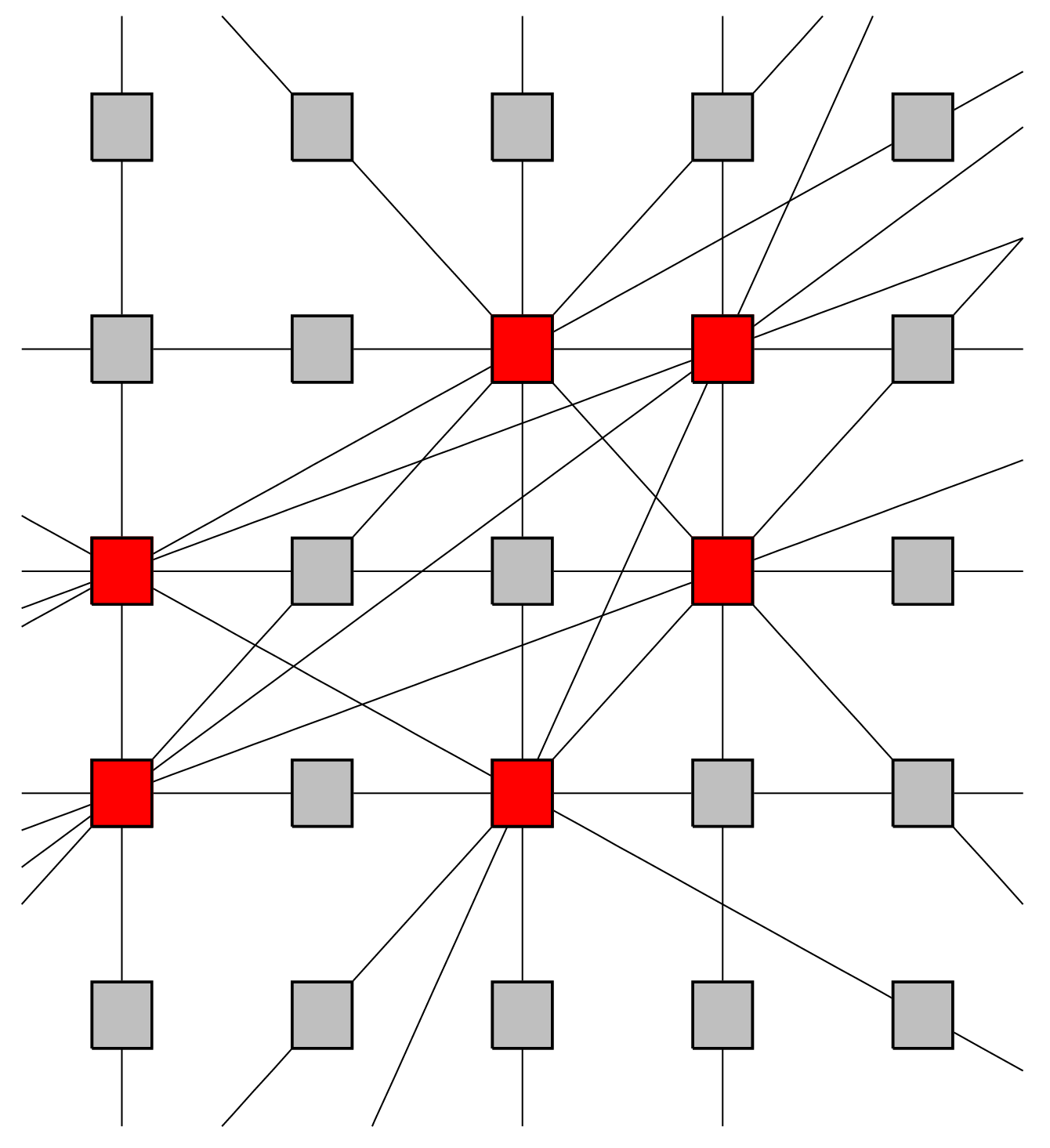}
\includegraphics[scale=0.3]{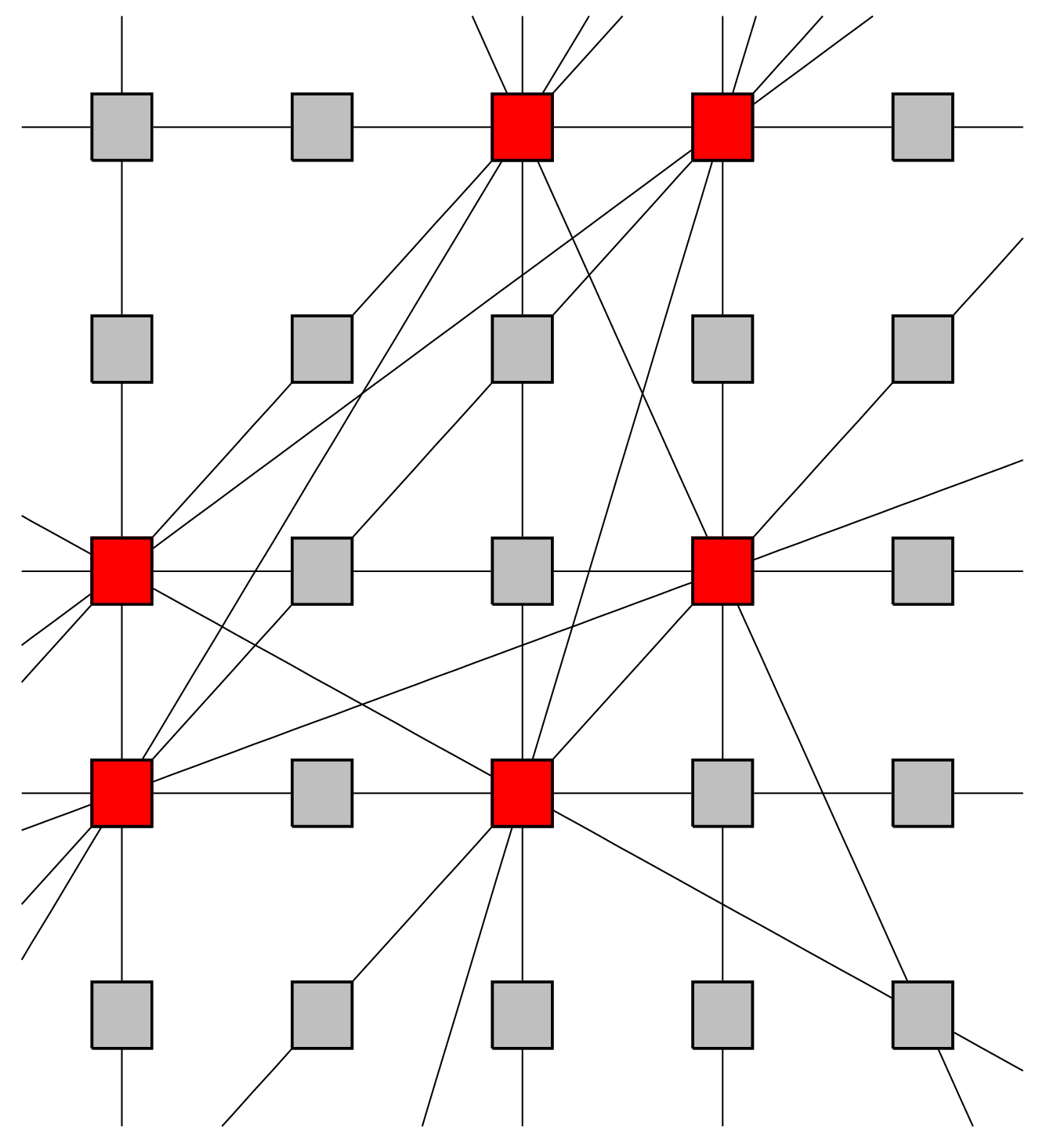}
\includegraphics[scale=0.3]{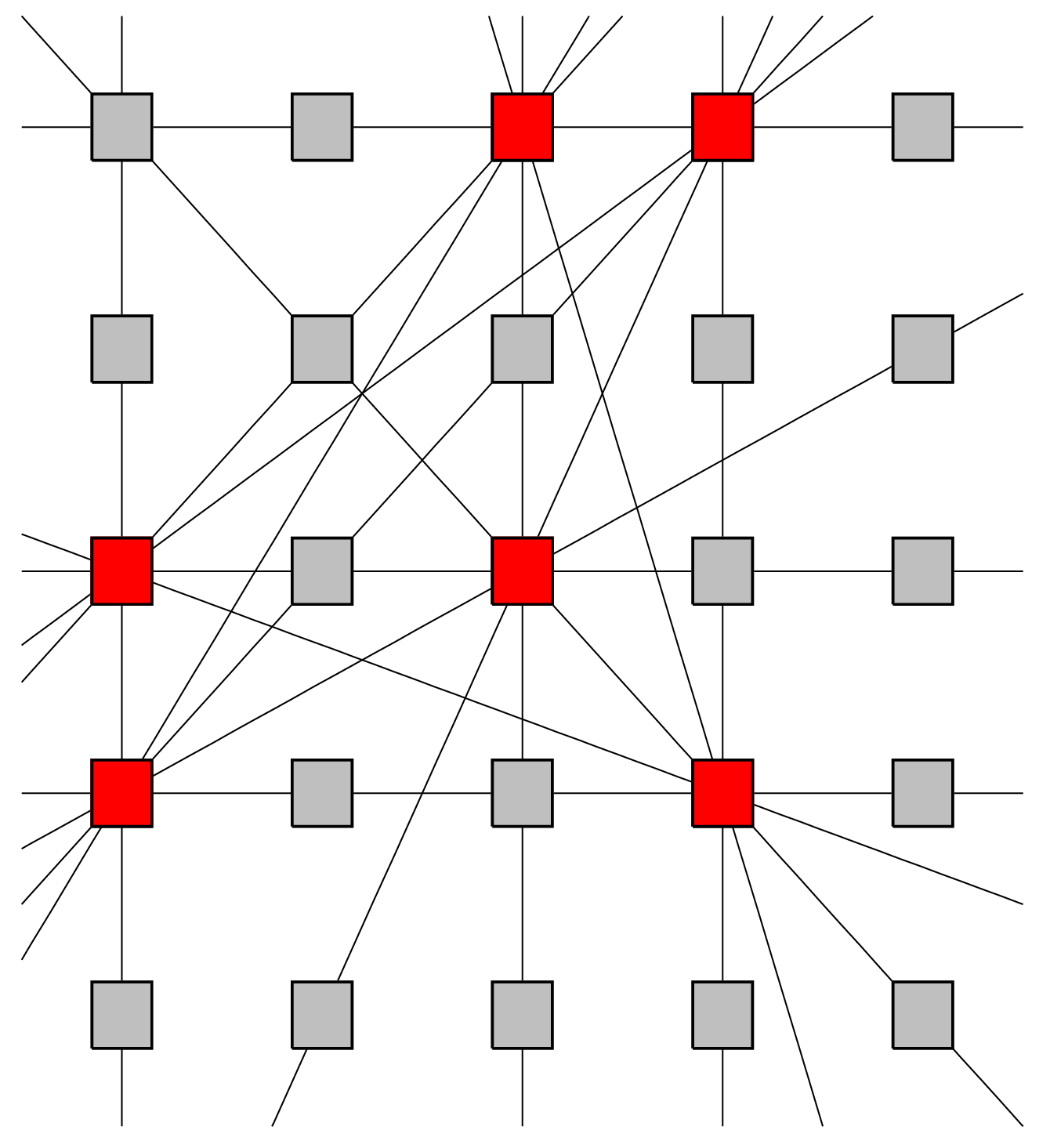}
\caption{$t(4)=6$ (continued).}
\label{fig.n4_4}
\end{figure}

\begin{figure}
\includegraphics[scale=0.3]{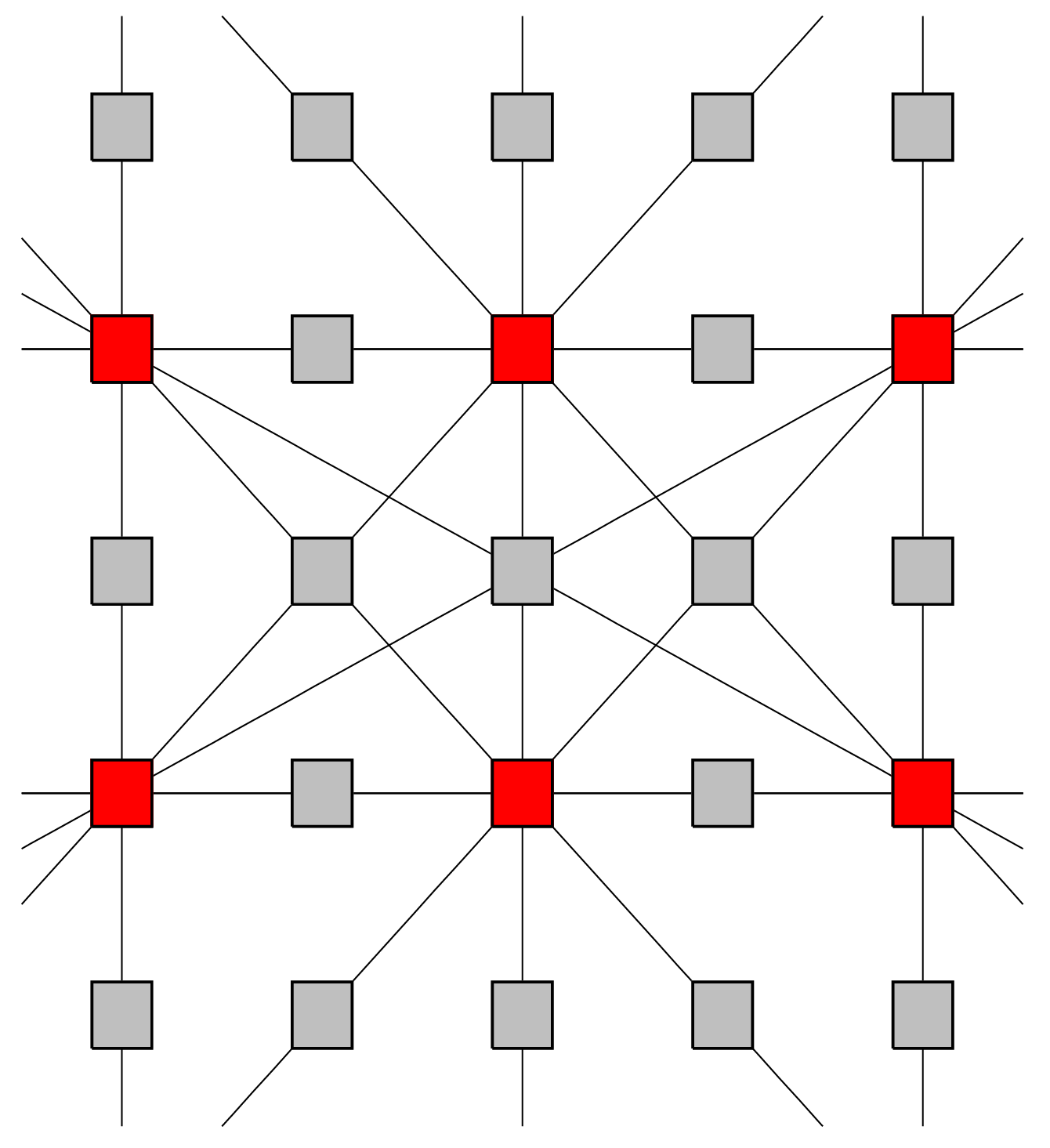}
\includegraphics[scale=0.3]{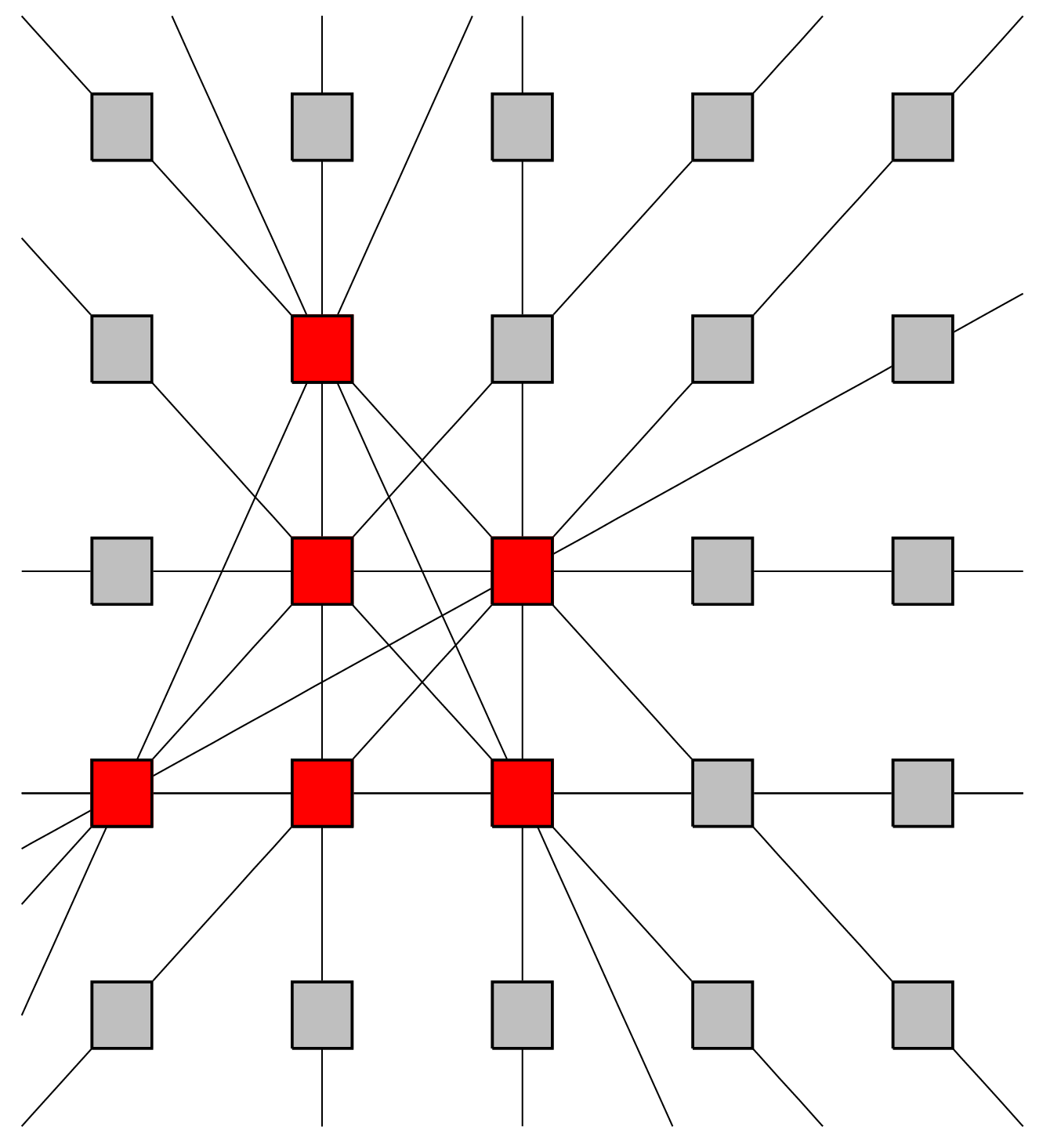}
\includegraphics[scale=0.3]{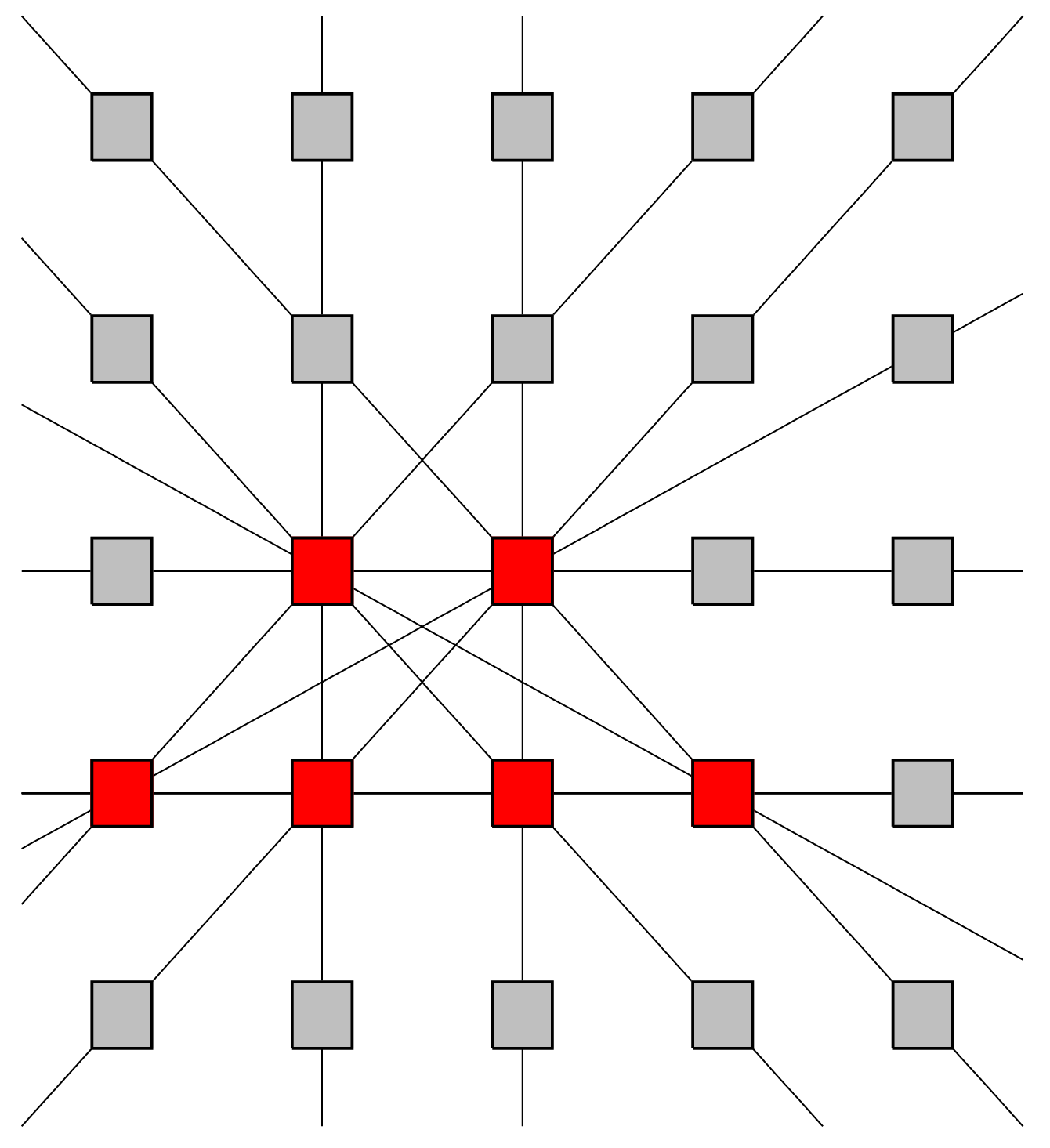}
\includegraphics[scale=0.3]{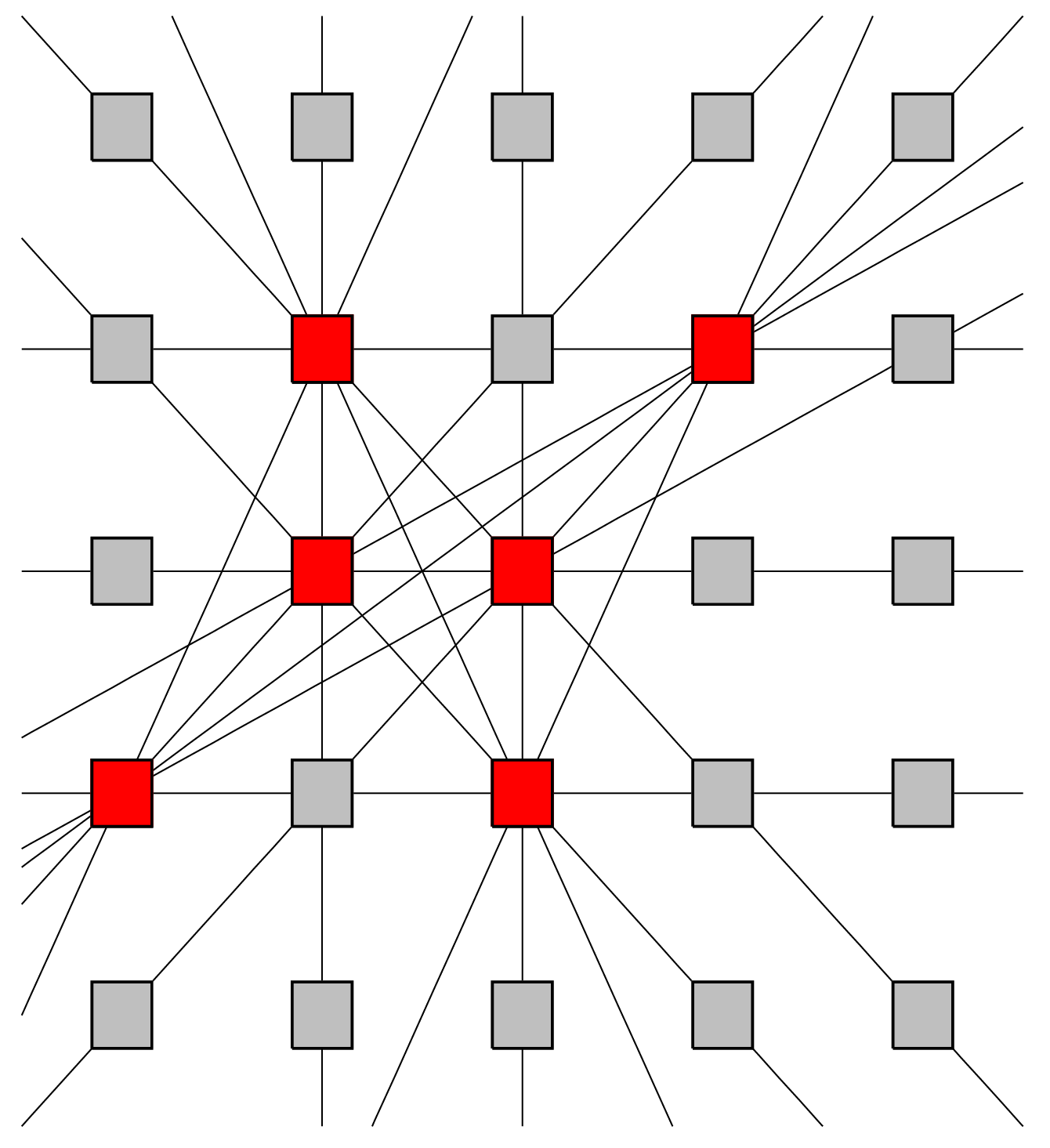}
\includegraphics[scale=0.3]{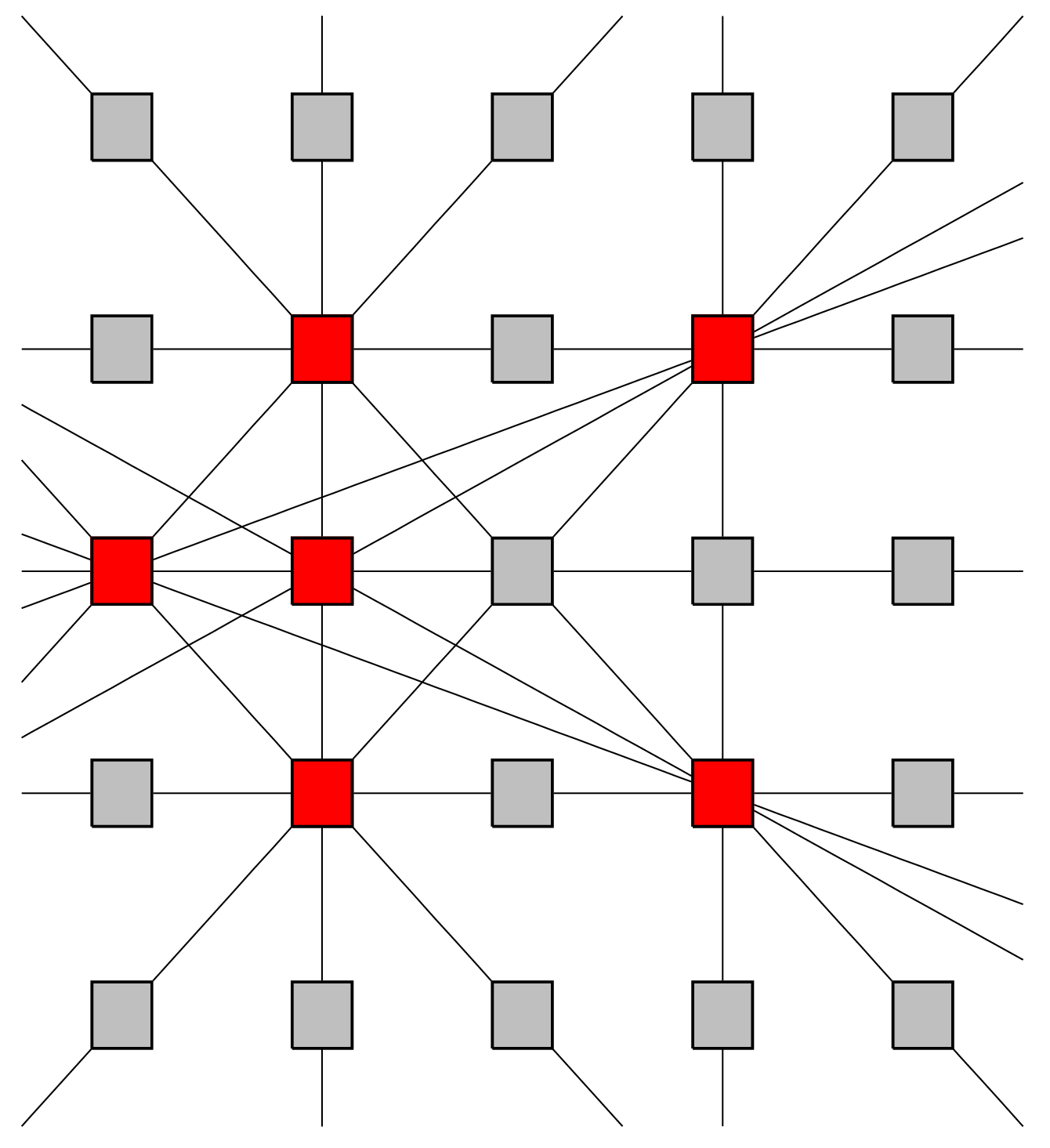}
\includegraphics[scale=0.3]{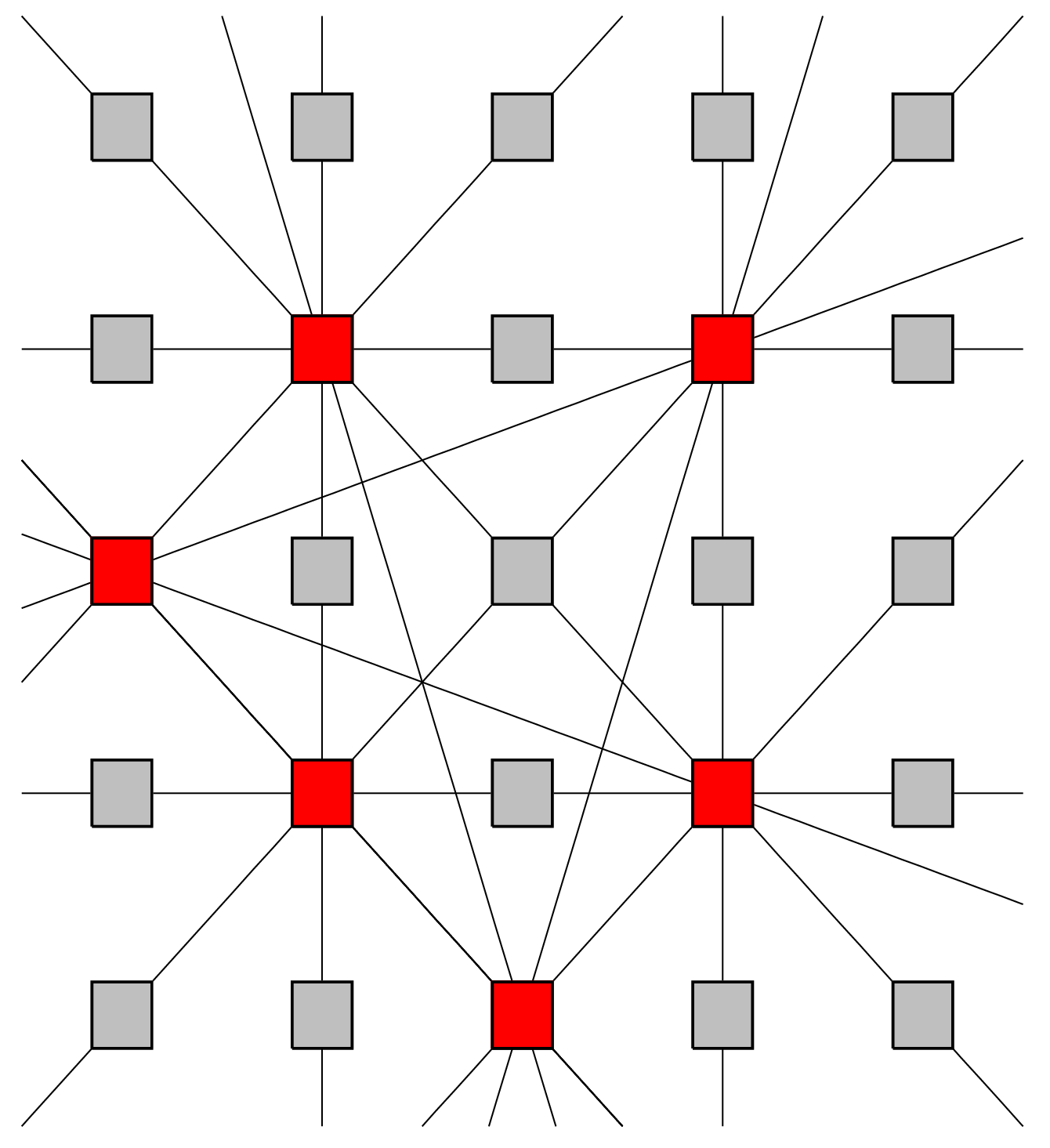}
\includegraphics[scale=0.3]{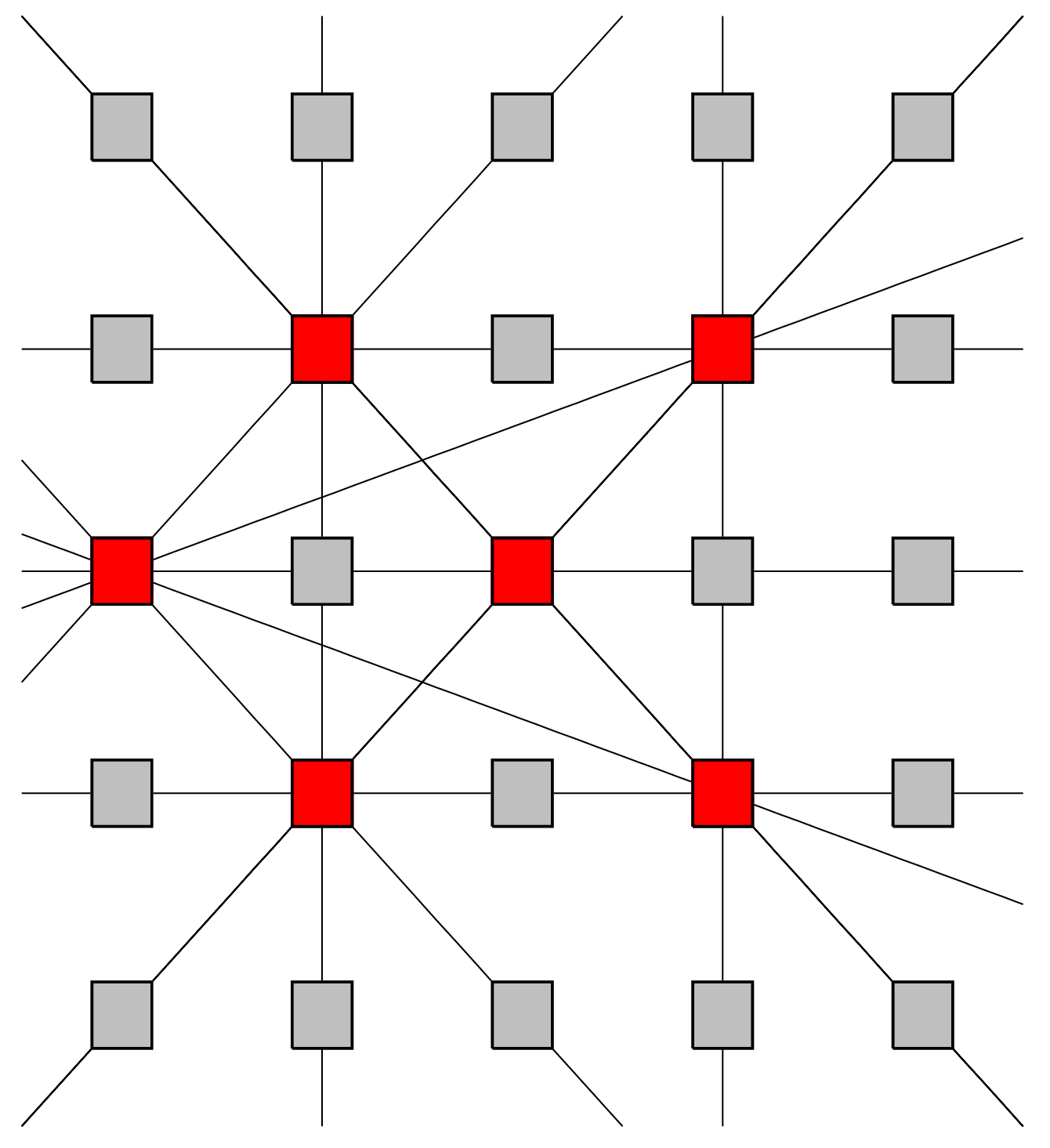}
\includegraphics[scale=0.3]{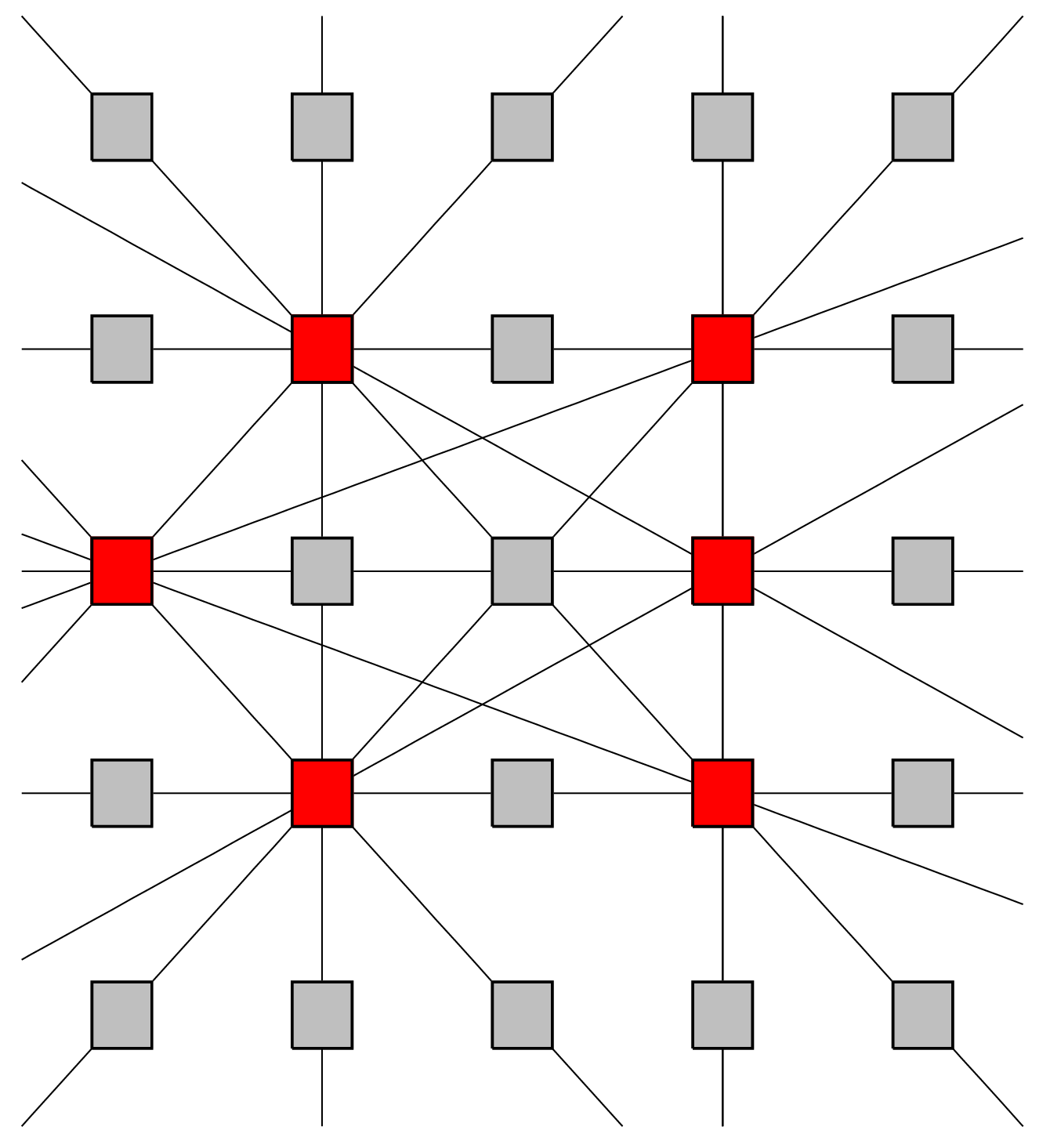}
\includegraphics[scale=0.3]{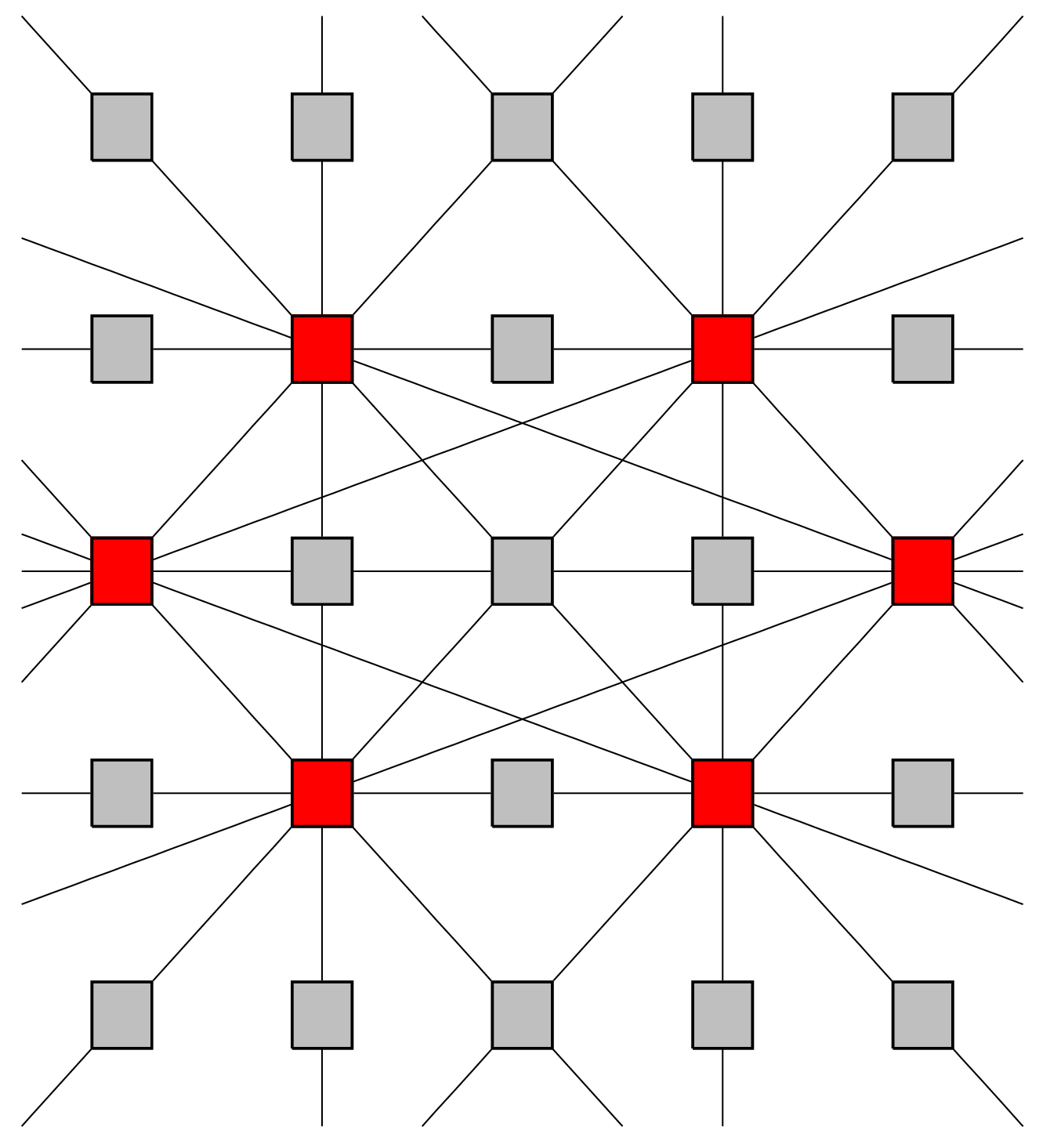}
\includegraphics[scale=0.3]{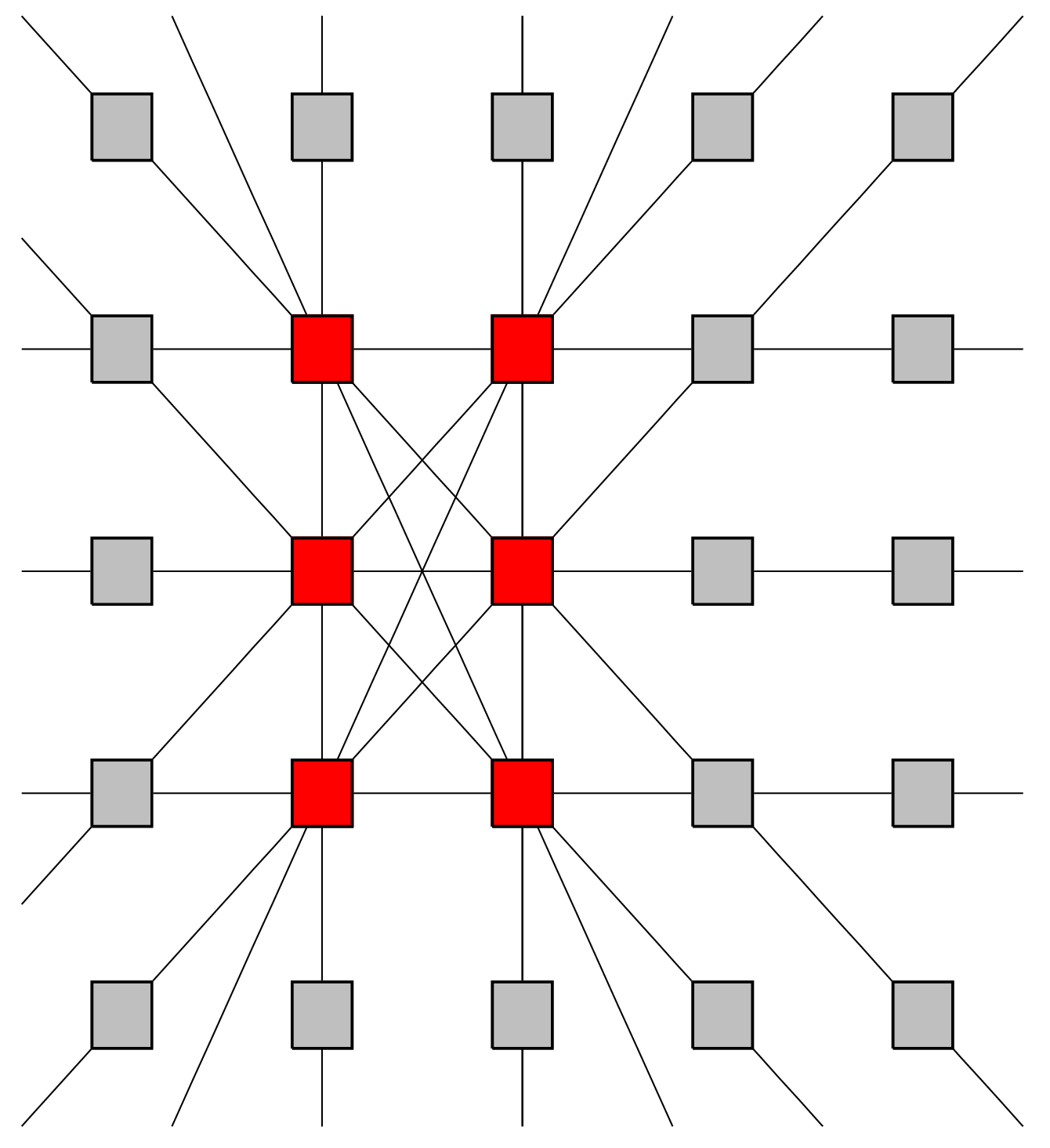}
\includegraphics[scale=0.3]{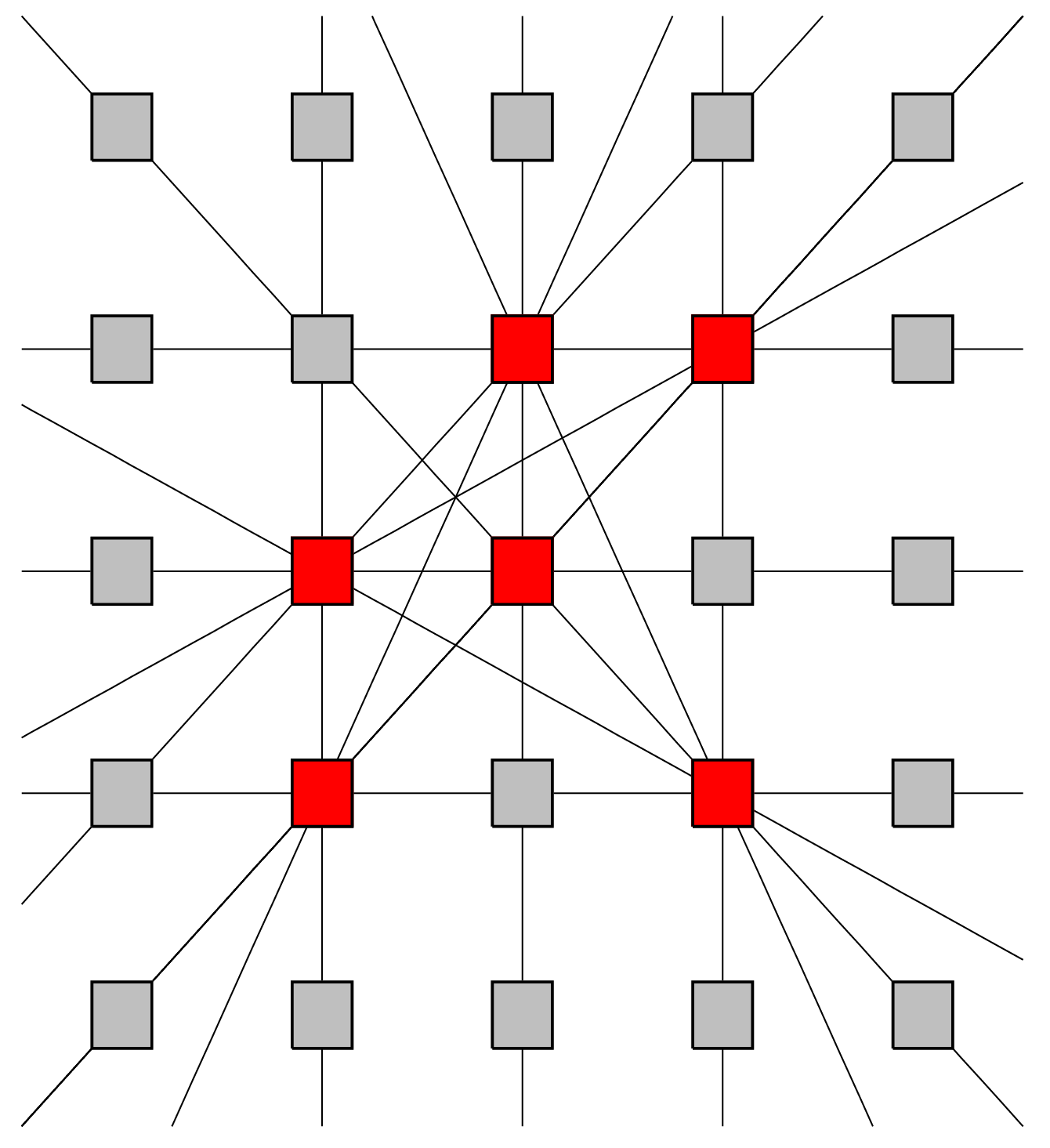}
\caption{$t(4)=6$ (continued).}
\label{fig.n4_l}
\end{figure}

59 incongruential choices of selecting $t(4)=6$ vertices are known:
Figures \ref{fig.n4_1}--\ref{fig.n4_l}.

\clearpage

\clearpage


\clearpage

\begin{figure}[h]
\includegraphics[scale=0.3]{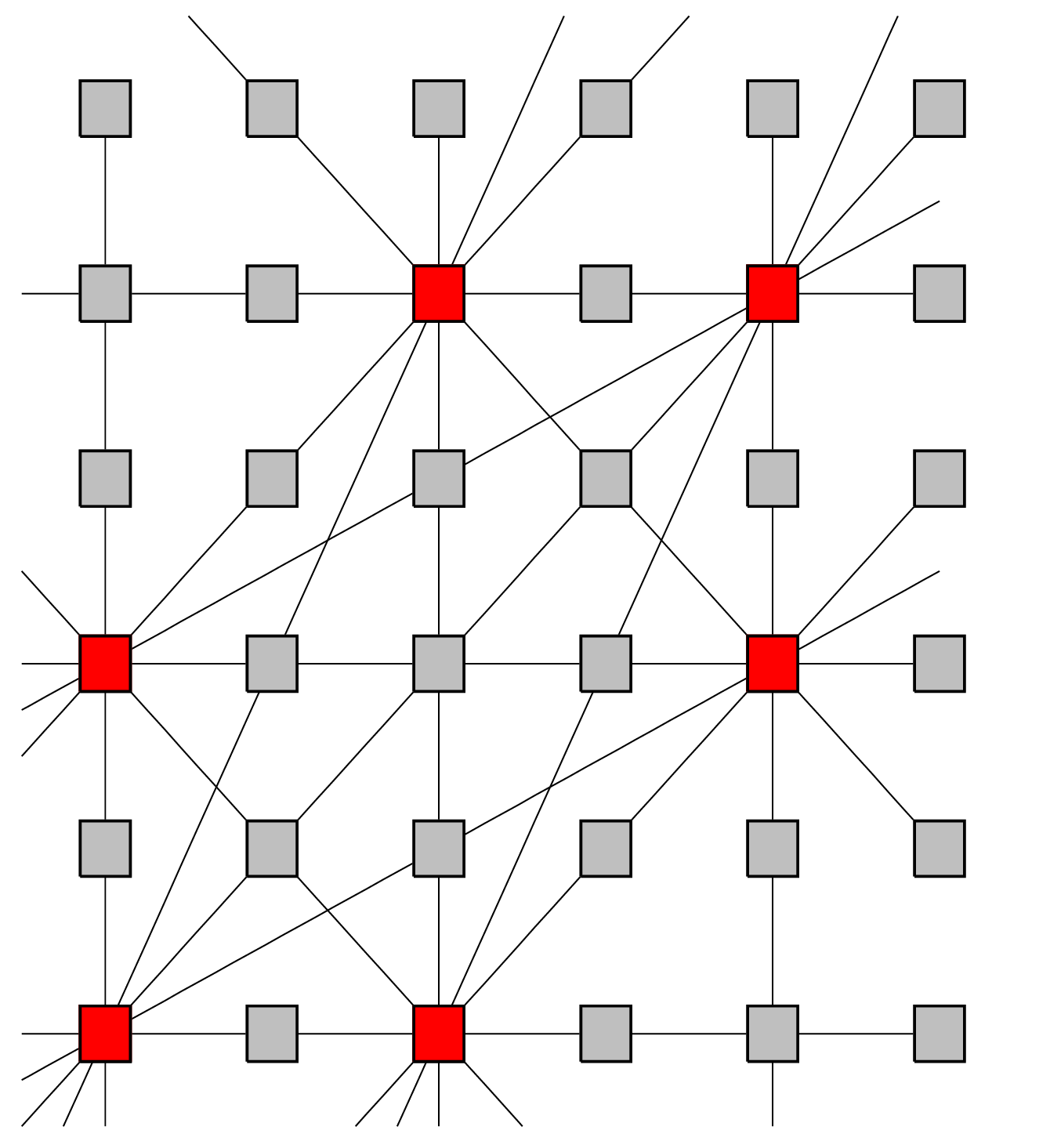}
\includegraphics[scale=0.3]{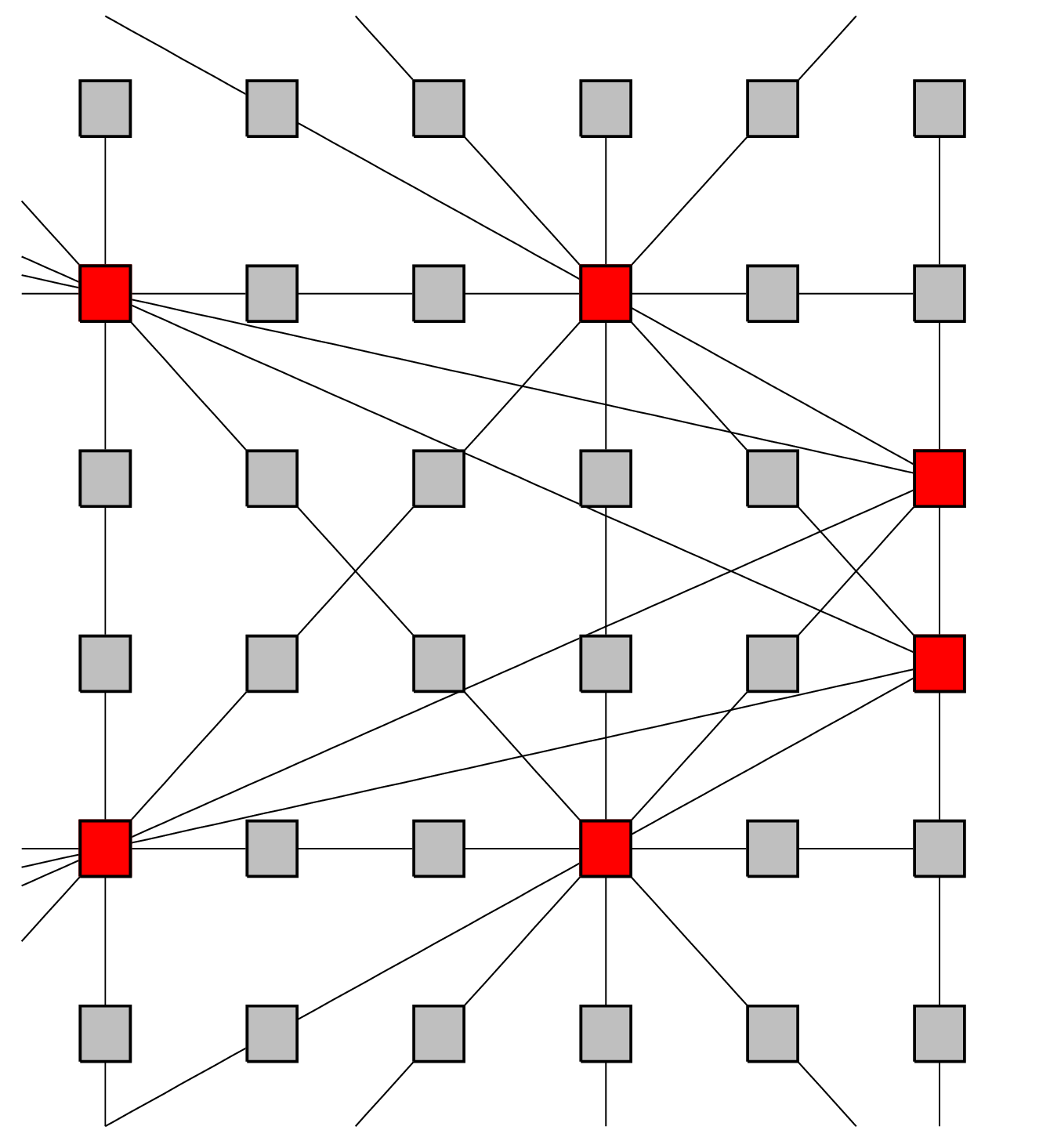}
\includegraphics[scale=0.3]{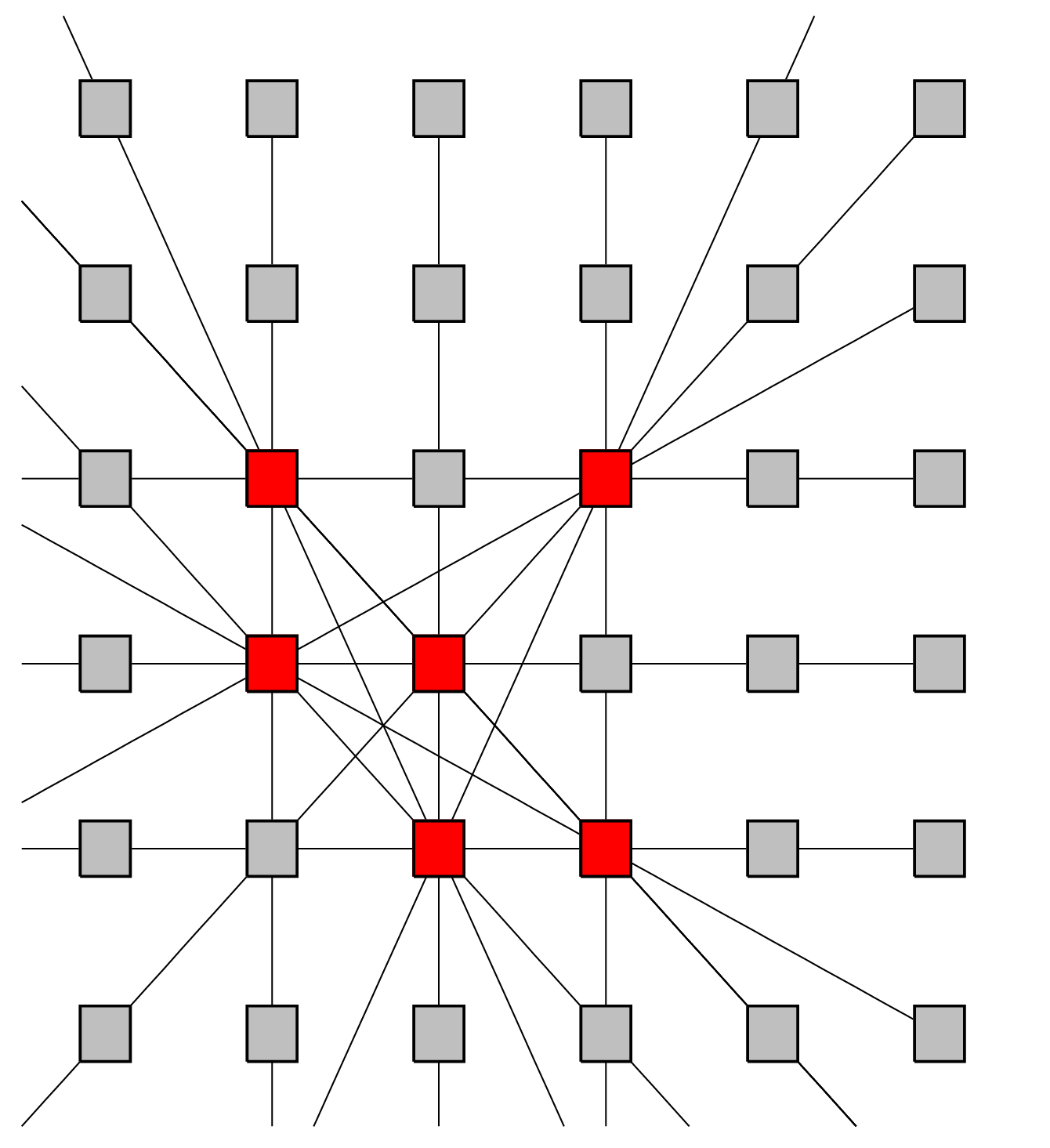}
\includegraphics[scale=0.3]{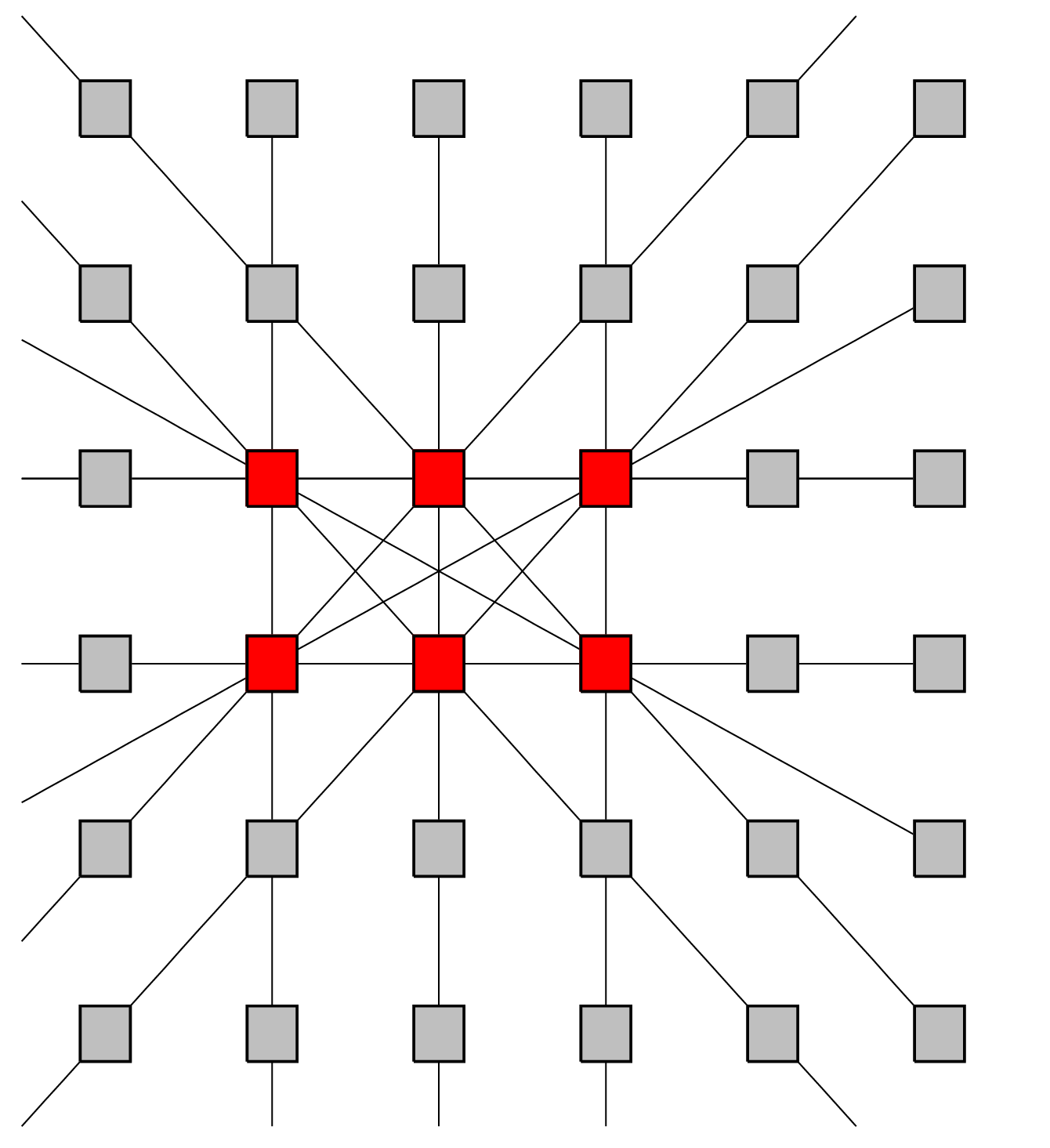}
\caption{$t(5)=6$.}
\label{fig.n5}
\end{figure}

At $n=5$, all 4 incongruential choices of coverage by $t(5)=6$ vertices are shown
in Fig.\ \ref{fig.n5}\@.
For three of these, the sublattice points fit already in the smaller 
$n=4$ lattice:
the upper left is found in the second row, third column
of Fig.\ \ref{fig.n4_3}.
The upper right is found in the second row, third column
of Fig.\ \ref{fig.n4_4}, in the
third row second column of 
Fig.\ \ref{fig.n4_4},
and
in the last row of Fig.\ \ref{fig.n4_l}.
The last is found in 
the third row, first column of Fig.\ \ref{fig.n4_4}
and in the last row of Fig.\ \ref{fig.n4_l}.

\clearpage

\begin{figure}[h]
\includegraphics[scale=0.3]{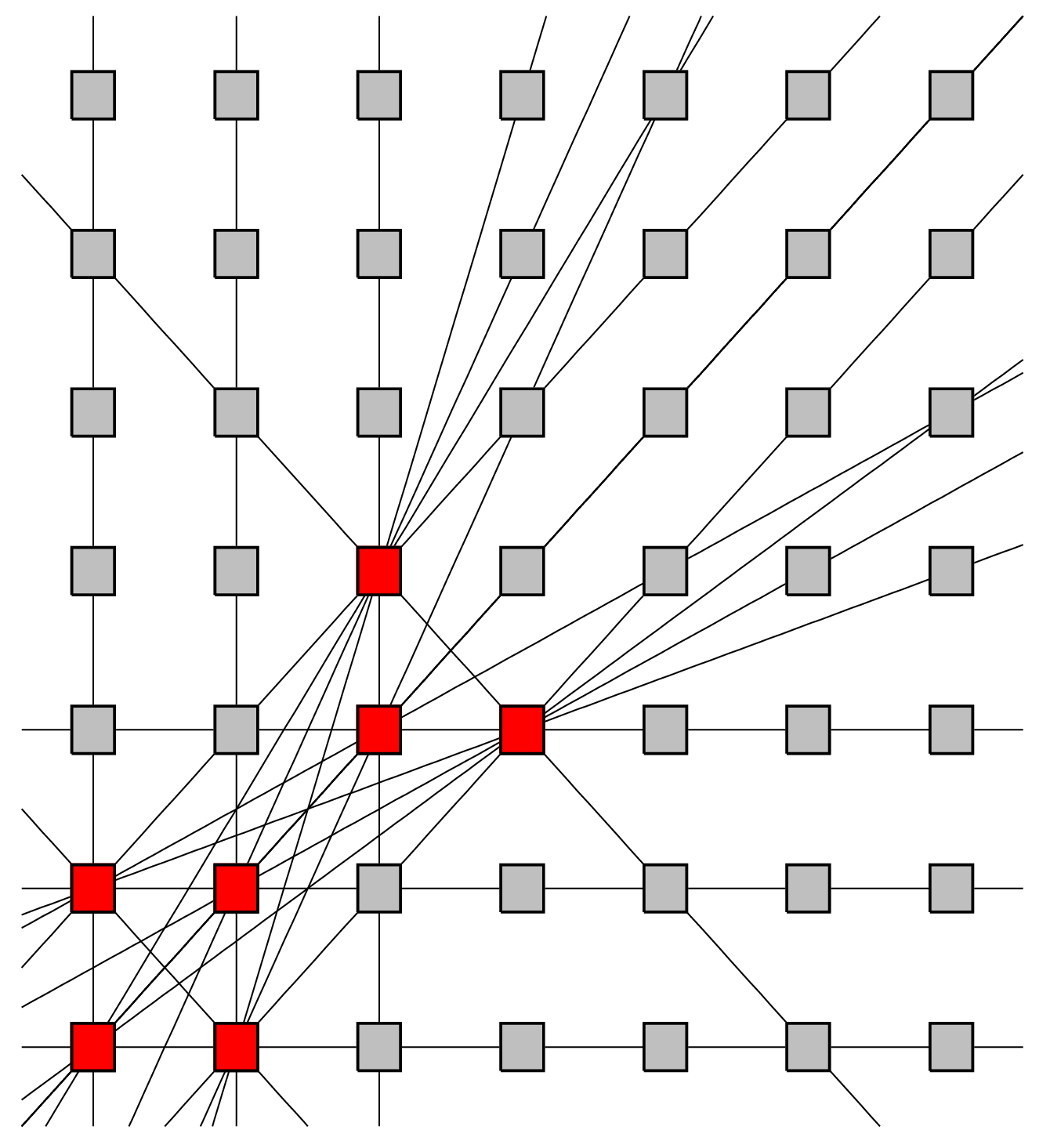}
\includegraphics[scale=0.3]{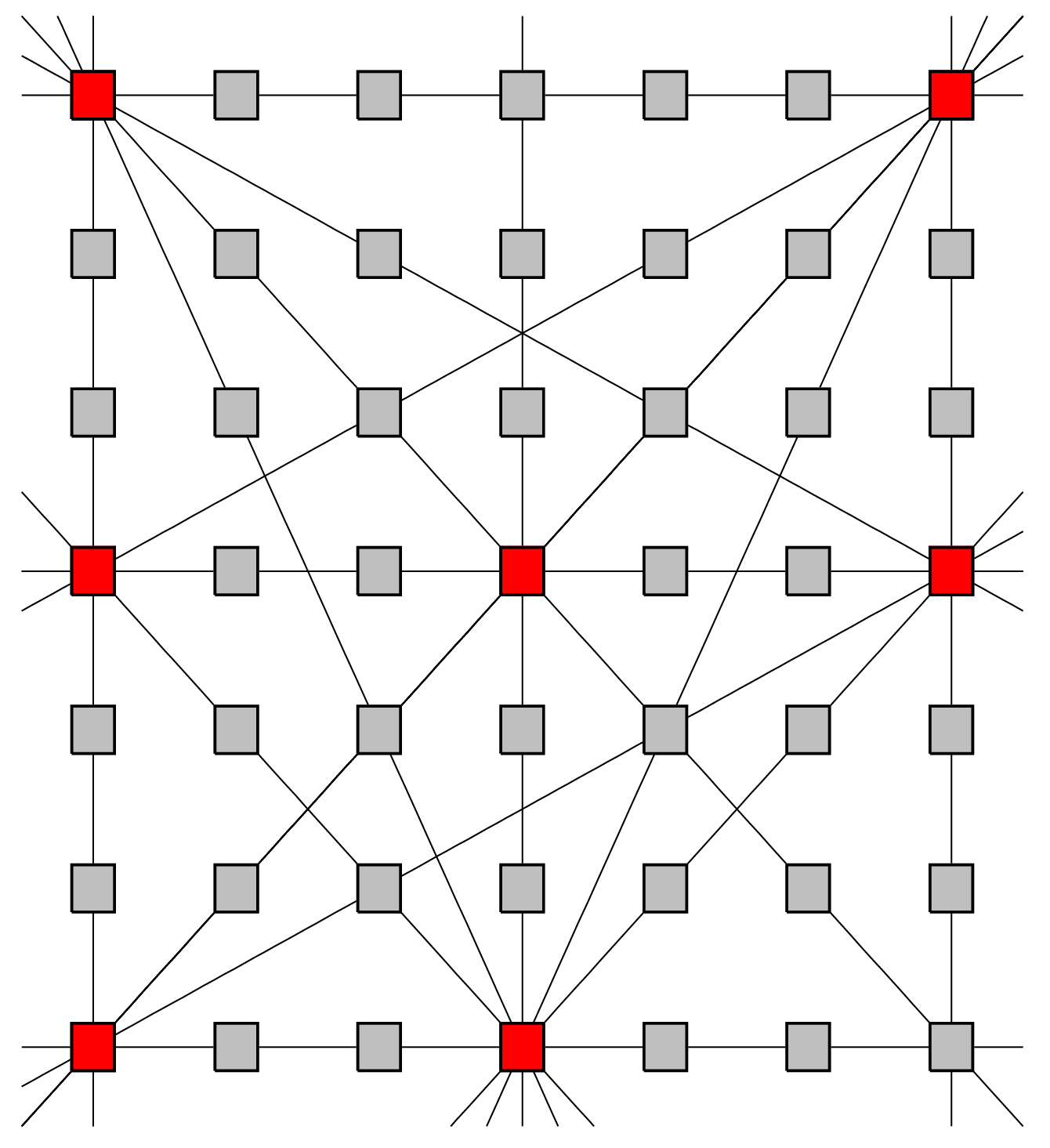}
\includegraphics[scale=0.3]{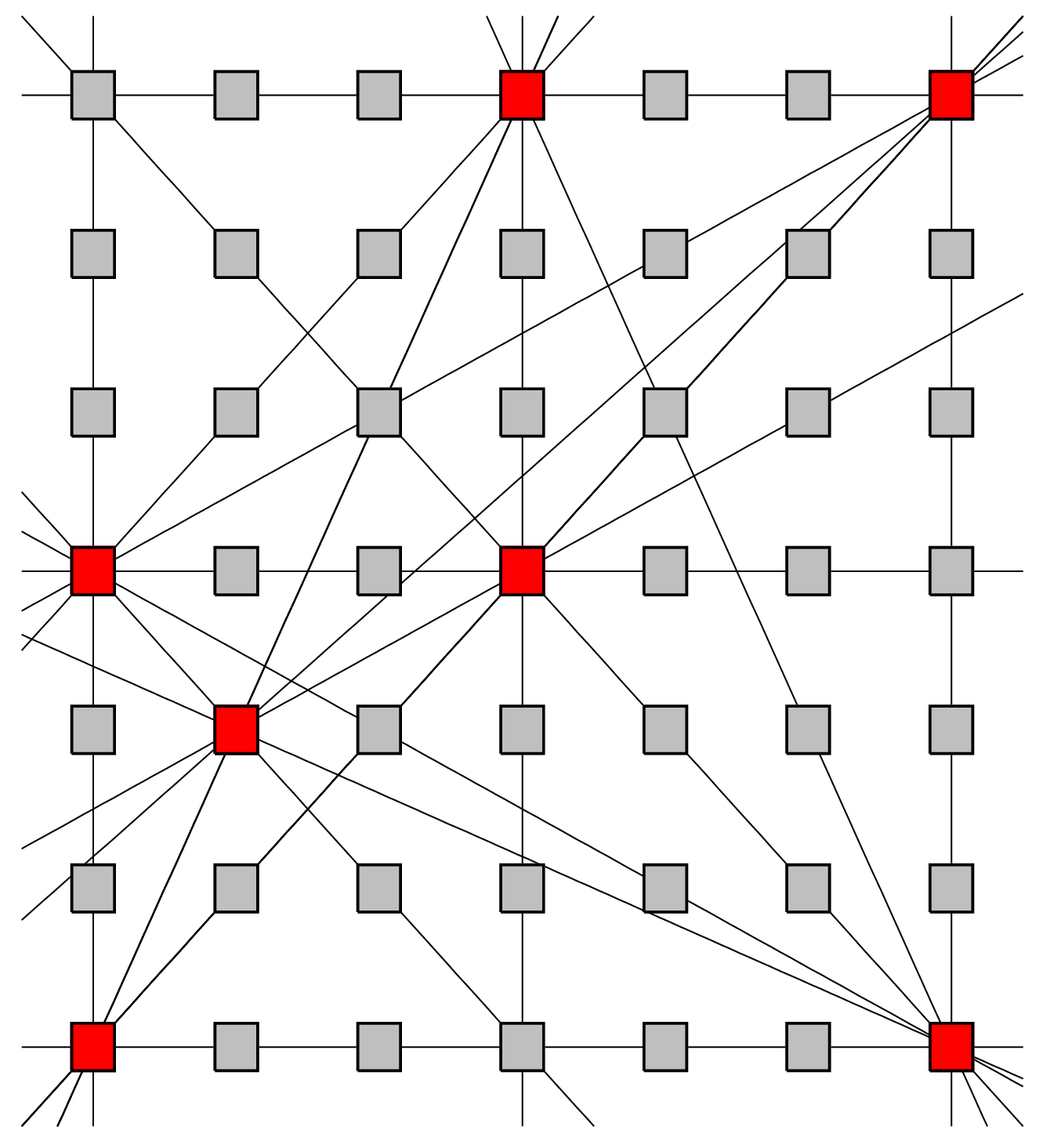}
\includegraphics[scale=0.3]{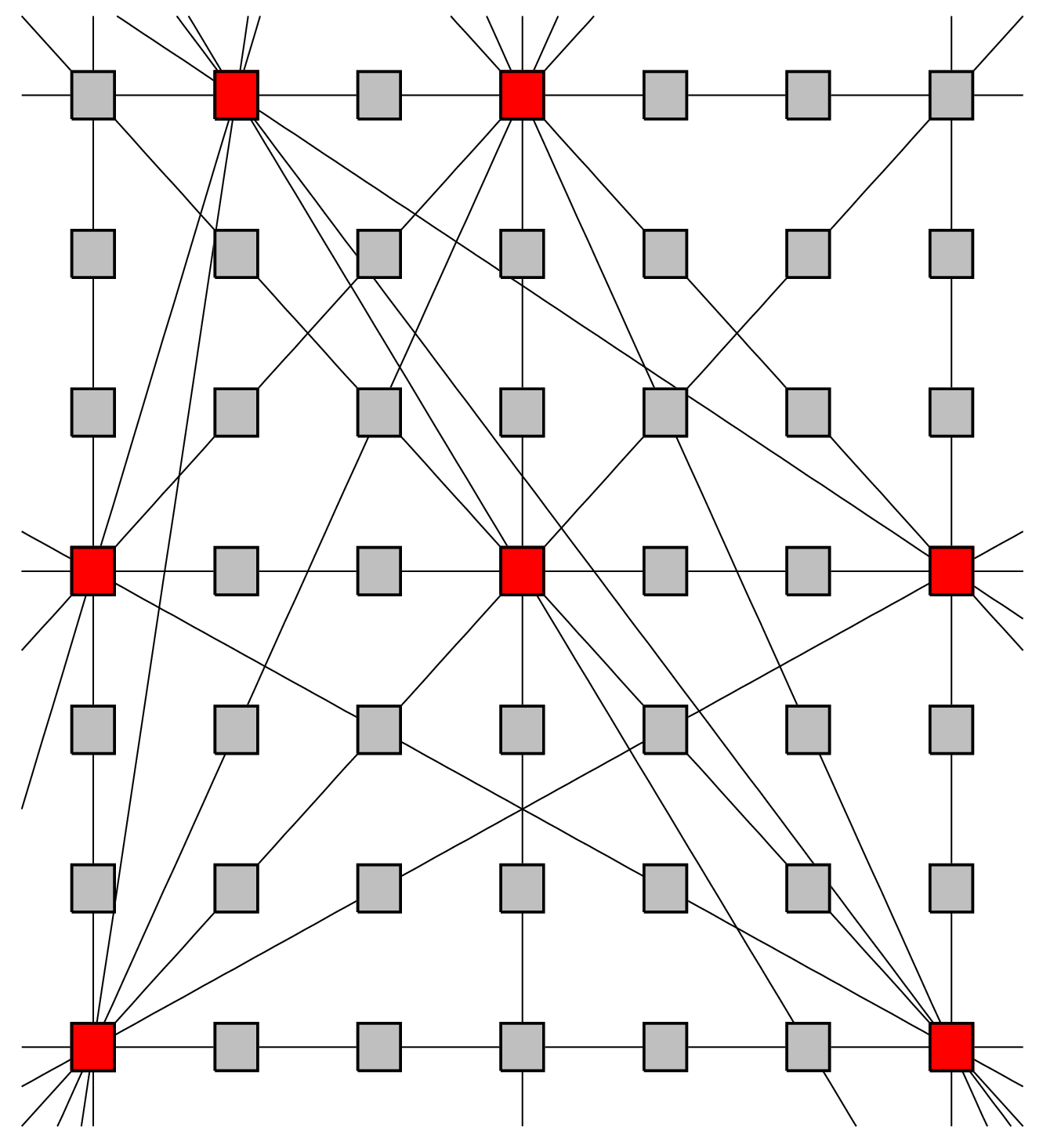}
\includegraphics[scale=0.3]{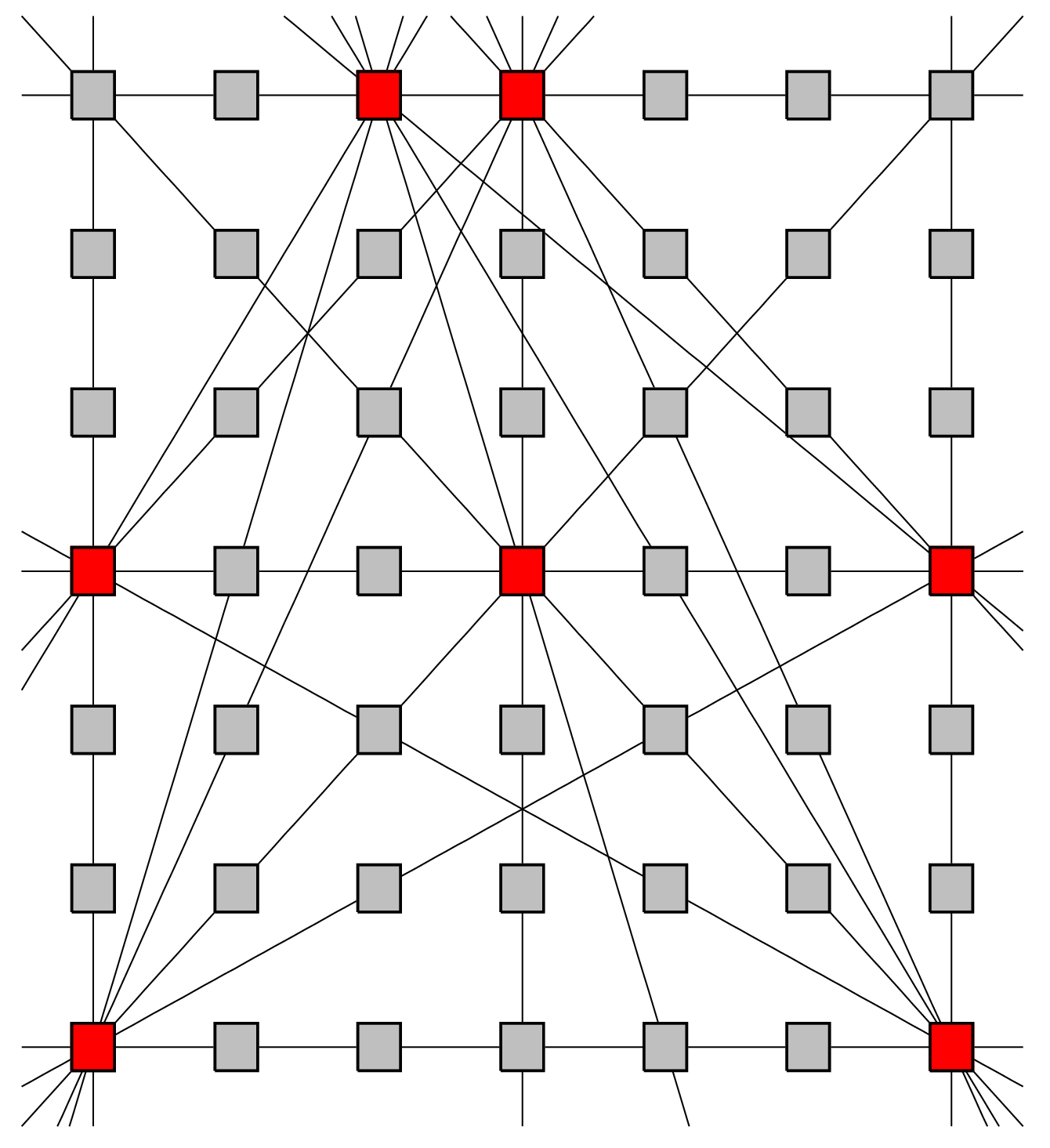}
\includegraphics[scale=0.3]{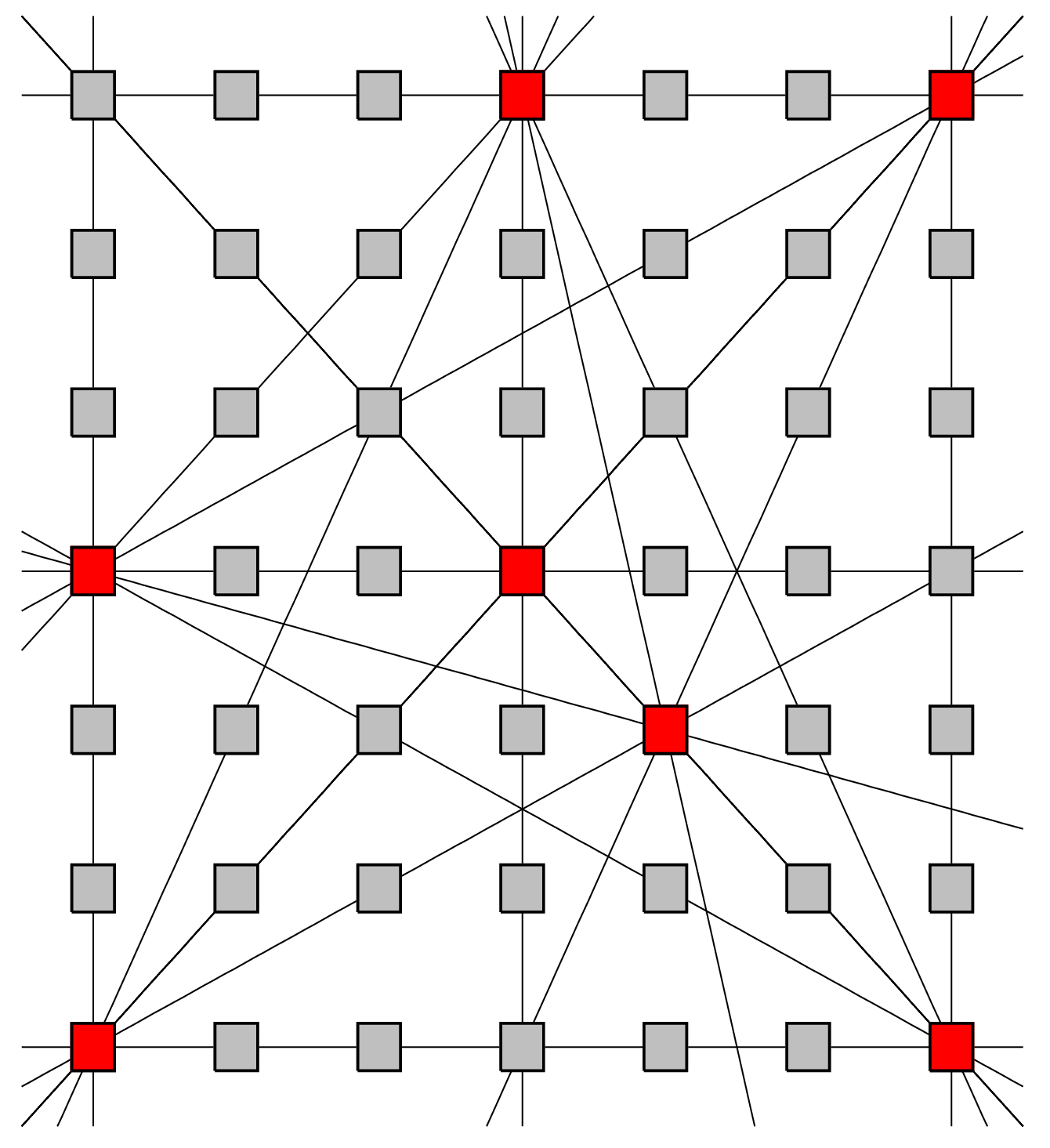}
\includegraphics[scale=0.3]{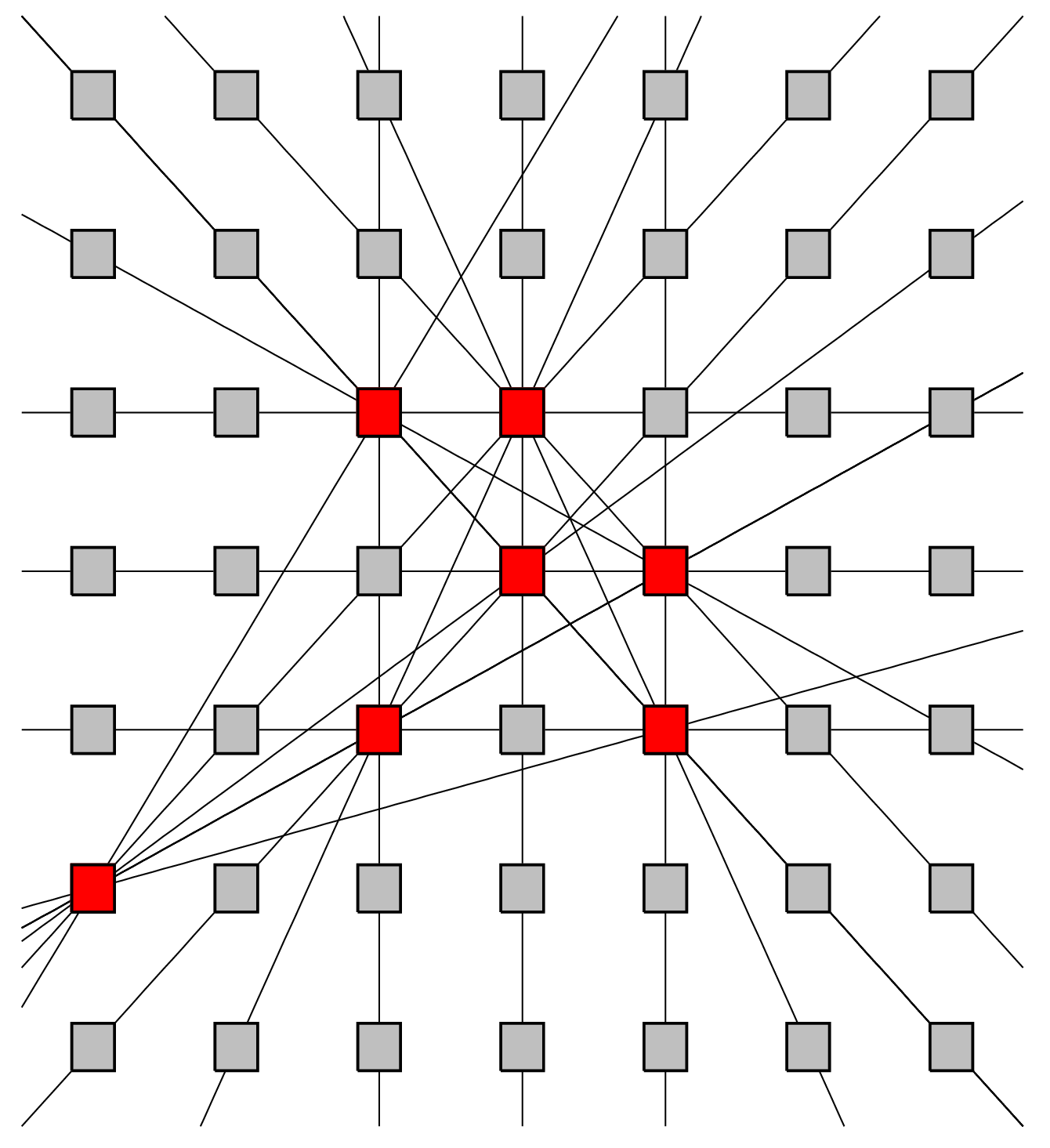}
\includegraphics[scale=0.3]{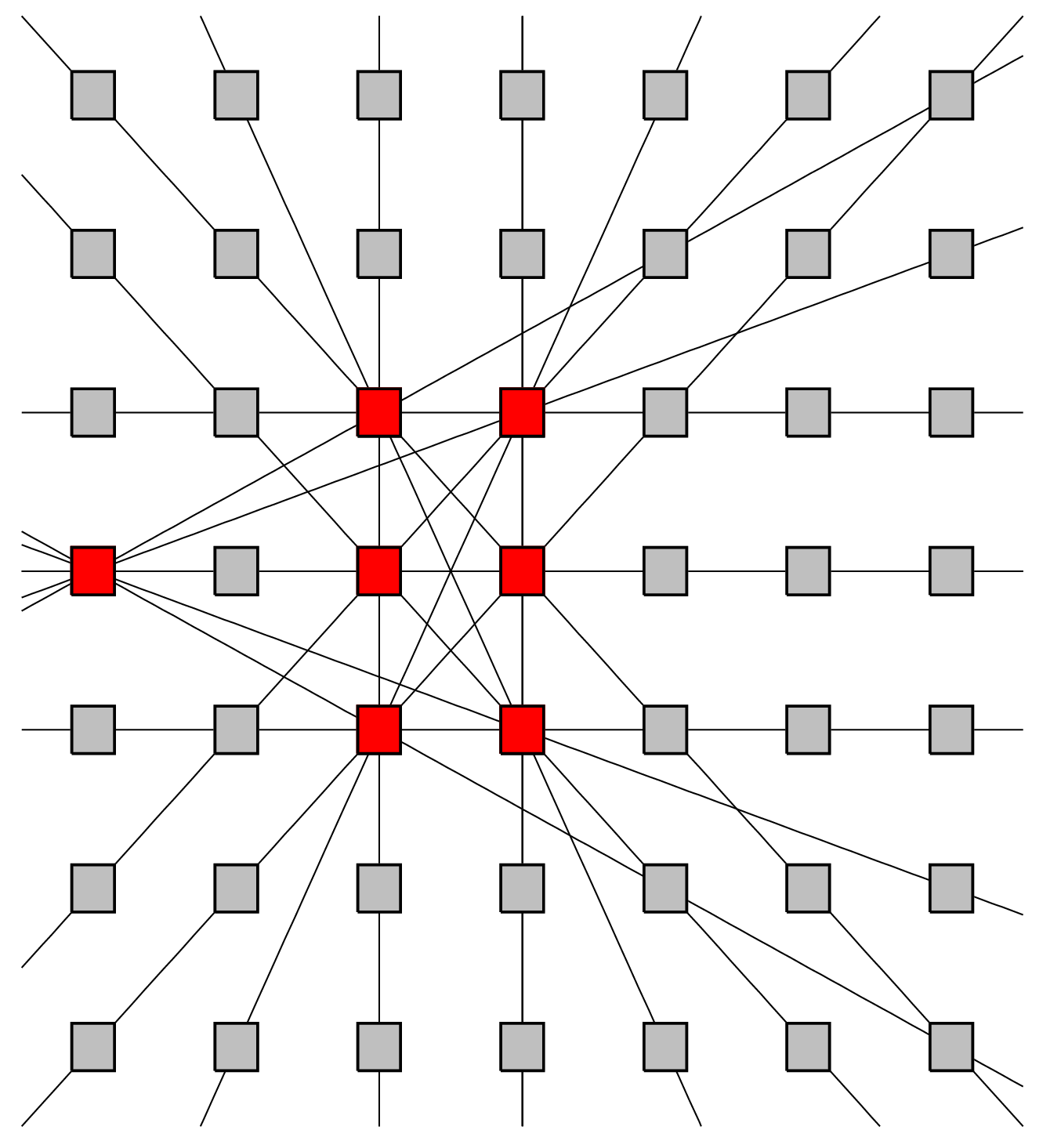}
\includegraphics[scale=0.3]{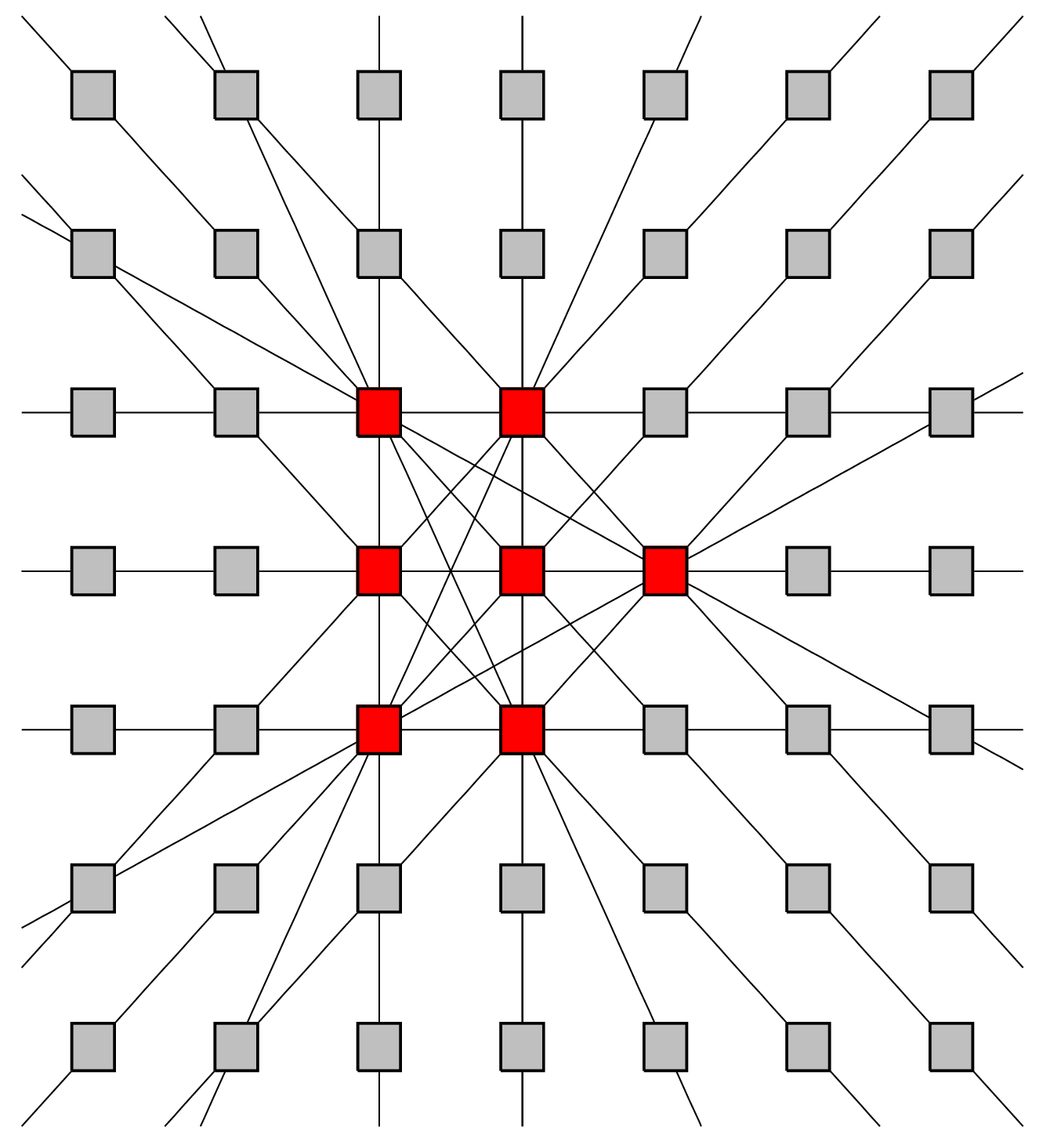}
\caption{$t(6)=7$. Two solutions have $m_d$ or $m_{d'}$ symmetry, two
others $m_x$ symmetry.}
\label{fig.cov6}
\end{figure}

Moving on to $n=6$, some (but perhaps not all) choices of coverage
by $t(6)=7$ vertices follow
in Figure \ref{fig.cov6}.

\section{Upper Bounds}\label{sec.boun}
\subsection{From Symmetric Design}

\begin{figure}[h]
\includegraphics[scale=0.3]{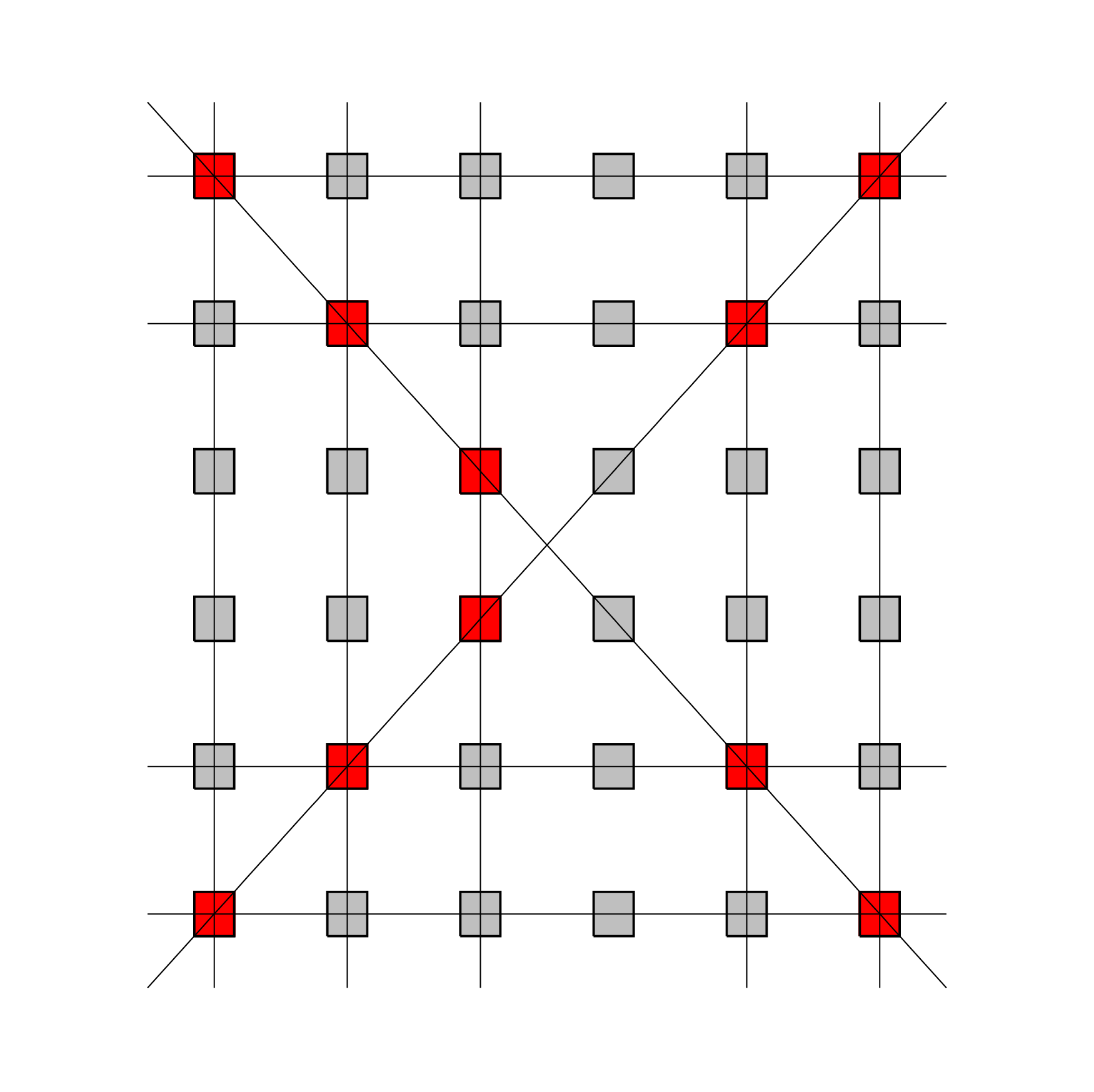}
\includegraphics[scale=0.3]{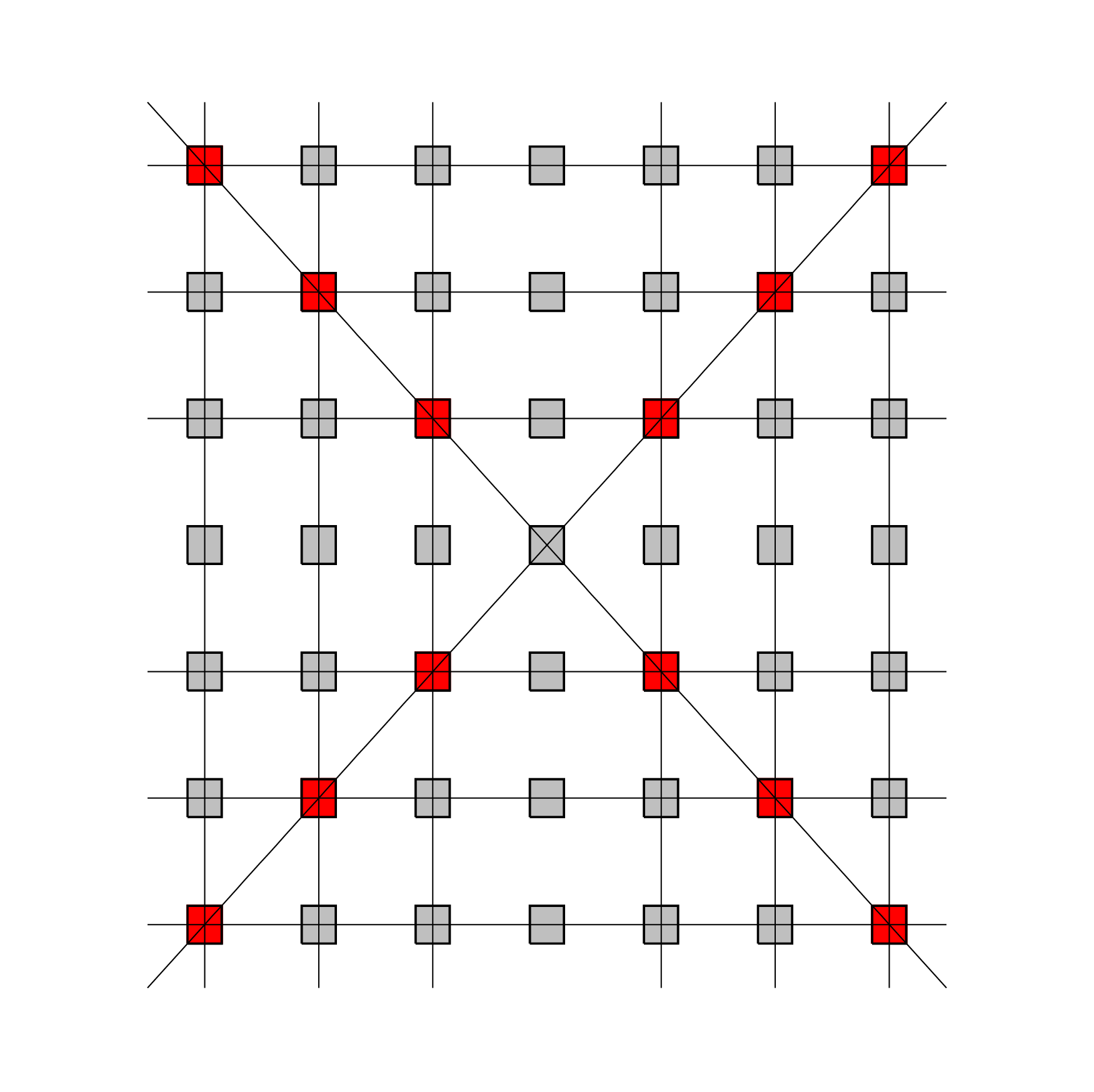}
\caption{Symmetric, generally non-optimum solutions for odd $n$ (left) or even $n$ (right).
Not all of the ${t \choose 2}$ base lines are shown,
only
the vertical, horizontal and diagonal ones which suffice to cover the lattice.}
\label{fig.symm}
\end{figure}

The upper bound 
\begin{equation}
t(n)\le 2n;\quad n\ge 2,
\label{eq.bouns}
\end{equation}
is found by selecting $2n$ of the vertices along the
two main diagonals of the lattice, with the exception of a column nearest
to the middle vertical if $n$ is odd, or with the exception of
the middle vertical if $n$ is even. Coverage is demonstrated
by Figure \ref{fig.symm}.

\subsection{From Recursion: Tapering Corners}
Once a covering for some $n$ is found, a covering for $n+1$ is obtained
by adding one row and one column to the lattice, and choosing three additional
sublattice vertices: (i) at the intersection of the new row and column, (ii)
one additional
vertex of the new row, and (iii) one additional vertex of the new
column.  These obviously add (at least) two new lines which cover the new row
and column:
\begin{equation}
t(n) \le t(n-1)+3.
\label{eq.rec}
\end{equation}
This is often a weaker bound than (\ref{eq.bouns}) because the growth
of (\ref{eq.bouns}) is only $t(n)\sim t(n-1)+2$. This kind of
relative growth estimate may result in better bounds if 
small $t(n-1)$ have been found by any of the other means.

\begin{figure}[hbt]
\includegraphics[scale=0.3]{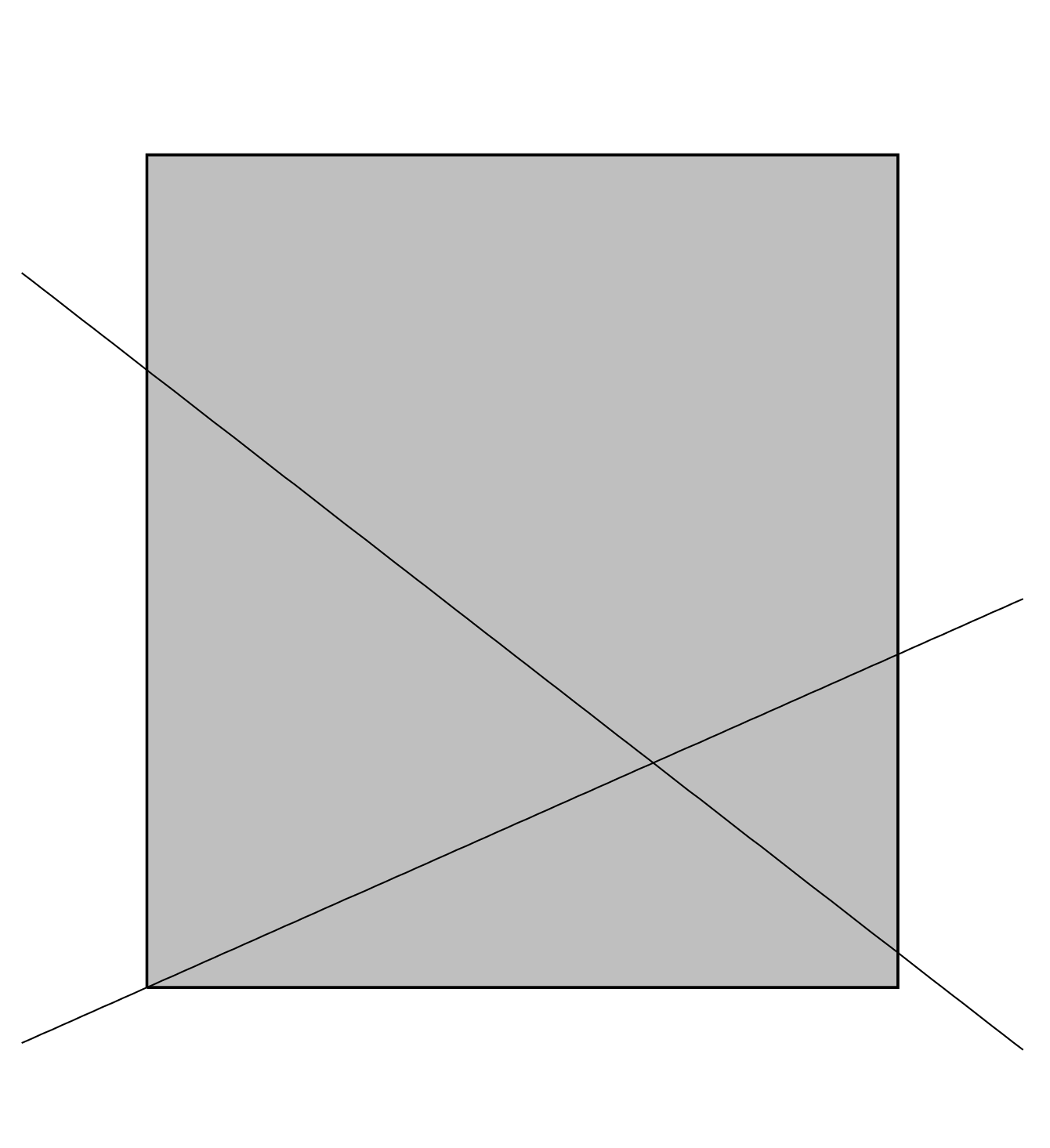}
\includegraphics[scale=0.3]{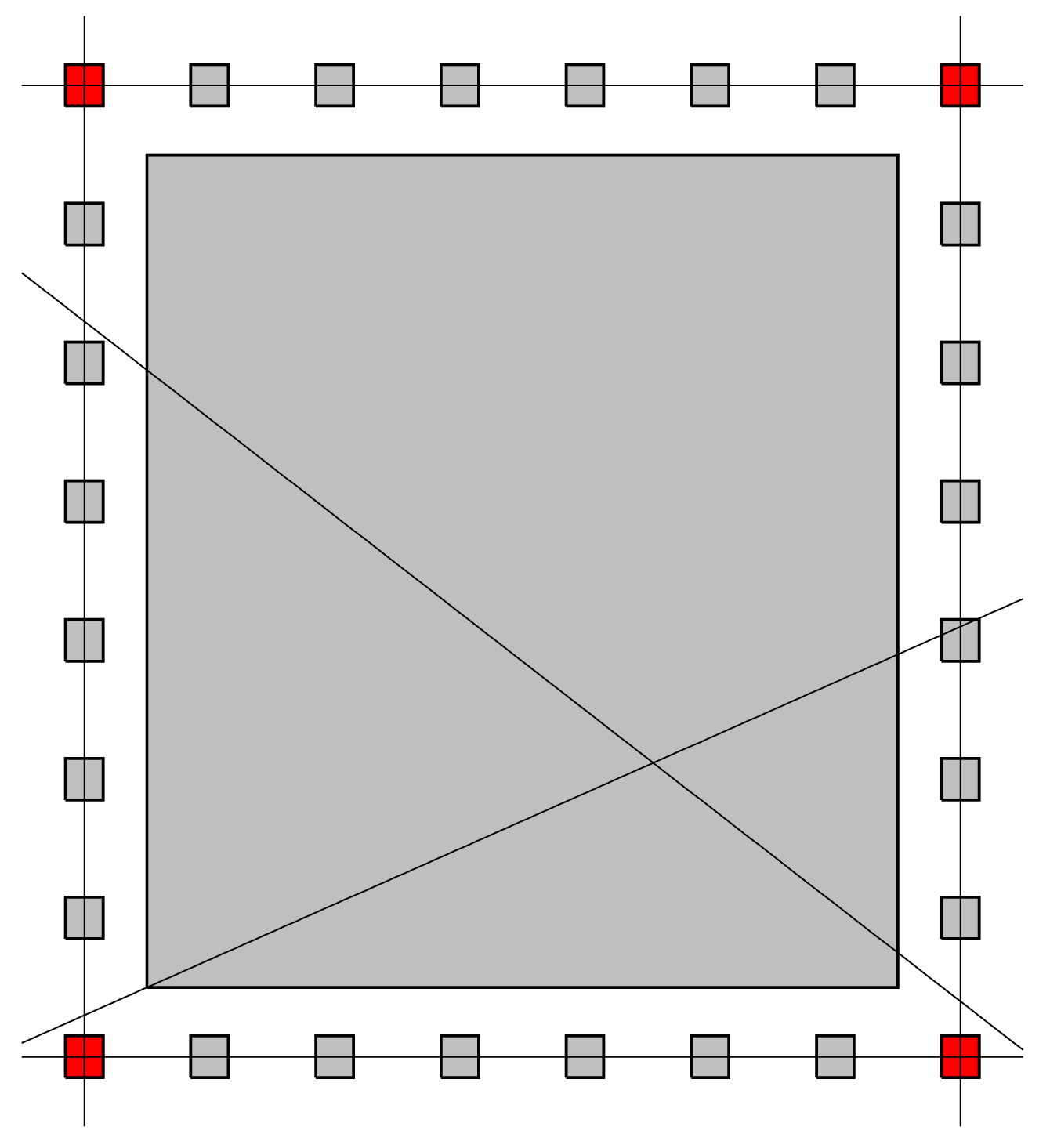}
\caption{Generation of solutions for $(n+1)\times (n+1)$ lattices
from $(n-1)\times (n-1)$ lattices associated with (\ref{eq.stein}).
The two random lines symbolize inheritance of base line coverage of
the smaller lattice (left) to the larger one (right).
}
\label{fig.stein}
\end{figure}

In an improved variant of this technique,
one may upgrade any solution for a $(n-1)\times (n-1)$ lattice
to a solution for a $(n+1)\times(n+1)$ lattice by placing four additional
sublattice vertices just outside the four corners of the smaller lattice
\cite{Steinerberger}.
Coverage of the two new rows and two new columns
that circumnavigate the four edges is ensured--see Fig.\ \ref{fig.stein}:
\begin{equation}
t(n) \le t(n-2)+4.
\label{eq.stein}
\end{equation}
Fig.\ \ref{fig.symm} could be bootstrapped by this synthesis.

Although it needs one more sublattice point than (\ref{eq.rec}), it
generally is advantageous since the index of $t$ has been advanced by 2,
not just by 1 as in (\ref{eq.rec}).

\subsection{Areal Subdivision (Tiling)}
If $n$ is odd, the lattice can be subdivided into
four square blocks (Northeast, Northwest, Southeast, Southwest) of
size $(n+1)/2\times (n+1)/2$ each. Covering each of these blocks
with $t[(n-1)/2]$ vertices ensures coverage of the
entire lattice. The same construction with partially overlapping
4 square blocks of size $(n/2+1)\times (n/2+1)$ works for even $n$, so
\begin{equation}
t(n)\le 4 t(\lfloor n/2\rfloor)
.
\end{equation}
Similar formulas can be worked out by subdivision into 9, 16, 25
etc square blocks for $n$ larger than some evident minimum.

Another aspect of this approach is that solutions equipped
with sublattice vertices at each of the 4 corners can be stacked horizontally
and vertically
with one column and one row overlap, exploiting
sublattice vertices points of these overlaps more than once. Taking the left
case of Figure \ref{fig.n3} as the building block, we arrive
at designs with $n=3,6,9,\ldots$ and the estimate
\begin{equation}
t(3i) \le (i+1)^2;\quad i=1,2,3,\ldots
\end{equation}
Taking tiles from one of the three figures at the bottom row of Fig.\ \ref{fig.n4_2}
or the top left row of Fig.\ \ref{fig.n4_3} we get
\begin{equation}
t(4i) \le (i+1)^2+2i^2;\quad i=1,2,3,\ldots
\end{equation}
where $(i+1)^2$ counts the sublattice vertices at all the block's edges and $2i^2$
the `interstitial' vertices in each of the $i^2$ blocks.
Use of the additional shared edge vertex of the tile of
the middle top row of Fig.\ \ref{fig.n4_3} leads to a slightly better estimate
\begin{equation}
t(4i) \le (i+1)^2+i(i+1);\quad i=1,2,3,\ldots
\end{equation}
The generic disadvantage of these estimates from subdivision
is that they grow with some second order polynomial of $n$, so eventually
they are not found to be competitive with the other approaches.

\subsection{From Central Star}
If one places one base line point at the lattice center ($n$ even)
or at a point closest to the center ($n$ odd), coverage of the lattice
is ensured if lines originating from there include the finite set of
inclinations $m=\Delta y/\Delta x$ needed to reach the other vertices
displaced by $\Delta x$ horizontally and $\Delta y$ vertically from
there
\cite{Faberarxiv05,FaberJNT124}.
Scanning only one quadrant from one such vintage vertex (say,
the Northeast, that is, under the constraint $\Delta y>0$, $\Delta x\ge 0$)
yields a number of directions
with `visible' vertices essentially counted by
sequence A049632 in the Online Encyclopedia of Integer Sequences \cite{EIS},
including
the central vertex itself. Moving this vertex near the middle of the
lattice takes advantage of the fact that the lines also cut backwards into the
opposite quadrant, reduces the size of the quadrant to monitor to
$\lfloor (n+3)/2 \rfloor \times \lfloor (n+3)/2 \rfloor$, and needs to 
scan one more quadrant to cover the full $(n+1)\times (n+1)$ lattice.
\begin{figure}[hbt]
\includegraphics[scale=0.35]{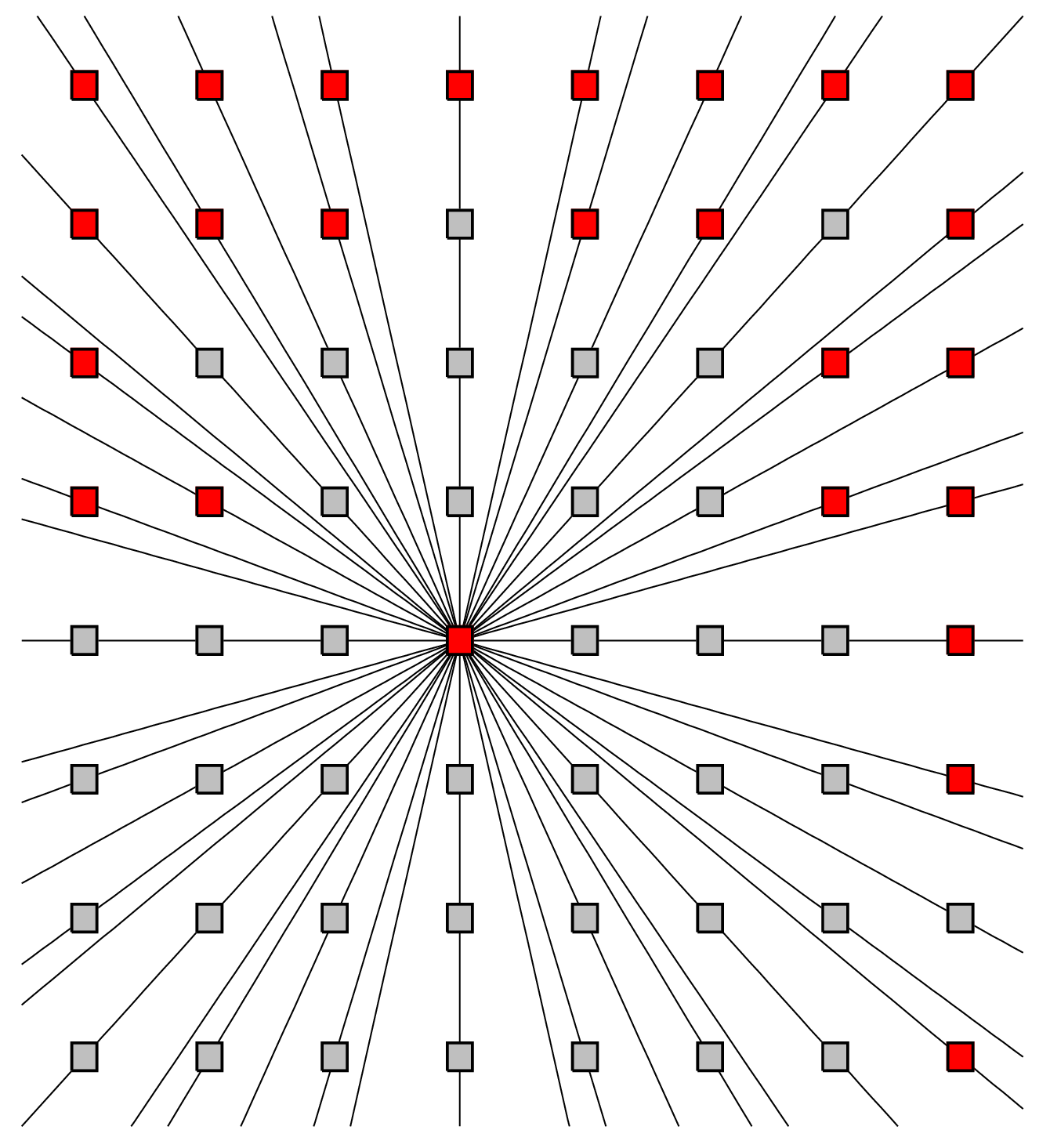}
\includegraphics[scale=0.35]{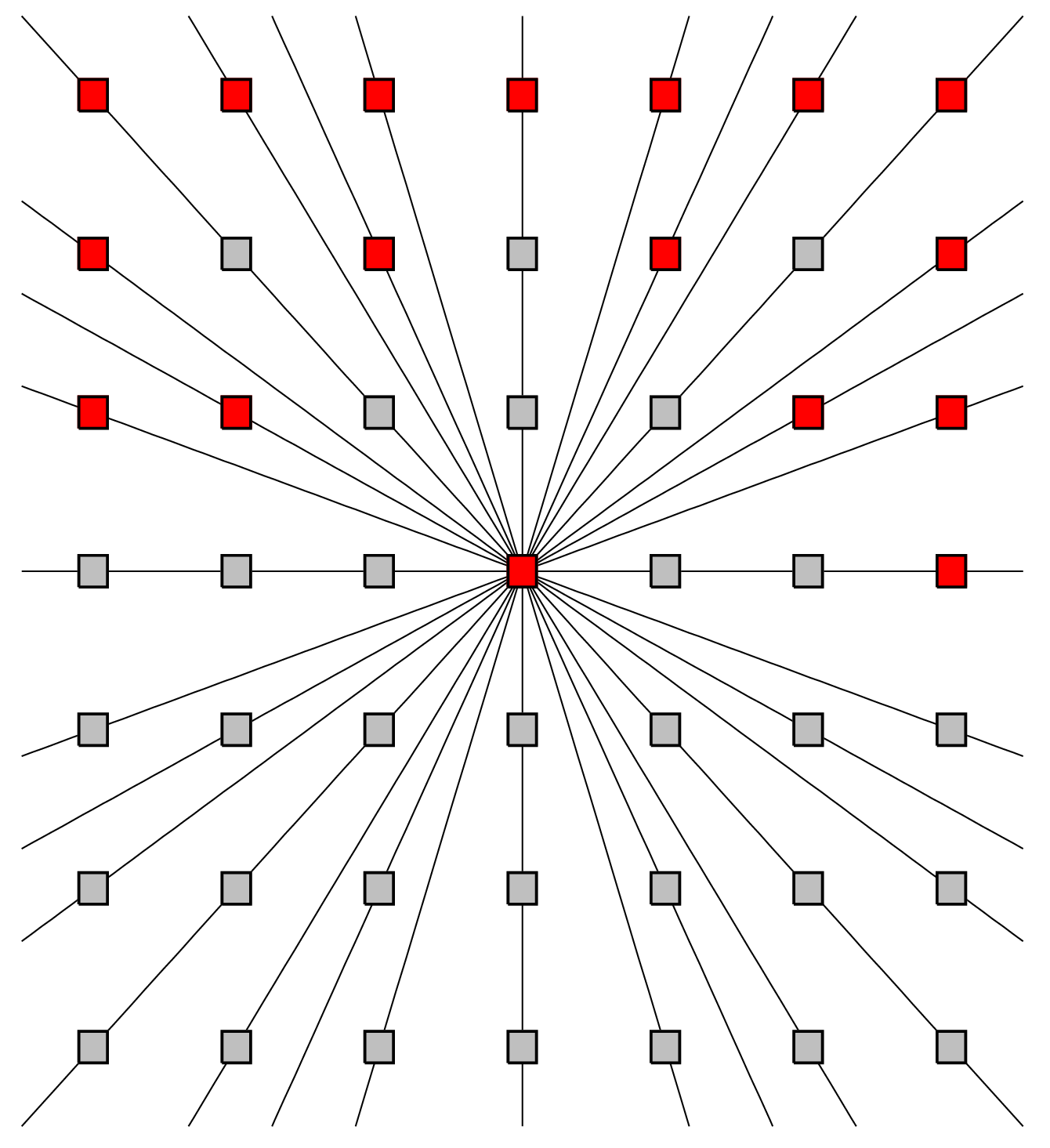}
\caption{
One central vertex and one more for each different value of
line ascensions $\Delta y/\Delta x$.
An example for odd $n$ shown at the left, for even $n$ at the right.
Not all ${t \choose 2}$ lines are shown, only those that connect
to the central vertex.}
\end{figure}
The number of reduced fractions $\Delta y/\Delta x$, plus one
for the center point, provides another upper bound,
\begin{equation}
t(n) \le 1+ 4\times A002088(\lfloor (n+1)/2 \rfloor),
\end{equation}
with reference to another sequence of Sloane's data base \cite{EIS}.

Obviously there is no improvement in comparison to (\ref{eq.bouns}) if no 
further pruning of vertices is possible.

\subsection{Monte Carlo Samples}
\begin{figure}[h]
\includegraphics[scale=0.3]{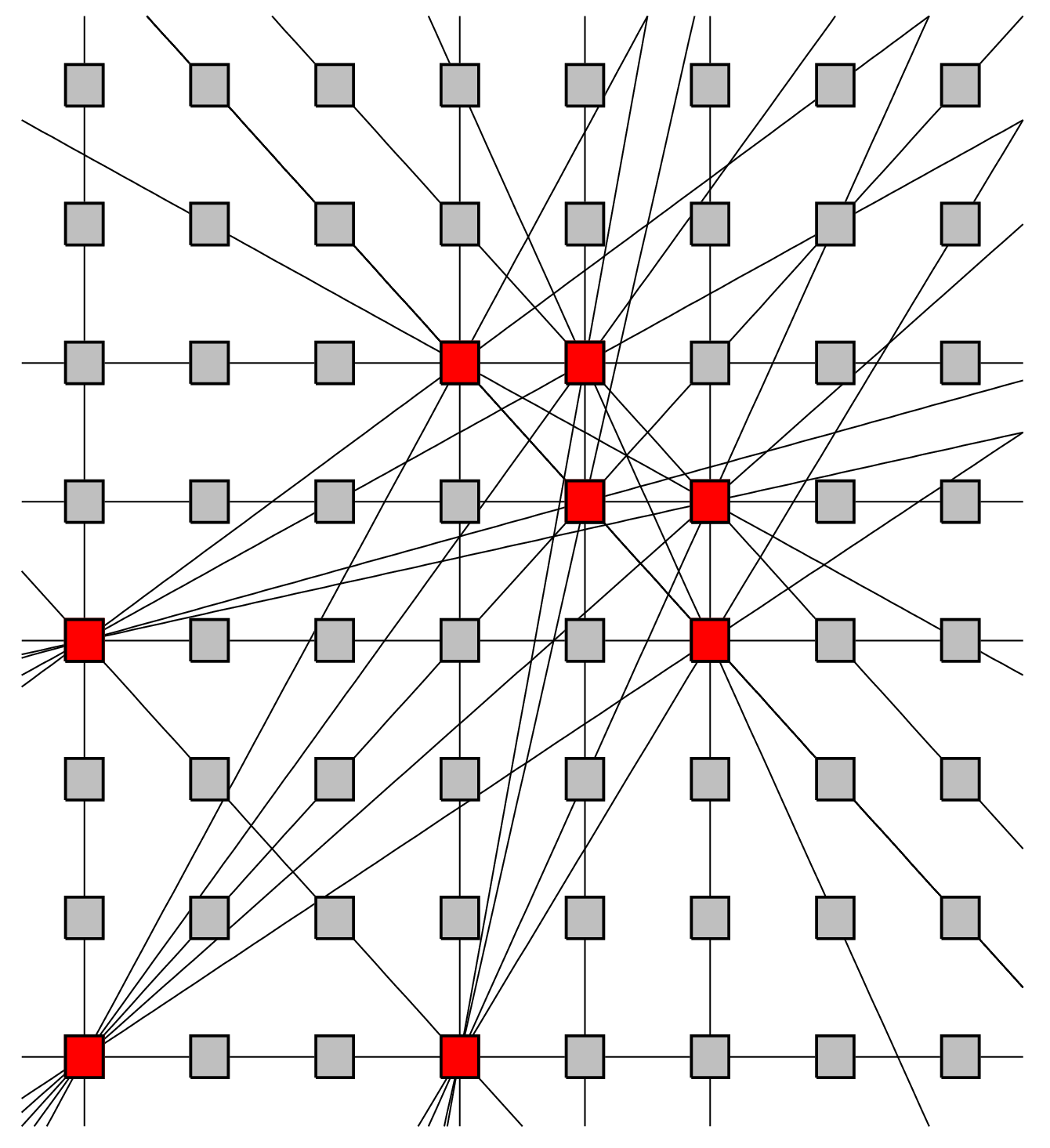}
\includegraphics[scale=0.3]{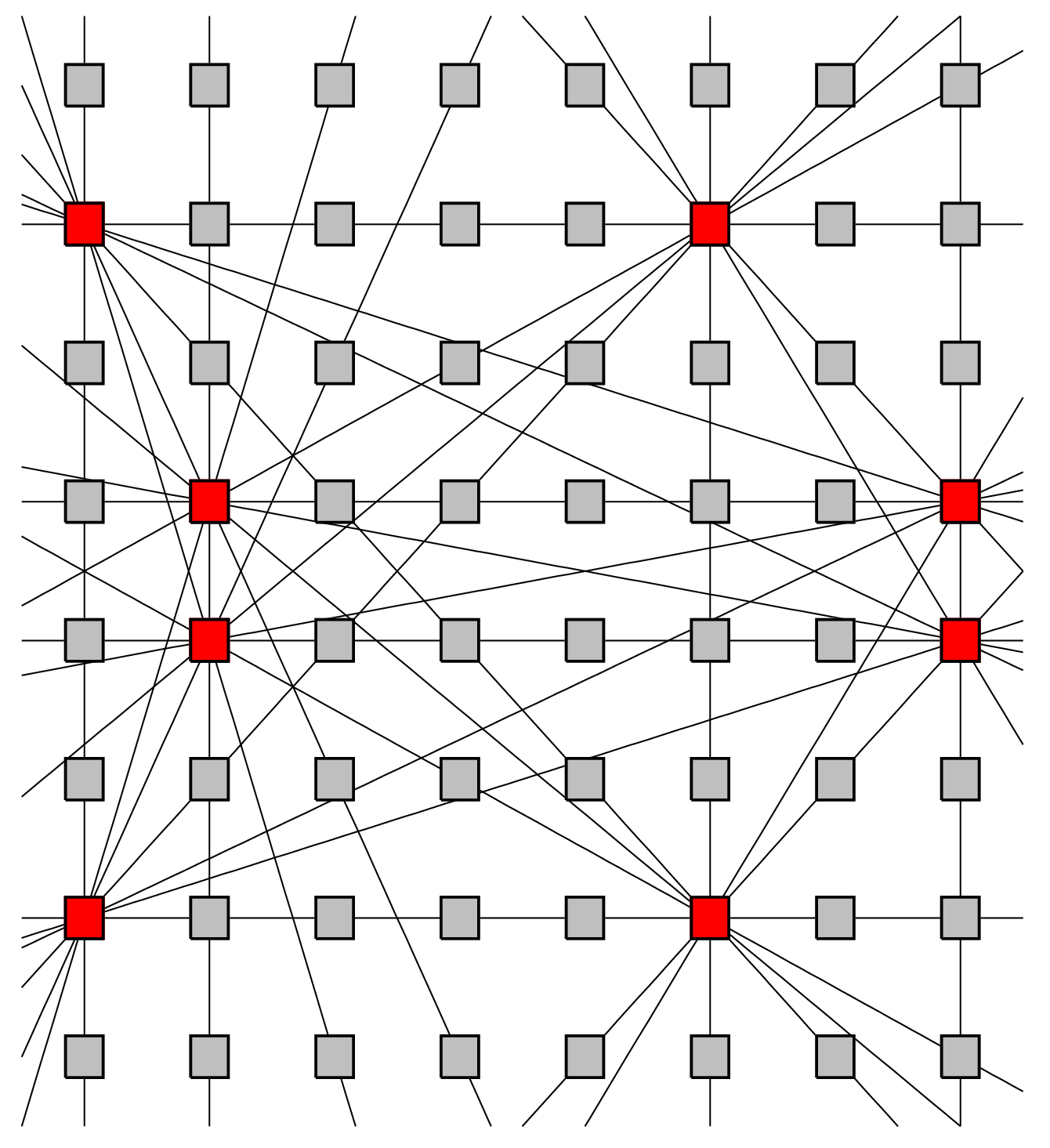}
\includegraphics[scale=0.3]{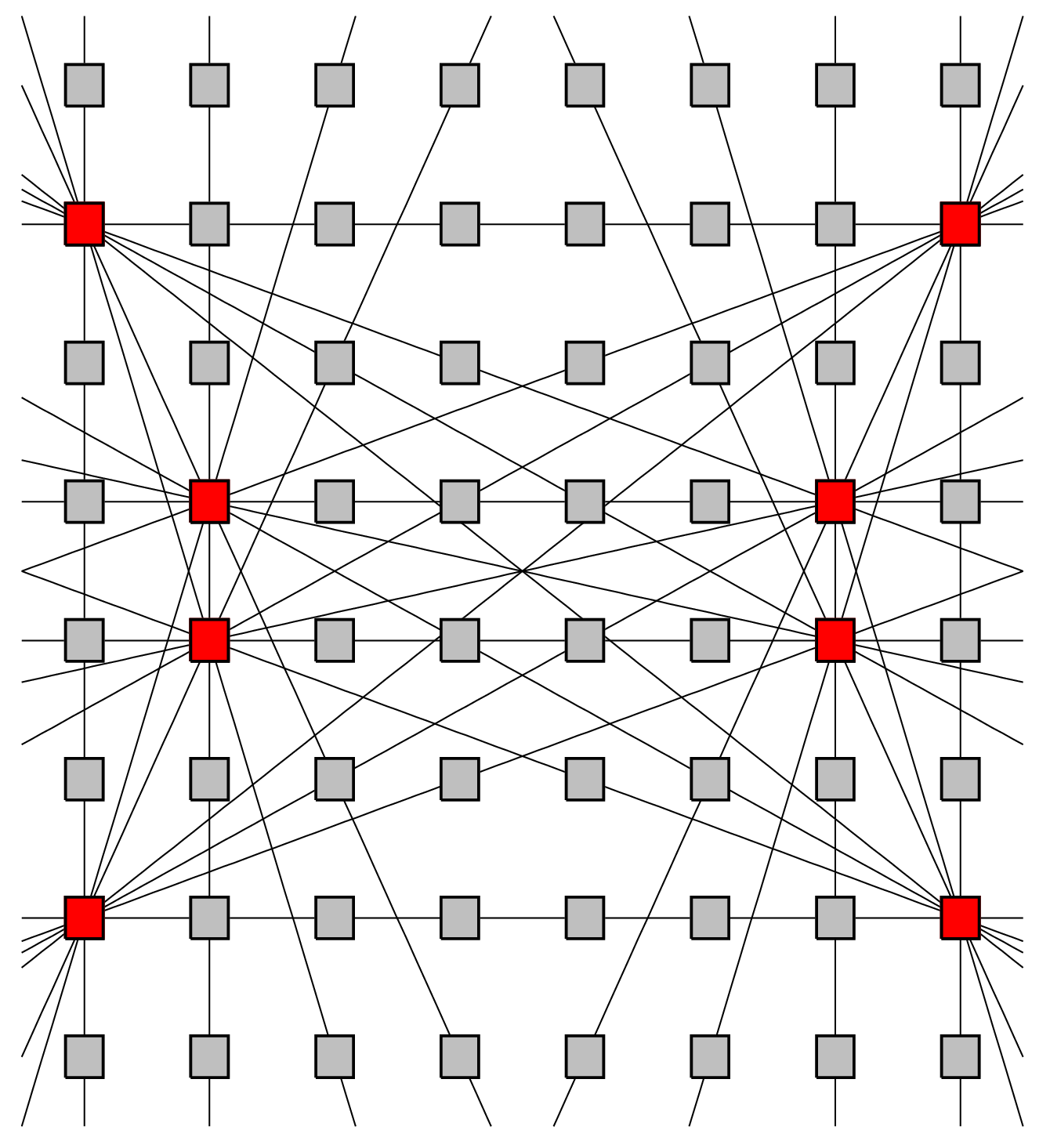}
\caption{
$t(7)\le 8$.
}
\label{fig.n7}
\end{figure}

\begin{figure}[h]
\includegraphics[scale=0.3]{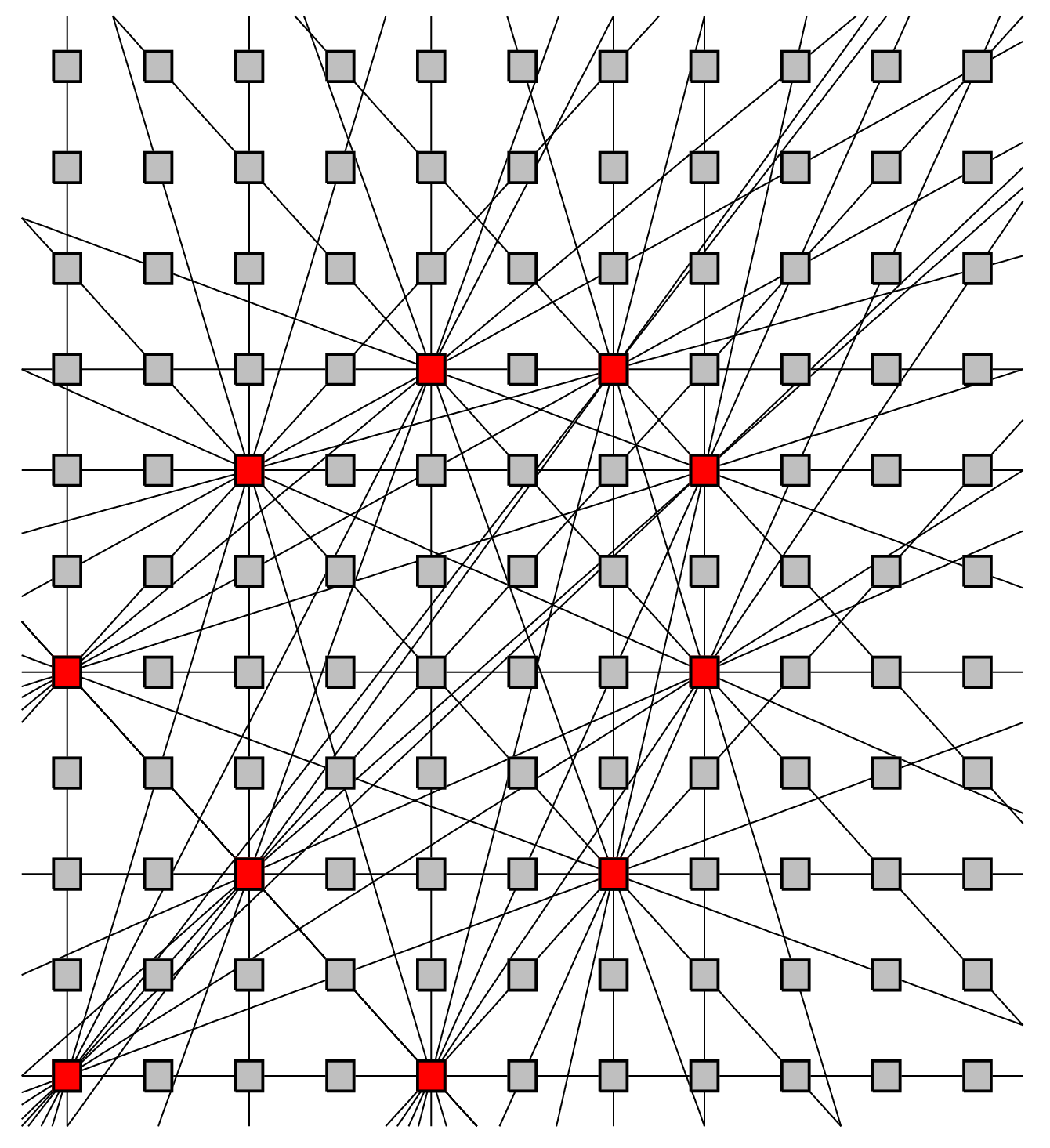}
\includegraphics[scale=0.3]{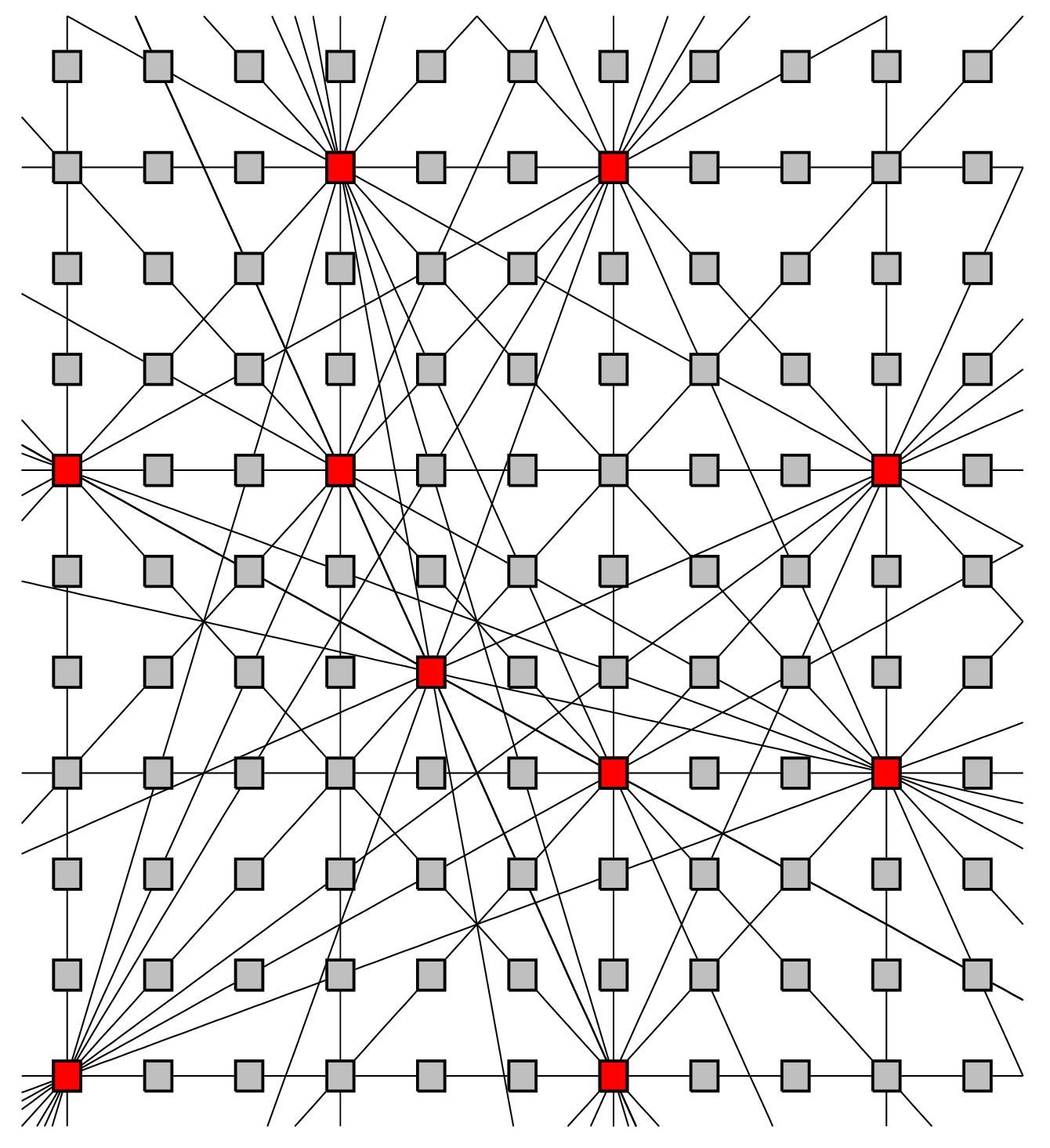}
\includegraphics[scale=0.3]{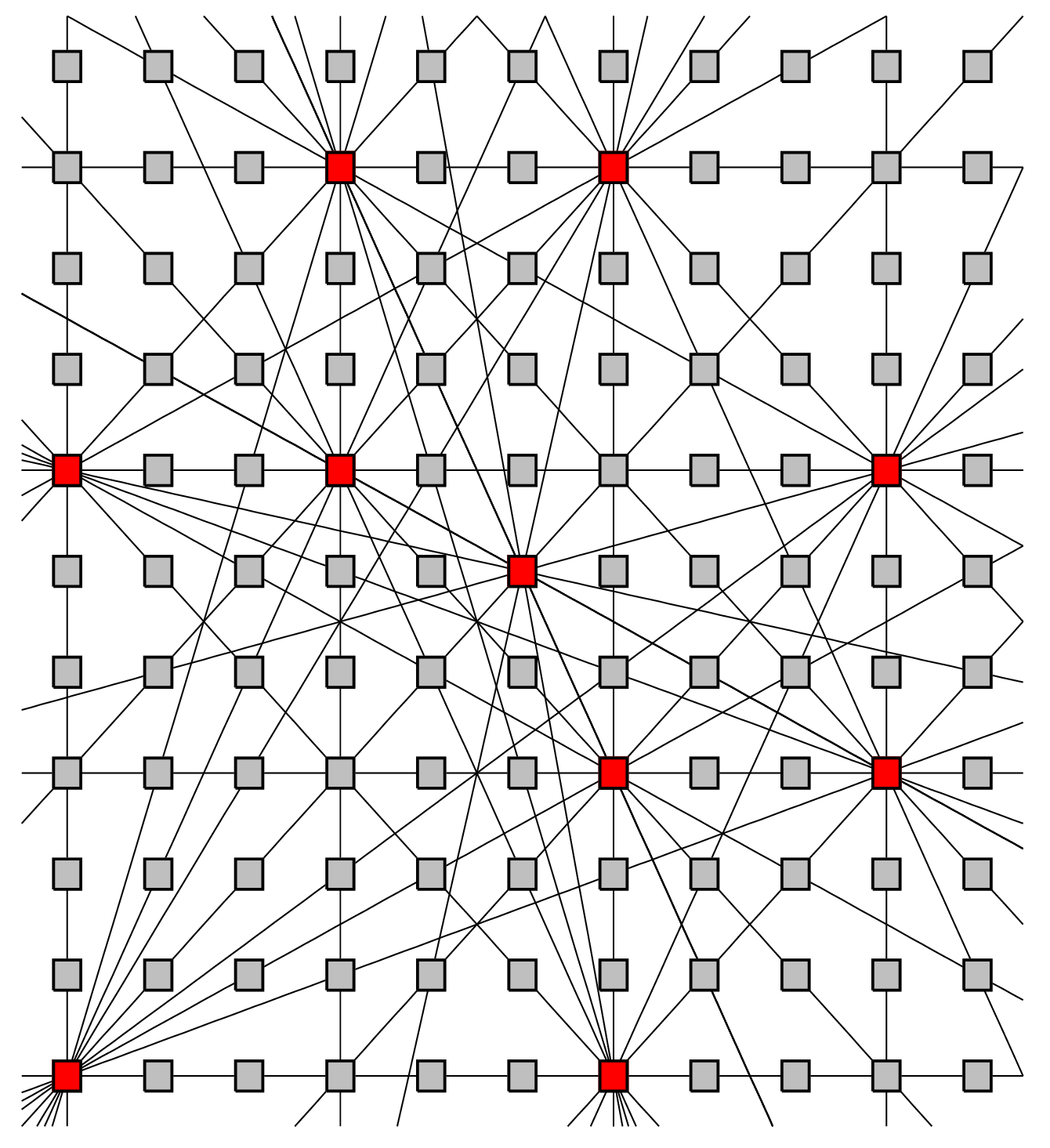}
\includegraphics[scale=0.3]{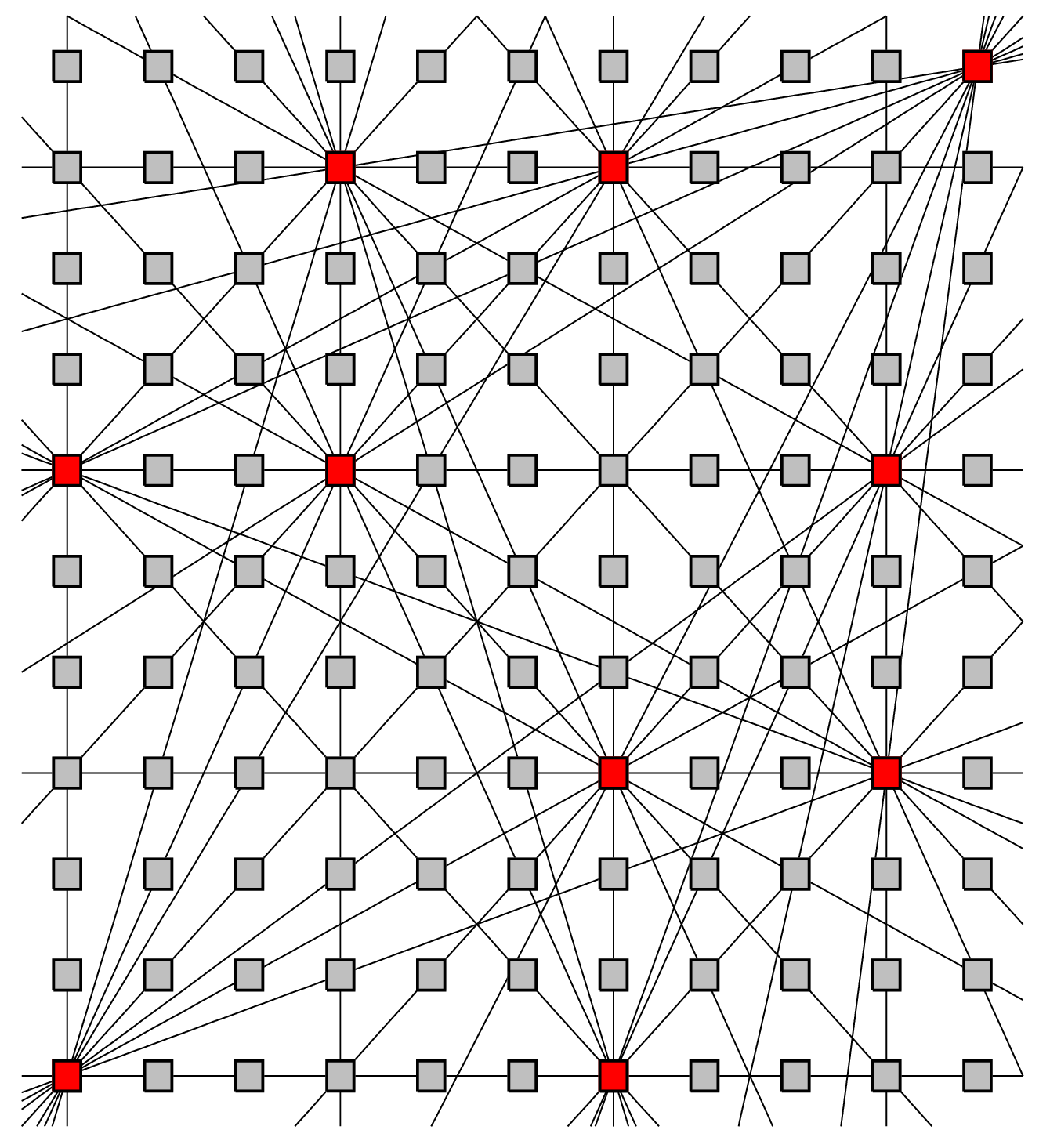}
\includegraphics[scale=0.3]{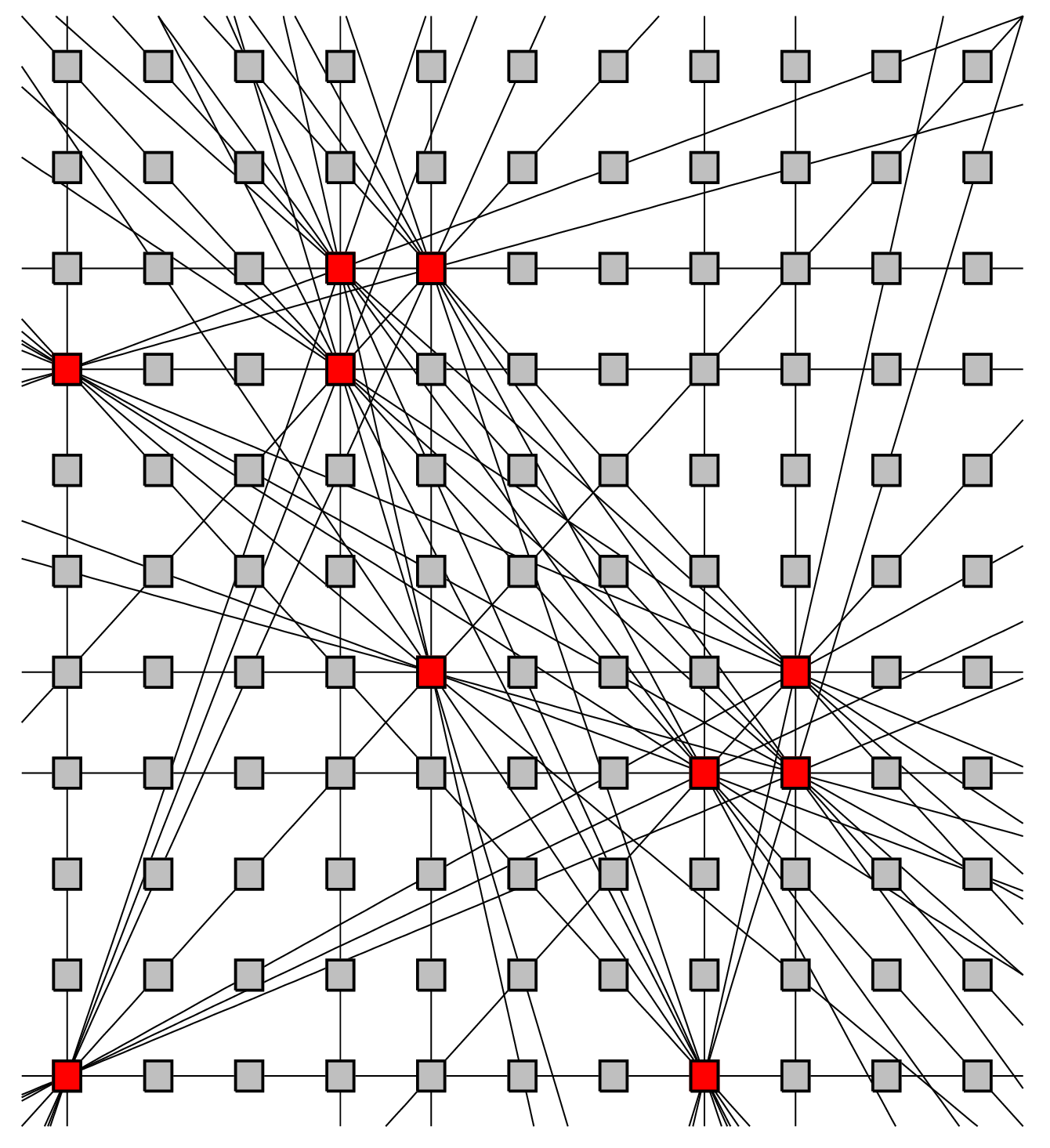}
\includegraphics[scale=0.3]{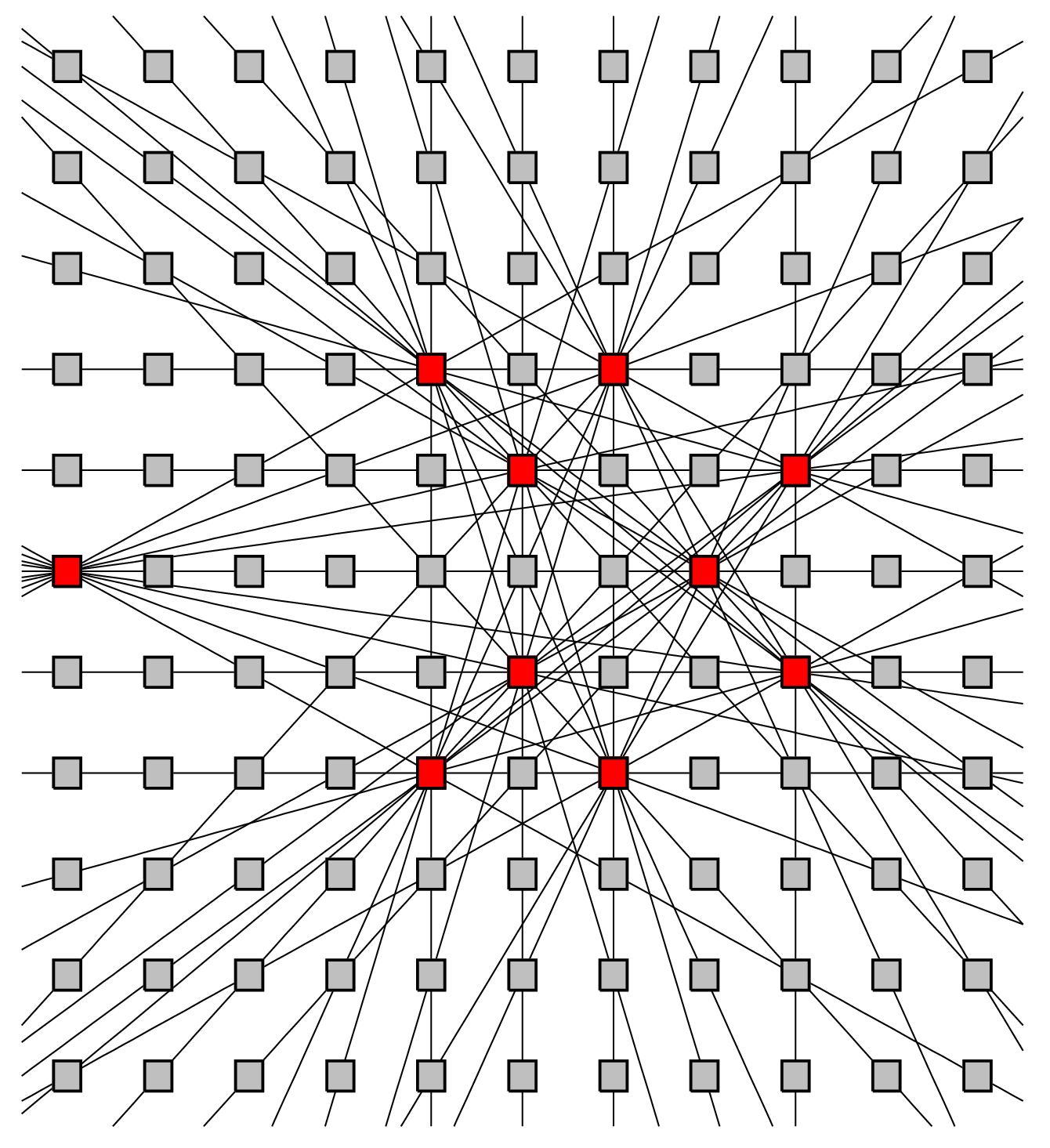}
\includegraphics[scale=0.3]{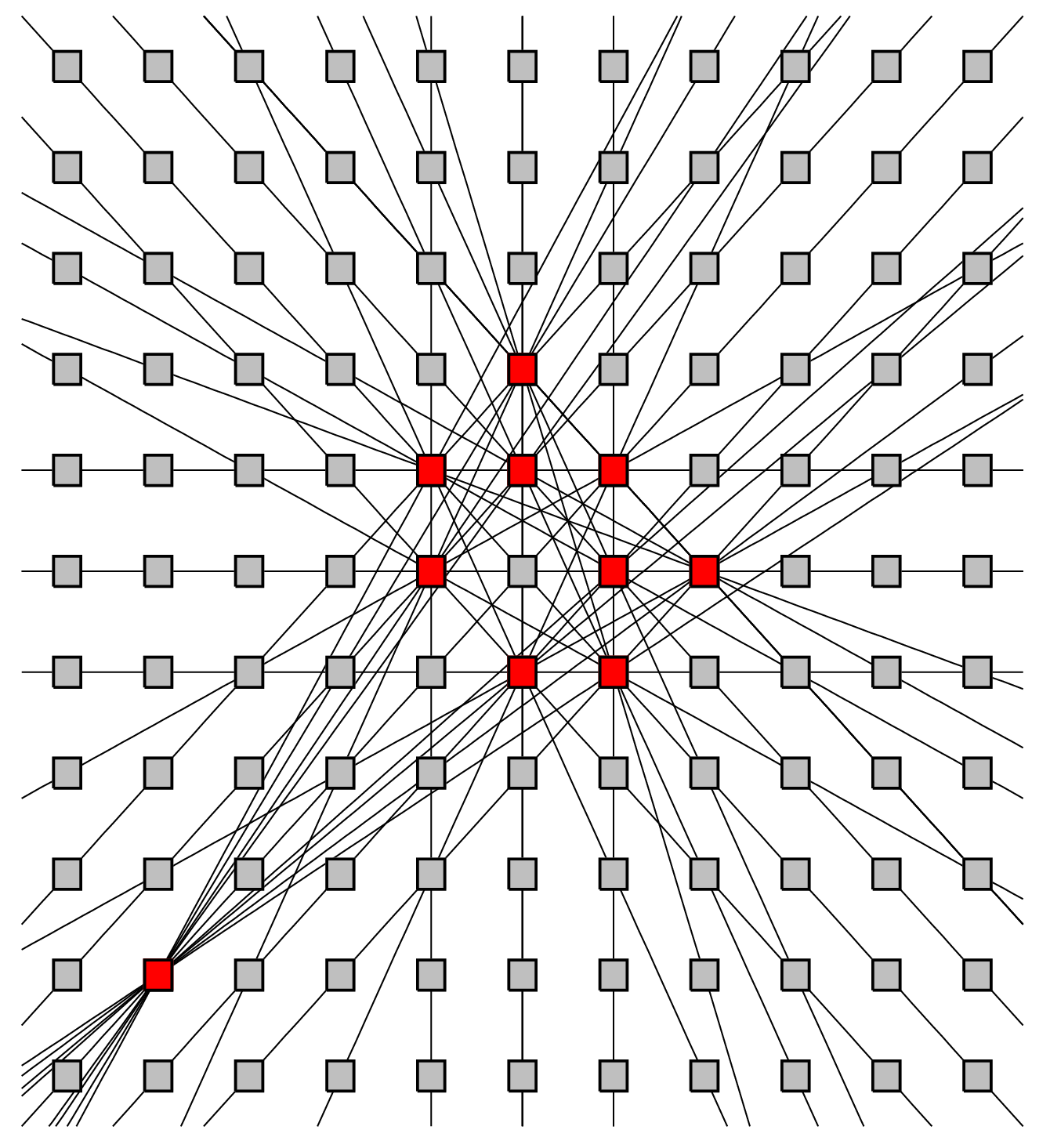}
\includegraphics[scale=0.3]{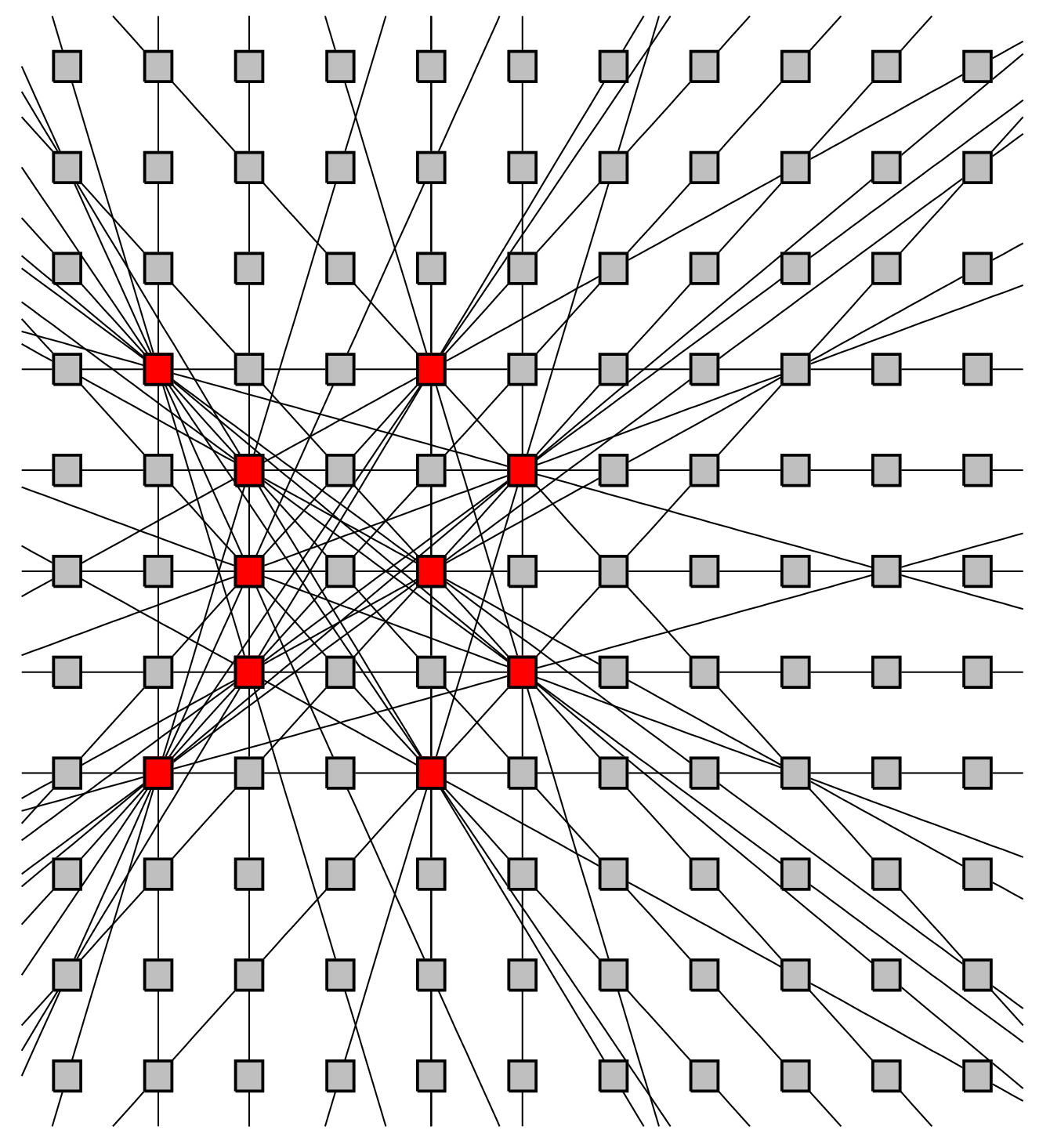}
\includegraphics[scale=0.3]{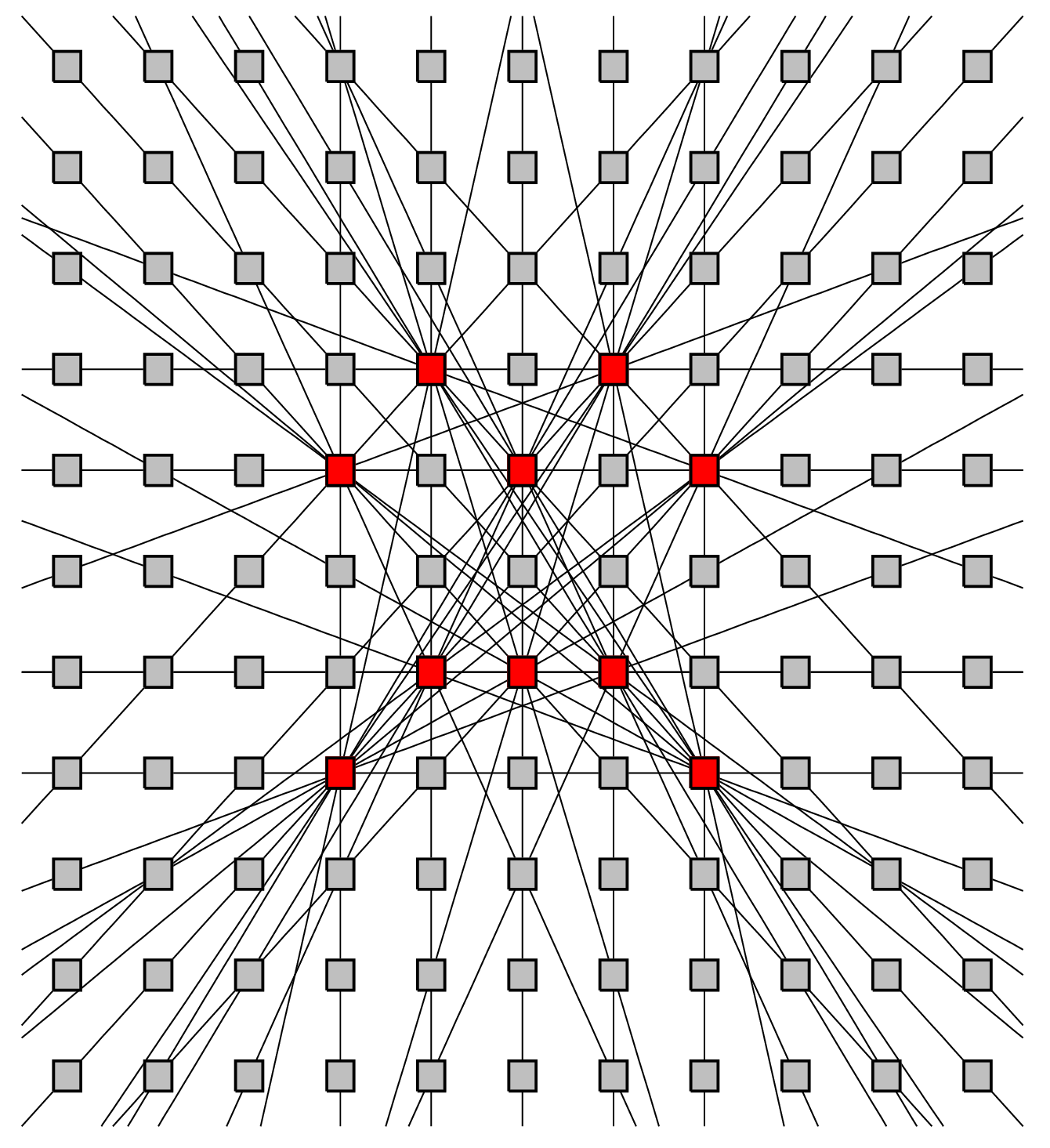}
\includegraphics[scale=0.3]{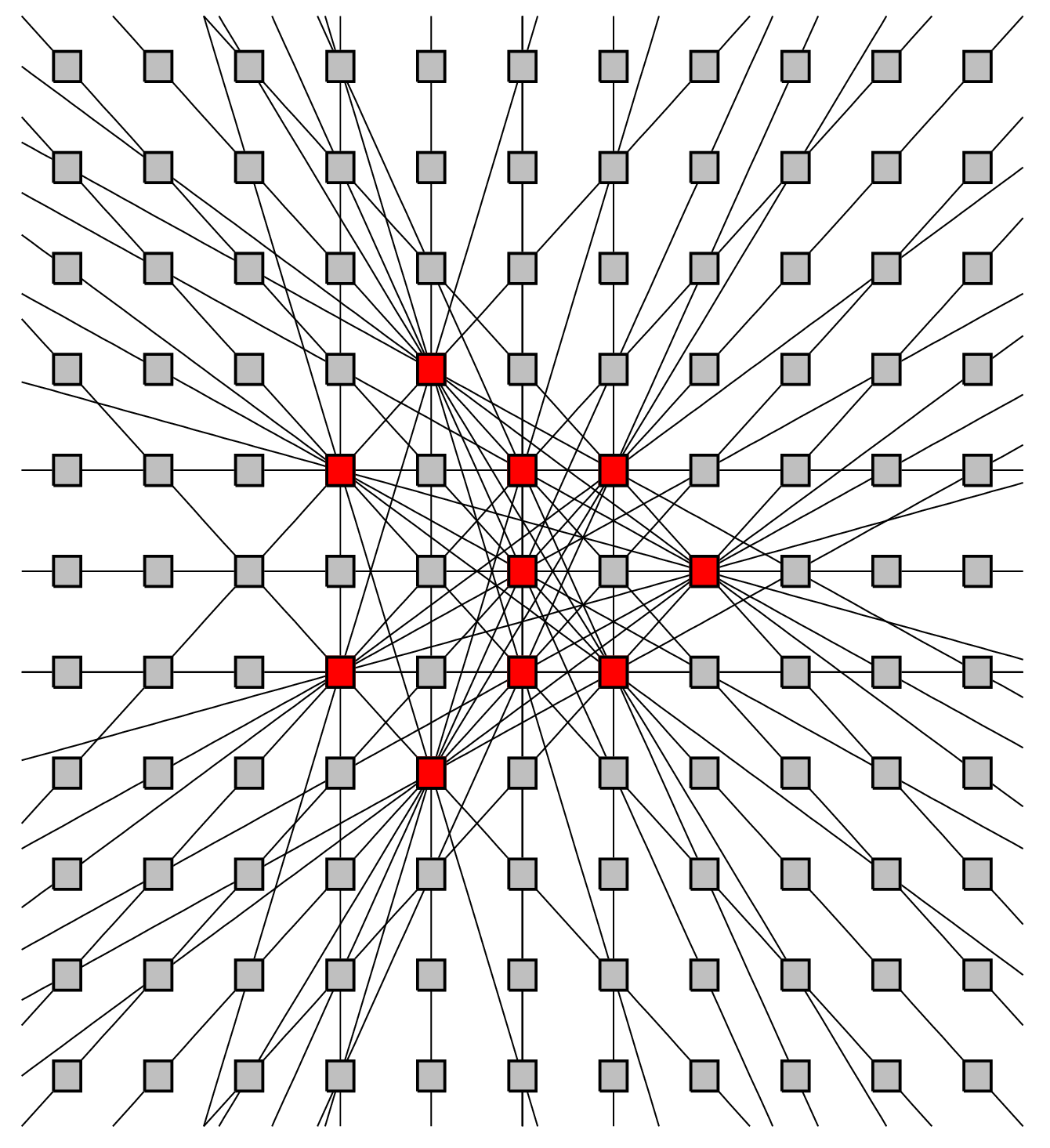}
\includegraphics[scale=0.3]{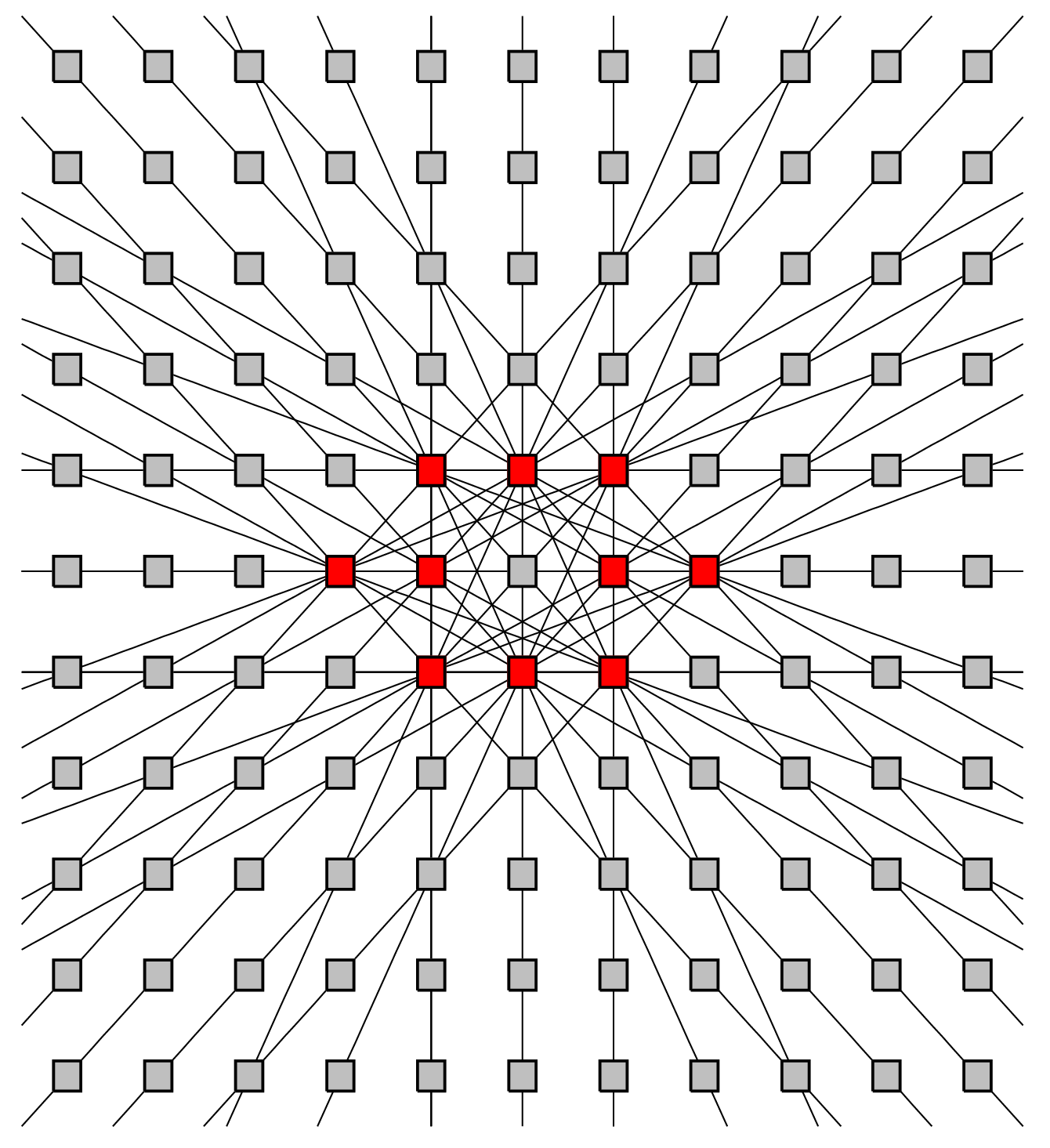}
\caption{
$t(10)\le 10$.
}
\label{fig.n10}
\end{figure}

\begin{figure}[h]
\includegraphics[scale=0.40]{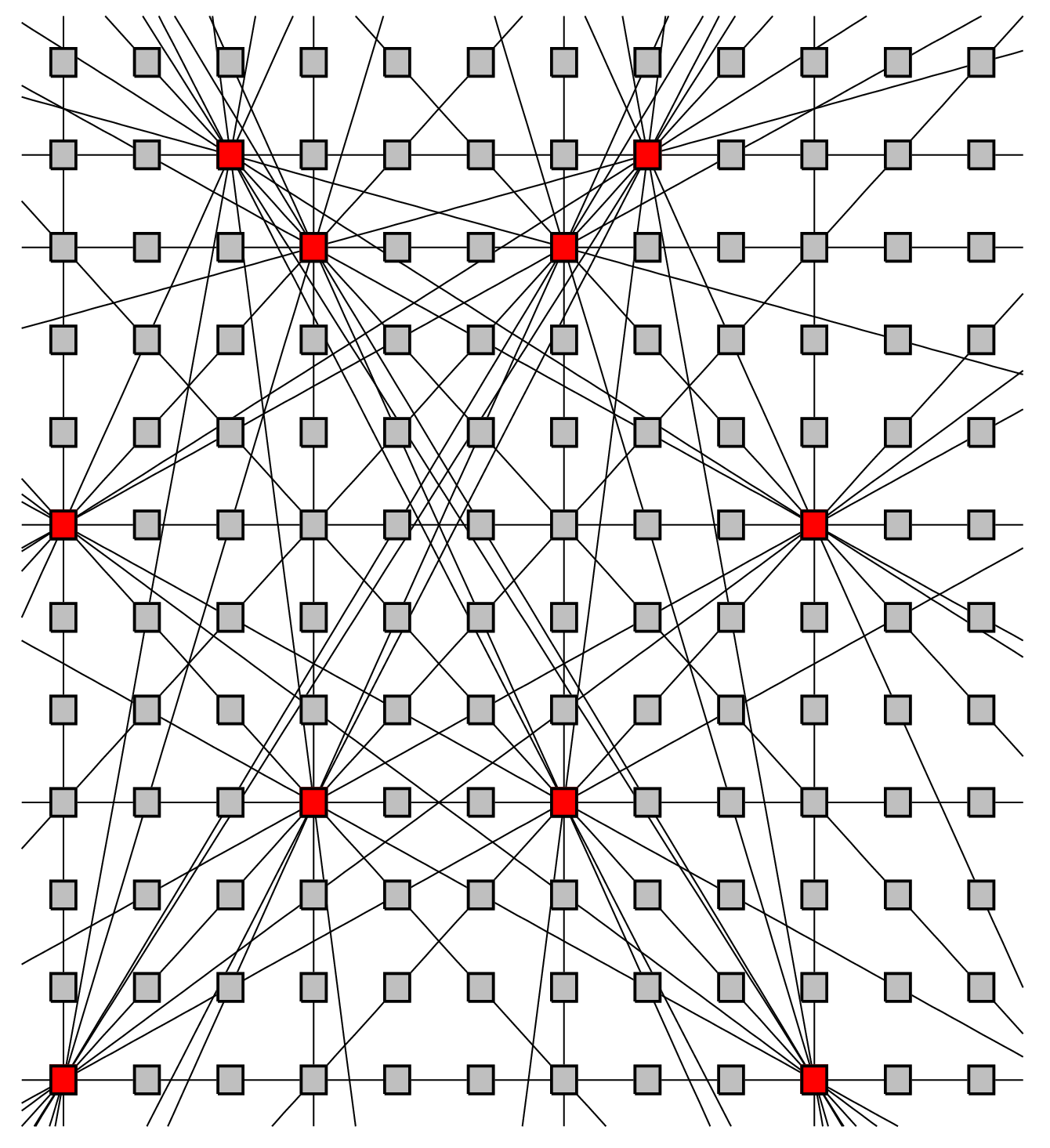}
\includegraphics[scale=0.40]{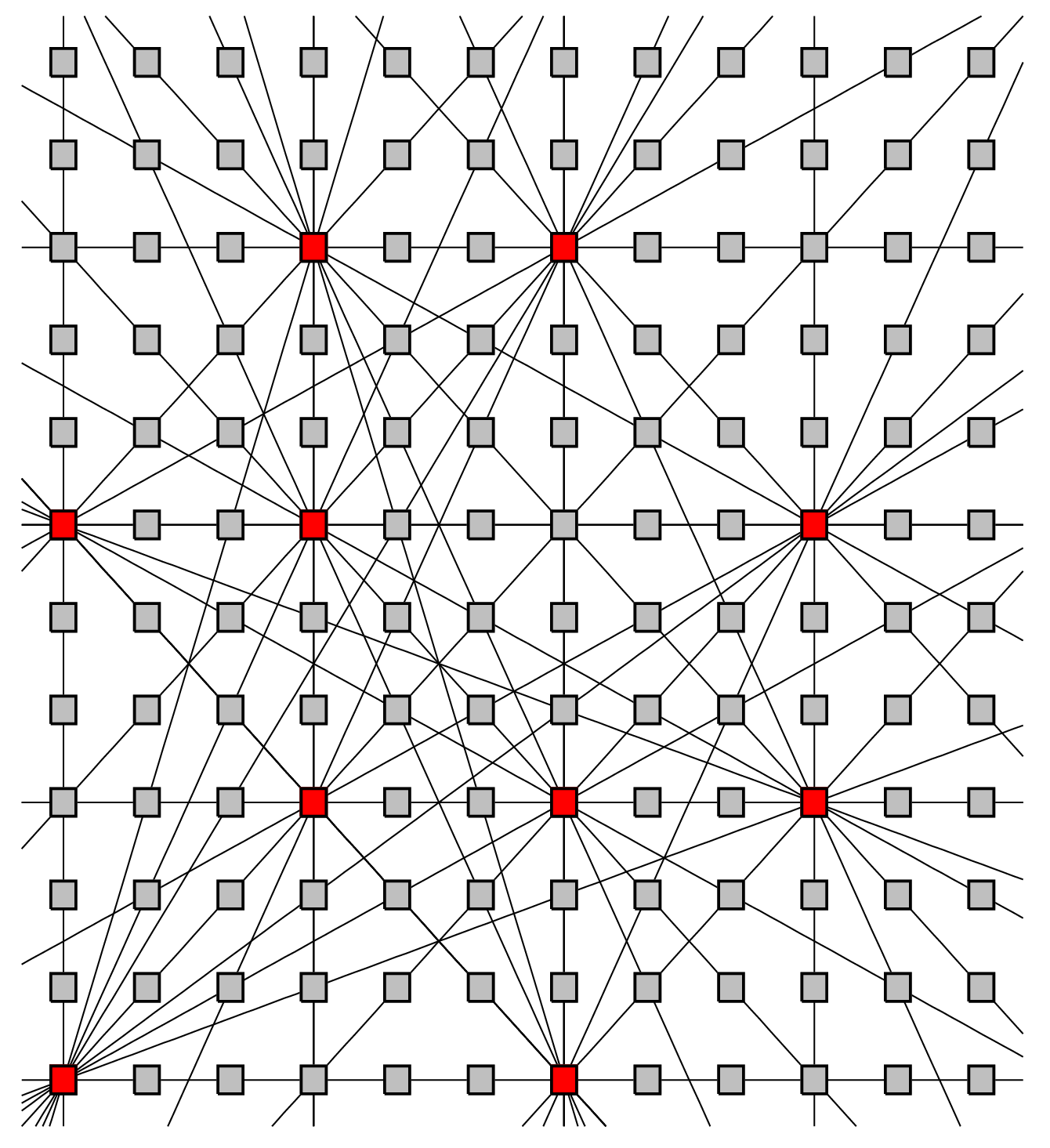}
\includegraphics[scale=0.40]{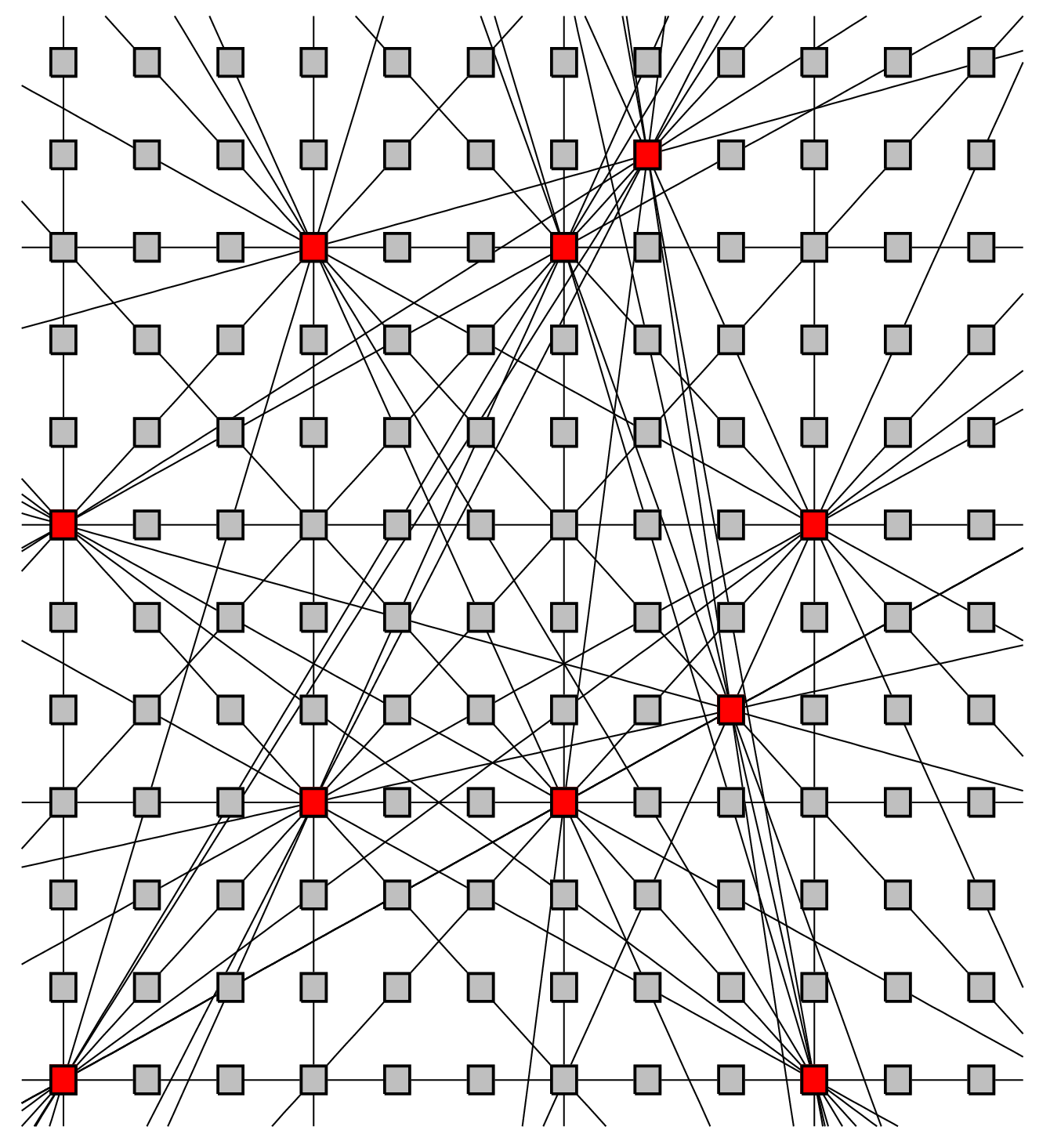}
\includegraphics[scale=0.40]{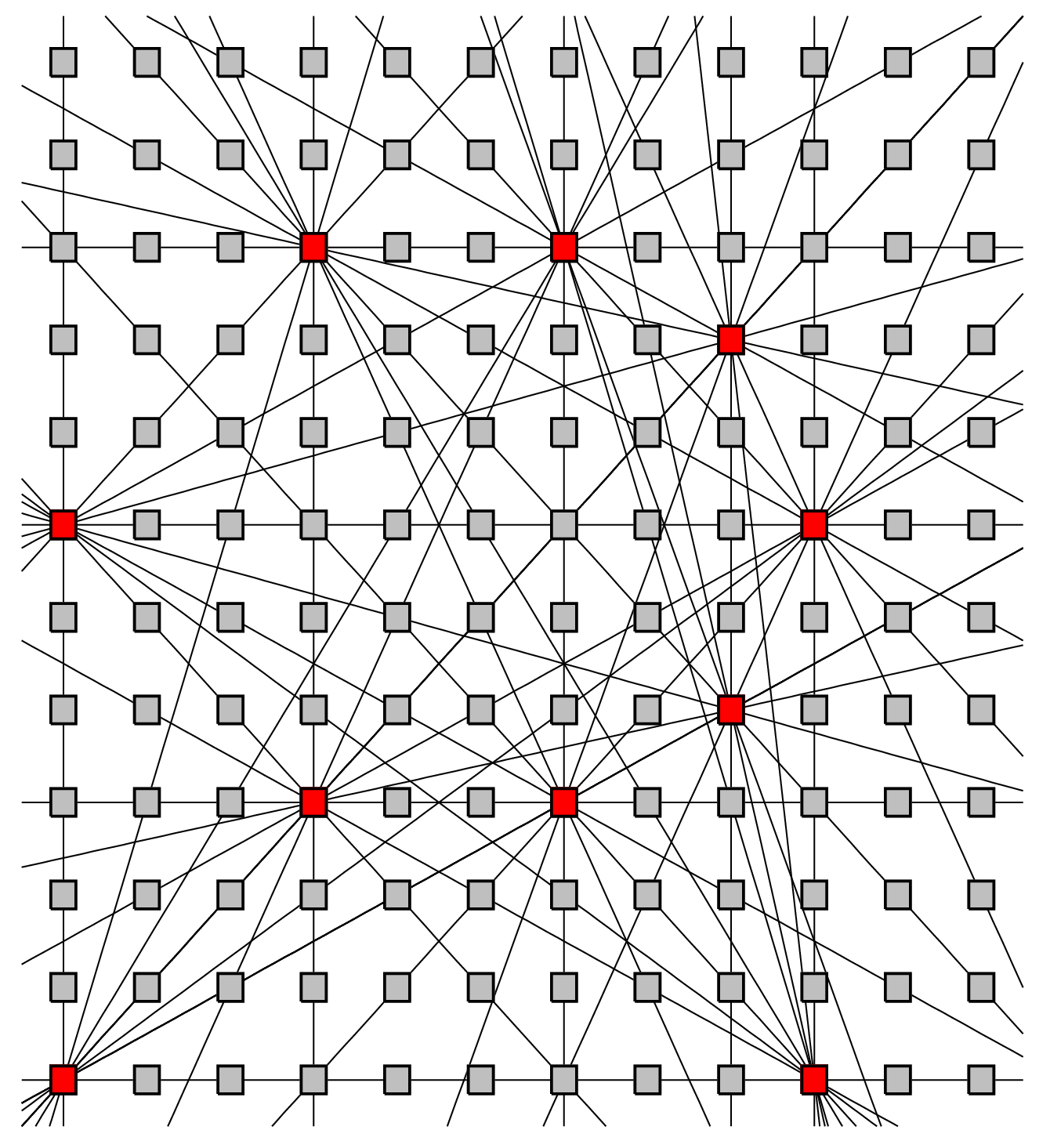}
\includegraphics[scale=0.40]{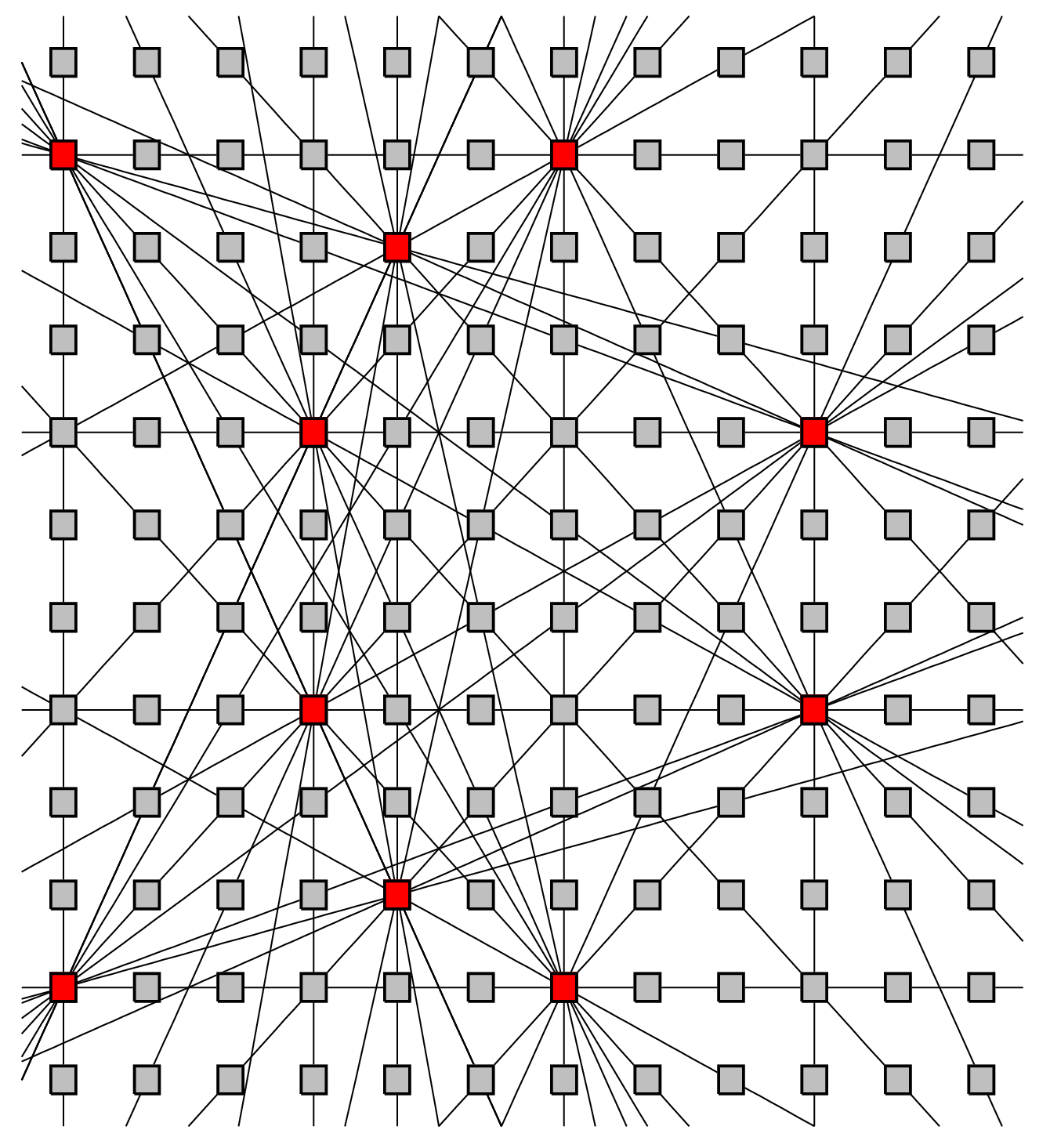}
\includegraphics[scale=0.40]{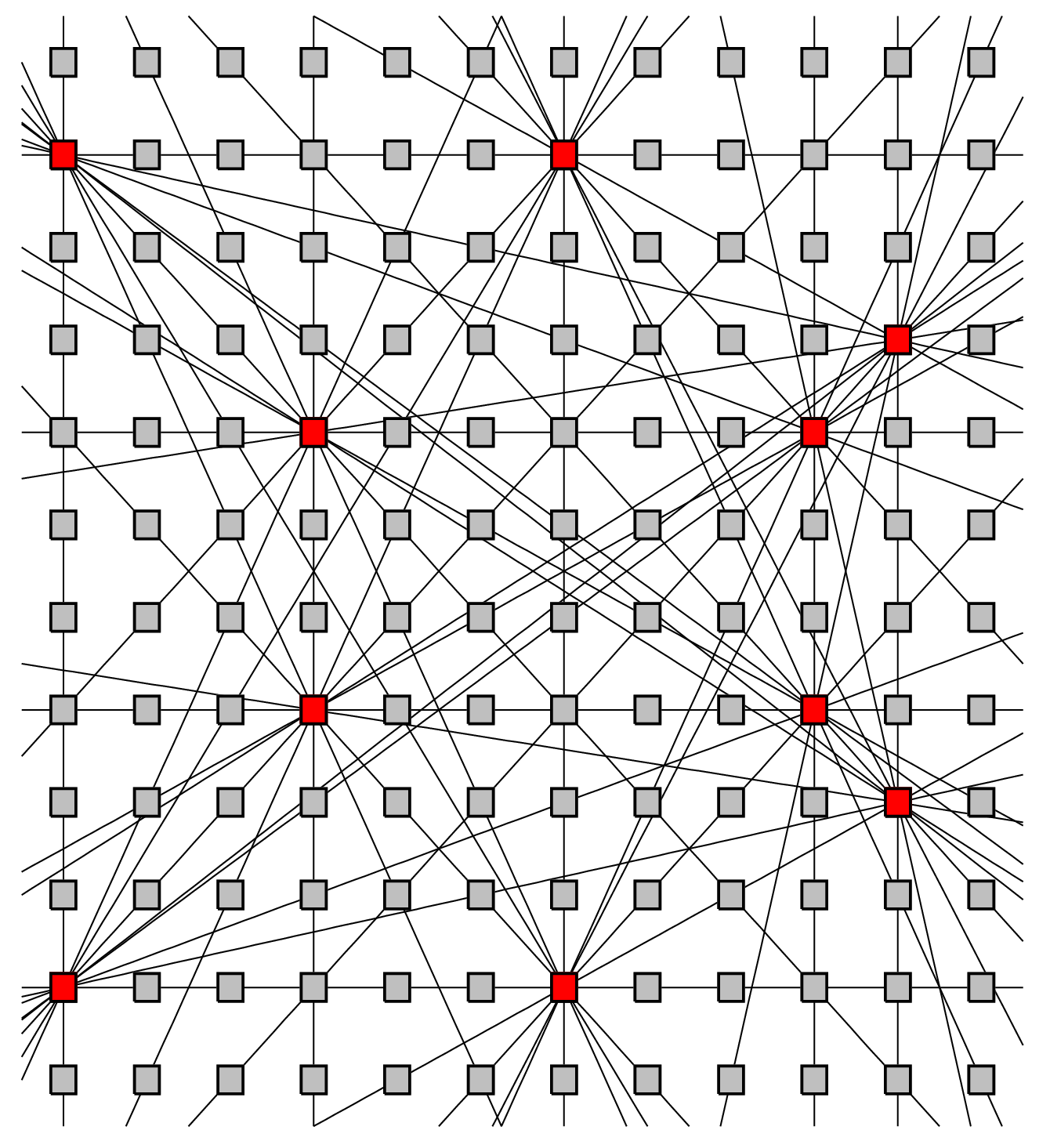}
\caption{
$t(11)\le 10$. Cutting off the top row and right column
generates solutions that might be added to Fig.\ \ref{fig.n10}.
}
\label{fig.n11}
\end{figure}

Distribution of red vertices randomly  as an alternative
to a massive complete search of all possible configurations leads to 
additional bounds, with examples in the final Figures \ref{fig.n7}--\ref{fig.n11}.
These are
results of two search strategies, one allowing for random placement
of the sublattice nodes, and a `symmetric' version which ensures that
they are placed as pairs on opposite sides of one of the two diagonals
or pairs on opposite sides of the horizontal or vertical through the center.
(This is inspired by the experience that at least one of the optimum
solutions of the smaller $n$ is of this class.)

\begin{figure}[h]
\includegraphics[scale=0.3]{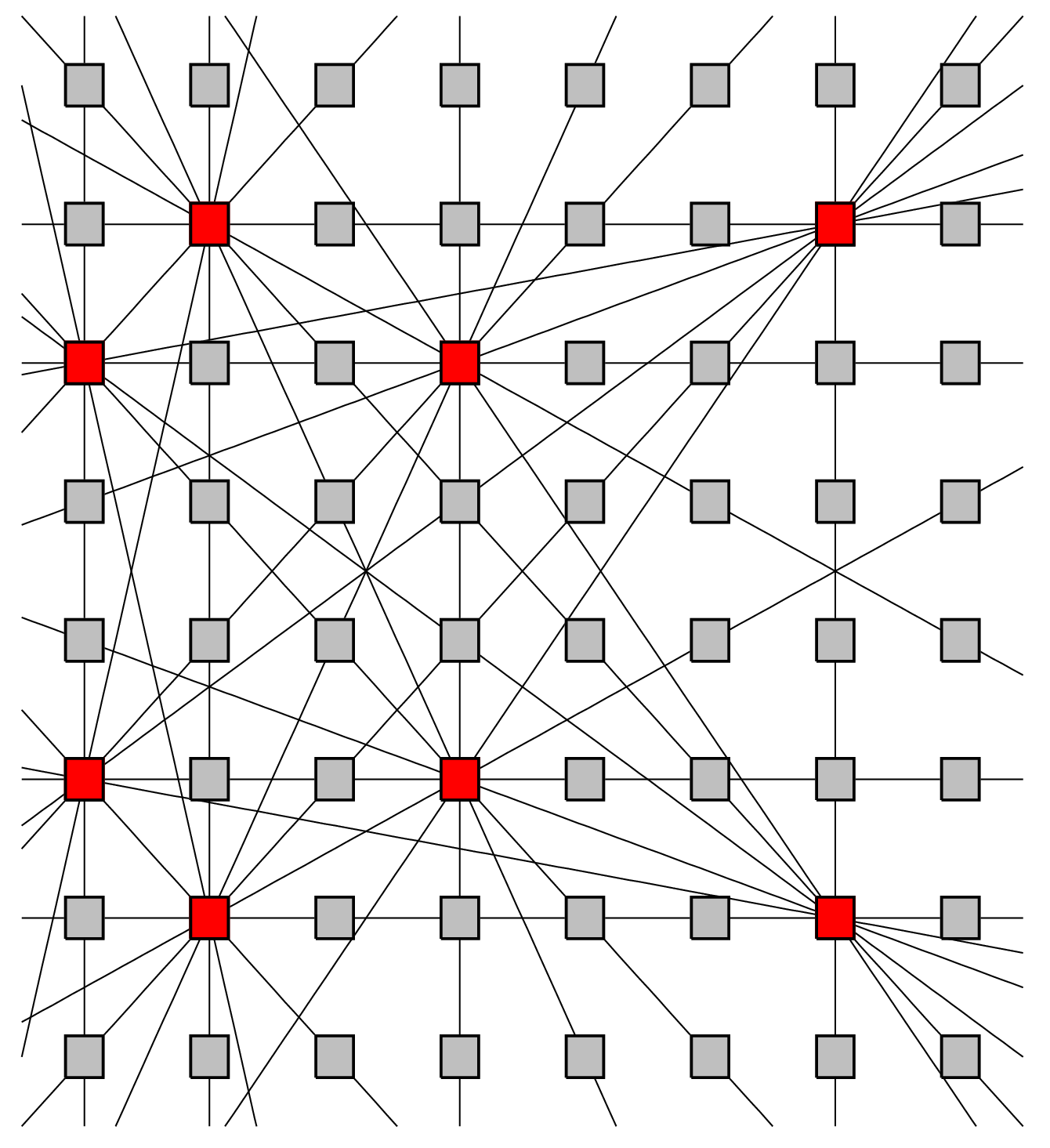}
\includegraphics[scale=0.3]{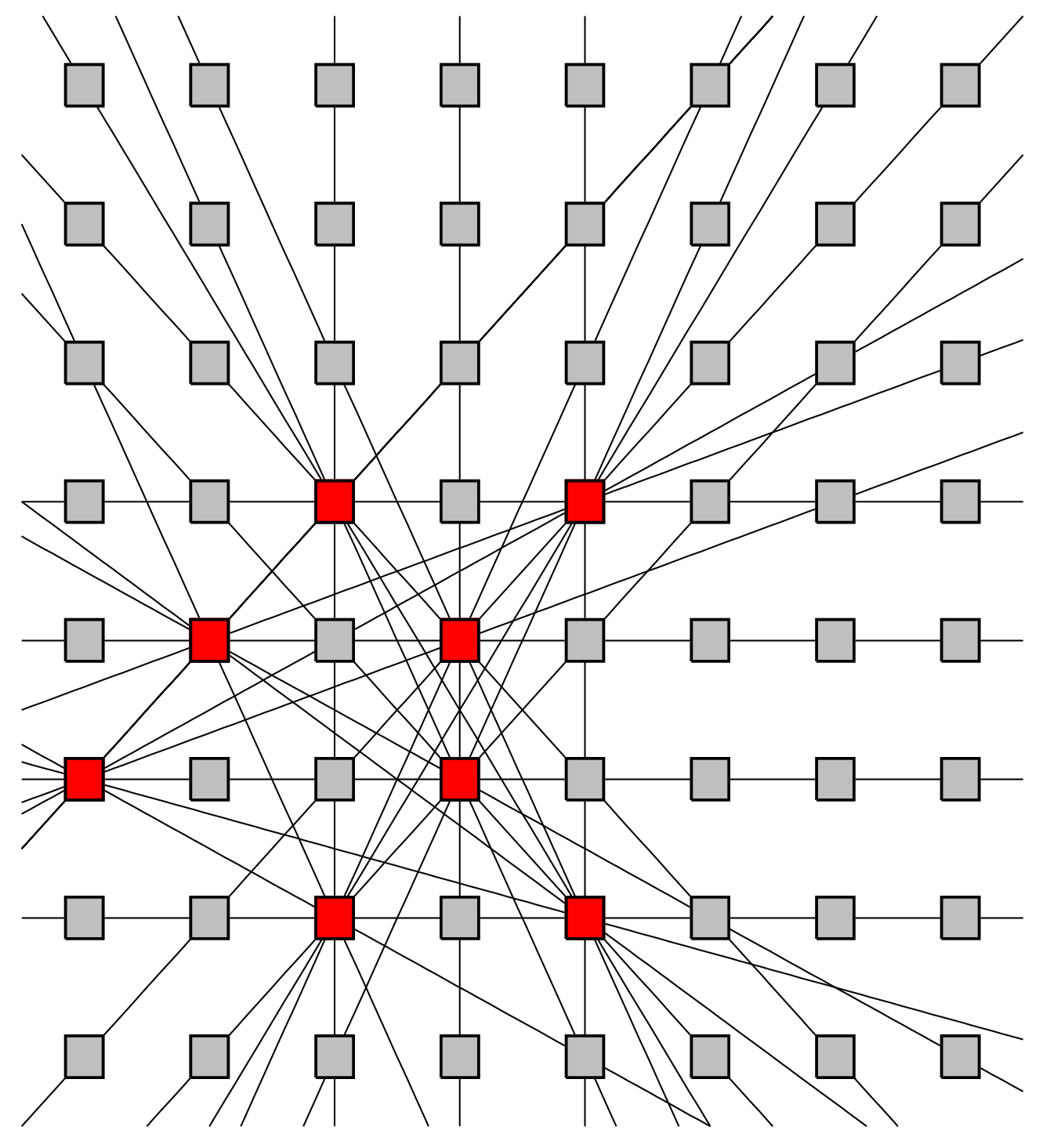}
\includegraphics[scale=0.3]{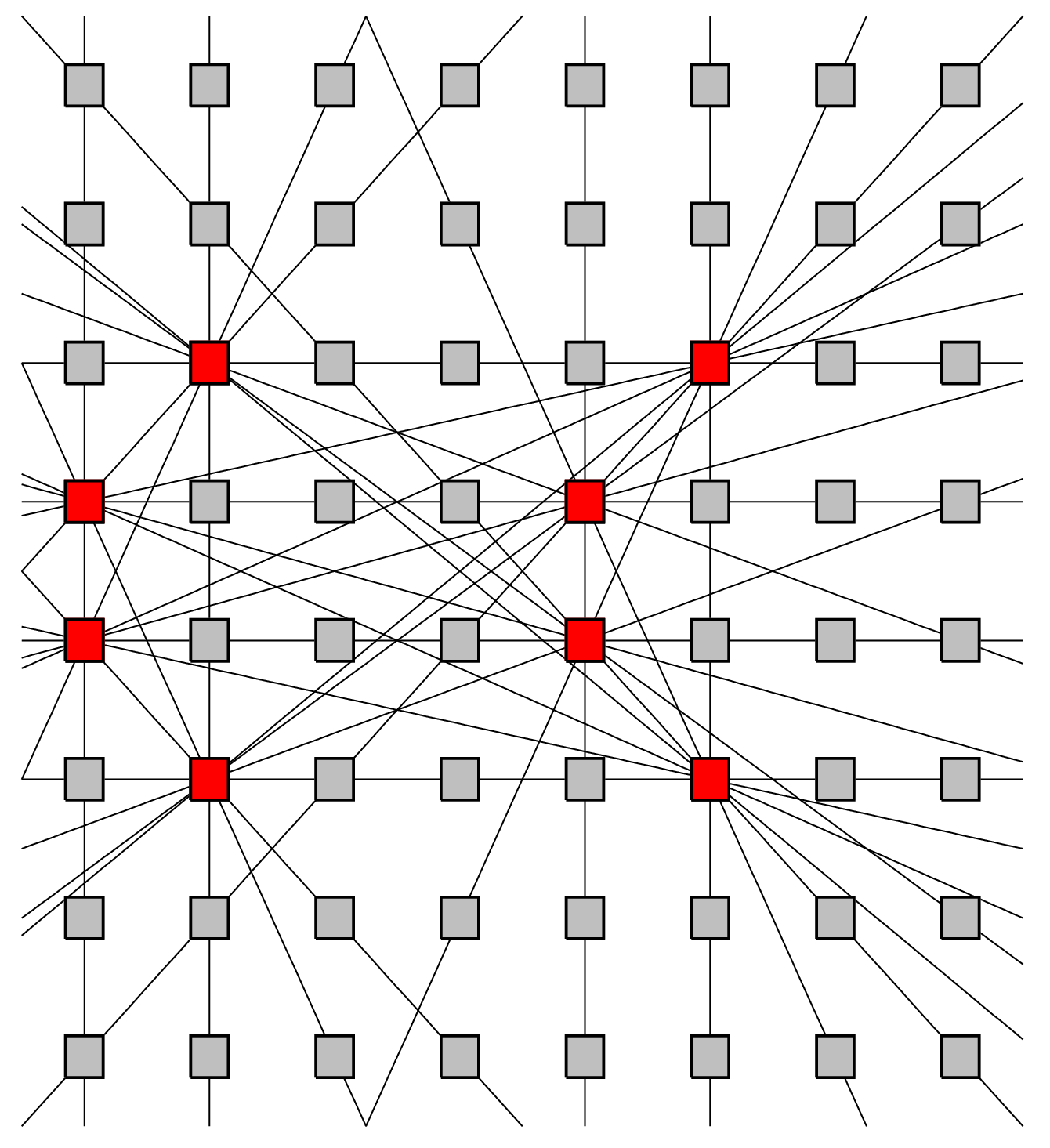}
\includegraphics[scale=0.3]{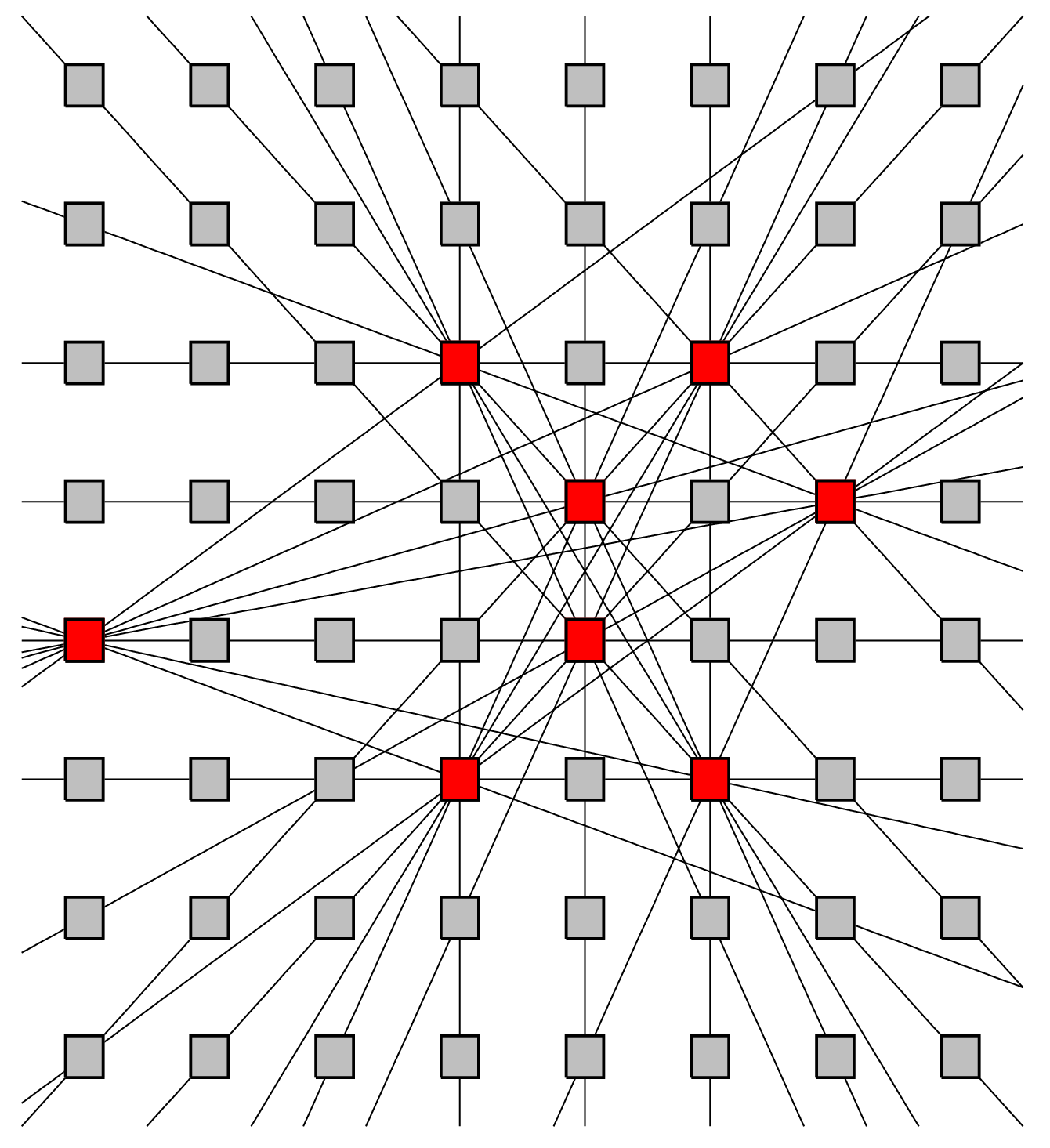}
\includegraphics[scale=0.3]{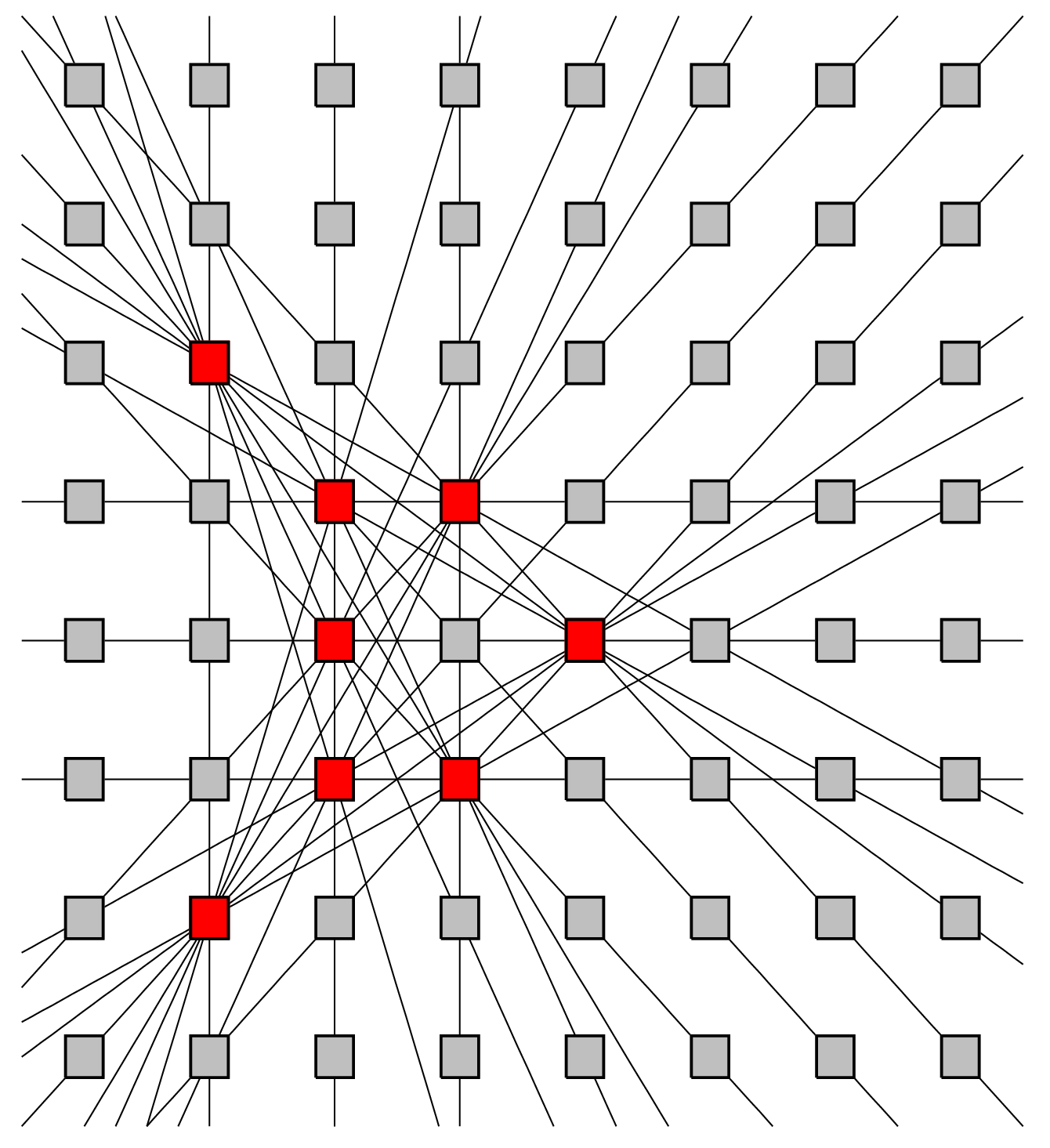}
\includegraphics[scale=0.3]{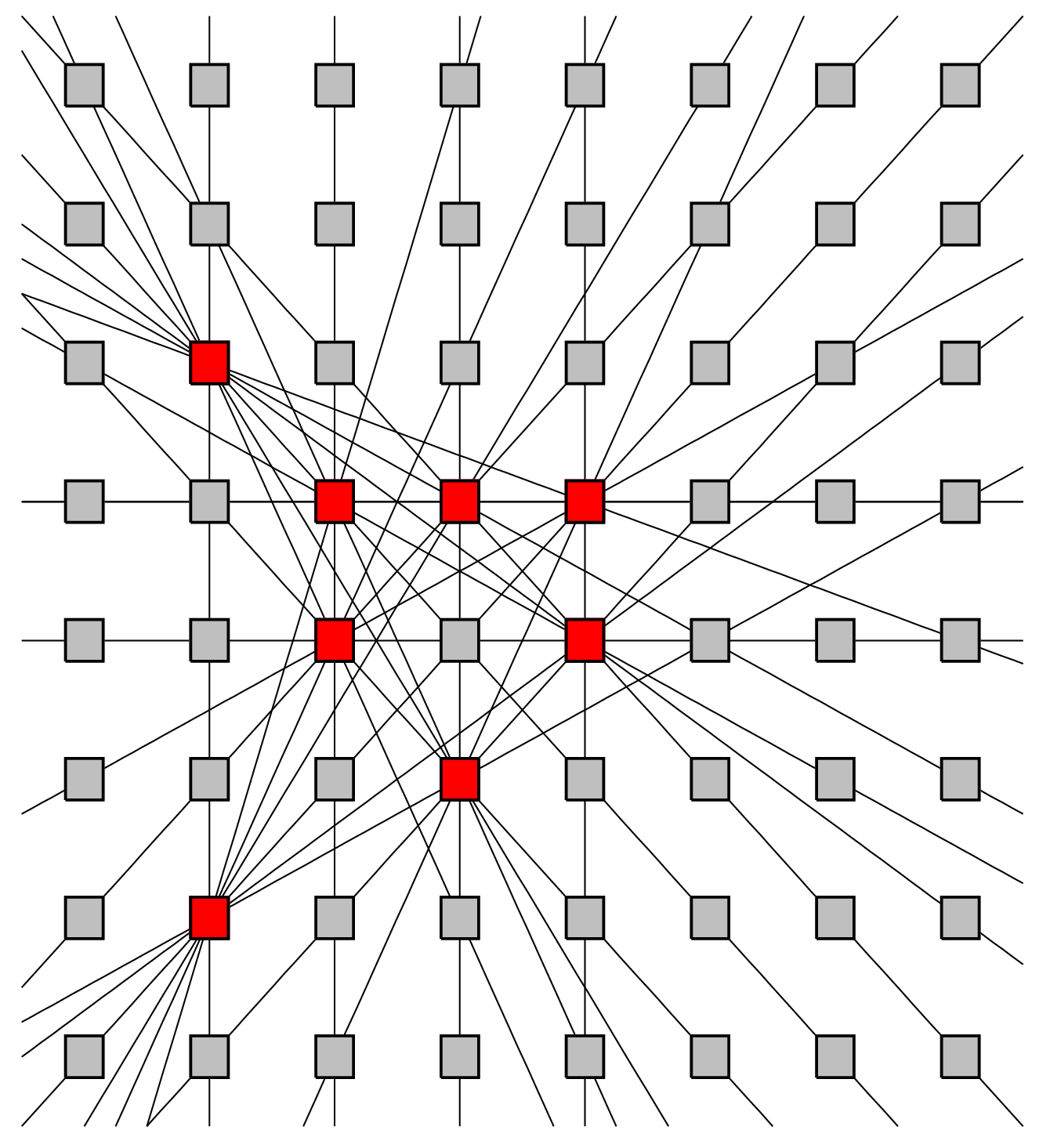}
\includegraphics[scale=0.3]{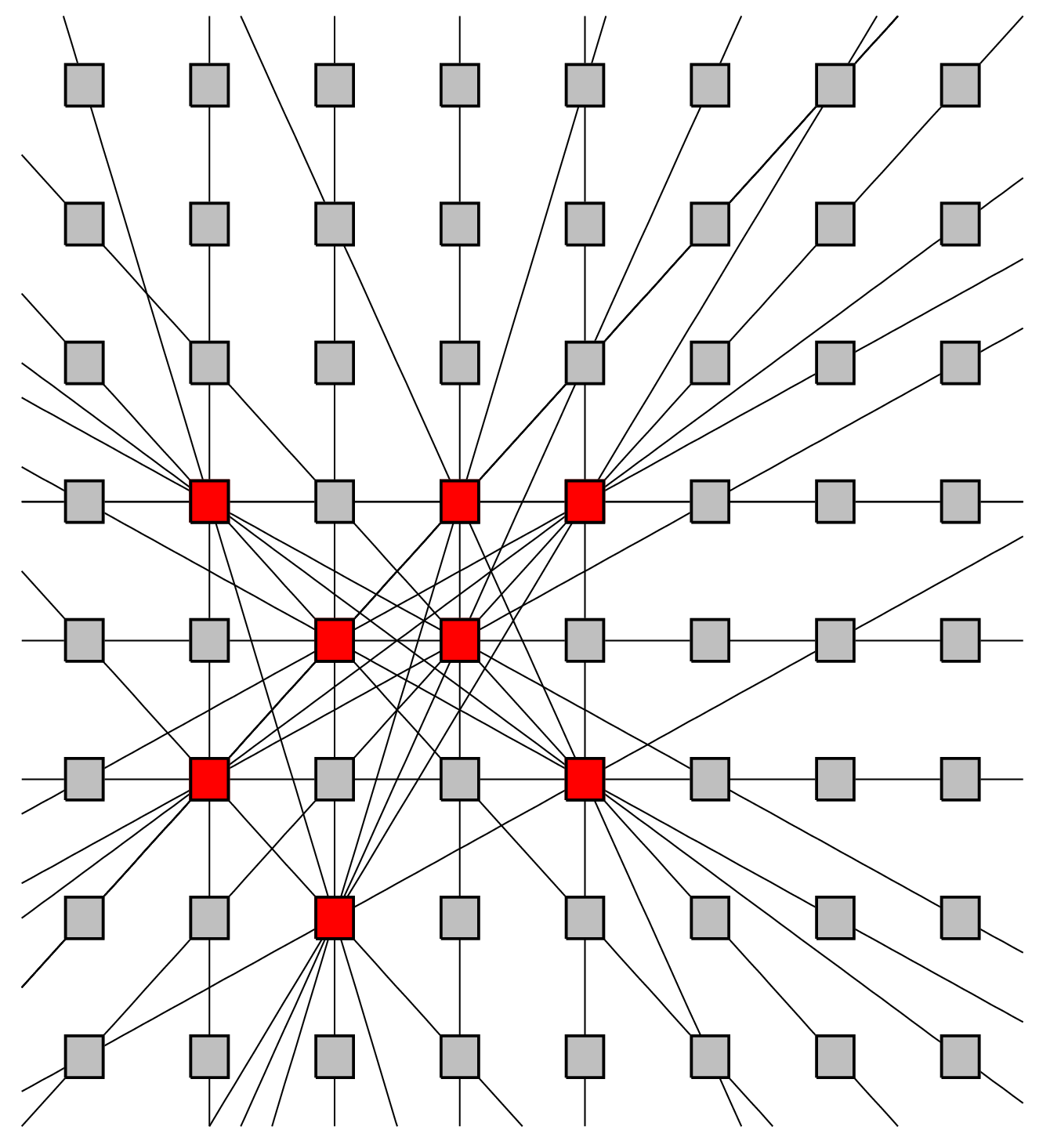}
\includegraphics[scale=0.3]{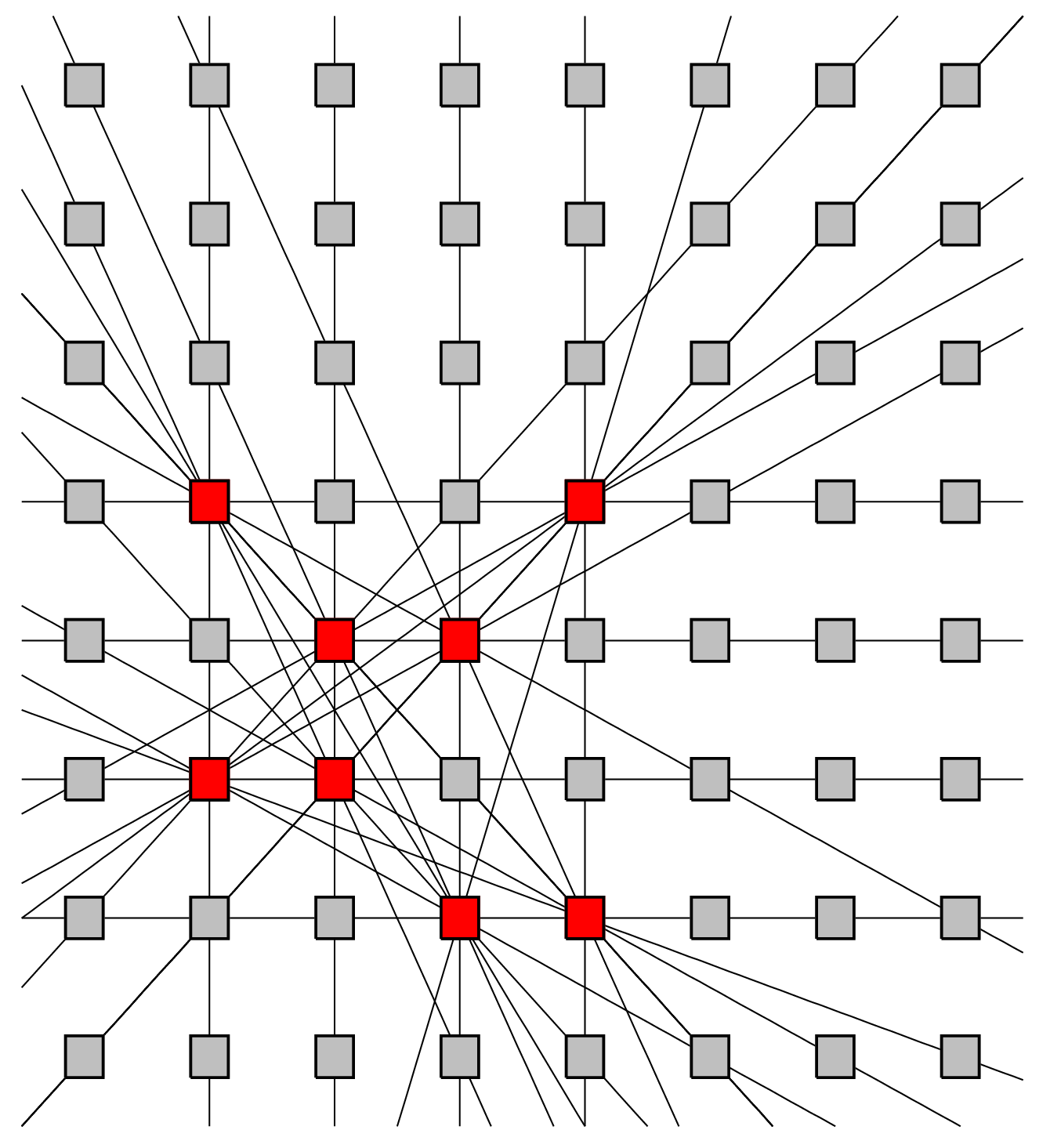}
\includegraphics[scale=0.3]{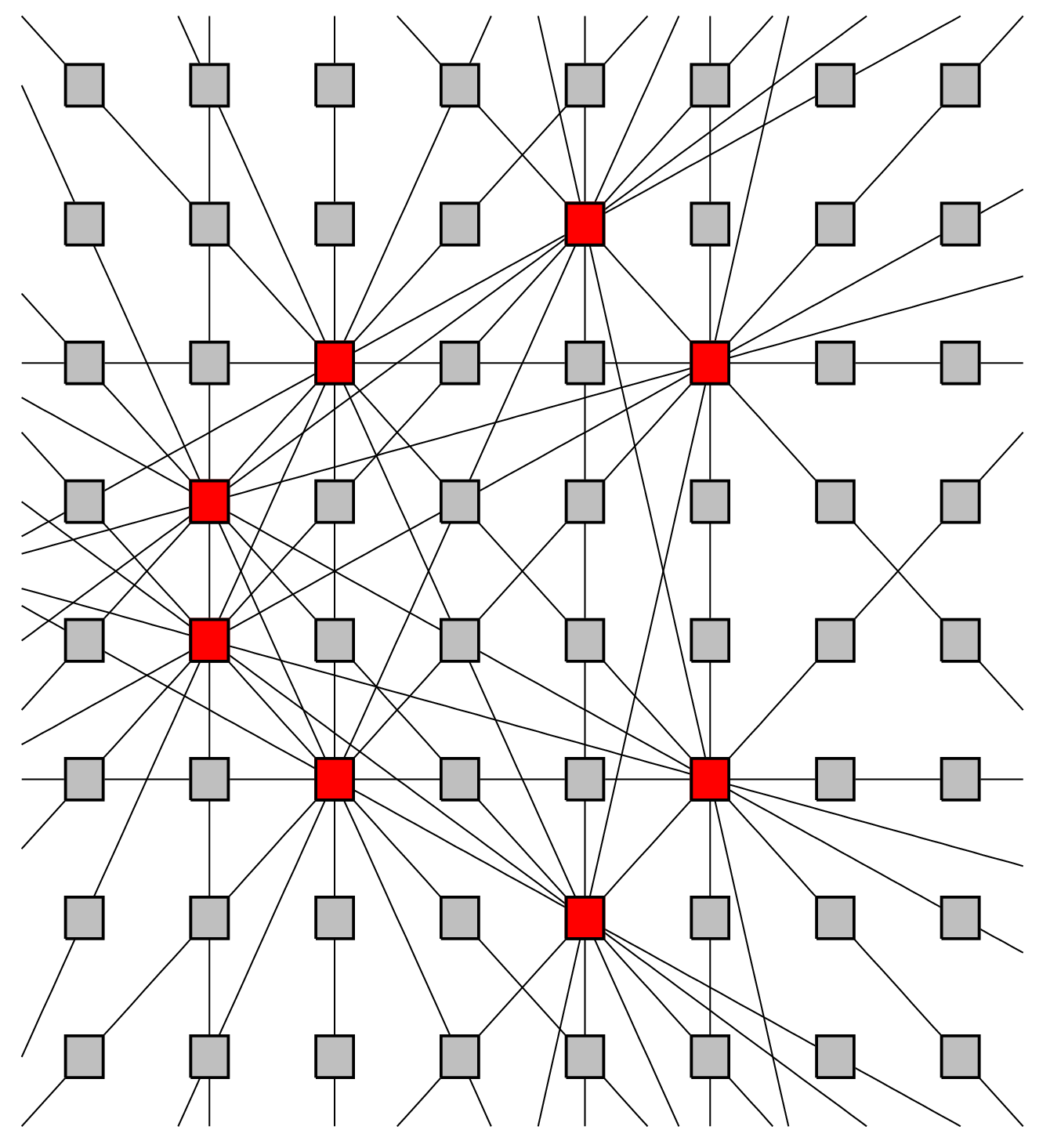}
\includegraphics[scale=0.3]{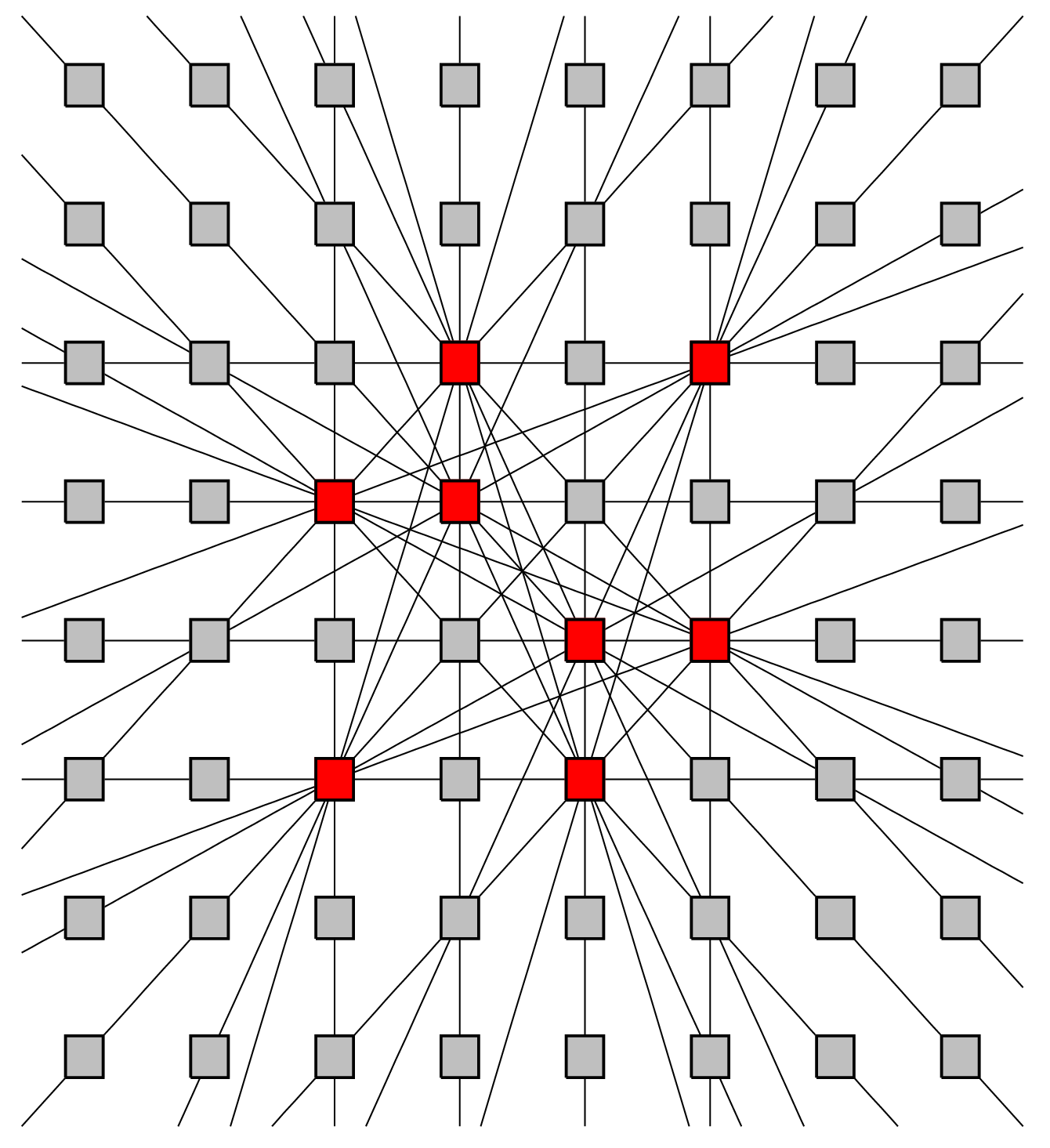}
\includegraphics[scale=0.3]{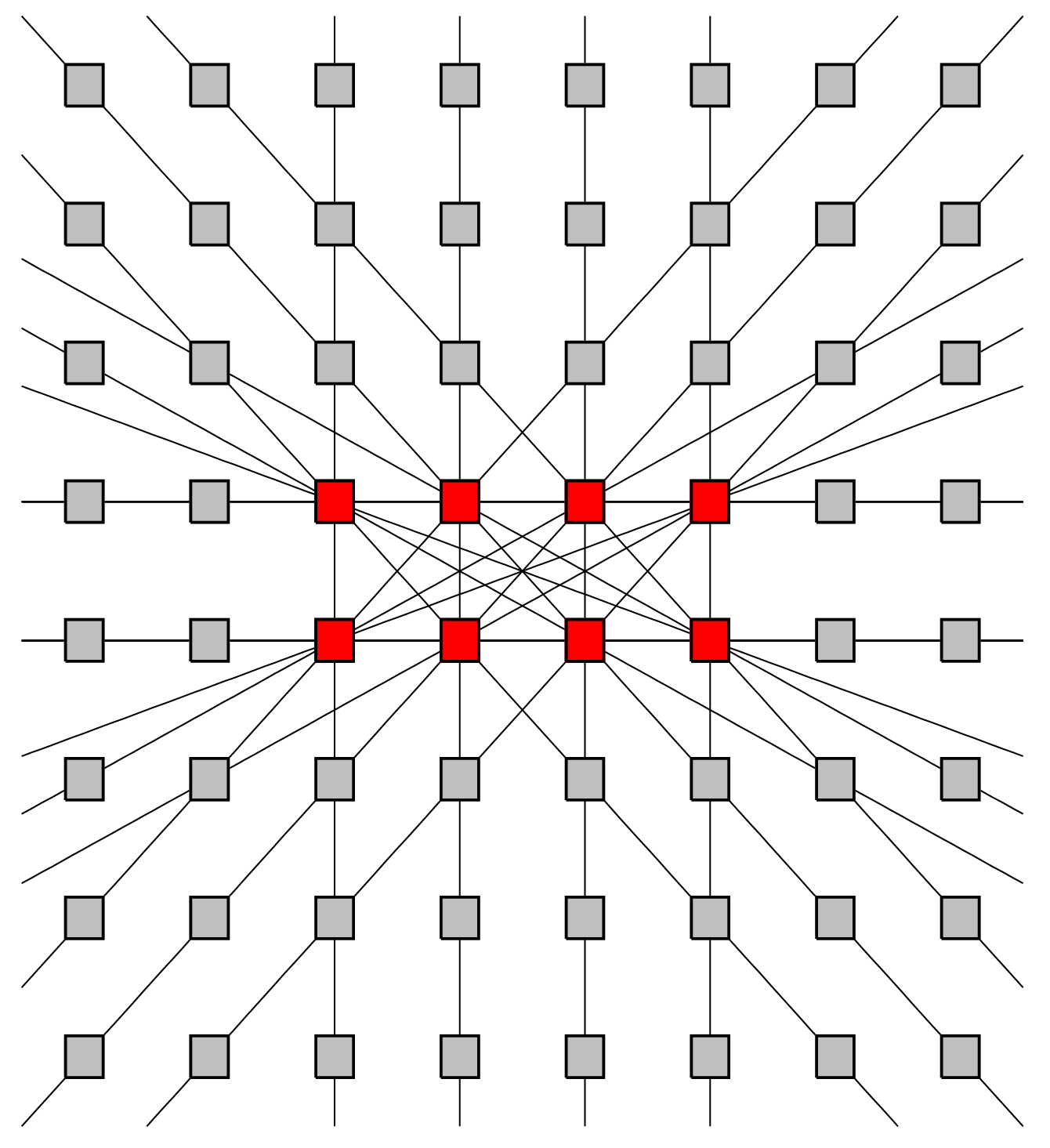}
\caption{
$t(7)\le 8$ (continued). 
}
\end{figure}

\begin{figure}[h]
\includegraphics[scale=0.32]{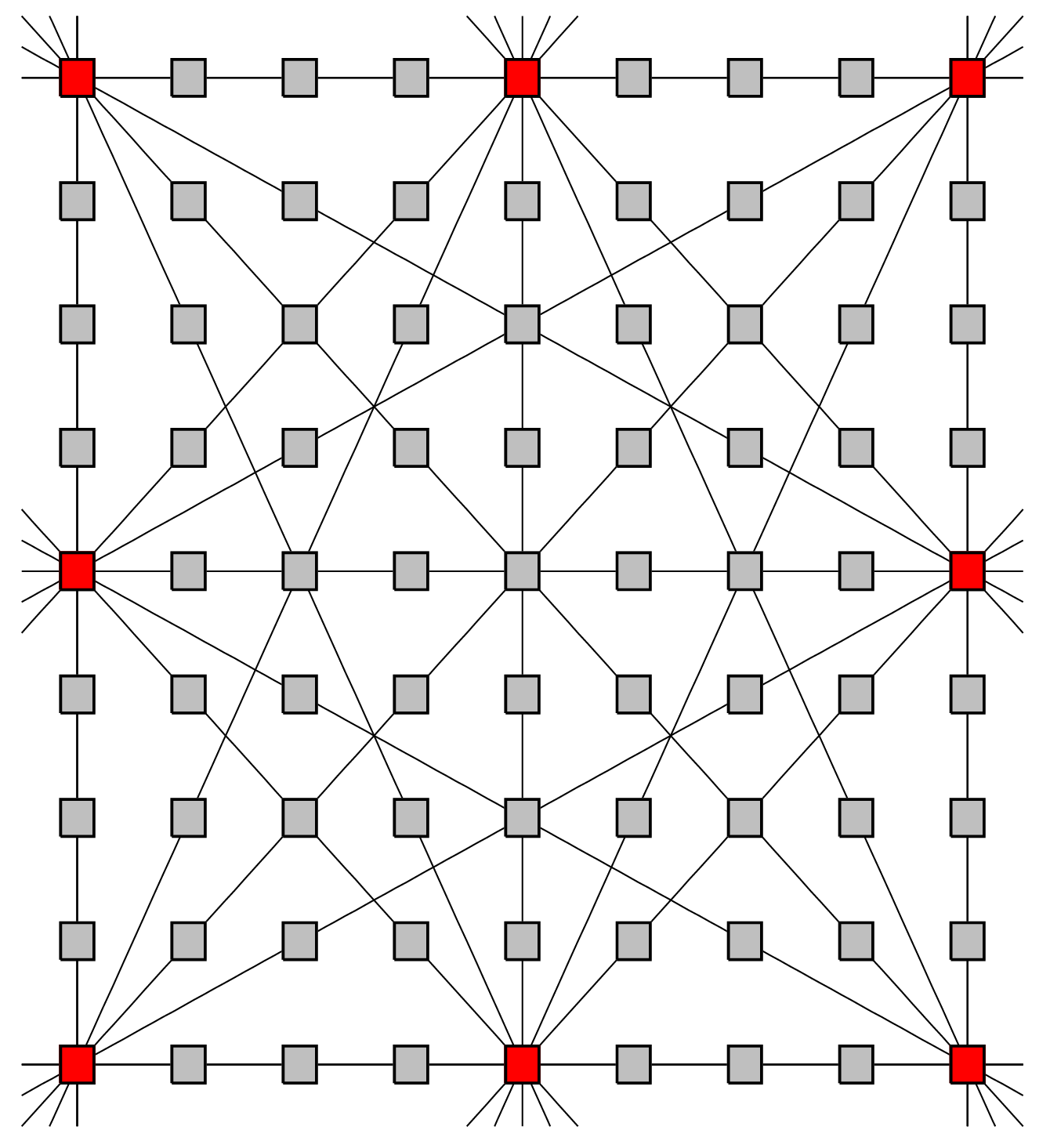}
\includegraphics[scale=0.32]{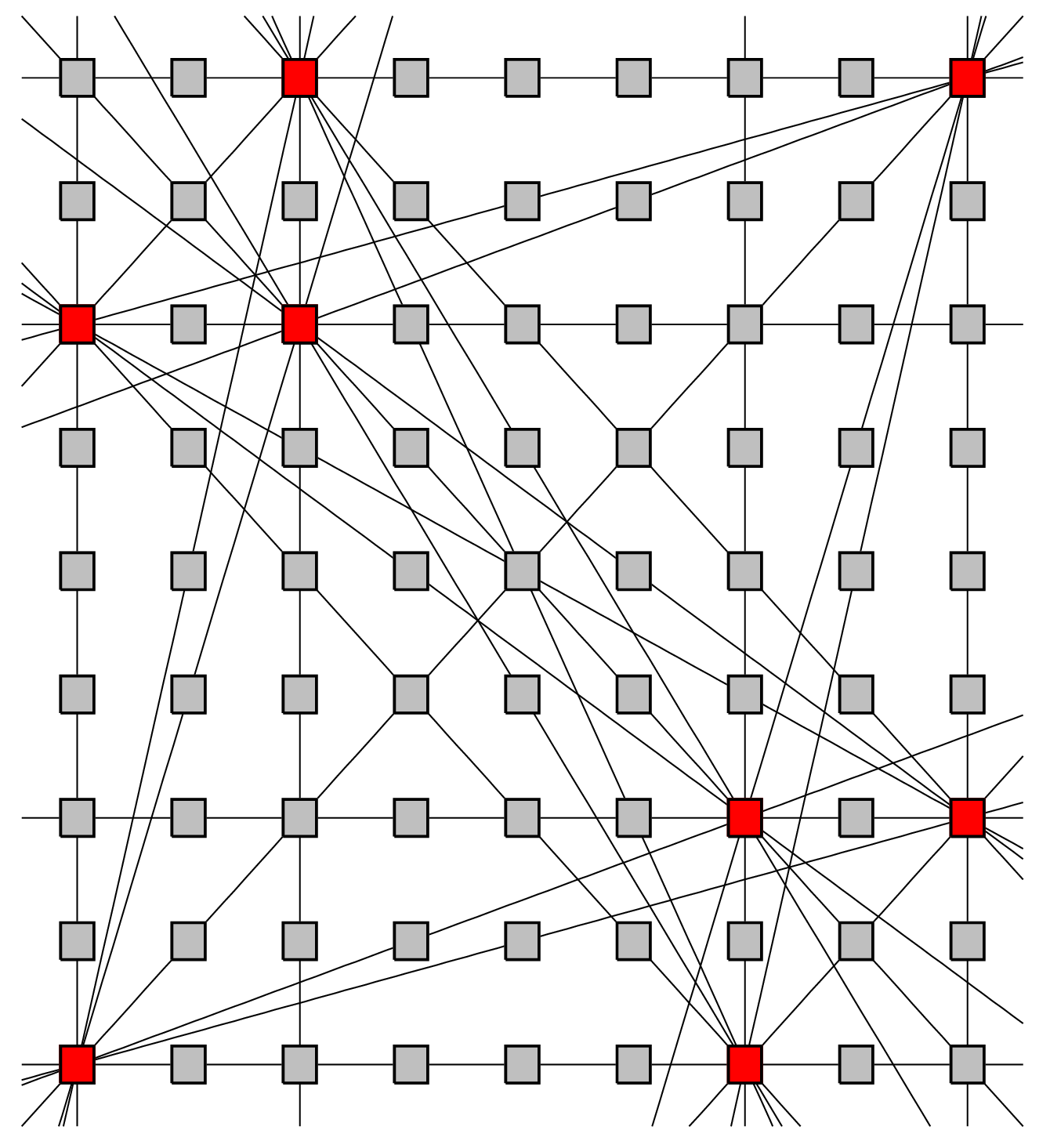}
\includegraphics[scale=0.32]{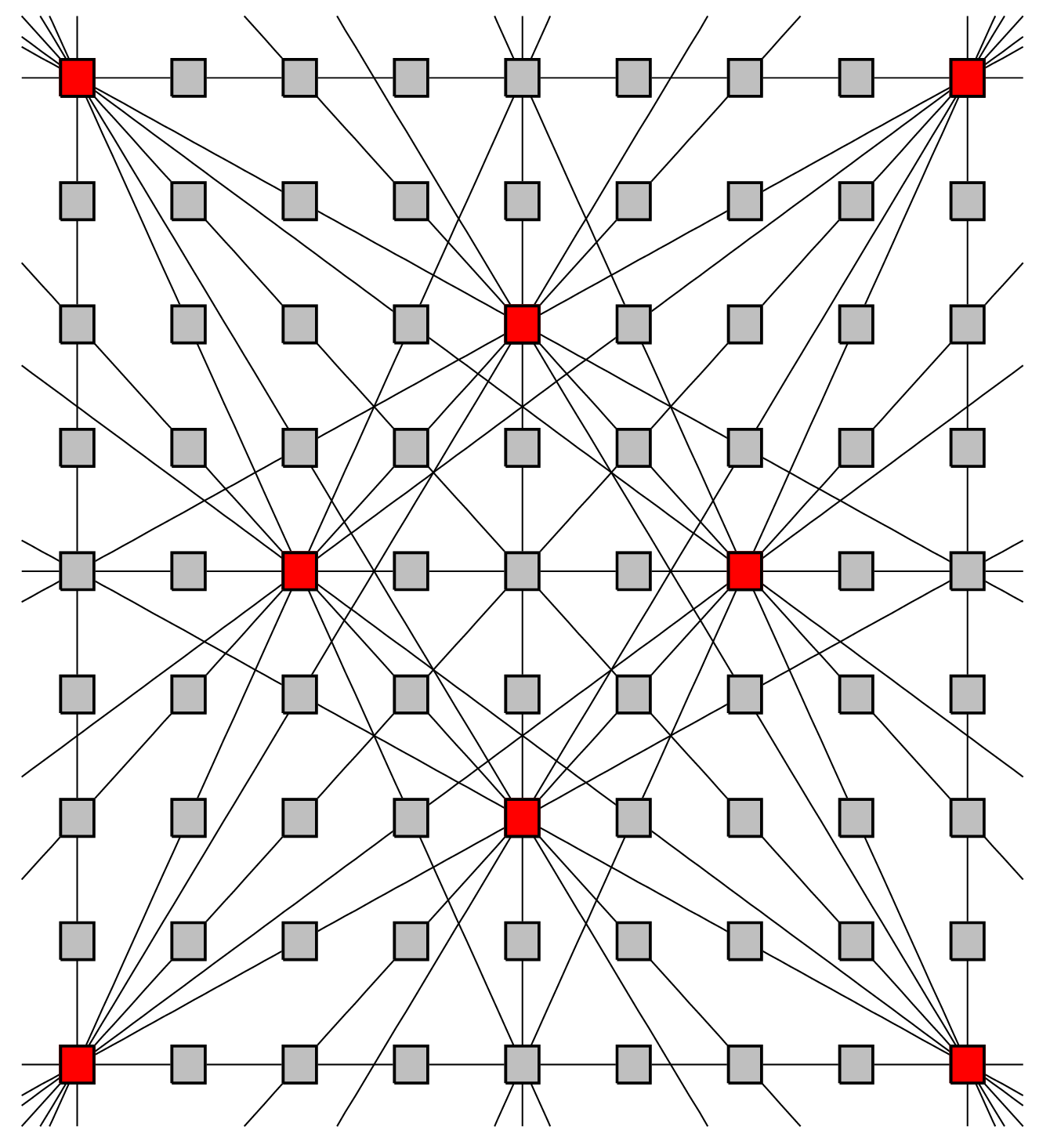}
\includegraphics[scale=0.32]{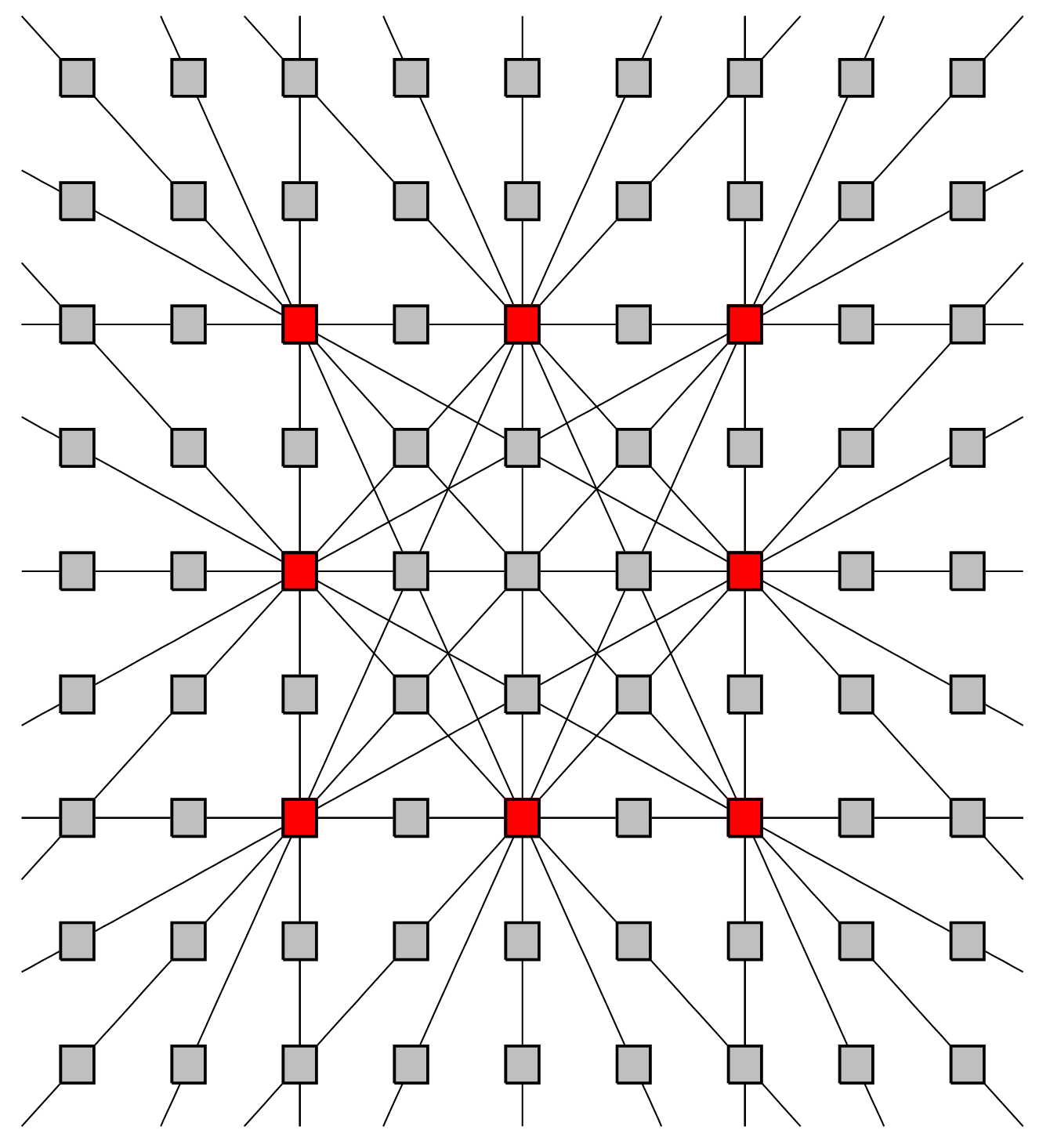}
\caption{
$t(8)\le 8$.
}
\end{figure}
\begin{figure}[h]
\includegraphics[scale=0.32]{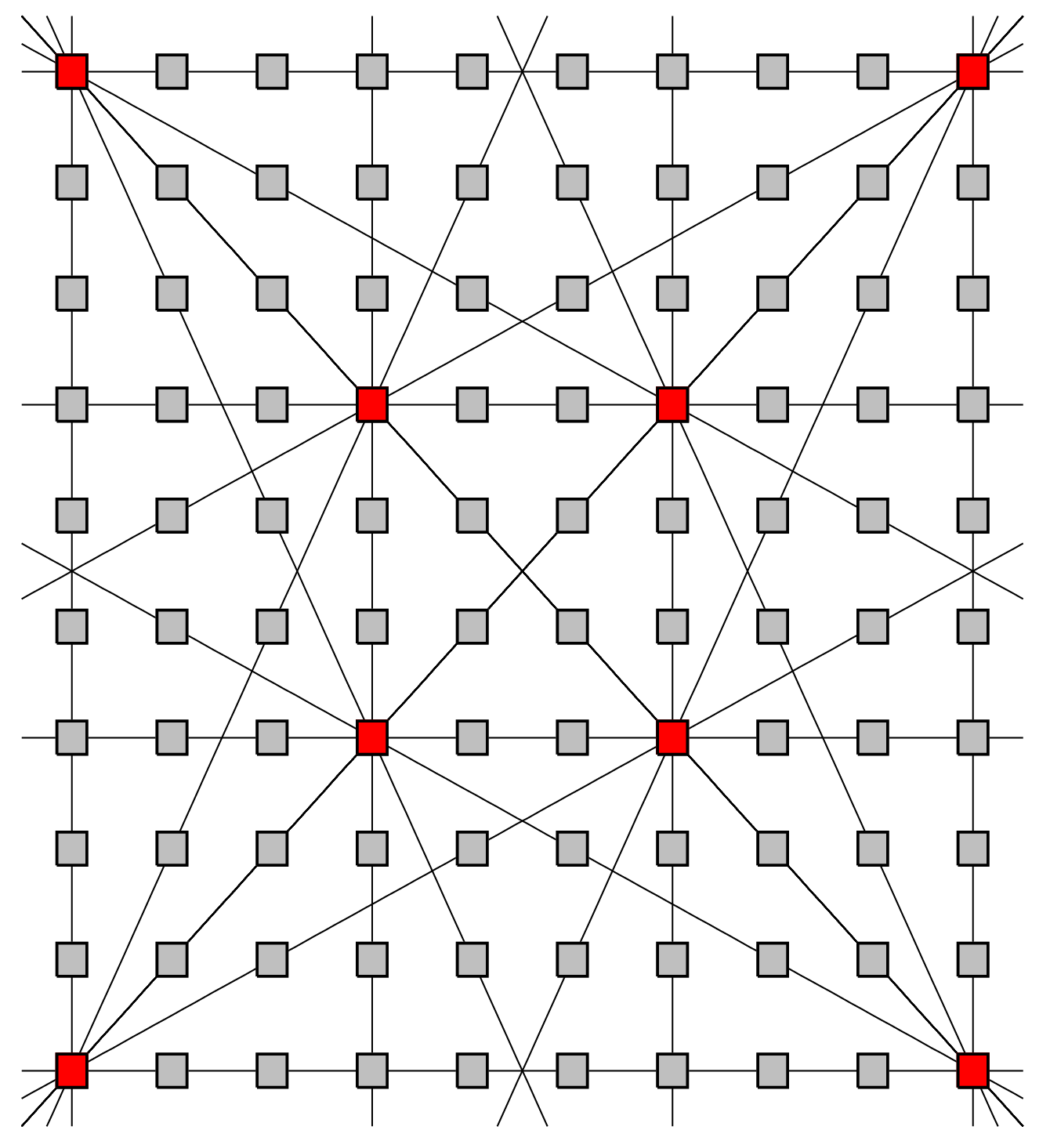}
\includegraphics[scale=0.32]{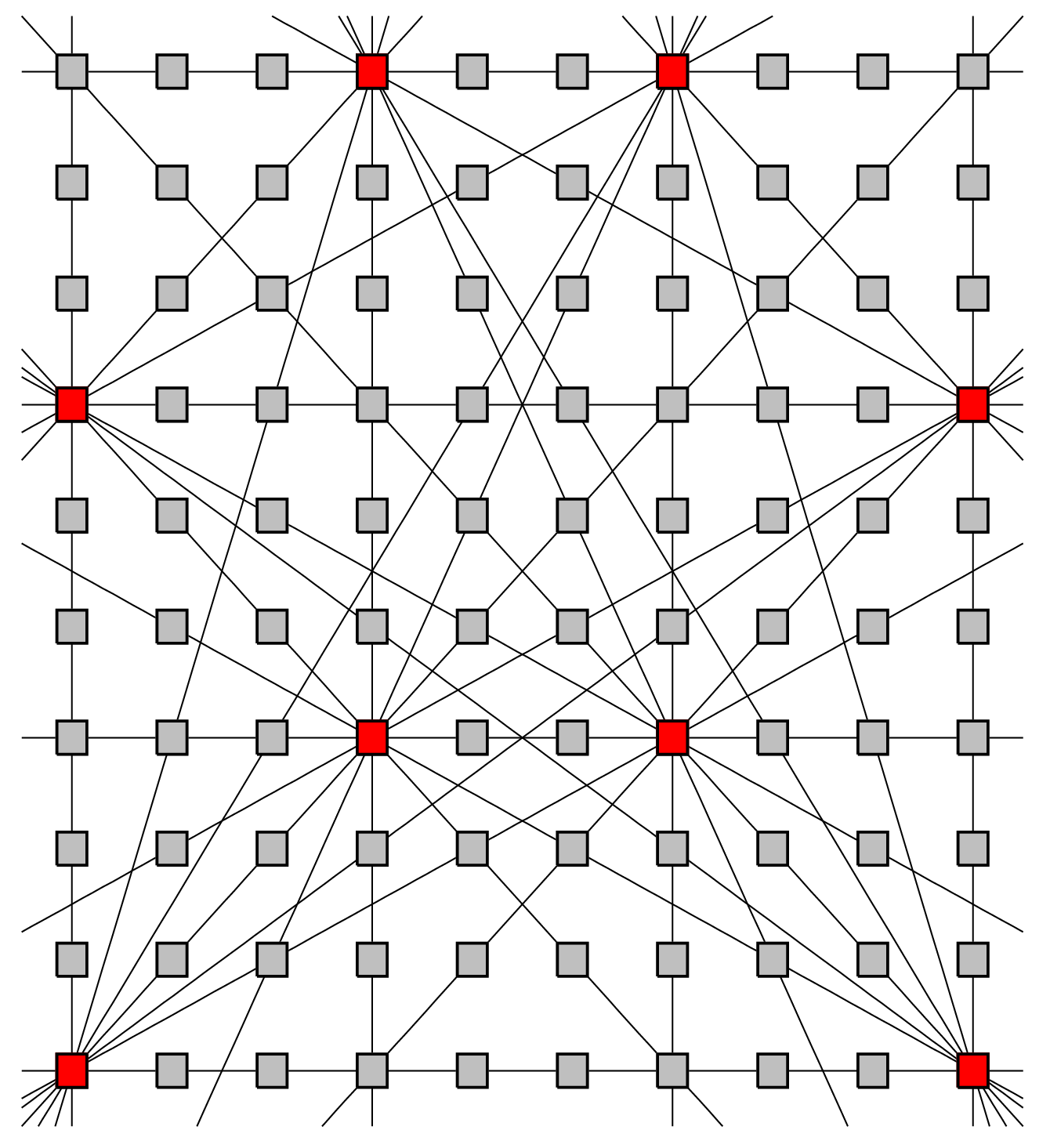}
\caption{
$t(9)\le 8$.
}
\end{figure}


\clearpage
One solution for each $n$ in the range $12$ to 36, then in larger steps up to $n=110$
is put into Table \ref{tab.highn}. This rendering of results saves space
in comparison to figurative presentations, but
symmetries of the sublattices are no longer evident.

\begin{longtable}{r|r|p{0.8\textwidth}}
$n$ & $t(n)$ & sublattice coordinates \\
\hline
\endhead
& & \hfill Table \ref{tab.highn}\ldots\\
\endfoot
\endlastfoot

 12 & $\le 11 $ &  (0,0) (0,3) (2,12) (3,6) (3,9) (6,3) (6,12) (9,6) (9,9) (12,0) (12,3)\\
 13 & $\le 12$ &  (1,1) (2,5) (2,8) (3,10) (5,2) (5,11) (7,6) (8,2) (8,11) (11,5) (11,8) (12,12)\\
 14 & $\le 13 $ &  (4,4) (4,8) (5,5) (5,7) (6,6) (6,7) (6,8) (7,5) (7,6) (7,8) (8,4) (8,6) (8,7)\\
 15 & $\le 13 $ &  (5,7) (5,8) (6,6) (6,7) (6,9) (7,7) (7,8) (7,9) (8,6) (8,7) (8,9) (9,7) (9,8)\\
 16 & $\le 14 $ &  (0,6) (0,9) (3,3) (3,12) (5,6) (6,0) (6,12) (6,15) (9,0) (9,15) (12,3) (12,12) (15,6) (15,9)\\
 17 & $\le 15 $ &  (0,7) (0,10) (3,4) (3,13) (6,1) (6,16) (9,1) (9,16) (12,4) (12,13) (13,8) (13,9) (14,8) (15,7) (15,10)\\
 18 & $\le 16 $ &  (0,4) (2,8) (2,11) (3,7) (4,4) (5,5) (5,14) (8,2) (8,17) (9,15) (11,2) (11,17) (14,5) (14,14) (17,8) (17,11)\\
 19 & $\le 16 $ &  (1,3) (1,16) (2,8) (2,11) (5,5) (5,14) (8,2) (8,17) (11,2) (11,17) (14,5) (14,14) (16,3) (16,16) (17,8) (17,11)\\
 20 & $\le 17 $ &  (0,0) (0,7) (0,20) (5,5) (5,10) (6,13) (7,7) (7,13) (10,5) (10,15) (13,6) (13,7) (13,13) (15,10) (15,15) (20,0) (20,20)\\
 21 & $\le 18 $ &  (0,0) (0,21) (1,6) (3,9) (3,12) (6,6) (6,15) (9,3) (9,18) (12,3) (12,9) (12,18) (15,6) (15,15) (18,9) (18,12) (21,0) (21,21)\\
 22 & $\le 19 $ &  (0,2) (0,16) (1,4) (1,9) (1,14) (6,4) (6,14) (6,19) (10,7) (10,11) (11,4) (11,14) (12,5) (16,4) (16,14) (16,19) (17,8) (21,9) (21,14)\\
 23 & $\le 21 $ &  (0,1) (0,11) (0,16) (5,1) (5,6) (5,21) (10,16) (10,21) (11,20) (13,4) (15,1) (15,6) (15,21) (17,12) (19,14) (19,23) (20,1) (20,6) (20,11) (20,16) (23,12)\\
 24 & $\le 21 $ &  (0,6) (0,15) (2,8) (2,18) (7,3) (7,8) (7,13) (7,18) (7,23) (10,14) (11,16) (12,3) (12,23) (13,16) (17,3) (17,8) (17,13) (17,18) (17,23) (22,8) (22,18)\\
 25 & $\le 21 $ &  (0,0) (0,25) (4,11) (5,10) (5,15) (5,20) (6,6) (10,5) (10,10) (10,20) (11,9) (15,5) (15,15) (15,20) (19,16) (19,19) (20,5) (20,10) (20,15) (25,0) (25,25)\\
 26 & $\le 23 $ &  (0,0) (0,15) (0,25) (5,10) (5,20) (6,7) (6,18) (10,5) (10,15) (10,25) (13,17) (14,14) (14,18) (15,5) (15,10) (15,15) (15,20) (20,0) (20,10) (20,20) (25,5) (25,15) (25,25)\\
 27 & $\le 23 $ &  (0,8) (0,18) (5,3) (5,8) (5,13) (5,18) (5,23) (10,3) (10,8) (10,18) (10,23) (15,3) (15,8) (15,18) (15,23) (20,3) (20,8) (20,13) (20,18) (20,23) (25,8) (25,13) (25,18)\\
 28 & $\le 24 $ &  (1,4) (1,9) (1,14) (1,19) (2,8) (2,16) (2,21) (6,4) (6,9) (6,24) (11,19) (11,24) (16,4) (16,9) (16,24) (18,15) (21,4) (21,9) (21,14) (21,19) (24,7) (24,21) (26,9) (26,19)\\
 29 & $\le 25 $ &  (0,4) (2,2) (2,12) (2,22) (7,12) (7,17) (7,22) (7,27) (12,2) (12,7) (12,22) (12,27) (17,2) (17,7) (17,12) (19,28) (22,7) (22,22) (22,27) (27,12) (27,17) (27,22) (28,4) (28,13) (28,19)\\
 30 & $\le 25 $ &  (0,2) (0,12) (4,3) (5,17) (5,22) (5,27) (10,2) (10,7) (10,22) (10,27) (15,2) (15,7) (15,22) (20,17) (20,22) (20,27) (25,2) (25,7) (25,12) (25,17) (25,22) (29,6) (29,24) (30,12) (30,22)\\
 31 & $\le 26 $ &  (2,7) (2,12) (2,17) (2,22) (7,2) (7,7) (7,12) (7,17) (7,22) (7,27) (12,2) (12,7) (12,22) (12,27) (17,2) (17,7) (17,12) (17,17) (17,22) (17,27) (22,7) (22,12) (22,17) (22,22) (27,12) (27,17)\\
 32 & $\le 27 $ &  (0,8) (0,13) (0,18) (0,23) (5,3) (5,18) (5,23) (5,28) (10,3) (10,8) (10,13) (11,19) (15,8) (15,23) (15,28) (20,8) (20,23) (20,28) (25,3) (25,8) (25,13) (25,28) (29,10) (30,13) (30,18) (30,23) (30,28)\\
 33 & $\le 28 $ &  (0,14) (0,19) (1,15) (5,9) (5,24) (8,13) (8,20) (9,21) (10,9) (10,14) (10,19) (10,24) (11,18) (12,26) (15,9) (15,19) (15,24) (20,4) (20,9) (20,29) (25,4) (25,19) (25,24) (25,29) (30,9) (30,14) (30,19) (30,24)\\
 34 & $\le 28 $ &  (0,5) (0,15) (0,20) (5,0) (5,25) (10,10) (10,15) (10,20) (10,25) (15,0) (15,5) (15,10) (15,25) (15,30) (19,22) (20,5) (20,10) (24,17) (25,5) (25,10) (25,15) (25,20) (25,25) (25,30) (30,10) (30,15) (30,20) (30,25)\\
 35 & $\le 29 $ &  (0,5) (1,1) (1,6) (1,11) (6,11) (6,16) (6,21) (6,26) (11,1) (11,6) (11,11) (11,16) (11,21) (11,26) (16,6) (16,11) (16,26) (16,31) (21,6) (21,11) (21,31) (26,6) (26,11) (26,16) (26,26) (26,31) (31,16) (31,21) (31,26)\\
 36 & $\le 30 $ &  (3,3) (3,13) (3,18) (3,28) (7,27) (8,13) (8,18) (8,23) (8,28) (10,11) (13,3) (13,28) (18,3) (18,8) (18,23) (18,28) (18,33) (20,29) (23,8) (23,13) (23,18) (23,23) (23,28) (23,33) (28,8) (28,13) (28,18) (28,23) (33,8) (33,23)\\
 40 & $\le 33 $ &  (2,10) (2,17) (2,24) (9,3) (9,10) (9,17) (9,24) (9,31) (10,9) (15,30) (16,3) (16,10) (16,17) (16,24) (16,31) (16,38) (23,3) (23,10) (23,17) (23,24) (23,31) (23,38) (24,30) (30,3) (30,10) (30,17) (30,24) (30,31) (36,9) (37,10) (37,17) (37,24) (37,31)\\
 50 & $\le 41 $ &  (0,6) (0,20) (0,34) (6,10) (7,13) (7,27) (7,34) (7,41) (7,48) (14,6) (14,13) (14,18) (14,20) (14,34) (14,41) (21,6) (21,13) (21,27) (21,34) (21,41) (21,48) (28,6) (28,13) (28,34) (28,41) (28,48) (35,6) (35,13) (35,20) (35,34) (35,41) (35,48) (40,45) (42,13) (42,20) (42,27) (42,34) (42,41) (49,6) (49,13) (49,20)\\
 60 & $\le 52 $ &  (0,12) (0,26) (0,40) (0,47) (1,26) (1,33) (7,5) (7,12) (7,47) (7,54) (14,5) (14,19) (14,26) (14,47) (14,54) (21,5) (21,12) (21,19) (21,26) (21,33) (21,40) (24,10) (28,19) (28,40) (28,54) (35,12) (35,19) (35,26) (35,33) (35,40) (35,47) (35,54) (42,5) (42,12) (42,19) (42,26) (42,33) (42,40) (42,47) (43,43) (44,38) (45,25) (49,19) (49,40) (49,47) (54,15) (56,5) (56,19) (56,40) (56,54) (59,8) (59,59)\\
 63 & $\le 55 $ &  (2,6) (2,52) (2,55) (5,3) (5,14) (5,25) (5,36) (5,47) (5,58) (11,34) (11,61) (13,20) (13,52) (14,4) (16,3) (16,14) (16,25) (16,36) (16,47) (16,58) (17,20) (20,28) (27,3) (27,14) (27,25) (27,36) (27,47) (27,58) (29,52) (34,41) (38,3) (38,14) (38,25) (38,36) (38,47) (38,58) (43,21) (45,2) (49,3) (49,14) (49,25) (49,36) (49,47) (49,58) (55,6) (57,8) (57,61) (59,49) (59,53) (60,3) (60,14) (60,25) (60,36) (60,47) (60,58)\\
 70 & $\le 63 $ &  (0,9) (6,3) (6,14) (6,25) (6,36) (6,47) (6,58) (6,69) (7,12) (8,62) (10,62) (11,58) (16,42) (17,3) (17,14) (17,25) (17,36) (17,47) (17,58) (17,69) (20,28) (22,43) (27,9) (28,3) (28,14) (28,25) (28,36) (28,47) (28,58) (28,66) (28,69) (29,37) (33,31) (35,1) (36,18) (39,3) (39,14) (39,25) (39,36) (39,47) (39,58) (39,69) (40,63) (47,19) (50,3) (50,14) (50,25) (50,36) (50,47) (50,58) (50,69) (56,52) (57,66) (58,65) (59,43) (60,12) (60,13) (61,3) (61,14) (61,25) (61,36) (61,47) (61,58)\\
 80 & $\le 69 $ &  (0,4) (0,51) (1,1) (1,12) (1,23) (1,34) (1,67) (5,7) (6,31) (6,42) (12,1) (12,23) (12,45) (12,56) (12,67) (12,78) (18,72) (20,18) (22,0) (23,1) (23,12) (23,23) (23,34) (23,45) (23,56) (23,67) (23,78) (29,4) (29,80) (34,12) (34,23) (34,34) (34,45) (34,56) (34,67) (34,78) (38,42) (38,74) (40,13) (45,1) (45,23) (45,34) (45,56) (45,67) (45,78) (49,31) (49,74) (56,12) (56,23) (56,34) (56,45) (56,56) (56,67) (56,78) (61,49) (67,1) (67,12) (67,23) (67,34) (67,56) (67,67) (67,78) (69,51) (69,80) (73,71) (78,12) (78,23) (78,45) (78,67)\\
 90 & $\le 80 $ &  (1,3) (1,59) (9,2) (9,65) (10,3) (10,72) (16,16) (16,23) (16,30) (16,37) (16,44) (16,51) (16,58) (16,65) (16,72) (19,13) (23,16) (23,23) (23,30) (23,37) (23,65) (23,67) (23,72) (23,79) (28,27) (28,67) (30,9) (30,16) (30,23) (30,30) (30,37) (30,44) (30,51) (30,58) (30,65) (37,9) (37,30) (37,51) (37,79) (40,23) (41,27) (41,58) (42,6) (42,26) (44,16) (44,54) (44,58) (44,65) (51,23) (51,30) (51,44) (51,65) (51,72) (51,79) (54,44) (56,26) (56,59) (58,9) (58,16) (58,44) (58,58) (65,23) (65,30) (65,44) (65,58) (65,65) (65,72) (65,79) (68,40) (72,16) (72,23) (72,58) (72,65) (72,72) (79,2) (79,23) (79,44) (79,51) (86,9) (86,44)\\
 101 & $\le 90 $ &  (0,5) (0,16) (0,27) (0,60) (0,71) (0,82) (2,7) (10,89) (11,5) (11,14) (11,27) (11,38) (11,60) (11,71) (11,82) (11,93) (14,88) (22,5) (22,11) (22,16) (22,27) (22,38) (22,49) (22,60) (22,71) (22,82) (22,93) (26,46) (30,4) (30,14) (32,64) (33,5) (33,16) (33,27) (33,49) (33,60) (33,71) (33,82) (33,93) (44,5) (44,16) (44,27) (44,38) (44,93) (49,24) (49,55) (50,54) (55,16) (55,27) (55,38) (55,49) (55,60) (55,71) (55,82) (55,93) (56,96) (62,100) (66,16) (66,27) (66,49) (66,71) (66,82) (66,93) (72,30) (74,68) (76,48) (77,5) (77,16) (77,27) (77,38) (77,49) (77,60) (77,82) (84,45) (85,89) (88,16) (88,27) (88,49) (88,71) (88,82) (88,93) (89,42) (90,38) (91,46) (91,54) (99,38) (99,49) (99,60) (99,71) (99,82)\\
 110 & $\le 100 $ &  (1,4) (1,91) (4,94) (6,22) (18,6) (20,4) (20,21) (38,73) (40,58) (41,54) (42,48) (44,44) (44,55) (46,10) (46,44) (46,46) (47,33) (47,47) (47,50) (48,28) (48,42) (48,45) (48,49) (49,48) (49,55) (50,24) (50,34) (50,45) (51,45) (51,48) (51,49) (52,48) (52,49) (52,51) (52,54) (53,43) (53,48) (53,50) (53,52) (53,53) (54,39) (54,43) (54,47) (54,49) (54,50) (54,51) (55,21) (55,39) (55,40) (55,42) (55,46) (55,49) (55,53) (56,45) (56,47) (56,53) (56,55) (57,43) (57,44) (57,47) (58,30) (58,46) (58,47) (58,49) (58,50) (58,53) (59,46) (59,50) (59,51) (59,52) (60,24) (60,34) (60,41) (60,49) (61,41) (61,50) (61,51) (61,62) (62,28) (62,40) (62,42) (63,42) (63,44) (63,47) (64,10) (64,46) (64,48) (65,42) (65,50) (66,45) (66,77) (69,64) (70,52) (70,76) (74,76) (80,7) (91,11) (92,6) (106,94) (109,91)\\

\hline
\caption{Examples of coverage of some larger arrays. The sublattice points'
Cartesian coordinates $(x,y)$ are listed in the range $0\le x,y\le n$.
}
\label{tab.highn}
\end{longtable}

\section{Summary}\label{sec.summ}

Square lattices of size $3\times 3$ or $4\times 4$ can be covered
by base lines spanned by sublattices of order $t=4$.
For square lattices of size $5\times 5$ or $6\times 6$,
sublattices must be at least of order $t=6$, for size $7\times 7$ of order $t=7$.
For all other $(n+1)\times(n+1)$ lattices studied, the
order stayed below an upper limit of $t(n)< (n+1)^{2/3}\log(n+1)$.

\bibliographystyle{amsplain}
\bibliography{all}

\end{document}